\documentclass[12pt]{amsart}
\usepackage{amscd,amssymb,graphics,a4wide,color,hyperref}
\usepackage{srcltx}

\usepackage{mathrsfs}
\input xy
\xyoption{all}
\numberwithin{equation}{section}

\newtheorem{theorem}{Theorem}[section]
\newtheorem{lemma}[theorem]{Lemma}
\newtheorem{proposition}[theorem]{Proposition}
\newtheorem{corollary}[theorem]{Corollary}

\theoremstyle{definition}
\newtheorem{definition}[theorem]{Definition}
\newtheorem{construction}[theorem]{Construction}
\newtheorem{example}[theorem]{Example}

\theoremstyle{remark}
\newtheorem{remark}[theorem]{Remark}

\newcommand{\HH} {\mathbb{H}}
\newcommand{\NN} {\mathbb{N}}
\newcommand{\ZZ} {\mathbb{Z}}
\newcommand{\QQ} {\mathbb{Q}}
\newcommand{\RR} {\mathbb{R}}

\renewcommand{\AA} {\mathbb{A}}
\newcommand{\GG} {\mathbb{G}}

\newcommand {\shAff} {\mathcal{A}\text{\textit{ff}}}

\newcommand {\shB}  {\mathcal{B}}

\newcommand {\shExt} {\mathcal{E} \!\text{\textit{xt}}}
\newcommand {\shE}  {\mathcal{E}}
\newcommand {\shF}  {\mathcal{F}}

\newcommand {\shHom} {\mathcal{H}\!\text{\textit{om}}}
\newcommand {\shI}  {\mathcal{I}}

\newcommand {\shK}  {\mathcal{K}}
\newcommand {\shL}  {\mathcal{L}}

\newcommand {\shQ}  {\mathcal{Q}}

\newcommand {\shT}  {\mathcal{T}}

\newcommand {\shV}  {\mathcal{V}}
\newcommand {\shX}  {\mathcal{X}}

\newcommand {\shZ}  {\mathcal{Z}}

\newcommand {\Aff}  {\operatorname{Aff}}

\newcommand {\bct} {\mathrm{bct}}

\newcommand {\C} {\mathscr{C}}

\newcommand {\cl}  {\operatorname{cl}}

\newcommand {\coker} {\operatorname{coker}}
\newcommand {\Conv} {\operatorname{Conv}}
\newcommand {\Cone} {\operatorname{Cone}}

\newcommand {\Diag} {\operatorname{Diag}}

\newcommand {\Di} {{\rm D}}
\newcommand {\di} {{\rm d}}

\newcommand {\dlog} {\operatorname{dlog}}
\newcommand {\Dlog} {\operatorname{Dlog}}

\newcommand {\dual} {\vee}

\newcommand {\eqand} {\quad\text{and}\quad}

\newcommand {\Ext}  {\operatorname{Ext}}

\newcommand {\GL}  {\operatorname{GL}}
\newcommand {\gp}  {{\operatorname{gp}}}

\newcommand {\Hom}  {\operatorname{Hom}}

\newcommand {\id}  {\operatorname{id}}
\newcommand {\im}  {\operatorname{im}}

\newcommand {\Int}  {\operatorname{Int}}

\renewcommand {\ker } {\operatorname{ker}}
\newcommand {\kk} {\Bbbk}

\newcommand {\lra}  {\longrightarrow}
\newcommand {\ls}  {\dagger}

\newcommand {\M} {\mathcal{M}}

\renewcommand{\O}  {\mathcal{O}}

\renewcommand{\P}  {\mathscr{P}}

\newcommand {\pre}  {\mathrm{pre}}

\newcommand {\rank} {\operatorname{rank}}

\newcommand {\Res}  {\operatorname{Res}}

\newcommand {\sing} {\mathrm{sing}}
\newcommand {\Sing} {\operatorname{Sing}}
\newcommand {\SL}  {\operatorname{SL}}
\newcommand {\Spec} {\operatorname{Spec}}
\newcommand {\Spf}  {\operatorname{Spf}}

\renewcommand {\top} {\mathrm{top}}
\newcommand{\Top}{\operatorname{Top}}
\newcommand {\Tor}  {\operatorname{Tor}}
\newcommand {\Tors}  {\operatorname{Tors}}
\newcommand {\tors}  {\operatorname{tors}}
\newcommand {\Tot}  {\operatorname{Tot}}

\renewcommand {\Vert} {\operatorname{Vert}}

\newcommand{\bbfamily}{\fontencoding{U}\fontfamily{bbold}\selectfont}
\newcommand{\textbb}[1]{{\bbfamily#1}}
\newcommand {\lfor} {\mbox{\textbb{[}}}
\newcommand {\rfor} {\mbox{\textbb{]}}}
\newcommand {\T} {\shT}
\newcommand {\X} {\shX}
\newcommand {\V} {\shV}

\newcommand {\Gm} {\GG_m}
\newcommand {\shLS} {\mathcal{LS}}
\newcommand {\lss} {{\bf ls}}

\def\mydate{\ifcase\month \or January\or February\or March\or
April\or May\or June\or July\or August\or September\or October\or 
November\or December\fi \space\number\day,\space\number\year}

\begin{document}
\def\mapright#1{\smash{
  \mathop{\longrightarrow}\limits^{#1}}}
\def\mapleft#1{\smash{
  \mathop{\longleftarrow}\limits^{#1}}}
\def\exact#1#2#3{0\to#1\to#2\to#3\to0}
\def\mapup#1{\Big\uparrow
   \rlap{$\vcenter{\hbox{$\scriptstyle#1$}}$}}
\def\mapdown#1{\Big\downarrow
   \rlap{$\vcenter{\hbox{$\scriptstyle#1$}}$}}
\def\dual#1{{#1}^{\scriptscriptstyle \vee}}
\def\invlim{\mathop{\rm lim}\limits_{\longleftarrow}}

\title[Mirror Symmetry via Logarithmic Degeneration Data II]{Mirror Symmetry
via Logarithmic Degeneration Data II}

\author{Mark Gross} \address{UCSD Mathematics,
9500 Gilman Drive, La Jolla, CA 92093-0112, USA}
\email{mgross@math.ucsd.edu}
\thanks{This work was partially supported by NSF grant 0505325 and
DFG priority programs ``Globale Methoden in der komplexen Geometrie''
and ``Globale Differentialgeometrie''.}

\author{Bernd Siebert} \address{Mathematisches Institut,
Albert-Ludwigs-Universit\"at, Eckerstra\ss e~1, 79104 Freiburg,
Germany}
\email{bernd.siebert@math.uni-freiburg.de}
\date{\today}
\maketitle
\tableofcontents
\bigskip

\section*{Introduction.}

This paper is a sequel to \cite{PartI}, where we laid the foundations
to a program for studying mirror symmetry using logarithmic
geometry. This program can be viewed as an algebro-geometric
version of the Strominger-Yau-Zaslow (SYZ) program \cite{SYZ}, 
and gives a way of
passing via affine geometry to the two sides of the mirror
symmetry picture.

We recall briefly the setup and goals of this program. We begin with
$B$ an integral affine manifold with singularities: there is an open
subset $B_0\subseteq B$ with an atlas with transition maps in
$\Aff(\ZZ^n)$, and $\Delta:=B\setminus B_0$ is codimension $\ge 2$
in $B$ (see \cite{PartI}, Def. 1.15). We also assume given a
multi-valued convex potential function $K$ on $B$. As currently
understood, one of the basic challenges in the  SYZ picture of
mirror symmetry is to associate a K\"ahler manifold $X$ to  $(B,K)$.
Mirror symmetry can then be accomplished by a Legendre transform:
there is a Legendre transform $(\check B,\check K)$ of the pair
$(B,K)$, and one expects that the K\"ahler manifold $\check X$
associated to $(\check B,\check K)$ is the mirror to $X$. One then
wants to construct a dictionary between things one would like to
compute on $X$, (such as Hodge numbers, periods, and  Gromov-Witten
invariants), and things one can compute on $B$. For example, one
expects that in this picture, Gromov-Witten invariants on $X$ and
periods on $\check X$ should both be expressible in terms of
tropical curves on $B$.

Let us review very quickly the relationship between affine and
complex manifolds. This has already been discussed in many places,
see for example
\cite{Leung},\cite{Announce},\cite{KS2},\cite{Fukaya}. There is a
local system $\Lambda \subset \T_{B_0}$ of integral constant vector
fields with respect to the integral affine structure, and a torus
bundle over $B_0$, $\T_{B_0}/\Lambda$. This manifold carries a
natural complex structure. One might furthermore consider a family
of complex manifolds obtained by rescaling the lattice,
$\T_{B_0}/\epsilon\Lambda$, with $\epsilon\rightarrow 0$
corresponding to the so-called ``large complex structure limit'' of
string theory.  One would like to compactify
$\T_{B_0}/\epsilon\Lambda$ to a complex manifold fibering over $B$.
While this can be done in many cases topologically, it is generally
impossible to do so in the complex category, and one must deform the
complex structure before a compactification can be performed. This
problem was first considered directly by Fukaya \cite{Fukaya}. An
explicit description in two dimensions was found by Kontsevich and
Soibelman by passing from the category of complex manifolds to the
category of rigid analytic spaces. In \cite{Smoothing},  we used the
program started in \cite{PartI} to solve this problem in all
dimensions, again in a slightly different  category. It is worth
noting that using the techniques in \cite{PartI} and known
deformation theory, one can in fact solve the two-dimensional
problem quite easily, but without an explicit description.

We outline the program as begun in \cite{PartI}. The basic idea is to
discretize the problem. We consider \emph{toric polyhedral decompositions}
$\P$ of $B$ (see \cite{PartI}, Definition 1.22). Essentially this is a
decomposition of $B$ into lattice polytopes, but there are some 
delicate conditions involving how these polytopes interact with the
singular locus $\Delta$ of $B$. In particular, no vertex of $\P$ is
contained in $\Delta$. By looking at $B$ and $\P$ in a neighbourhood
of a vertex $v$ of $\P$, one obtains a rational polyhedral fan $\Sigma_v$,
and hence a toric variety $X_v$. These toric varieties can then be
glued together along toric strata
using the combinatorics dictated by $\P$. This gluing can be modified
by equivariant automorphisms of the strata, and so a gluing is
specified by some additional data, which we called
\emph{open gluing data} in \cite{PartI}, Definition 2.25. A choice of
open gluing data then specifies a scheme (or algebraic space)
$X_0(B,\P,s)$, which is a union of toric varieties.

We then want to consider certain sorts of degenerations of Calabi-Yau
varieties of the form $\X\rightarrow S$ over a base scheme $S$,
with some fibre isomorphic to $X_0(B,\P,s)$. The basic idea then is
to pass between $B$ and genuine Calabi-Yau varieties by way of these
degenerations. 

There are two principal problems, in general, with such a degeneration
approach. First,
it is very important to note that this scheme $X_0(B,\P,s)$ does not
contain enough information by itself to carry out mirror symmetry.
There may be many different smoothings of $X_0(B,\P,s)$ to Calabi-Yau
varieties, even with different Hodge numbers. Thus, at first glance,
studying mirror symmetry via degenerations seems far-fetched.

The key innovation of \cite{PartI} is the observation that one should
consider $X_0(B,\P,s)$ as a \emph{log scheme} of Illusie-Fontaine
and Kato (see \cite{PartI}, \S 3.1 and references therein
for an introduction to log schemes as
needed for this program). This extra structure is exactly what is needed
to extract useful information from $X_0(B,\P,s)$, such as moduli,
Hodge numbers, and eventually, periods and Gromov-Witten invariants.

The second principal problem is that not all degenerations of Calabi-Yau
varieties are degenerations to varieties of the type $X_0(B,\P,s)$. 
There are many examples of degenerations which are, however. For example, 
all complete intersections in toric varieties have degenerations of this sort.
More generally, one might conjecture that \emph{all} large complex structure
limit degenerations will be birational to a degeneration of the form
discussed here. In any event, our belief is that the class of Calabi-Yau
varieties for which our degeneration approach will be useful will turn out
to be a far broader class than that of complete intersections in toric 
varieties. We leave this second problem for future work.

Returning to the program at hand, we recall it was shown in 
\cite{PartI} that $(B,\P)$
contains more information than just the scheme $X_0(B,\P,s)$: it
also says precisely how to put a log structure on $X_0(B,\P,s)$,
which we write as $X_0(B,\P,s)^{\dagger}$,
along with a suitable
morphism $X_0(B,\P,s)^{\dagger}\rightarrow\Spec\kk^{\dagger}$
log smooth away from a set $Z\subset X_0(B,\P,s)$ which is codimension
2 and doesn't contain any toric stratum of $X_0(B,\P,s)$. Very roughly
$Z$ corresponds to $\Delta\subset B$.
(Here $\Spec \kk^{\dagger}$ denotes
the standard log point over the field $\kk$, which we always
take to be an algebraically closed field of characteristic zero, 
and the daggers represent log schemes). Once we have done this, 
we obtain so-called log Calabi-Yau spaces which behave very much
like non-singular Calabi-Yau varieties. Conversely, giving the log
Calabi-Yau space
structure on $X_0(B,\P,s)$ gives enough information to reconstruct
$(B,\P)$; the scheme $X_0(B,\P,s)$ by itself does not determine the
affine structure on $B$.

These log Calabi-Yau spaces can be viewed as
playing an intermediate role between the affine manifold $B$ and a
complex manifold or non-singular variety. Without this extra log structure,
there would be no way to extract useful information from $X_0(B,\P,s)$
and much of the information contained in $B$ would be lost.

A large part of \cite{PartI}
was devoted to classifying log Calabi-Yau spaces arising from $(B,\P)$.
With suitable hypotheses, namely that $(B,\P)$ be \emph{simple}
(see \cite{PartI}, Definition 1.60), one finds the set of all log Calabi-Yau
spaces arising from $(B,\P)$ is $H^1(B,
i_*\Lambda\otimes\kk^{\times})$, where $i:B_0\hookrightarrow B$ is the
inclusion. As we shall see in this paper, this moduli space has the
expected dimension given by the first cohomology of the tangent bundle
of a smoothing of $X_0(B,\P,s)^{\dagger}$ (if it exists). 

More importantly, in \cite{Smoothing}, 
we show, with some more general hypotheses on $(B,\P)$,
that a log Calabi-Yau space $X_0(B,\P,s)^{\dagger}\rightarrow\Spec\kk^{\dagger}$
can be put into a formal family ${\frak X}$ over $\kk\lfor t\rfor$. Assuming the
family can be polarized, Grothendieck existence can be applied
and one obtains a scheme $\X$ over $\Spec \kk\lfor t\rfor$ which is a flat
deformation of $X_0(B,\P,s)$. This accomplished one of the original aims
of the program, allowing us to associate a scheme to $B$. Again, the
log structure plays a vital role, providing the ``initial conditions''
to canonically determine the family ${\frak X}$.

Furthermore, in \cite{PartI} we devoloped the discrete Legendre transform.
Given a multi-valued strictly convex integral piecewise linear function
$\varphi$ on $B$, one can construct from the triple $(B,\P,\varphi)$
the Legendre dual $(\check B,\check\P,\check\varphi)$. This then should
yield the mirror family. In \cite{GrossBB}, the first author
checked that this agrees with Batyrev-Borisov duality.

Thus there remains the problem of building a dictionary between geometric
objects
on $B$ and geometric
objects on $X_0(B,\P,s)$ and $\X$.

This paper starts this process. A great deal of the difficulty of this
project is that the necessary theory of log structures has not been
developed sufficiently in the literature for the types of log schemes
that we have to deal with. As a result, both in \cite{PartI} and this
paper, some space has to be devoted to essentially foundational
issues concerning log structures. 

In this paper, the goal is two-fold: first, we wish to
compute Dolbeault cohomology groups
of log Calabi-Yau spaces and their smoothings. This
will enable us to verify the usual exchange of Hodge numbers under mirror
symmetry. Second, we wish to lay the groundwork for calculation of
periods of the family of Calabi-Yaus $\X\rightarrow\Spec k\lfor t\rfor$. This
calculation will be carried out in a future paper, but many of the results
in this paper will be needed. In fact, while these calculations
are not visible, the results of this paper were essential for carrying
out certain period calculations needed to make interesting new
enumerative predictions of mirror symmetry which were stated and verified
in the article \cite{TropVert}. 

To explain the results of this paper in further detail, we give an
outline of the paper. In \S 1, we review the definition of logarithmic
derivations and logarithmic differentials and give local calculations
in a set-up suited for our needs. The only difficulty that arises here
is the presence of the set $Z\subset X_0(B,\P,s)$ where
the log structure breaks down. In logarithmic geometry parlance, 
the log structure fails to be \emph{coherent} at these points. One
consequence is that the sheaf of logarithmic
differentials $\Omega^1_{X_0(B,\P,s)^{\dagger}/\kk^{\dagger}}$ fails to
be coherent at $Z$. As a result, it turns out the proper sheaf to look at
is $j_*\Omega^1_{X_0(B,\P,s)^{\dagger}/\kk^{\dagger}}$, where
$j:X_0(B,\P,s)\setminus Z\hookrightarrow X_0(B,\P,s)$ is the inclusion.
This is similar to the usual definition of Danilov differentials.
\S 1 then is devoted to local calculations of these various sheaves.

\S 2 develops deformation theory for log Calabi-Yau spaces. Results 
currently available do not apply directly, see e.g. \cite{F.Kato}, again
because of the lack of coherence. Unfortunately, we did not find
a clean way of dealing with this problem. In general, there can be
many unpleasant deformations of non-coherent log schemes. We restrict
these deformations by only allowing deformations with a very specific
sort of local model, and we call these \emph{divisorial deformations}.
We then prove that this deformation theory is controlled, as one would
hope, by the sheaf of logarithmic derivations. 

In a sense, \S 2 is a relic of our original hope, when we started this
program, of proving the existence of smoothings of log Calabi-Yau spaces
by proving a version of the Bogomolov-Tian-Todorov theorem, \cite{Tian},
\cite{Todorov}, as had
already been done for normal crossings Calabi-Yau varieties
in \cite{KN}. However, there are technical problems in trying
to prove such a result, and it resisted all our efforts. Instead,
we constructed explicit smoothings in \cite{Smoothing}, and this is much more
valuable anyway as it allows an explicit calculation of periods.
It is important to know that those explicit smoothings
are divisorial deformations, and to do so, we need the deformation theory
results of \S 2. In this way we learn that the results of this
paper, notably the basechange result of \S 4, apply to the smoothings
constructed in \cite{Smoothing}. Other than this point, the results
of \S 2.2 are not needed elsewhere in the paper.

\S 3 is the heart of this paper. We investigate the Dolbeault cohomology groups
\[
H^q(X_0(B,\P,s),j_*\Omega^p_{X_0(B,\P,s)^{\dagger}/\kk^{\dagger}}).
\]
Our main results here follow from a lengthy explicit computation:

\begin{theorem}
Under a hypothesis slightly stronger than simplicity (see Theorem 
\ref{bigtheorem}) we have
\begin{enumerate}
\item
\[
H^q(X_0(B,\P,s),j_*\Omega^p_{X_0(B,\P,s)^{\dagger}/\kk^{\dagger}})
\cong H^q(B,i_*{\bigwedge}^p\check\Lambda\otimes\kk),
\]
where $\check\Lambda$ is the dual local system to $\Lambda$.
\item (Hodge decomposition) The algebraic log de Rham cohomology
groups satisfy the Hodge decomposition, i.e.
\[
\HH^r(X_0(B,\P,s),j_*\Omega^{\bullet}_{X_0(B,\P,s)^{\dagger}/\kk^{\dagger}})
\cong\bigoplus_{p+q=r} H^q(B,i_*{\bigwedge}^p\check\Lambda\otimes
\kk).
\]
\item If in addition the holonomy of $B_0\subseteq B$ is contained
in $\ZZ^n\rtimes \SL_n(\ZZ)$ rather than just $\ZZ^n\rtimes
\GL_n(\ZZ)$, then
\[
H^p(X_0(B,\P,s),\Theta_{X_0(B,\P,s)^{\dagger}/\kk^{\dagger}})
\cong H^p(B,i_*\Lambda\otimes\kk).
\]
By results of \S 2, under the milder assumption of simplicity,
the tangent space of the log deformation functor
is 
\[
H^1(X_0(B,\P,s),\Theta_{X_0(B,\P,s)^{\dagger}/\kk^{\dagger}}),
\]
so this is $H^1(B,i_*\Lambda\otimes\kk)$. This fits with the description
of the moduli space of log Calabi-Yau spaces with dual intersection
complex $B$, which is $H^1(B,i_*\Lambda\otimes\kk^{\times})$.
\end{enumerate}
\end{theorem}

The extra hypothesis, over and above simplicity, is important. Essentially,
this hypothesis
says that the mirror variety is non-singular, rather than an orbifold.
As is well-known \cite{Bat},
to consider mirror pairs in dimension $\ge 4$, one
needs to include orbifold singularities, as not all Gorenstein abelian quotient
singularities in dimension $\ge 4$ have crepant resolutions. 
In the orbifold
context, one needs to consider \emph{stringy} Hodge numbers \cite{BatStringy}.
Calculation of stringy Hodge numbers would take us too far afield in
this paper, and hence the extra hypothesis. However, stringy Hodge numbers
are not necessary in dimension three, so the results of this paper are
complete in this case. For exploration of the relationship between
the calculations in this paper and stringy Hodge numbers in dimension
4, see the paper \cite{Ruddat} of Ruddat.

This calculation of the log Dolbeault groups demonstrates that
the mirror duality proposed in \cite{PartI} in fact interchanges Hodge
numbers. Indeed, consider a Legendre dual pair $(B,\P,\varphi)$
and $(\check B,\check\P,\check\varphi)$. Then $\Lambda^B\cong
\check\Lambda^{\check B}$ (see \cite{PartI}, Proposition~1.50,(1)), 
where the superscripts denote which affine
structure is being used to define the sheaf. Moreover, if 
$\dim B=n$ and the holonomy of $B$ is contained in $\ZZ^n\rtimes
\SL_n(\ZZ)$, then 
$\bigwedge^q\check\Lambda^B\cong
\bigwedge^{n-q}\Lambda^{B}$ and hence
$\bigwedge^q\check\Lambda^B\cong
\bigwedge^{n-q}\check\Lambda^{\check B}$. Thus the isomorphism
$H^p(B,i_*\bigwedge^q\check\Lambda^B)\cong H^p(\check B,i_*\bigwedge^{n-q}
\check\Lambda^{\check B})$ gives the usual exchange of Hodge numbers
on the level of log Calabi-Yau spaces.

To relate these to the usual Hodge numbers of a smoothing, in
\S 4 we prove a base-change theorem, which tells us that the log de Rham
groups of $X_0(B,\P,s)$ coincide with the ordinary algebraic de Rham
groups of a smoothing. Again, with the extra hypothesis, the same
holds for the Dolbeault groups. This demonstrates we have defined the
log Dolbeault groups correctly, and also demonstrates that when one
has non-singular Calabi-Yau varieties on both sides of the picture,
the Hodge numbers are exchanged, as expected. This implies, for example,
by the results of \cite{GrossBB}, the usual interchange of Hodge numbers
in the context of Batyrev-Borisov duality \cite{BB}. In 
particular, this gives a new way of computing Hodge numbers in this case;
it is not at all clear from a combinatorial viewpoint why these computations
should give the same answer, but of course they must.

\S 5 calculates the Gauss-Manin connection in two different contexts.
In \S 5.1, we consider the following situation. If 
\[
\X\rightarrow\Spec\kk\lfor t\rfor
\]
is a one-parameter deformation of 
\[
X_0(B,\P,s)^{\dagger}\rightarrow
\Spec\kk^{\dagger},
\]
then one would ideally like to calculate the
Gauss-Manin connection on this family. Once this is done, one can
write a family of holomorphic $n$-forms in terms of flat sections
of the Gauss-Manin connection; this describes the periods of
the holomorphic $n$-forms which yields the $B$-model predictions
of mirror symmetry. This must wait for further
work, but one can easily
calculate the monodromy of the system of flat sections of the Gauss-Manin
connection as the exponential of the residue of the connection. 
If $\nabla$ is the Gauss-Manin connection, this residue is
\begin{align*}
\Res_0(\nabla)=\nabla_{t\partial/\partial t}|_{t=0}:
\HH^r(X_0(B,\P,s),& j_*\Omega^{\bullet}_{X_0(B,\P,s)^{\dagger}/\kk^{\dagger}})
\\
&\rightarrow
\HH^r(X_0(B,\P,s),j_*\Omega^{\bullet}_{X_0(B,\P,s)^{\dagger}/\kk^{\dagger}}).
\end{align*}

In \S \ref{monodromysection}, this residue is calculated explicitly. 
It is given on the level of the summands of the
Hodge decomposition by, for $r=p+q$, the map
\begin{eqnarray*}
H^q(B,i_*{\bigwedge}^p\check\Lambda)&\rightarrow&H^{q+1}(B,
i_*{\bigwedge}^{p-1}\check\Lambda)\\
\alpha&\mapsto&c_B\cup\alpha
\end{eqnarray*}
where $c_B\in H^1(B,i_*\Lambda)$ is the \emph{radiance obstruction}
of $B$, an invariant of affine manifolds introduced in \cite{GH}.
This description of monodromy coincides with that given in \cite{SlagI}
for the SYZ picture.

In \S 5.2, we consider a universal family of log Calabi-Yau spaces
over the moduli space $S=H^1(B,i_*\Lambda\otimes\kk^{\times})$, and
consider the Gauss-Manin connection of this family. In this way,
we obtain a variation of mixed Hodge structures over $S$. This material
is largely included for future applications to $B$-model calculations
of enumerative predictions.

\emph{Acknowledgements:} We would like to thank Arthur Ogus and Helge 
Ruddat for useful conversations, as well as the referee.

\section*{Notational summary}

By the very nature of the subject this paper involves a considerable
amount of notation, which to a large part has been introduced in
\cite{PartI}.
Although we  must assume the reader has some familiarity with \cite{PartI}, 
we will give references to notation in \cite{PartI} whenever it is first
used in this paper. In addition, here we will survey the
basic ideas and notation of \cite{PartI}. 
While it would take up too much space to be
completely self-contained, this section attempts to
make this paper more accessible.

The setup for the bulk of this paper is as follows. We assume given
an $n$-dimensional toric log Calabi-Yau space $X_0^\ls:=
(X_0,\M_{X_0})$ over an algebraically closed field $\kk$ of
characteristic $0$, as defined in \cite{PartI}, Definition~4.3. One
should think of $X_0$ as the central fibre of what we called a
\emph{toric degeneration of Calabi-Yau varieties} (\cite{PartI},
Definition~4.1). This is a Calabi-Yau variety $\X$ over a discrete
valuation ring with closed fibre a union of toric varieties that
mutually intersect in toric strata, and $\M_{X_0}$ is then the log
structure defined by the embedding $X_0\subseteq \X$. One of the
central results of \cite{PartI} (Theorem~5.4) shows that, under
suitable hypotheses, toric log
Calabi-Yau spaces are equivalent to discrete data $(B,\P)$ and a
cohomology class $s\in H^1(B,i_*\Lambda\otimes\GG_m)$ called
lifted gluing data. Here $B$ is
an integral affine manifold of real dimension $n$, with
singularities along a real codimension two subset $\Delta$, together
with a polyhedral decomposition $\P$ consisting of integral lattice
polytopes. Moreover, $i:B\setminus\Delta\to B$ is the inclusion,
$\Lambda$ is the sheaf of integral tangent vectors on
$B\setminus\Delta$, and $\check\Lambda$ is the dual local system. 

As a rule, general elements of $\P$ are denoted
$\tau$, while $\omega$ usually denotes edges, $\rho$ cells of
codimension one and $\sigma$ maximal cells. Vertices of $\P$ are
written $v$ or $w$. In \cite{PartI}, we allowed elements of $\P$
to have self-intersections: For example, if $B=\RR/\ZZ$, we can
take $\P=\{\{0\},B\}$, i.e., we view $B$ as a line segment
with opposite ends identified. Because of this decision, in order
to be clear about how elements of $\P$ attach, we viewed
$\P$ as a category, with morphisms $\omega\rightarrow\tau$
indicating a choice of inclusion of a face $\omega$ in the 
``normalization'' of $\tau$. To simplify notation, we shall assume
here that cells do not self-intersect, the general case being
straightforward. Most of the time, this will not make a difference,
and we will write $\omega\rightarrow\tau$ and $\omega\subseteq\tau$
interchangeably, writing the latter when it makes notation cleaner.

The integral affine structure on the underlying
topological space $B$ is defined by the lattice polytopes together
with the structure of a fan $\Sigma_v$ at each vertex $v$. In
particular, $\Delta$ is disjoint from vertices and the interiors of
maximal cells of $\P$. The toric variety $X_v$ associated to
$\Sigma_v$ is one of the irreducible components of $X_0$. A higher
dimensional cell $\tau\in\P$ labels a lower dimensional toric
stratum $X_\tau= \bigcap_{v\in\tau} X_v \subseteq X_0$, and $q_\tau:
X_\tau\to X_0$ denotes the inclusion.
Then $X_0$ is the categorical
limit of the strata with attaching maps
\[
F_{\tau_1,\tau_2}: X_{\tau_2}\lra X_{\tau_1},
\]
for any two cells $\tau_1\subseteq \tau_2$. 

The meaning of the description of the cells of $\P$ as lattice
polytopes is that they
define the discrete data of the log structure, given by
$\overline\M_{X_0} =\M_{X_0}/\O_{X_0}^\times$, along with a section
$\rho$ for the morphism to the standard log point. The cohomology
class $s$ takes care of the moduli of gluings of the irreducible
components and the choice of log structure with given $\overline
\M_{X_0}$. In particular, the attaching maps $F_{\tau_1,\tau_2}$
depend on $s$, and hence were denoted $F_{S,s}(\tau_1\to \tau_2)$ in
\cite{PartI}. In this paper we suppress $s$ and the base scheme
$S$, which is $\Spec\kk$ except in \S5, to simplify the notation.

Local monodromy of the affine structure around $\Delta$ leads to
singularities of the log structure on $X_0$, and this correspondence
is quite important throughout. The local monodromy of $\Lambda$
is completely
determined by pairs $\omega\subseteq\rho$ where $\omega$ is an edge
of the codimension one face $\rho$. This data determines a closed loop
passing from one vertex $v^+$ of $\omega$ into one of the two maximal
cells separated by $\rho$ and back to $v^+$ via $v^-$ and the other
maximal cell bounding $\rho$. Parallel transport along this loop
leads to monodromy in $\Lambda$ of the form
\[
m\lra m+ \kappa_{\omega\rho} \langle \check d_\rho,m\rangle
d_\omega,\quad m\in\Lambda_{v^+},
\]
with $\kappa_{\omega\rho}\in\ZZ$, $d_\omega$ a primitive tangent
vector of $\omega$ and $\check d_\rho\in\check \Lambda_{v^+}$ a
primitive linear form vanishing on tangent vectors of $\rho$, see
\S1.5 in \cite{PartI}, where the symbol $n_{\omega\rightarrow\rho}
=\kappa_{\omega\rho}\check d_{\rho}$
was used. Signs can be fixed in such a way that cases
arising from actual degenerations fulfill $\kappa_{\omega\rho}\ge
0$. These are called \emph{positive} (\cite{PartI},
Definition~1.54), and positivity of $(B,\P)$ is an assumption
throughout. More generally, if $\omega\in\P$ is an edge and
$\sigma^{\pm}$ are two maximal cells containing $\omega$, monodromy around 
a loop
starting at a vertex $v^+$ of $\omega$ into $\sigma^+$ to $v^-$
and then via $\sigma^-$ back to $v^+$ can be written as
\[
m\lra m+\langle n^{\sigma^+\sigma^-}_{\omega},m\rangle d_{\omega},
\quad m\in\Lambda_{v^+},
\]
for some $n^{\sigma^+\sigma^-}_{\omega}\in \check\Lambda_{v^+}$.
Dually, if $\rho\in\P$ is a codimension one cell contained in
two maximal cells $\sigma^{\pm}$, and $v^{\pm}$ are two vertices
of $\rho$, then monodromy around a loop starting at $v^+$
into $\sigma^+$ to $v^-$ and then via $\sigma^-$ back to $v^+$
can be written as
\[
m\lra m+\langle \check d_{\rho},m\rangle m^{\rho}_{v^+v^-},
\quad m\in\Lambda_{v^+},
\]
for some $m^{\rho}_{v^+v^-}\in\Lambda_{v^+}$.

We also assume a kind of indecomposability of the local
monodromy, formalized in the concept of \emph{simplicity}
(see \cite{PartI}, Definition~1.60 for the full definition). 
A necessary condition is
$\kappa_{\omega\rho}\in\{0,1\}$ for all $\omega, \rho$. Under this
latter condition, local monodromy is determined completely
combinatorially by the pairs $(\omega,\rho)$ with
$\kappa_{\omega,\rho}=1$. In any case, for each $\tau\in\P$
simplicity allows one to capture the local monodromy information around
some $\tau\in\P$ in terms of two collections of
subsets
\[
\Omega_1,\ldots,\Omega_p\subseteq \big\{\omega\subseteq \tau
\,\big|\, \dim\omega=1\big\},\quad
R_1,\ldots,R_p\subseteq \big\{\rho\supseteq \tau\,\big|\,
\dim\rho=n-1 \big\},
\]
of edges and facets with the same behaviour with respect to local
monodromy
satisfying several properties, including $\kappa_{\omega\rho}=1$
if and only if there exists an $i$ such that $\omega\in\Omega_i$,
$\rho\in R_i$.
There are then corresponding monodromy polytopes, well-defined
after fixing a vertex $v\in\tau$ and a maximal cell $\sigma$
containing $\tau$ and any $\omega_i\in\Omega_i$, $\rho_i\in R_i$:
\begin{eqnarray*}
\Delta_i&=&\Conv\big\{ m^{\rho_i}_{vv'}\,\big|\,
v'\in\tau\big\}\subseteq \Lambda_{\tau,\RR},\\
\check\Delta_i&=&\Conv\big\{ n^{\sigma\sigma'}_{\omega_i}\,\big|\,
\tau\subseteq\sigma'\in\P_{\max}\big\}\subseteq \check \Lambda_{\tau,\RR}.
\end{eqnarray*}
In a sense, $\Delta_i$ captures the $i$-th part of ``inner
monodromy'' of $\tau$, fixing $\rho_i$ containing $\tau$ and letting
$\omega\subseteq \tau$ vary, while $\check \Delta_i$ captures the
corresponding outer part. The partitioning into $p$ parts is
motivated by the study of the complete intersection case, see
\cite{GrossBB}. Simplicity is then defined by requiring the convex
hulls of $\bigcup_{i=1}^p \Delta_i\times\{e_i\}$ and of
$\bigcup_{i=1}^p \check\Delta_i\times\{e_i\}$ to be elementary
simplices.

The monodromy polytopes are related to singularities of the
log structure as follows. The charts for the log structure from
\cite{PartI} use the affine cover of $X_0$ given by the sets
\[
V(\sigma)= X_0\setminus \bigcup_{\tau\cap\sigma=\emptyset} X_\tau,
\]
for $\sigma\in\P$ maximal cells. Now $V(\sigma)$ is canonically the
boundary divisor of the affine toric variety $\Spec\kk[P]$, where
$P$ is the monoid of integral points of $C(\sigma)^\vee$, the dual
of the cone generated by $\sigma\times\{1\}$ in
$\Lambda_{\sigma,\RR}\oplus\RR$, see Construction~2.15 in \cite{PartI}.
This embedding suggests a log structure on $V(\sigma)$ over the
standard log point $\Spec \kk^\ls$, but these log structures are not
locally isomorphic unless the local monodromy is trivial. The class
$s\in H^1(B, i_*\Lambda\otimes_\ZZ \GG_m)$ provides the necessary
local changes to the standard log structures on $V(\sigma)$ to allow
a consistent gluing. Explicitly, if $\omega\subseteq\sigma$ is an edge
of integral length $e$ and endpoints $v^\pm$, then a chart for the
standard log structure in a neighbourhood of
the big cell of the codimension one
stratum $X_\omega=X_{v^+}\cap X_{v^-}$ (an irreducible component of
$(X_0)_\sing$) is given by an equation $xy=t^e$. Here $t$ is the
deformation parameter and $x,y$ monomials vanishing along
$X_\omega\cap V(\sigma)$. The log structure modified by $s$
reads
\[
xy=f(w)t^e
\]
for some function $f$ depending only on monomials $w$ not vanishing identically
along $X_\omega$. Along the zero locus of $f$ the induced log
structure is singular (not \emph{fine}). The global meaning of $f$
is as a section of a coherent sheaf $\shLS^+_{\pre,X_0}$ on
$(X_0)_\sing$. The zero locus of this section is the
$(n-2)$-dimensional singular locus $Z\subseteq (X_0)_\sing$ of the log
structure of $X_0^\ls$, see Theorems~3.22 and 3.27 in \cite{PartI}.

The meaning of the monodromy polytopes with respect to $Z$ is that
$p$ is the number of irreducible components $Z_1,\ldots,Z_p$ of
$Z\cap X_\tau$, the log structure locally along $Z_i$
is determined by $\Delta_i$, and $\check \Delta_i$ is the Newton
polytope of $Z_i$ on the big cell of $X_\tau$.
Moreover, if the convex hull of $\bigcup_i \Delta_i\times
\{e_i\}$ is an elementary simplex there is still a toric model
defining a chart for the log structure locally along $Z$. This is
shown in \S2.1 of the present paper. We end this outline of the
results of \cite{PartI} with a list of relevant notation for easy
reference. 

\bigskip

{\small
\begin{tabbing}
\centerline{\bf Relevant standard notation from \cite{PartI}}\\[1ex]
$M$, $M_\RR$\hspace{20ex}
	\=free abelian group of rank $n$, $M_\RR=M\otimes_\ZZ\RR$\\
$N$, $N_\RR$
	\>the dual groups $N=\Hom(M,\ZZ)$, $N_\RR=N\otimes_\ZZ\RR$\\
$B$
	\>$n$-dimensional integral affine manifold with singularities\\
$\Delta\subseteq B$, $i:B\setminus\Delta \to B$
	\>discriminant locus and the inclusion of its complement\\
$\Lambda_\RR, \Lambda$
	\>local system of flat (integral) tangent vector fields\\
$\check\Lambda_\RR, \check\Lambda$
	\>duals to $\Lambda_\RR$ and $\Lambda$\\
$\shAff(B,\RR)$, $\shAff(B,\ZZ)$
	\>sheaf of continuous (integral) affine functions
	on $B\setminus\Delta$\\
$\P$
	\>polyhedral decomposition of $B$\\
$v$, $w$
	\>vertices of $\P$\\
$\rho$, $\sigma$, $\tau$, $\omega$
	\>cells of $\P$\\
$\Int\tau$
	\>relative interior of $\tau\in\P$\\
$\Lambda_{\tau,\RR}$, $\Lambda_\tau$
	\>space of (integral) tangent vector fields on $\tau\in\P$\\
$\Sigma_v$
	\>fan in $\Lambda_{\RR,v}$ induced from $\P$\\
	\>for $\tau\in\P$ generalizes to the fan $\Sigma_\tau$ in
	$\Lambda_{\tau,\RR}$\\
$\check \Sigma_\tau$
	\>normal fan of $\tau$ in $\check\Lambda_{\tau,\RR}$\\
$d_\omega$
	\>for $\dim\omega=1$, generator of $\Lambda_\omega$\\
$\check d_\rho$
	\>for $\dim\rho=n-1$, generator of
	$\Lambda_\rho^\perp\subseteq \check\Lambda_v$\\
$W_{\tau_1\to\tau_2}$
	\>open star in the barycentric subdivision of $\P$ of the interior\\
	\>of the edge connecting the barycenters of $\tau_1$ and $\tau_2$\\
$\kappa_{\omega\rho}$
	\>integer determining monodromy around loop given by
	$\omega\subseteq\rho$\\
$\Omega_i$, $\Delta_i$
	\>for given $\tau\in\P$, set of edges and the associated monodromy
	polytope\\
$R_i$, $\check \Delta_i$
	\>dual data to $\Omega_i$, $\Delta_i$\\
$C(\tau)$
	\>cone over $\tau\times\{1\}$ in $\Lambda_{\tau,\RR}\oplus\RR$
          for $\tau\in\P$\\
$s$
	\>element of $H^1(B,i_*\Lambda\otimes\GG_m)$; represents moduli
	of gluing\\
	\>of the strata and of the log structure\\
$X_0^\ls= (X_0,\M_{X_0})$
	\>the toric log Calabi-Yau space associated to $(B,\P)$\\
	\>and $s\in H^1(B, i_*\Lambda\otimes\GG_m)$.\\
$Z\subset (X_0)_\sing$
	\>singular locus of the log structure; $\dim Z= n-2$\\
$\overline\M_{X_0}$
	\>$\M_{X_0}/\O_{X_0}^\times$, the discrete part of the
	log structure\\
$\Spec\kk ^\ls$
	\>the standard log point $(\Spec \kk,\kk^\times\times \NN)$\\
$q_\tau: X_\tau\to X_0$
	\>inclusion of the toric stratum of $X_0$
	isomorphic to $X(\Sigma_\tau)$\\
$F_{\tau_1,\tau_2}$
	\>for $\tau_1\subseteq \tau_2$, the attaching map
	$X_{\tau_2}\to X_{\tau_1}$\\
$V(\tau)$
	\>canonical open affine neighbourhood in $X_0$ of the big cell of
	$X_\tau$\\
$V_{\tau_1\rightarrow\tau_2}$
        \>toric stratum of $V(\tau_2)$ corresponding to $\tau_1$, equal to
          $X_{\tau_1}\cap V(\tau_2)$\\
$\rho$
	\>section of $\overline\M_{X_0}$ defining the morphism to
	$\Spec\kk^\ls$ locally.
\end{tabbing}
}

\section{Derivations and differentials}

Let $\pi:X^\ls=(X,\M_X)\to S^\ls= (S,\M_S)$ be a morphism of
logarithmic spaces. Here $\M_X$ is a sheaf of monoids on $X$ and
the dagger is always used to denote logarithmic spaces.
See \cite{PartI}, \S 3.1 for an introduction
to log schemes as needed here, and for further references.

\begin{definition}\label{log derivation}
A \emph{log derivation} on $X^\ls$ over $S^\ls$ with values in an
$\O_X$-modules $\shE$ is a pair $(\Di,\Dlog)$, where $\Di: \O_X\to
\shE$ is an ordinary derivation of $X/S$ and $\Dlog: \M^\gp_X\to
\shE$ is a homomorphism of abelian sheaves with
$\Dlog\circ\pi^\#=0$; these fulfill the following compatibility
condition
\begin{eqnarray}\label{compatibility}
\Di\big(\alpha_X(m)\big)=\alpha_X(m)\cdot \Dlog(m),
\end{eqnarray}
for all $m\in\M_X$, where $\alpha_X:\M_X\rightarrow\O_X$ is the
log structure.

We denote by $\Theta_{X^{\dagger}/S^{\dagger}}$ the sheaf of log derivations 
of $X^{\dagger}$ over $S^{\dagger}$ with values in $\O_X$. \qed
\end{definition}

\begin{remark} 
\label{Thetaextend}
Suppose $X$ has no embedded components, $\M_X$ has no sections
with support contained in a codimension two subset, and the $\O_X$-module
$\shE$ is $S_2$, so that sections extend across codimension two subsets.
Then if $Z\subseteq X$ is a codimension $\ge 2$ closed subset of $X$,
the module of log derivations on $X^{\dagger}$ with values in $\shE$ is the same
as the module of log derivations on $X^{\dagger}
\setminus Z$ with values in $\shE$.
We will use this fact freely in what follows.
\end{remark}

In many cases a log derivation $(\Di,\Dlog)$ is already determined by
$\Di$.

\begin{proposition}
\label{Thetainclude}
Assume that $\M_X=\pi^*(\M_S)$ holds on an open, dense subset
$U\subset X$, and that $\shE$ has no sections with support in $X
\setminus U$. Then the forgetful map
\[
(\Di,\Dlog) \longmapsto \Di
\]
from the sheaf of log derivations on $X^\ls/S^\ls$ with values in
$\shE$ to the sheaf of usual derivations on $X/S$ with values in
$\shE$ is injective.
\end{proposition}
\proof
On $U$ each $m\in \M_X$ may be written as $h\cdot \pi^\#(n)$
for $h\in\O_X^\times$ and $n\in \M_S$. Hence $\Dlog(m)$ is
determined by $\Di$ via Equation~(\ref{compatibility}). Thus if $\Di=0$ then
$\Dlog|_U=0$, which under the assumption on $\shE$ implies
$\Dlog=0$. 
\qed
\medskip

We may thus often think of log derivations as usual derivations with
certain vanishing behaviour determined by the log structure:

\begin{example} 
\label{preserveideal}
Let $Y$ be a normal integral scheme over $\kk$,
and $X\subseteq Y$ a reduced Weil divisor. Endow $Y$ with the divisorial log
structure $\M_{(Y,X)}:=j_*\O_{Y\setminus X}^{\times}\cap\O_Y$,
where $j:Y\setminus X\hookrightarrow Y$ is the inclusion. 
Then $\Theta_{Y^{\dagger}/\kk}$ consists of
the usual derivations of $Y$ which preserve the ideal of $X$. Indeed,
if $\Di$ is a log derivation and $f\in \shI_{X/Y}$, then at the generic
point $\eta$ of an irreducible component of $X$, we can write
$f=f'\cdot t^p$ for $t$ a generator of $\shI_{X/Y}$ at $\eta$, $p>0$,
and $f'$ a regular function. Then $t$ defines an element
of $\M_{(Y,X)}$ in a neighbourhood of $\eta$, so 
$\Di f=t^p \Di f'+pf' t^{p-1}\Di t=t^p(\Di f'+pf'\Dlog t)$
is in $\shI_{X/Y}$ in a neighbourhood of $\eta$. Thus
$\Di f$ vanishes along every component of $Y$, so is in $\shI_{X/Y}$.

Conversely if $\Di$ is an ordinary derivation preserving $\shI_{X/Y}$,
then for $f\in \M_{(Y,X)}$, we define $\Dlog f$ as ${\Di f\over f}$;
that this is a regular function is immediately checked again as
above at the generic points of $\eta$. \qed
\end{example}

\medskip

The sheaf of log derivations is well-understood for 
log-smooth morphisms of fine log schemes.
However, the type of log schemes that arise in \cite{PartI}
are not in general fine: they fail to be \emph{coherent} in the sense
of Ogus \cite{ogus}, i.e. there may not be local charts of the form
$P\rightarrow \M_X$ with $P$ a finitely generated monoid. However, the
examples we wish to consider are \emph{relatively coherent} \cite{ogus}. We will
not use the full formalism of \cite{ogus} 
here as all the cases we need fall into
a narrow class of examples, which we shall now introduce.

Let $P$ be a toric monoid, $F\subseteq P$ a face. Set $Y=\Spec \kk[P]$,
and let $j:U\hookrightarrow Y$ 
be the largest open subset where $z^p$ is invertible for all $p\in F$. 
Let $X:=Y\setminus U$.
There are two natural log structures on $Y$. 
The first is given by $\M_Y=\M_{(Y,X)}=
j_*(\O_U^{\times})\cap\O_Y$. 
The other one is induced by the chart $P\rightarrow
\kk[P]$, which is a fine log structure, which we write as
$\M_{\tilde Y}$. Note this can also be described as $\M_{\tilde Y}=
\tilde j_*(\O_{\tilde U}^{\times})\cap\O_Y$, where $\tilde j:
\tilde U\hookrightarrow Y$ is the inclusion of the big torus orbit
$\tilde U$ of $Y$.
There is an obvious inclusion $\M_{Y}
\subseteq\M_{\tilde Y}$, which is relatively coherent in the
language of \cite{ogus}. Such log structures still have good properties.
We write $Y^{\dagger}$ and $\tilde Y^{\dagger}$ for the two log structures
respectively. See Example~\ref{ODPexample} for a standard example.

Now let $M=\ZZ^n$, $M_{\RR}=M\otimes_{\ZZ}\RR$,
$N=\Hom_{\ZZ}(M,\ZZ)$, $N_{\RR}=N\otimes_{\ZZ}\RR$, and
suppose $P$ is given by $\sigma^{\vee}\cap N$ for a strictly
convex rational polyhedral cone $\sigma$ in $M_{\RR}$. 
Write 
\[
X=Y\setminus U=\bigcup_{i=1}^s X_i,
\]
where the $X_i$'s are the toric divisors of $Y$ contained in $X$.
Let \[
D=\bigcup_{j=1}^t D_j
\]
be the union of toric divisors of $Y$ not 
contained in $X$.
We take primitive generators of extremal rays of $\sigma$
to be $v_1,\ldots,v_{s+t}$,
with $v_1,\ldots,v_s$ corresponding to $X_1,\ldots,X_s$, and $v_{s+1},\ldots,
v_{s+t}$ corresponding to $D_1,\ldots,D_t$. For ease of notation, we sometimes
write 
\[
\hbox{$w_j=v_{s+j}$ for $1\le j\le t$.}
\]
Let $P_1,\ldots,P_s$, $Q_1,\ldots,Q_t$
be the facets (maximal proper faces)
of $P$ corresponding to $v_1,\ldots,v_s$, $w_1,\ldots,w_t$ 
respectively. Note that $X_i=\Spec \kk[P_i]$, the $P_i$'s are the
facets of $P$ not containing $F$, and the face $F\subseteq P$ is given by
\[
F=\langle w_1,\ldots,w_t\rangle^{\perp}\cap P.
\]

\begin{proposition}
\label{Thetatoric}
In the above situation,
$\Gamma(Y,\Theta_{Y^{\dagger}/\kk})$ splits into
$P^{\gp}$-homogeneous pieces
\[
\bigoplus_{p\in  P^{\gp}} z^p(\Theta_{Y^{\dagger}/\kk})_p,
\]
where
\[
(\Theta_{Y^{\dagger}/\kk})_p
=
\begin{cases}
M\otimes_{\ZZ} \kk&\hbox{if $p\in P$,}\\
\ZZ v_i\otimes_{\ZZ} \kk&\hbox{if there exists an $i$, $s+1\le i\le s+t$,}\\
&\hbox{with
$\langle v_i,p\rangle=-1$,
$\langle v_j,p\rangle\ge 0$ for $j\not=i$,}\\
0&\hbox{otherwise}.
\end{cases}
\]
We write an element $m\in(\Theta_{Y^{\dagger}/\kk})_p$
as $\partial_m$, and $z^p\partial_m$ acts on the monomial
$z^q$ by
\[
z^p\partial_m z^q=\langle m,q\rangle z^{p+q}.
\]
\end{proposition}

\begin{remark} Note that in the case that
$F=P$, we have $\Theta_{Y^{\dagger}/\kk}=M\otimes_\ZZ \O_Y$. This is
the standard case where $X$ is the toric boundary of $Y$, in which case
$Y^{\dagger}$ is log smooth over $\kk$. Otherwise, $\Theta_{Y^{\dagger}/\kk}$
may not be locally free.
\end{remark}

\noindent
\emph{Proof of Proposition \ref{Thetatoric}.}
We first identify ordinary derivations on $Y$ (which is the special case
$F=\{0\}$). As any derivation
on $Y$ restricts to a derivation on the big torus orbit of $Y$,
\[
\Gamma(Y,\Theta_{Y/\kk})\subseteq \bigoplus_{p\in P^{\gp}} z^p(M\otimes_{\ZZ}
\kk),
\]
where $m\in M\otimes_{\ZZ} \kk$ corresponds to the derivation $\partial_m$.
Furthermore, the torus $\Spec \kk[P^{\gp}]$ acts on $Y$,
$\Theta_{Y/\kk}$, and $\Theta_{Y^{\dagger}/\kk}$, so $\Gamma(Y,\Theta_{Y/\kk})$
and $\Gamma(Y,\Theta_{Y^{\dagger}/\kk})$ both decompose into $P^{\gp}$-homogeneous
pieces. Thus we need to determine for each $p\in P^{\gp}$ for which
$m\in M\otimes \kk$ is $z^p\partial_m$ a derivation on $Y$, i.e.
when is $z^p\partial_mz^q\in \kk[P]$ for all $q\in P$.
But $\langle m,q\rangle z^{p+q}$ is regular if and only if $p+q\in P$
or $\langle m,q\rangle=0$.

So first suppose $z^p\partial_m$ is an ordinary derivation on $Y$. There are
two cases. \emph{Case 1.} $m^{\perp}$ does not contain a
facet of $P$. In this case, for each $i$ we can find a $q\in P$
such that $\langle v_i,q\rangle=0$, $\langle m,q\rangle\not=0$. Then
$p+q\in P$, so $0\le\langle v_i,p+q\rangle=\langle v_i,p\rangle$.
Thus $p\in P$. \emph{Case 2}. $m$ is proportional to $v_i$ for some
$i$. Then for each $j\not=i$, the same argument shows $\langle v_j,p\rangle
\ge 0$, while if $q\in P$ with $\langle v_i,q\rangle=1$, then 
$0\le \langle v_i, p+q\rangle=\langle v_i,p\rangle +1$.

Reversing this argument, we see that
\[
\Gamma(Y,\Theta_{Y/\kk})=\bigoplus_{p\in P^{\gp}} z^p(\Theta_{Y/\kk})_p
\]
where
\[
(\Theta_{Y/\kk})_p=\begin{cases}
M\otimes_{\ZZ} \kk&\hbox{if $p\in P$,}\\
\ZZ v_i\otimes_{\ZZ} \kk&\hbox{if there exists an $i$ with
$\langle v_i,p\rangle=-1$,
$\langle v_j,p\rangle\ge 0$ for $j\not=i$,}\\
0&\hbox{otherwise}.
\end{cases}
\]

By Proposition~\ref{Thetainclude},
$(\Theta_{Y^{\dagger}/\kk})_p\subseteq (\Theta_{Y/\kk})_p$.
So now let's consider log derivations.
The ideal of $X$ (with the reduced induced scheme structure)
is generated by $P\setminus (P_1\cup\cdots \cup P_s)$.
Let $z^p\partial_m\in\Gamma(Y,\Theta_{Y/\kk})$.
Then $z^p\partial_m$ certainly preserves the ideal
if $p\in P$. On the other hand, if $p\not\in P$, then we can take
$m=v_i$ for some $i$ with  $\langle v_i,p\rangle=-1$. 
For any given $i$, we can find a $q\in P\setminus (P_1\cup\cdots\cup P_s)$
with $\langle v_i,q\rangle=1$, $\langle v_j,q\rangle \ge 0$
for all $j\not=i$, and then $z^p\partial_{v_i} z^q=\langle v_i,q\rangle
z^{p+q}$, but $\langle v_i,p+q\rangle=0$. Thus if $1\le i\le s$, we would
have $p+q\in P\cap v_i^{\perp}=P_i$. So in this case $z^p\partial_{v_i}$
preserves the ideal of $X$ if and only if $i\not\in \{1,\ldots,s\}$, i.e.
$s+1\le i\le s+t$. This gives the desired result. \qed

\begin{corollary} 
\label{Thetatoriccor}
In the situation of Proposition \ref{Thetatoric}, let $S=\Spec \kk[\NN]$
with the log structure defined by the obvious chart $\NN\rightarrow
\kk[\NN]$, and let $\rho\in P$ lie in the
interior of the face $F$, so that $\Spec\kk[P]/(z^{\rho})$ yields a scheme
with reduction $X$. Then  $z^{\rho}$ induces
a log morphism $Y^{\dagger}\rightarrow S^{\dagger}$, and
\[
\Gamma(Y,\Theta_{Y^{\dagger}/S^{\dagger}})=
\bigoplus_{p\in  P^{\gp}} z^p(\Theta_{Y^{\dagger}/S^{\dagger}})_p,
\]
where
\[
(\Theta_{Y^{\dagger}/S^{\dagger}})_p
=
\begin{cases}
\rho^{\perp}\otimes_{\ZZ} \kk&\hbox{if $p\in P$,}\\
\ZZ v_i\otimes_{\ZZ} \kk&\hbox{if there exists an $i$, $s+1\le i\le s+t$,}\\
&\hbox{with
$\langle v_i,p\rangle=-1$,
$\langle v_j,p\rangle\ge 0$ for $j\not=i$,}\\
0&\hbox{otherwise}.
\end{cases}
\]
\end{corollary}

\proof This follows by imposing the condition that an element of
$\Theta_{Y^{\dagger}/S^{\dagger}}$ must annihilate $z^{\rho}$.
\qed

\medskip

The following result shows we obtain the same description for
log derivations on suitable thickenings of $X$; this is essentially
a base-change result for derivations, but does not follow immediately
from generalities.

\begin{proposition} 
\label{ThetaXl}
Under the same hypotheses as Corollary \ref{Thetatoriccor},
assume in addition that $Y$ is Gorenstein and that $X=\Spec\kk[P]/(z^{\rho})$ 
is reduced.
Let $\X_k=\Spec \kk[P]/(z^{(k+1)\rho})$, with the induced log structure
from $Y^{\dagger}$. Then $\Gamma(\X_k,\Theta_{\X_k^{\dagger}/\kk})$
splits into $P^{\gp}$-homogeneous pieces
\[
\bigoplus_{p\in P^{\gp}} z^p\left(\Theta_{\X_k^{\dagger}/\kk}\right)_p,
\]
where
$\left(\Theta_{\X_k^{\dagger}/\kk}\right)_p=0$ if there does not exist an $i$,
$1\le i\le s$, such that $0\le \langle v_i,p\rangle\le k$;
otherwise
\[
\left(\Theta_{\X_k^{\dagger}/\kk}\right)_p=\left(\Theta_{Y^{\dagger}/\kk}\right)_p.
\]

In addition, let $A_k=\kk[t]/(t^{k+1})$, with natural map $\Spec A_k
\rightarrow S$. Pull back the log structure $S^{\dagger}$ on $S$ to 
$\Spec A_k$ to yield the log scheme $\Spec A_k^{\dagger}$. Then
$\Gamma(\X_k,\Theta_{\X_k^{\dagger}/A_k^{\dagger}})$
splits into $P^{\gp}$-homogeneous pieces
\[
\bigoplus_{p\in P^{\gp}} z^p\left(\Theta_{\X_k^{\dagger}/A_k^{\dagger}}\right)_p,
\]
where
$\left(\Theta_{\X_k^{\dagger}/A_k^{\dagger}}\right)_p=0$ 
if there does not exist an $i$,
$1\le i\le s$, such that $0\le \langle v_i,p\rangle\le k$;
otherwise
\[
\left(\Theta_{\X_k^{\dagger}/A_k^{\dagger}}
\right)_p=\left(\Theta_{Y^{\dagger}/S^{\dagger}}\right)_p.
\]

\end{proposition}

\proof We first observe there is a restriction map $\Theta_{Y^{\dagger}/\kk}
\rightarrow\Theta_{\X_k^{\dagger}/\kk}$. Indeed, given a log derivation $(\Di,\Dlog)$ 
on $Y$, for a function $f$ on $\X_k^{\dagger}$, we define $\Di|_{\X_k}f=(\Di\tilde f)|_{\X_k}$,
where $\tilde f$ is an extension of $f$ to $Y$. If $\tilde f,\tilde f'$ are two such
extensions, $\tilde f-\tilde f'=h\cdot z^{(k+1)\rho}$ for some function $h$,
and then
\[
\Di(\tilde f-\tilde f')=z^{(k+1)\rho}\Di h+h\cdot (k+1)z^{k\rho}\Di(z^{\rho}).
\]
Since $\Di(z^{\rho})$ is proportional to $z^{\rho}$, this is in the ideal of $\X_k$, hence
vanishes on the restriction to $\X_k$.

We also need to restrict $\Dlog$: this is important as the hypotheses of Proposition
\ref{Thetainclude} don't hold for $\X^{\dagger}_k$, so $\Theta_{\X_k^{\dagger}/\kk}$
is not contained in $\Theta_{\X_k/\kk}$. (For example, $\Theta_{\kk^{\dagger}/
\kk}$ has a non-zero element defined by $D=0$ and $\Dlog(n)=n$ for 
$n\in\NN$.)
But an element of $\M_{\X_k}$ is an equivalence class of pairs $(f,h)$ where 
$f\in \M_Y$, $h\in \O_{\X_k}^{\times}$, and $(f,h)\sim (f\cdot g, h\cdot g|_{\X_k}^{-1})$
for $g\in\O_Y^{\times}$. Then we can define
\[
\Dlog|_{\X_k}(f,h)=(\Dlog f)|_{\X_k}+h^{-1}\Di|_{\X_k}(h).
\]
This is easily checked to be well-defined.

Now all monomials in $(k+1)\rho+P$ restrict to zero on $\X_k$. Since
$\langle v_i,\rho\rangle=1$ for $1\le i\le s$ and $\langle w_j,\rho\rangle
=0$, one sees that $p\in P$
is in $(k+1)\rho+P$ if and only if $\langle v_i,p\rangle \ge k+1$
for all $1\le i\le s$. In case $z^p$ is one of the monomials with poles
along some component of $D$
occuring in the description of $\Theta_{Y^{\dagger}/\kk}$ of Proposition
\ref{Thetatoric}, then the same
condition tests to see if $z^p|_{\X_k}=0$, as the poles allowed are not
along any component of $X$. Thus we will obtain the desired description
of $\Theta_{\X_k^{\dagger}/\kk}$ if we show the restriction map is surjective,
and $z^p\partial_m$ restricts to zero on $\X_k$ if and only if 
$z^p$ restricts to zero.

To show this latter fact, note that $Y^{\dagger}\setminus D$ is log smooth
over $\kk$ as the log structure on $Y^{\dagger}\setminus D$ is
defined by the entire toric boundary. Thus $(\Theta_{Y^{\dagger}/\kk})|_{Y\setminus
D}$ is locally free, and it is easy to see that 
\[
(\Theta_{\X_k^{\dagger}/\kk})|_{\X_k\setminus D}
\cong (\Theta_{Y^{\dagger}/\kk})|_{Y\setminus D}\otimes \O_{\X_k\setminus D}.
\]
Now the restriction map $\Gamma(Y,\Theta_{Y^{\dagger}/\kk})\rightarrow\Gamma(Y\setminus
D,\Theta_{Y^{\dagger}/\kk})$ is injective, so a derivation of the form
$z^p\partial_m$ restricts to zero on $\X_k\setminus D$ if and only if $z^p$
restricts to zero on $\X_k\setminus D$. This proves there are no derivations
unexpectedly restricting to zero.

To show surjectivity of the restriction map
$\Theta_{Y^{\dagger}/\kk}\rightarrow\Theta_{\X_k^{\dagger}/\kk}$,
we make use of the Gorenstein condition. This is equivalent
to the existence of a $\rho_K\in P$ such that $\langle v_i,\rho_K\rangle =1$ for
$1\le i\le s+t$; this represents the canonical divisor on $Y$. 
Thus we can set 
\[
\bar\sigma=\{m\in\sigma|\langle m,\rho_K\rangle=1\}.
\]
This
polytope is the convex hull of $v_1,\ldots,v_{s+t}$.
Set $\tilde\rho=\rho_K-\rho$;
because $X$ is reduced, this takes non-negative values on all the $v_i$'s and hence
$\tilde\rho\in P$. In fact, $z^{\tilde\rho}=0$ defines the divisor $D$. 
Then $\Gamma(Y\setminus D,
\Theta_{Y^{\dagger}/\kk})$ is the localization of 
$\Gamma(Y,\Theta_{Y^{\dagger}/\kk})$ 
at $z^{\tilde\rho}$. Since $\X_k$ has no
embedded components as $Y$ is Cohen-Macaulay, 
there can be no log derivations with support in $\X_k\cap D$,
and hence $\Gamma(\X_k,\Theta_{\X_k^{\dagger}/\kk})$ injects into
$\Gamma(\X_k\setminus D,\Theta_{\X_k^{\dagger}/\kk})$. In addition,
$\Gamma(\X_k\setminus D,
\Theta_{\X_k^{\dagger}/\kk})$ is generated by derivations of the form
$z^{\tilde n\tilde\rho+n\rho+p}\partial_m$, where $\tilde n$ is any integer, $0\le n\le k$,
and $p\in P\setminus(\rho+P)$. 

Now given such a derivation $z^{\tilde n\tilde\rho+n\rho+p}\partial_m$, 
set $I=\{i|p\in P_i\}$. Let us assume this derivation on $\X_k\setminus D$
extends to a derivation on $\X_k$, and see what restriction on $\tilde n
\tilde\rho+n\rho+p$ we find.
First assume $m$ is not proportional to $w_j$ for any $j$. 
Take an $i\in I$. Look at the vertices $v$ of $\bar\sigma$ connected to $v_i$
by an edge. There are two cases.
If $v=v_j$ for some $1\le j\le s$, then $\langle v_j,
\tilde 
n\tilde\rho+p\rangle=\langle v_j,p\rangle\ge 0$. If $v=w_j$ for some $j$, then
$X_i\cap D_j$ is a divisor on $X_i$. We can choose $q\in Q_j$ such that $\langle m,q\rangle
\not=0$ and $\langle v_i,q\rangle\not=0$ for all $1\le i\le s$. Let $D_q=D\cap
\{z^q=0\}$. Note $D_j\not\subseteq D_q$.
Then $z^q$ is a function on $Y\setminus D_q$ which
vanishes only along $X\setminus D_q$, and hence $z^q$ is in $\Gamma(Y\setminus D_q,\M_Y)$.
Thus the pair $(z^q,1)$ represents an element of $\Gamma(\X_k\setminus D_q,\M_{\X_k})$.
Computing, we see that
\[
z^{\tilde n\tilde\rho+n \rho+p}\partial_m\log (z^q,1)=
\langle m,q\rangle z^{\tilde n\tilde\rho+n\rho+p}\not=0
\]
on $\X_k\setminus D_q$. If $z^{\tilde n\tilde\rho+n \rho+p}\partial_m$
were then a log derivation on all of $\X_k$, the function $z^{\tilde n\tilde\rho+n \rho+p}$
would not be allowed to have a pole on $\X_k\setminus D_q$. But the support of this
scheme includes a dense open subset of $X_i\cap D_j$, so the absence of poles implies
\[
0\le\langle w_j,\tilde n\tilde\rho+n\rho+p\rangle =\langle w_j,\tilde n\tilde\rho+p\rangle.
\]
Thus we learn that $\langle v, \tilde n\tilde\rho+p\rangle\ge 0$ 
for all vertices $v$
of $\bar\sigma$ adjacent to $v_i$. As $\langle v_i,\tilde n\tilde\rho+p\rangle =0$,
it then follows that $\langle v,\tilde n\tilde\rho+p\rangle\ge 0$ for all
vertices of $\bar\sigma$, hence $\tilde n\tilde\rho+p\in P$, hence
$\tilde n\tilde\rho+n\rho+p\in P$.

If $m$ is proportional to $w_j$ for some $j$, we can assume $m=w_j$. The same
argument shows that if $i\in I$, $v$ adjacent to $v_i$ in $\bar\sigma$,
then $\langle v,\tilde n\tilde\rho+p\rangle\ge 0$ unless $v=w_j$. Suppose
$v=w_j$ is adjacent to $v_i$. Then
$P_i\cap Q_j$ is a maximal proper face of $P_i$, and we can find a $q\in P_i$
such that $\langle w_j,q\rangle =1$. Then
\[
z^{\tilde n\tilde\rho+n\rho+p}\partial_{w_j} z^q=z^{\tilde n\tilde\rho+n \rho+p+q},
\]
and again if this derivation extends to $\X_k$, the function on the right hand side
cannot have a pole along $X_i\cap D_j$, i.e.
\[
0\le \langle w_j,\tilde n\tilde\rho+n\rho+p+q\rangle
=\langle w_j,\tilde n\tilde\rho+p\rangle+1.
\]
Thus either $\langle w_j,\tilde n\tilde\rho+p\rangle\ge 0$, and we finish
the argument as before to show $\tilde n\tilde\rho+n\rho+p\in P$, or else
$\langle w_j,\tilde n\tilde\rho+p\rangle =-1$, and then $(\tilde n+1)\tilde
\rho+p\in P$. In this latter case, 
$\bar\sigma\cap ((\tilde n+1)\tilde\rho+p)^{\perp}$
is a face of $\bar\sigma$ spanned by $S_1=\{v_i|i\in I\}$ and $S_2=\{w_k|
\langle w_k,\tilde n\tilde\rho+p\rangle=-1\}$.
Noting that the convex hulls of $S_1$ and $S_2$ are contained in
the planes $\langle\cdot,\rho\rangle=1$ and $\langle\cdot,\tilde\rho\rangle=1$
respectively, it follows that if $w_k\in S_2$, then in the convex hull
of $S_1\cup S_2$, $w_k$ is adjacent to some $v_i\in S_1$. But if $w_k\not=
w_j$, $w_k\in S_2$, the above argument applied with initial choice of
vertex $v_i$ shows that
$\langle w_k,\tilde n\tilde\rho+p\rangle
\ge 0$, a contradiction. Thus $\langle v_i, \tilde n\tilde\rho+n\rho+p\rangle\ge 0$
for all $1\le i\le s+t$ with $v_i\not=w_j$, and $\langle w_j,\tilde n\tilde\rho+n\rho+p\rangle
=-1$. This shows that $z^{\tilde n\tilde\rho+n\rho+p}\partial_{w_j}$ is in 
$(\Theta_{Y^{\dagger}/k})_{\tilde n\tilde\rho+n\rho+p}$, so the restriction
map is surjective, as desired. 

The second statement of the Proposition follows immediately from the first.
\qed

\bigskip

There exists also a universal log-derivation $(\di,\dlog)$ (see \cite{K.Kato}
for the coherent case):

\begin{lemma}
Given a morphism $\pi:X^{\dagger}\rightarrow S^{\dagger}$ of log schemes, let
\[
\Omega^1_{X^\ls/S^\ls}=
\big(\Omega^1_{X/S}\oplus (\O_X \otimes_\ZZ \M^\gp_X)\big)
\big/ \shK,
\]
with $\shK$ the $\O_X$-module generated by
\[
(\di\alpha_X(m), -\alpha_X(m)\otimes m),\eqand (0,1\otimes\pi^*(n)),
\]
for $m\in\M_X$, $n\in\M_S$. Then the pair $(\di,\dlog)$ of natural
maps
\[
\di: \O_X\stackrel{d}{\lra}\Omega^1_{X/S}\lra
\Omega^1_{X^\ls/S^\ls},\quad
\dlog: \M^\gp_X\stackrel{1\otimes\,\cdot}{\lra}
\O_X \otimes \M_X^{\gp}\lra \Omega^1_{X^\ls/S^\ls},
\]
is a universal log derivation.
\end{lemma}
\proof
We verify the universal property. Let $(\Di,\Dlog)$
be a log derivation with values in the coherent $\O_X$-module $\shE$:
\[
\Di:\O_X\lra \shE,\quad
\Dlog: \M^\gp_X\lra \shE.
\]
By the universal property of $\Omega^1_{X/S}$ there is a unique
morphism $\varphi: \Omega^1_{X/S}\to \shE$ fulfilling
\[
\Di= \varphi\circ \di.
\]
Define
\[
\Phi: \Omega^1_{X/S}\oplus (\O_X \otimes_\ZZ \M^\gp_X)
\lra \shE,\quad \Phi(\gamma, h\otimes m)= \varphi(\gamma)+h\cdot
\Dlog(m).
\]
This descends to the quotient by $\shK$ because
\[
\varphi\big(\di\alpha_X(m)\big) - \alpha_X(m)\cdot \Dlog(m)=0,
\quad \Dlog(\pi^\#(n))=0,
\]
for every $m\in \M_X$, $n\in\M_S$. Uniqueness follows since
$\Omega^1_{X/S}\oplus (\O_X \otimes_\ZZ \M^\gp_X)$ is generated as
$\O_X$-modules by $\Omega^1_{X/S}$ and by $1\otimes \M_X^\gp$. On
these subsets $\Phi$ is defined by $\varphi$ and by $\Dlog$
respectively.
\qed
\medskip

The $\O_X$-module $\Omega^1_{X^\ls/S^\ls}$ is the module of \emph{log
differentials}. It is coherent for fine log structures if $\pi$ is
locally of finite type. If $\pi$ is log smooth then
$\Omega^1_{X^\ls/S^\ls}$ is locally free (\cite{K.Kato}, Proposition~3.10).
We define
\[
\Omega^r_{X^{\dagger}/S^{\dagger}}={\bigwedge}^r \Omega^1_{X^{\dagger}/
S^{\dagger}}.
\]

\begin{remark}
If $\alpha:P\to \O_U$ is a chart for the log structure on $X$, then
in the formula for $\Omega^1_{X^\ls/S^\ls}$ one may replace
$\M^\gp_X$ by $P^\gp$ and $\alpha_X$ by $\alpha$. In fact, any
$h\in \O_X^\times$ gives a relation
\[
\big(\di h, -h\otimes \alpha_X^{-1}(h)\big)\in \shK.
\]
Therefore, for any $m\in\M_X^\gp$ the log differential
$(0,1\otimes(\alpha_X^{-1}(h)\cdot m))$ may be written as $h^{-1}(\di
h,0\otimes 1) +(0,1\otimes m)$, which is the sum of an ordinary
differential and a log differential involving only $m$.
\qed
\end{remark}

\begin{example} 
\label{ODPexample}
For the relatively coherent examples we wish
to consider, the sheaf of log differentials is poorly behaved
at points where the log structure is not coherent. For example,
take $P\subseteq\ZZ^3$ generated by $(1,0,0)$, $(-1,0,1)$, 
$(0,-1,1)$ and $\rho=(0,1,0)$. If these generators correspond
to variables $x,y,w$ and $t$ respectively, then $\kk[P]
\cong \kk[x,y,w,t]/(xy-wt)$ and $\kk[P]/(z^{\rho})
\cong \kk[x,y,w,t]/(xy,t)$. Let $Y=\Spec \kk[P]$, $X=\Spec \kk[P]/(z^{\rho})$
as usual, with the log structure on $X$ induced by the inclusion 
$X\subseteq Y$. Note any function on a neighbourhood $0\in Y$
with zero locus contained in $X$ is of the form
$f\cdot t^l$, where $f$ is invertible. Then a section of
$\M_X$ in a neighbourhood of $0$ is necessarily induced by such
a function: we write $(ft^l)|_X$ for the corresponding section
of $\M_X$. But in $\Omega^1_{X^{\dagger}/\kk^{\dagger}}$,
$(0,1\otimes (ft^l)|_X)=(f|_X^{-1}df|_X,0)$. Thus we see in
fact that $\Gamma(X,\Omega^1_{X^{\dagger}/\kk^{\dagger}})=
\Gamma(X,\Omega^1_{X/\kk})$. But if $\Omega^1_{X^{\dagger}/\kk^{\dagger}}$
were quasi-coherent, we would then have $\Omega^1_{X^{\dagger}/\kk^{\dagger}}
=\Omega^1_{X/\kk}$, which is not the case. \qed
\end{example}

As a consequence, it is not natural to use this universal sheaf.
Instead, we will use the push-forward of the sheaf of log differentials
on the log smooth part of $X$, as in the following proposition.

\begin{proposition}
\label{OmegaXl} 
In the situation of Proposition \ref{ThetaXl},
let $j:Z:=D\cap X_{\Sing}\hookrightarrow |\X_k|=|X|$ be the inclusion.
(Here $|\X_k|$ denotes the underlying topological space.)
Then $\Gamma(\X_k\setminus Z,
\Omega^r_{\X_k^{\dagger}/\kk})$ is naturally a $P$-module with decomposition into
$P$-homogeneous pieces given as follows:
\[
\Gamma(\X_k\setminus Z,\Omega^r_{\X^{\dagger}_k/\kk})
=\bigoplus_{p\in P\setminus ((k+1)\rho+P)}
{\bigwedge}^r\left(\bigcap_{\{j|p\in Q_j\}} Q_j^{\gp}\right)
\otimes_{\ZZ} \kk.
\]
Here $an_1\wedge\cdots\wedge
n_r$, for $a\in \kk$, $n_i\in P^{\gp}$, in the summand
of degree $p$ corresponds to the restriction of $az^p\dlog n_1\wedge
\cdots\wedge \dlog n_r \in\Gamma(Y\setminus Z,\Omega^r_{Y^{\dagger}/\kk})$ to $\X_k$.
\end{proposition}

\proof 
It is clear that $\Omega^1_{Y^{\dagger}/\kk}|_{Y\setminus Z}\subseteq 
\Omega^1_{\tilde Y^{\dagger}/\kk}|_{Y\setminus Z}$. (See the discussion
before Proposition~\ref{Thetatoric} for the definition of $\tilde Y^{\dagger}$.)
Note $\tilde Y^{\dagger}$ is log smooth over $\Spec \kk$.
Denoting
$\tilde\X_k^{\dagger}$ the restriction of the log structure 
$\tilde Y^{\dagger}$to $\X_k$, we obtain
\[
\Gamma(\X_k\setminus Z,\Omega^r_{\X_k^{\dagger}/\kk})
\subseteq\Gamma(\X_k\setminus Z,\Omega^r_{\tilde\X_k^{\dagger}/\kk})
=\Gamma(\X_k,\Omega^r_{\tilde\X_k^{\dagger}/\kk}).
\]
On the right-hand side we deal with a coherent sheaf on an affine space and 
hence there is a surjection
\[
\kappa:\Gamma(Y,\Omega^r_{\tilde Y^{\dagger}/\kk})\mapright{}\Gamma(\X_k,
\Omega^r_{\tilde\X_k^{\dagger}/\kk}).
\]
Now $\Omega^r_{\tilde Y^{\dagger}/\kk}$ is generated by $\dlog m$ for
$m\in\bigwedge^r P^{\gp}$ and hence
\[
\Gamma(Y,\Omega^r_{\tilde Y^{\dagger}/\kk})
=\bigoplus_{p\in P} z^p\left({\bigwedge}^r P^{\gp}\right)\otimes_{\ZZ} \kk.
\]
This description exhibits the $P$-grading.
Letting $\shI=(z^{(k+1)\rho})$ be the ideal sheaf of $\X_k\subseteq  Y$
we see
\begin{eqnarray*}
\ker(\kappa)&=&\Gamma(Y,\shI\Omega^r_{\tilde Y^{\dagger}/\kk})=
\Gamma(Y,\shI)\cdot \Gamma(\Omega^r_{\tilde Y^{\dagger}/\kk})\\
&=&\bigoplus_{p\in (k+1)\rho+P} z^p\left({\bigwedge}^r P^{\gp}\right)
\otimes_{\ZZ} \kk,
\end{eqnarray*}
and in turn 
\[
\Gamma(\X_k,\Omega^r_{\tilde\X_k^{\dagger}/\kk})=\bigoplus_{p\in P\setminus((k+1)\rho
+P)} z^p\left({\bigwedge}^rP^{\gp}\right)\otimes_{\ZZ} \kk.
\]
Another way to view the $P$-grading is as weights with respect to the
action of the algebraic torus $\Spec \kk[P^{\gp}]$. As this action 
respects the inclusions $X\subseteq Y$ and $D\subseteq Y$ it induces
an action on $\Gamma(Y,\Omega^r_{Y^{\dagger}/\kk})
\subseteq \Gamma(Y,\Omega^r_{\tilde Y^{\dagger}/\kk})$. From this
it is clear that for each $p\in P$ there exists a $\kk$-vector subspace
$V_p^r\subseteq \bigwedge^r P^{\gp}\otimes_{\ZZ} \kk$ such that 
\[
\Gamma(\X_k\setminus Z,
\Omega^r_{\X_k^{\dagger}/\kk})=\bigoplus_{p\in P\setminus ((k+1)\rho+P)} 
z^p V_p^r. 
\]
To finish the proof
it remains to describe $V_p^r$ for $p\in P\setminus ((k+1)\rho+P)$.
An element of $z^p(\bigwedge^r P^{\gp})$ is in 
$\Gamma(\X_k\setminus Z,\Omega^r_{\X_k^{\dagger}/\kk})$ if and only if the
contraction of it with any element of $\Gamma(\X_k\setminus Z,
\Theta_{\X_k^{\dagger}/\kk})=\Gamma(\X_k,\Theta_{\X_k^{\dagger}/\kk})$ 
(described in Propositions~\ref{ThetaXl} and \ref{Thetatoric})
is in
$\Gamma(\X_k\setminus Z,\Omega^{r-1}_{\X^{\dagger}_k/\kk})$. Thus we can
compute $V_p^r$ by induction. The result is obvious for $r=0$.
Assume now the formula for $V_p^{r-1}$ stated in the proposition
is correct, and let $n\in
\bigwedge^r P^{\gp}\otimes \kk$. 
Suppose $n\in V^r_p$, and in addition suppose $p\in Q_j$ for some $j$. Since 
$p\in P\setminus ((k+1)\rho+P)$, there must be a 
$v_i$ connected to $w_j$ by an edge in $\bar\sigma$ (as defined in the proof
of Proposition~\ref{ThetaXl}) such that $\langle v_i,p\rangle
\le k$. 
(Indeed, $\bar\sigma$ is the convex hull of two polytopes, $\bar\sigma_0$
and $\bar\sigma_1$ on which $\rho$ takes the values $0$ and $1$ respectively.
Then $\langle p,\cdot\rangle=0$ defines a supporting hyperplane for
$\bar\sigma_0$. Let $l$ be the minimal value $p$ takes on $\bar\sigma_1$.
By assumption, $l\le k$.
Then $p-l\rho$ defines a supporting hyperplane of $\bar\sigma$ which
contains at least one vertex of $\bar\sigma_1$; 
one of these vertices will be connected by an edge of $\bar\sigma$
to $w_j$, say $v_i$. Then $\langle v_i,p\rangle=l\le k$.)
We can then find a $q_j\in P^{\gp}$ such that $\langle w_j,q_j\rangle
=\langle v_i,q_j\rangle=-1$ and $\langle v_k,q_j\rangle \ge 0$ for
$v_k\not\in\{w_j,v_i\}$. Indeed, we just need $q_j+\rho+\tilde\rho$ to
define a supporting
hyperplane for the edge joining $w_j$, $v_i$. Replacing $q_j$
with $q_j+\rho$, we obtain a $q_j$ with $\langle w_j,q_j\rangle =-1$,
$\langle v_i,q_j\rangle=0$, $\langle v_k,q_j\rangle\ge 0$ for
all $v_k\not\in\{w_j,v_i\}$. Then $\langle v_i,q_j+p\rangle\le k$, so 
$z^{q_j+p}$ does not vanish on $\X_k$. Now
\[
\iota(z^{q_j}\partial_{w_j})(z^p\dlog n)= z^{p+q_j}\dlog \iota(w_j)n
\in\Gamma(\X_k\setminus Z,\Omega^{r-1}_{\X_k^{\dagger}/\kk})
\]
if and only if $\iota(w_j)n\in V^{r-1}_{p+q_j}$. Since $p\in Q_j$,
$\langle w_j,p+q_j\rangle=\langle w_j,q_j\rangle=-1$, so $p+q_j\not\in P$ and $V^{r-1}_{p+q_j}=0$.
Hence if $p\in Q_j$, we must have $\iota(w_j)n=0$, which is the case if and
only if $n\in \bigwedge^r Q_j^{\gp}\otimes \kk$. Thus we see
\[
n\in \bigcap_{\{j| p\in Q_j\}} {\bigwedge}^r Q_j^{\gp}\otimes \kk
={\bigwedge}^r\bigcap_{\{j| p\in Q_j\}} Q_j^{\gp}\otimes \kk.
\]
(This equality is easily checked, say, for the intersection of
two subspaces of any vector space.) Note this statement is true
even if $p\not\in Q_j$ for any $j$, as then it is vacuous.

Conversely, if $n\in \bigwedge^r\bigcap_{\{j|p\in Q_j\}} Q_j^{\gp}\otimes \kk$,
then we need to check 
\[
z^p\dlog n\in\Gamma(\X_k\setminus Z,\Omega^r_{\X_k^{\dagger}/\kk}).
\]
If $q\in P$, $m\in P^{\gp}\otimes \kk$, then
\[
\iota(m)n\in {\bigwedge}^{r-1}\bigcap_{\{j|p\in Q_j\}} Q_j^{\gp}\otimes \kk
\subseteq {\bigwedge}^{r-1}\bigcap_{\{j| p+q\in Q_j\}} Q_j^{\gp}\otimes \kk
=V_{p+q}^{r-1},
\]
so $\iota(z^q\partial_m)z^p\dlog n\in\Gamma(\X_k\setminus Z,
\Omega^{r-1}_{\X_k^{\dagger}/\kk})$. 
The other case to check is $q\not\in P$, so there exists a $j$ with
$\langle w_j,q\rangle=-1$,
$\langle v_i,q\rangle\ge 0$ for all $v_i\not=w_j$, $1\le i\le s+t$; then
\[
\iota(z^q\partial_{w_j})z^p\dlog n=
\begin{cases} 0 &\hbox{if $p\in Q_j$,}\\
z^{p+q}\dlog \iota(w_j) n&\hbox{if $p\not\in Q_j$.}
\end{cases}
\]
In the latter case, $\langle w_j,p+q\rangle\ge 0$ so $p+q\in P$,
and 
\[
\iota(w_j)n\in {\bigwedge}^{r-1} \left(Q_j^{\gp}\cap \bigcap_{\{l|p\in Q_l\}}
Q_l^{\gp}\right)\otimes \kk
\subseteq {\bigwedge}^{r-1} \bigcap_{\{l|p+q\in Q_l\}} Q_l^{\gp}\otimes \kk
=V^{r-1}_{p+q}.
\]
\qed

We then immediately obtain

\begin{corollary}
\label{OmegaXldagger}
In the situation of Proposition \ref{OmegaXl},
$\Gamma(\X_k\setminus Z,\Omega^r_{\X_k^{\dagger}/A_k^{\dagger}})$
is naturally a $P$-module with decomposition into $P^{\gp}$-homogeneous
pieces
\[
\Gamma(\X_k\setminus Z,\Omega^r_{\X_k^{\dagger}/A_k^{\dagger}})
=
\bigoplus_{p\in P\setminus ((k+1)\rho+P)}
{\bigwedge}^r\bigg(\bigcap_{\{j|p\in Q_j\}} Q_j^{\gp}/\ZZ\rho\bigg)
\otimes_{\ZZ} \kk.
\]
\end{corollary}

\bigskip


\section{Log Calabi-Yau spaces: local structure and deformation theory}

We now return to the type of log spaces considered in \cite{PartI}. There,
we defined the notion of a log Calabi-Yau space. In particular, we
showed how to construct log Calabi-Yau spaces from integral affine
manifolds with singularities with toric polyhedral decompositions 
$(B,\P)$ which are positive and simple (see \cite{PartI}, Definition~1.60). Our
goal in this section are two-fold. First, we wish to show that the
log Calabi-Yau spaces we constructed in \cite{PartI} in fact are locally
of the form considered in \S1; we can then use the results of \S1 to
describe their sheaves of differentials, something we will do in \S3.
Second, we wish to understand the deformation theory of log schemes
with local structure as described in \S1. We will again need to
make detailed use of the local descriptions of these log structures.

\subsection{Local structure}

\begin{construction}
\label{keyP}
We will begin by refining the somewhat general examples considered in
the previous section. We fix the following data. Let $M'$ be a lattice,
$N'$ the dual lattice, and set $M=M'\oplus \ZZ^{q+1}$, $N$ the dual lattice.
We write $e_0,\ldots,e_q$ for the standard basis of $\ZZ^{q+1}$, and we
identify these with $(0,e_0),\ldots,(0,e_q)$ in $M$. Thus we can write a
general element of $M$ as $m+\sum a_i e_i$ for $m\in M'$. Similarly, we
write $e_0^*,\ldots,e_q^*$ for the dual basis, which we view as elements of
$N$. 

Fix a convex lattice
polytope $\tau\subseteq M'_{\RR}$ with $\dim\tau=\dim M'_{\RR}$, 
with normal fan 
$\check\Sigma_{\tau}$ living in $N'_{\RR}$ (see \cite{PartI},
Definition 1.38 for our convention concerning the normal fan). We obtain a cone
$C'(\tau)\subseteq M'_{\RR}\oplus\RR$, $C'(\tau)=
\{(rm,r)|r\ge 0,m\in\tau\}$, and a monoid $P'=
\dual{C'(\tau)}\cap (N'\oplus\ZZ)$. Define $\rho'\in P'$ to be
given by the projection $M'\oplus\ZZ\rightarrow \ZZ$. We set
\[
V'(\tau)=\Spec \kk[P']/(z^{\rho'})=\Spec \kk[\partial P']
\]
(c.f.\ \cite{PartI}, Definition 2.13). Here $\partial P'$ is
the monoid consisting of elements of the boundary of $P'$ and $\infty$,
with $p+p'$ defined to be $p+p'$ if $p+p'$ lies in the boundary of
$P'$ and $\infty$ otherwise.
As in \cite{PartI}, we identify $\partial P'$ as a set
with $N'\cup\{\infty\}$ via projection to $N'$. We always use the 
convention that $z^{\infty}=0$.

Let $\check\psi_1,\ldots,\check\psi_q$
be integral piecewise linear functions on $\check\Sigma_{\tau}$
whose Newton polytopes are
$\Delta_1,\ldots,\Delta_q\subseteq M_{\RR}'$, i.e.
\[
\check\psi_i(n)=-\inf\{\langle n,m\rangle | m\in\Delta_i\}.
\]
Similarly, let $\check\psi_0$ have Newton polytope $\tau$, i.e.
\[
\check\psi_0(n)=-\inf\{\langle n,m\rangle | m\in\tau\}.
\]
For convenience of notation, we set $\Delta_0:=\tau$.

Given this data, we can define a monoid $P\subseteq N$ given by
\[
P=\bigg\{n+\sum_{i=0}^q a_ie_i^*\,\bigg|\,\hbox{$n\in N'$ and $a_i\ge\check\psi_i(n)$ for
$0\le i\le q$}\bigg\}.
\]
Set $Y=\Spec \kk[P]$ as in \S 1. Note that
$P=\dual{K}\cap N$ where $K$ is the cone in $M_{\RR}$ 
generated by
\[
\bigcup_{i=0}^q(\Delta_i\times \{e_i\}).
\]
In particular, we see $Y$ is Gorenstein because $\rho_K
=\sum_{i=0}^q e_i^*$ takes the value $1$ on each primitive integral
generator of an extremal ray of $K$. Letting $X=\Spec \kk[P]/(z^{\rho})$
as usual with $\rho:=e_0^*$, we describe $X$ explicitly by defining 
\[
Q=\bigg\{n+\sum_{i=0}^q a_ie_i^*\in P\,\bigg|\,a_0=\check\psi_0(n)\bigg\}
\cup\{\infty\}
\]
with addition on $Q$ defined by
\[
q_1+q_2=
\begin{cases}
q_1+q_2& \hbox{if $q_1+q_2\in Q$}\\
\infty&\hbox{otherwise.}
\end{cases}
\]
Then $Q\setminus\{\infty\}$ is, as a set, $P\setminus(\rho+P)$, so
it is clear that $X=\Spec \kk[Q]$. Note that $Q\cong\partial P'
\oplus\NN^q$, via
\[
(n,a_0,\ldots,a_p)\mapsto (n,0,a_1-\check\psi_1(n),\ldots,
a_q-\check\psi_q(n)).
\]
Thus $X=V'(\tau)\times\AA^q$. 

We define subschemes $Z_i$ of $X$
by their ideals, for $1\le i\le q$,
with $I_{Z_i/X}$ generated by the set of monomials
\[
\{z^{e_i^*}\}\cup 
\left\{z^p\bigg|{\hbox{$p=n+\sum a_je_j^*$ such that there exists a unique}
\atop\hbox{vertex $w$ of $\Delta_i$ such that $\langle n,w\rangle=-\check
\psi_i(n)$}}\right\}.
\]
The effect of the right-hand set is to select those irreducible components
of the singular locus of $X$ corresponding to edges of $\Delta_i$,
and $z^{e_i^*}$ defines a closed subscheme of this set of components.
Set
\[
Z=\bigcup_{i=1}^q Z_i.
\]
This will be the locus where the log structure on $X$ fails to
be coherent. 
Let 
\[
u_i:= z^{e_i^*}
\]
for $1\le i\le q$. For any vertex $v$ of $\tau$, denote by
$\Vert_i(v)$ the vertex of $\Delta_i$ which represents the function
$-\check\psi_i$ restricted to the maximal cone $\check v$ of $\check\Sigma_{\tau}$
corresponding to $v$. 
For every edge $\omega\subseteq\tau$, choose a primitive
generator $d_{\omega}$ of the tangent space of $\omega$, and
let $v^{\pm}_{\omega}$
be the two vertices of $\omega$, labelled so that 
$d_{\omega}$ points from $v^+_{\omega}$ to $v^-_{\omega}$
as in \cite{PartI}, \S 1.5. Set
\[
\Omega_i=\{\omega\subseteq\tau| \hbox{$\dim\omega=1$ and 
$\Vert_i(v^+_{\omega})\not=\Vert_i(v^-_{\omega})$}\}.
\]
(This notation is compatible with that in the definition of simplicity,
\cite{PartI}, Definition 1.60.) It is then easy to see that
\[
Z_i=\{u_i=0\}\cap \bigcup_{\omega\in \Omega_i} V_{\omega}.
\]
Here for $\omega\subseteq\tau$ any face, we define $V_{\omega}\subseteq
X$ to be the
closed toric stratum of $Y$ defined by the face
of $K$ generated by $\omega\times\{e_0\}$.

Similarly we define $V'_{\omega}$, for any face $\omega\subseteq\tau$,
to be the closed stratum of $V'(\tau)$ corresponding to $C'(\omega)
\subseteq C'(\tau)$.
\qed
\end{construction}

\bigskip

For future use, we note

\begin{proposition}
\label{smoothcase}
Let $P$ be as in Construction \ref{keyP}, determined by data
$\tau,\Delta_1,\ldots,\Delta_q$. Then the generic fibre of 
$f:\Spec \kk[P]\rightarrow\Spec \kk[\NN]$ induced by $\rho=e_0^*$
is non-singular if and only if 
\[
\Conv\bigg(\bigcup_{i=1}^q\Delta_i\times\{e_i\}\bigg)
\]
is a standard simplex.

If 
\[
\Conv\bigg(\bigcup_{i=1}^q\Delta_i\times\{e_i\}\bigg)
\]
is an elementary simplex (i.e. its only integral points are its
vertices) then the generic fibre of $f$ has codimension at least four
Gorenstein quotient singularities.
\end{proposition}

\proof Since $K=\Cone(\bigcup_{i=0}^q \Delta_i\times\{e_i\})$ is the
cone defining the monoid $P$ via $P=\dual{K}\cap N$, it is
a standard fact of toric geometry that the generic fibre of $f$
is then defined similarly by the cone 
\[
K\cap\rho^{\perp} =\Cone\bigg(\bigcup_{i=1}^q\Delta_i\times\{e_i\}\bigg).
\]
This cone is the standard simplicial cone corresponding to affine space
if and only if
$\Conv(\bigcup_{i=1}^q\Delta_i\times\{e_i\})$ is the standard simplex.

The second part is then obvious, given that an elementary simplex is
standard if it is dimension $\le 2$.
\qed

\bigskip

We now turn to the global situation considered in \cite{PartI}. 
We will need to consider certain
sorts of \'etale neighbourhoods of log schemes:

\begin{definition}
\label{strictetale}
A morphism $\phi:X^{\dagger}\rightarrow Y^{\dagger}$ is \emph{strict
\'etale} if it is \'etale as a morphism of schemes and is strict,
i.e. the log structure on $X^{\dagger}$ is the same as the pull-back
of the log structure on $Y^{\dagger}$.
\end{definition}

\begin{remark}
This is consistent with standard notation: if a morphism of coherent 
log structures is log \'etale and strict, then it is also \'etale
in the scheme theoretic sense. Thus in the coherent case
the above definition fits with the standard one, and we are extending this
to the non-coherent case. 

Strict \'etale morphisms have the following standard
property of \'etale morphisms: If $Y^{\dagger}$
is a log scheme, and $Y_0^{\dagger}$ is a closed subscheme
of $Y^{\dagger}$ defined by a nilpotent sheaf of ideals 
with the induced log structure on $Y_0$, then 
there is an equivalence between the
categories of strict \'etale $Y^{\dagger}$-schemes and strict \'etale 
$Y_0^{\dagger}$-schemes. Indeed, $X\mapsto X_0=X\times_Y Y_0$ gives an 
equivalence of categories between \'etale $Y$-schemes and \'etale
$Y_0$-schemes (\cite{Milne}, Chap. I, Theorem 3.23), 
and to obtain the log structures on 
$X$ or $X_0$, one just pulls back the log structure on $Y$ or $Y_0$.

In particular, if we have a strict \'etale morphism $X_0^{\dagger}
\rightarrow Y_0^{\dagger}$ and a thickening $Y^{\dagger}$ of
$Y_0^{\dagger}$, we can talk about pulling back this thickening to
$X_0^{\dagger}$ giving $X^{\dagger}$.

Note also that if $f:X^{\dagger}\rightarrow Y^{\dagger}$ is a strict
\'etale morphism over $\Spec \kk^{\dagger}$, then $\Theta_{X^{\dagger}/\kk}
=f^*\Theta_{Y^{\dagger}/\kk}$ and $\Theta_{X^{\dagger}/\kk^{\dagger}}
=f^*\Theta_{Y^{\dagger}/\kk^{\dagger}}$, as is easily checked. \qed
\end{remark}

We now turn to the basic setup we will consider in the remainder
of this paper, applying the above results to the log spaces constructed
in \cite{PartI}. In what follows, 
let $B$ be an integral affine manifold with singularities with
a toric polyhedral decomposition $\P$, (\cite{PartI}, Definitions 1.15, 1.22) 
and suppose $(B,\P)$ is
positive and simple (\cite{PartI}, Definitions~1.54 and 1.60). 
Then a choice of lifted open gluing data $s\in H^1(B,i_*\Lambda
\otimes \kk^{\times})$ (\cite{PartI}, Definition 5.1)
determines a toric log Calabi-Yau space (\cite{PartI}, Definition 4.3)
$X_0(B,\P,s)^{\dagger}$, 
with a morphism of log schemes $X_0(B,\P,s)^{\dagger}\rightarrow\Spec \kk^{\dagger}$
which is log smooth off of a codimension two set $Z$ as demonstrated in
\cite{PartI}, Theorem 5.2. We will fix $(B,\P,s)$ now, and write $X_0^{\dagger}$
instead of $X_0(B,\P,s)^{\dagger}$.
We wish to describe local
models for $X_0$ at points of $Z$; because of the simplicity of $(B,\P)$,
the singularities of the log structure will be well-behaved, and essentially
will be of the sort considered above.

We remark here that in \cite{PartI}, Remark 4.5, we mentioned we never really
need to worry about the precise log structure of $X_0$ along
$Z$. In particular, the definition of an isomorphism of
toric log Calabi-Yau spaces (in \cite{PartI}, Definition 4.3) ignores the
log structure on $Z$, and we will continue this policy below. We will
never check isomorphisms of the log structure on the locus where the
log structures fail to be coherent. 
This won't affect any calculations below. Indeed, 
all of our log schemes $X^{\dagger}$ are $S_2$, with $\M_X$
having no sections with support in codimension two, and thus by
Remark \ref{Thetaextend}, modules of log derivations with values
in $\O_X$ are insensitive to the precise log structure on the 
singular locus $Z$. Furthermore, sheaves of log differentials are only
considered off of $Z$.

\begin{remark}
\label{simplicitydata}
We shall almost always assume $(B,\P)$ is positive and simple in
this paper (see \cite{PartI}, Definition~1.60). Thus for a cell
$\tau\in\P$ with $0<\dim\tau<\dim B$, we always obtain associated to
$\tau$ the data $\Omega_1,\ldots,\Omega_p$, $R_1,\ldots,R_p$,
$\Delta_1,\ldots,\Delta_p$ and
$\check\Delta_1,\ldots,\check\Delta_p$, with $\Delta_i\subseteq
\Lambda_{\tau,\RR}$ and $\check\Delta_i\subseteq
\Lambda_{\tau,\RR}^{\perp}$ elementary simplices. 
($\Lambda_{\tau,\RR}$ is the tangent space to $\tau$ in $B$: see
\cite{PartI}, Definition~1.31.)  We call this data \emph{simplicity
data} for $\tau$.
\end{remark}

\begin{theorem}
\label{localmodel}
Given $(B,\P)$ positive and simple and $s$ lifted open gluing data,
let $\bar x\rightarrow Z\subseteq X_0$ be a closed geometric point. 
Then there exists data $\tau$, $\check\psi_1,\ldots,
\check\psi_q$ as in Construction \ref{keyP}
defining a monoid $P$, and an element
$\rho\in P$, hence log spaces $Y^{\dagger}$, $X^{\dagger}\rightarrow
\Spec \kk^{\dagger}$ as in \S 1, 
such that there is a diagram over $\Spec \kk^{\dagger}$ 
\begin{equation}
\label{etalediagram}
\xymatrix{
&V^{\dagger}\ar[ld]\ar[rd]^{\phi}&\\
X_0^{\dagger}&&X^{\dagger}
}
\end{equation}
with both maps strict \'etale and $V^{\dagger}$
an \'etale neighbourhood of $\bar x$.
\end{theorem}

\proof 
By \cite{PartI}, Lemma 2.29, for each $\tau\in\P$ there is a map
$q_{\tau}:X_{\tau}\rightarrow X_0$,
where $X_{\tau}$ is the toric variety defined in \cite{PartI},
Definition 2.7; this is the normalization of the stratum of
$X_0$ corresponding to $\tau$.
There exists a unique
$\tau\in\P$ with $\bar x\in q_{\tau}(X_{\tau}\setminus\partial X_{\tau})$,
and since $\bar x\in Z$, we have 
$0<\dim\tau<\dim B$. We then obtain simplicity data associated to
$\tau$ as in Remark \ref{simplicitydata}.
Set
$M'=\Lambda_{\tau}$, $M'_{\RR}=
\Lambda_{\tau,\RR}$, and identify $\tau$ as a polytope in $\Lambda_{\tau,\RR}$,
well-defined up to translation. (Without loss of generality we won't
worry about self-intersections of $\tau$.) 
By \cite{PartI}, Corollary 5.8, $q_{\tau}^{-1}(Z)\setminus\partial X_{\tau}=
Z_1^{\tau}\cup\cdots\cup Z_p^{\tau}$
where each $Z_i^{\tau}$ 
is an irreducible and smooth
divisor on $\Int(X_{\tau}):=X_{\tau}\setminus\partial
X_{\tau}$, with Newton
polytope given by $\check\Delta_i$. 
Furthermore, it is not difficult to see these divisors meet transversally,
using the fact that simplicity implies
the Newton polytopes $\check\Delta_i$ are elementary
simplices and their tangent space $T_{\check\Delta_i}$ in 
$\Lambda_{\tau,\RR}^{\perp}$ form an interior direct sum (see \cite{PartI}, 
Remark 1.61, (4)).
Assume that $\bar x
\in q_{\tau}(Z_i^{\tau})$, for $1\le i\le r$ for some $r\le p$ (reorder the indices
if necessary) and redefine $\Delta_i=\{0\}$ for $i>r$. In addition,
set $\Delta_i=\{0\}$ 
for $p< i \le\dim B-\dim\tau$. This provides data
$\tau$, $\Delta_1,\ldots,\Delta_q$ for $q=\dim B-\dim\tau$, and the 
$\Delta_i$'s define piecewise linear functions $\check\psi_i$ on the normal
fan $\check\Sigma_{\tau}$ in $N'_{\RR}$ (see \cite{PartI}, Remark~1.59:
there is a typo in that Remark; $\varphi_{\rho}$ should be $\check\psi_{\rho}$.)
This data produces monoids $P'\subseteq N'$ with $P'=\dual{C'(\tau)}\cap N'$,
a monoid $P\subseteq N=N'\oplus \ZZ^{q+1}$ as constructed from this
data in Construction \ref{keyP}, $\rho\in P$ given by $\rho=e_0^*$, yielding
$Y=\Spec\kk[P]$ and $X=\Spec\kk[P]/(z^{\rho})$. 
We recall $X\cong \Spec \kk[\partial P'\oplus\NN^q]$.
We will now construct the diagram (\ref{etalediagram}).

First, pick some $g:\tau\rightarrow\sigma\in\P_{\max}$, so we obtain
as in \cite{PartI}, Construction 2.15
an open set $V(\tau)\subseteq V(\sigma)$.
As observed in \cite{PartI}, Construction 2.15, $V(\tau)\cong
\Spec \kk[\partial P']\times \Gm^q$. The factor $\Gm^q$ can be identified
with $\Int(X_{\tau})$. In particular,
there is a unique point $0\in \Spec \kk[\partial P']$ whose ideal is generated
by $\{z^p| p\in\partial P', p\not=0\}$, and we identify the stratum
$0\times \Gm^q$ with $\Int(X_{\tau})$ also.
Then there is an \'etale
map $p_{\sigma}:V(\sigma)\rightarrow X_0$ by the construction
of $X_0$: see \cite{PartI}, Definition 2.28. 
Furthermore $p_{\sigma}$ is one-to-one on the interior of any
toric stratum, so in particular is one-to-one on $\Int(X_{\tau})$.
So there exists a unique point $\bar v\in \Int(X_{\tau})$ with 
$p_{\sigma}(\bar v)=\bar x$.
In particular, $V(\tau)$ is an \'etale neighbourhood of $\bar x$.

We pull back the
log structure on $X_0$ to $V(\sigma)$: along with the induced
morphism to $\Spec \kk^{\dagger}$, this structure is then determined by 
a section $f_{\sigma}\in \Gamma(V(\sigma),\shLS^+_{\pre,V(\sigma)})$,
(see \cite{PartI}, Definition 4.21)
or equivalently,
\[
(f_{\sigma,e})_{e:\omega\rightarrow\sigma}
\in
\Gamma\bigg(V(\sigma),\bigoplus_{e:\omega
\rightarrow\sigma\atop \dim\omega=1} \O_{V_e}\bigg).
\]
(Recall \cite{PartI}, Definition~2.22 that for $e:\omega\rightarrow
\sigma$, $V_e$ is the open affine subset of $X_{\omega}$ determined
by $\sigma$.)
Furthermore, it is argued in the proof of \cite{PartI}, Theorem 5.2, (2),
that for any factorization $\omega\rightarrow\tau'\mapright{g'}
\sigma$ of $e$, $f_{\sigma,e}|_{V_{g'}}$ is completely determined 
by $s$ and the Newton polytope of $f_{\sigma,e}|_{V_{g'}}$. 

We first claim the functions $f_{\sigma,g\circ e}$ for $e\in\Omega_i$ glue.
Indeed, let $e_1,e_2\in\Omega_i$, with $e_i:\omega_i\rightarrow\tau$,
and suppose we have a diagram
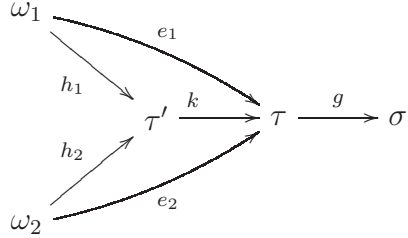
\[
\xymatrix@C=30pt
{\omega_1\ar@/^/[rrd]^{e_1}\ar[rd]_{h_1}&&&\\
&\tau'\ar[r]^(.3){k}&\tau\ar[r]^g&\sigma\\
\omega_2\ar@/_/[rru]_{e_2}\ar[ru]^{h_2}&&&\\ 
}
\]
Let $\Omega_1',
\ldots,\Omega'_{p'}$, $\Delta_1',\ldots,\Delta'_{p'}$ etc. be the simplicity data
for $\tau'$. Then as there exists $l:\tau\rightarrow\rho\in\P$ with
$l\in R_i$ such that $\kappa_{\omega_1,\rho},\kappa_{\omega_2,\rho}
\not=0$, (see \cite{PartI}, Remark 1.61, (3)) 
it follows that $h_1,h_2$ are both in the
same $\Omega_j'$ for some $j$ by \cite{PartI}, Definition 1.60, (1). 
Thus $f_{\sigma,g\circ e_1}|_{V_{g\circ k}}$ and
$f_{\sigma,g\circ e_2}|_{V_{g\circ k}}$ both have the same Newton polytope, namely
$\check\Delta'_j$, and hence must coincide, since as remarked
in the previous paragraph, these functions are completely determined
by $s$ and their Newton polytope. 

This gives rise, for each $i\le r$, to a function $f_i$ on 
the closed subscheme
$\bigcup_{e\in\Omega_i} V_e\subseteq V(\tau)$. We can then extend
each $f_i$ to a function on $V(\tau)$, which we also call $f_i$. 

Next, choose 
coordinates $z_1,\ldots,z_q$
on $\Int(X_{\tau})\cong\Gm^q$, and pull these back to functions on
$V(\tau)=\Spec \kk[\partial P']\times \Gm^q$.
Recalling from the first paragraph of this proof that
 $Z_1^{\tau},\ldots,Z_r^{\tau}$ intersect transversally,
one can find a subset $\{i_1,\ldots,i_r\}\subseteq \{1,\ldots,q\}$
such that $\det (\partial f_i/\partial z_{i_j})_{1\le i,j\le r}
\not=0$ at $\bar v$. By reordering the indices,
we can assume $\{i_1,\ldots,i_r\}=\{1,\ldots,r\}$, and set $f_i:=z_i-z_i(\bar v)$
for $r+1\le i\le q$, where $z_i(\bar v)$ denotes the value of the function
$z_i$ at the point $\bar v$.

Now consider the section $f'_{\sigma}\in \Gamma(V(\tau),\shLS^+_{\pre,V(\sigma)})$
defined by
\[
f'_{\sigma,e}=\begin{cases} 1&\hbox{if $e\not\in \bigcup_{i=1}^r \Omega_i$;}\\
f_i|_{V_e}&\hbox{if $e\in\Omega_i$, $1\le i\le r$.}
\end{cases}
\]
Then essentially repeating the argument of the proof of \cite{PartI}, Theorem
5.2, (2), one sees that $f'_{\sigma}$ satisfies the multiplicative
condition of \cite{PartI}, Theorem 3.22. As $f_{\sigma}$ itself
satisfies this multiplicative condition, since $f_{\sigma}$ defines a log
smooth structure away from the singular locus $p_{\sigma}^{-1}(Z)$,
we see that
the section $a_{\sigma}=f_{\sigma}/f'_{\sigma}
\in\Gamma(V(\tau),\shLS^+_{\pre,V(\sigma)})$
also satisfies this multiplicative condition. Note
\[
a_{\sigma,e}=\begin{cases} f_{\sigma,e}
&\hbox{if $e\not\in \bigcup_{i=1}^r \Omega_i$;}\\
1&\hbox{if $e\in\Omega_i$, $1\le i\le r$.}
\end{cases}
\]
Also $a_{\sigma,e}$ is always invertible at $\bar v$, and so $a_{\sigma}$ is
a section of $\shLS_{V(\sigma)}$ (see \cite{PartI}, Definition~3.19)
on some Zariski open neighbourhood
$V\subseteq V(\tau)$ of $\bar x$. We can assume that 
$\det (\partial f_i/\partial z_j)_{1\le i,j\le r}$ 
is also invertible on $V$. Then, possibly shrinking $V$
more, the section $a_{\sigma}$ determines, via \cite{PartI}, Theorem 3.22,
a chart
\[
P'\rightarrow \Gamma(V,\O_V)
\]
of the form
\[
p\mapsto h_pz^p
\]
where $h_p$ is some invertible function on $V(p)\cap V$. 
(Here, $V(p)$ is the closure of the open subset of $V(\tau)$
where $z^p\not=0$.)
This is \emph{not}
a chart defining the log structure $V^{\dagger}$ on $V$. 

We now define a morphism $\phi:V\rightarrow X=\Spec \kk[\partial P'\oplus
\NN^q]$ by
\begin{eqnarray*}
\phi^*(z^p)&=&h_pz^p\quad\hbox{for $p\in\partial P'$}\\
\phi^*(u_i)&=&f_i\quad\hbox{for $1\le i\le q$}
\end{eqnarray*}
where $u_i$ is the monomial in $\kk[\partial P'\oplus \NN^q]$
corresponding to the $i$th standard basis vector of $\NN^q$. 

\medskip

\noindent
{\it Claim}: Possibly after shrinking $V$, $\phi$ is \'etale.

\proof We prove this via a slight modification of the argument
in \cite{PartI}, Proposition 3.20, which essentially dealt with the case when
$\bar x\not\in Z$. First rewrite $\kk[\partial P']$ as
$\kk[{\mathbf X}]/({\mathbf G})$, where ${\mathbf X}$ denotes a collection
of variables $\{X_k\}$, with $X_k=z^{p_k}$ for a set of generators $\{p_k\}$
of $\partial P'$. The relation ideal $({\mathbf G})$ is generated by
a finite number of polynomials $G_{\lambda}$ of the form $\prod X_k^{a_k}\
-\prod X_k^{b_k}$ with $\sum a_kp_k=\sum b_kp_k$, or $\prod X_k^{a_k}$ with
$\sum a_kp_k=\infty$. For each $k$, choose an extension $h_k$ of
$h_{p_k}$ to $V$; after shrinking $V$ if necessary we can assume
$h_k$ is invertible on $V$. Again, by shrinking $V$ we can assume
$V\subseteq V(\tau)$ is given by $G\not=0$ for some function $G$
on $V(\tau)$. Let $A=\kk[\partial P'\oplus \NN^q]=\kk[{\mathbf X},u_1,\ldots,
u_q]/({\mathbf G})$. Then the morphism $\phi:V\rightarrow \Spec A$ is given
by the map of rings induced by the inclusion of $A$ in the polynomial ring
over $A$,
\[
A\rightarrow A[{\mathbf S},z_1^{\pm 1},\ldots,z_q^{\pm 1},G^{-1}]/I
\]
where the variables ${\mathbf S}=\{S_k\}$ are in one-to-one correspondence
with the variables $\{X_k\}$, and the ideal $I$ is generated by
\begin{eqnarray*}
{\mathbf X}-{\mathbf h}'\cdot {\mathbf S}&&\\
u_i-f_i({\mathbf S},{\mathbf z})&\quad&1\le i\le q
\end{eqnarray*}
where ${\mathbf h}'$ denotes
${\mathbf h}=\{h_k\}$ with the $X_k$'s replaced by $S_k$'s. According to
standard results concerning \'etale morphisms (see e.g. \cite{Milne}, I. 3.16),
a map of affine schemes induced by a map of rings $C\rightarrow C[Z_1,\ldots,Z_n]/
(F_1,\ldots,F_n)$ is \'etale if and only if $\det (\partial F_i/\partial Z_j)$
is a unit in $C[Z_1,\ldots,Z_n]/(F_1,\ldots,F_n)$. Furthermore, localizing
a ring $C$ at an element $a$ is the same thing as $C[t]/(at-1)$, so the same
result holds if we further localize $C[Z_1,\ldots,Z_n]$
before dividing out by $(F_1,\ldots,F_n)$ with $F_i$ in the localized ring.
This is our situation, and in this case the Jacobian is
\[
\begin{pmatrix}
-\Diag({\mathbf h}')-(\partial {\mathbf h}'/\partial {\mathbf S})
\cdot {\mathbf S}&
-(\partial {\mathbf h'}/\partial {\mathbf z})\cdot {\mathbf S}\\
-\partial{\mathbf f}/\partial{\mathbf S}&-\partial{\mathbf f}/\partial
{\mathbf z}
\end{pmatrix}
\]
As $S_k=0$ for all
$k$ at $\bar v$, we see that at $\bar v$ this matrix is
\[
\begin{pmatrix}
-\Diag({\mathbf h}')&{\mathbf 0}\\
-\partial{\mathbf f}/\partial{\mathbf S}&-\partial{\mathbf f}/\partial
{\mathbf z}
\end{pmatrix}
\]
By our choice of the $f_i$'s, the determinant of this is non-zero
at $\bar v$. Thus, possibly further shrinking $V$, we get an \'etale
morphism. \qed

\smallskip

We now only need to show the log structure $X^{\dagger}\rightarrow
\Spec \kk^{\dagger}$ induced by $X\subseteq Y$ pulls back to the given
structure $V^{\dagger}\rightarrow\Spec \kk^{\dagger}$. Since these structures
are determined by sections of $\shLS^+_{\pre,X}$ and $\shLS^+_{\pre,V}$
respectively it is actually enough to show this on
an open subset $X^o\subseteq X$ which is dense on each codimension one
stratum of $X$. We define $Y^o$ by localizing $Y$ at $\prod_{i=1}^q u_i$, and
similarly define $X^o$, so $X^o=X\cap Y^o$. The diagram
\[
\begin{matrix}
Y&\hookleftarrow&Y^o\\
\mapup{}&&\mapup{}\\
X&\hookleftarrow&X^o
\end{matrix}
\]
corresponds to the diagram of rings 
\[
\begin{matrix}
\kk[P]&\hookrightarrow&\kk[P'\oplus\ZZ^q]\\
\mapdown{}&&\mapdown{}\\
\kk[Q]&\hookrightarrow&\kk[\partial P'\oplus \ZZ^q]
\end{matrix}
\]
Here the horizontal inclusions are just
\[
(n,a_0,\ldots,a_q)\mapsto (n,a_0,\ldots,a_q)
\]
and the vertical maps are quotients by $\rho$. 
The chart on $X^o$ for the log structure given by the inclusion
$X^o\rightarrow Y^o$ is clearly just
\[
P'\ni p\mapsto z^p.
\]
However, the isomorphism of $\Spec \kk[\partial P'
\oplus\NN^q]$ with $\Spec \kk[Q]$ is given by the map of rings
$\kk[Q]\rightarrow \kk[\partial P'\oplus\NN^q]$ given by
\[
z^p\mapsto z^p\prod_{i=1}^q u_i^{-\check\psi_i(p)},
\]
as described in Construction \ref{keyP}. Thus the chart for the log
structure on $X^o$ \emph{as a subset of} $\Spec \kk[\partial P'\oplus\NN^q]$
is
\[
P'\ni p\mapsto z^p\prod_{i=1}^q  u_i^{-\check\psi_i(p)}.
\]
We then pull back this chart to $V$, getting 
\[
P'\ni p\mapsto h_pz^p\prod_{i=1}^q f_i^{-\check\psi_i(p)}=h_pz^p\prod_{i=1}^r
f_i^{-\check\psi_i(p)},
\]
the latter equality since $\check\psi_i(p)=0$ for $i>r$.
Set $f_p=h_p\prod_{i=1}^r f_i^{-\check\psi_i(p)}$, $f'_p=f_p/h_p$.
We now have to compare this chart with the log structure on $V$.
We do this by calculating the section of $\shLS_V$ determined by this
chart, using the construction of \cite{PartI}, Theorem~3.22, which associates
to the data $(f_p)_{p\in \partial P'}$ a section of $\shLS_V$. As this
map is multiplicative, and by construction $f_p=f'_ph_p$ and 
$(h_p)_{p\in \partial P'}$ maps to the section $a_{\sigma}$
of $\shLS_V$, it is enough to show that $(f'_p)_{p\in\partial P'}$ maps
to the section $f'_{\sigma}$ of $\shLS_V$.

Consider $e:\omega\rightarrow\tau$ with $\omega$ one-dimensional, 
with endpoints $e^{\pm}_{\omega}:v^{\pm}_{\omega}\rightarrow\omega$ as usual.
On the maximal cones $\check v_{\omega}^{\pm}$
of $\check\Sigma_{\tau}$, suppose $\check\psi_i$
is given by $-m_i^{\pm}$. By definition of $\check\psi_i$, $m_i^{\pm}=
\Vert_i(v^{\pm}_{\omega})$
are vertices of $\Delta_i$. In order to show the construction
of \cite{PartI}, Theorem 3.22 gives the section $f'_{\sigma}$, we have to show
for $n\in N'$ we have on $V_e$
\[
(d_{\omega}\otimes f'_{\sigma,e})(n)={\prod_{i=1}^r f_i^{\langle m_i^-,n\rangle}
\over 
\prod_{i=1}^r f_i^{\langle m_i^+,n\rangle}}
=\prod_{i=1}^r f_i^{\langle m_i^--m_i^+,n\rangle}.
\]
But $m_i^--m_i^+$ is always proportional to $d_{\omega}$, and by the fact that
$T_{\Delta_1},\ldots,T_{\Delta_r}$ are transversal (as follows
from the simplicity assumption, \cite{PartI}, Remark 1.61, (4)) and the fact
that $\Delta_i$ are elementary simplices, we have
\[
m_i^--m_i^+=\begin{cases}
0&e\not\in\Omega_i,\\
d_{\omega}&e\in\Omega_i.
\end{cases}
\]
Thus on $V_e$
\[
\prod_{i=1}^r f_i^{\langle m_i^--m_i^+,n\rangle}=f_i^{\langle d_{\omega},
n\rangle}
=(f'_{\sigma,e})^{\langle d_{\omega},n\rangle}=(d_{\omega}\otimes f'_{\sigma,e})(n)
\]
if $e\in\Omega_i$, as desired; if there is no such $i$, then this is $1$,
again as desired. \qed

\subsection{Deformation theory}

We will now develop deformation theory for log Calabi-Yau spaces
$f:X^{\dagger}\rightarrow \Spec \kk^{\dagger}$. This is foundational
material, and it is possible that there is a much more general context
in which to do deformation theory of relatively coherent log structures.
Here, however, we restrict attention to log structures whose local models
have been described in the previous section. We note, however, that the
results of this section, while of interest by themselves, are only needed
in order to be able to apply the results of Theorem 
\ref{hypercohomology base change} to the degenerations of Calabi-Yau
manifolds constructed in \cite{Smoothing}. In particular, Corollary 
\ref{smoothingisdivisorial} shows that those degenerations are of
the correct sort to apply the results of \S 4. So this section is not
necessary for the understanding of the remainder of the paper.
 
The dificulty is that
as $f$ is not log smooth along $Z$, we can't apply standard log smooth
deformation theory as in
\cite{F.Kato},\cite{K.Kato} immediately. As an abstract log space over $\Spec \kk^{\dagger}$, 
$X^{\dagger}$ has many deformations which are fairly perverse and which
we don't wish to allow. So we consider the following restricted deformation
problem. 

Fix once and for all the $\kk$-algebra $R:=\kk\lfor \NN\rfor=\kk\lfor t\rfor$.
$\Spec R$ carries the log structure induced by the natural chart
$\NN\rightarrow R$, $n\rightarrow t^n$. We write this log scheme
as $\Spec R^{\dagger}$. Thus for any $R$-algebra $A$,
$\Spec A$ carries the pull-back log structure, which we write as 
$\Spec A^{\dagger}$. We denote by $\C_{R}$ the category of
Artin local $R$-algebras with residue field $\kk$. As usual,
if $A,B\in Ob(\C_{R})$ with $B=A/I$, we say $A$ is a 
\emph{small extension} of $B$ if $m_AI=0$.

\begin{definition}
\label{divisorialdef}
Let $\X_\kk^{\dagger}$ be a toric log Calabi-Yau space over $\Spec \kk^{\dagger}$, 
with positive and simple dual intersection complex $(B,\P)$, and let
$A\in Ob(\C_{R})$. Then a \emph{divisorial 
log deformation} of $\X_\kk^{\dagger}$ over $\Spec A^{\dagger}$ is data
$f_A:\X_A^{\dagger}\rightarrow \Spec A^{\dagger}$ together with an isomorphism
$\X_A^{\dagger}\times_{\Spec A^{\dagger}} \Spec \kk^{\dagger}\cong
\X_{\kk}^{\dagger}$ over $\Spec \kk^{\dagger}$ such that
\begin{enumerate}
\item $f_A$ is flat as a morphism of schemes, and $f_A|_{\X_A\setminus Z}$
is log smooth.
\item For every closed geometric point $\bar x\in Z$, let $P,Y,X$ be the 
data of Theorem \ref{localmodel} giving a diagram (\ref{etalediagram})
over $\Spec \kk^{\dagger}$. Let $X_A^{\dagger}=Y^{\dagger}\times_{\Spec\kk[\NN]^{\dagger}}
\Spec A^{\dagger}$. Then there exists a diagram over $\Spec A^{\dagger}$
\begin{equation}
\label{etalediagram2}
\xymatrix{
&\V_A^{\dagger}\ar[ld]\ar[rd]^{\phi'}&\\
\X_A^{\dagger}&&X_A^{\dagger}
}
\end{equation}
with both maps strict \'etale.
\end{enumerate}
\end{definition}

\begin{example} 
\label{standarddivisorial}
The standard behaviour of a divisorial deformation in codimension
two can be described as follows. Given $B,\P$ and $s$ as usual,
let $\omega\in\P$ be an edge of affine length~$l$. Let
$\shQ_{\omega}=\Lambda_y/\Lambda_{\omega}$ for any $y\in
\Int(\omega)\setminus\Delta$, as in \cite{PartI}, Definition~1.38.
Then $\omega$ defines an \'etale open neighbourhood $V(\omega)$
of $X_0(B,\P,s)$, which can be written as 
\[
V(\omega)\cong\Spec\kk[\shQ_{\omega}^{\vee}][x,y]/(xy).
\]
The log structure and log morphism to $\Spec \kk^{\dagger}$ is determined
by $f_{\omega}\in\Gamma(V(\omega),\shLS^+_{\pre,V(\omega)})=\kk[\shQ_{\omega}^{\vee}]$,
and it is not difficult to describe the log structure:
We have a morphism
\[
U=\Spec\kk[\shQ_{\omega}^{\vee}][x,y,t]/(xy-f_{\omega}t^l)\rightarrow \Spec\kk[t]
=\AA^1
\]
given by $t\mapsto t$, and $V(\omega)$ is the fibre over the origin.
The inclusion $V(\omega)\subseteq U$ induces a log structure on $U$,
with $\M_U$ given as usual by the regular functions invertible away from
$V(\omega)$. The restriction of this log structure to $V(\omega)$
is the desired log structure. If $\AA^1$ is given the log structure induced
by $0\in \AA^1$ (or equivalently by the chart $\NN\mapsto \kk[t]$,
$n\mapsto t^n$), and $\Spec\kk[t]/(t^{k+1})$ is given a log structure
via the same chart, then
\[
U^{\dagger}\times_{(\AA^1)^{\dagger}} \Spec \kk[t]/(t^{k+1})^{\dagger}
\rightarrow
\Spec \kk[t]/(t^{k+1})^{\dagger}
\]
gives a $k$th order divisorial deformation of $V(\omega)^{\dagger}
\rightarrow\Spec\kk^{\dagger}$ (given by $k=0$). To check this,
it is enough to check that $U\rightarrow\AA^1$ is \'etale locally
equivalent to $Y\rightarrow\AA^1$ as given in Theorem~\ref{localmodel},
which is easily done as in the proof of that Theorem.
\end{example}

\begin{example}
Let us return to Example \ref{ODPexample}, and show there
are many non-divisorial deformations of $X^{\dagger}/\kk^{\dagger}$.
Consider $X[\epsilon]=\Spec\kk[P][\epsilon]/(z^{\rho},\epsilon^2)$
as a scheme over $\Spec\kk[\epsilon]/(\epsilon^2)$. Now the log
structure on $X\setminus\{0\}$ is fine and can be described
using charts on open sets $U_x$, $U_y$ and $U_w$ of $X$, the sets
where $x\not=0$, $y\not=0$ and $w\not=0$ respectively. Then the charts
$\NN\rightarrow\O_{U_x}$ and $\NN\rightarrow\O_{U_y}$ given by
\[
n\mapsto\begin{cases} 1 & n=0\\
0&n\not=0
\end{cases}
\]
and $\NN^2\rightarrow\O_{U_w}$ given by
\[
(n_1,n_2)\mapsto w^{-n_1}x^{n_1}y^{n_2}
\]
are easily seen to define isomorphic log structures on overlaps. Furthermore,
these charts also define the log morphisms to $\Spec \kk^{\dagger}$, via
$\NN \ni 1\mapsto 1\in\NN$ for the first two charts and $\NN\ni 1\mapsto 
(1,1)\in\NN^2$ for the third chart. There are unique identifications
of the three log structures on overlaps which identify these morphisms,
and these identifications
define the log structure $X^{\dagger}$ on $X\setminus\{0\}$.
Furthermore $\M_X=i_*\M_{X\setminus\{0\}}$, where $i:X\setminus\{0\}
\hookrightarrow X$ is the inclusion. We have now described the log morphism
$X^{\dagger}\rightarrow\Spec\kk^{\dagger}$ in terms of charts.

We can then lift this log structure and morphism to $X[\epsilon]$ by keeping the
same charts on $U_x$ and $U_y$, but lifting $\NN^2\rightarrow \O_{U_w}$
to $\NN^2\rightarrow\O_{U_w}[\epsilon]$ given by
\[
(n_1,n_2)\mapsto (w+f(w,w^{-1})\epsilon)^{-n_1}x^{n_1}y^{n_2},
\]
for $f(w,w^{-1})\in k[w,w^{-1}]$. Furthermore, 
$\Spec \kk[\epsilon]/(\epsilon^2)$ has a log structure 
induced by the $R$-algebra structure $t\mapsto 0$.
The same maps of monoids $\NN\rightarrow\NN$ and $\NN\rightarrow\NN^2$
as above then induce log morphisms $U_x[\epsilon]^{\dagger},U_y[\epsilon]^{\dagger},U_w[\epsilon]^{\dagger}\rightarrow
\Spec\kk[\epsilon]/(\epsilon^2)^{\dagger}$.
These log structures and morphisms glue canonically. 
Note that applying an automorphism
of $X[\epsilon]$ fixing $X$ can remove any positive powers of $w$ from
$f$, but not negative powers. Thus this gives an infinite dimensional
family of possible ``exotic'' deformations of the log structure. It
follows from Lemma \ref{claim} below that none of these are divisorial.
\qed
\end{example}

\begin{lemma}
\label{logauto}
Let $A'\in Ob(\C_{R})$ and let $\X_{A'}^{\dagger}\rightarrow
\Spec A'^{\dagger}$ be a morphism of log schemes which
is flat as a morphism of schemes.
If $A\in Ob(\C_{R})$ with a surjective homomorphism $A'\rightarrow A$
and kernel $I$ with $I^2=0$, then the set of log
automorphisms of $\X_{A'}^{\dagger}$ fixing $\X_A^{\dagger}=\X_{A'}^{\dagger}
\times_{(A')^{\dagger}} A^{\dagger}$ is in one-to-one
correspondence with the log derivations of $\X_{A'}^{\dagger}$ 
over $\Spec (A')^{\dagger}$
with values in $\O_{\X_{A'}}\otimes_{A'} I$.

If $A'$ is a small extension of $A$, $\X_{A'}^{\dagger}\rightarrow\Spec
{A'}^{\dagger}$ log smooth off of a codimension two closed subset $Z$
of $\X_{A'}$, $\X_{A'}$ has no embedded components and is $S_2$, and
$\M_{\X_{A'}}$ has no sections with support in $Z$,
then this space of log derivations is
$\Gamma(\X_\kk,\Theta_{\X^{\dagger}_{\kk}/\kk^{\dagger}}\otimes I)$, 
where $\X_\kk^{\dagger}
=\X_{A'}^{\dagger}\times_{{A'}^{\dagger}} \kk^{\dagger}$.
\end{lemma}

\proof This is standard. An automorphism is induced by an automorphism
of sheaves of rings $\psi:\O_{\X_{A'}}\rightarrow \O_{\X_{A'}}$
and an automorphism $\Psi:\M_{\X_{A'}}\rightarrow\M_{\X_{A'}}$ compatible
with $\psi$, which are the identity
after restricting to $\X_A$. Then 
$x\mapsto \Di x=\psi(x)-x\in \O_{\X_{A'}}\otimes_{A'} I$ is an ordinary derivation, while
if $m\in \M_{\X_{A'}}$, then $\Psi(m)=h\cdot m$ for $h\in 1+\O_{\X_{A'}}
\otimes_{A'} I$ as $\Psi$ induces the identity on $\overline{\M}_{\X_A}
=\overline{\M}_{\X_{A'}}$.
Then we define $\Dlog(m)=h-1$. It is easy to check $(\Di,\Dlog)$ yields 
a log derivation.
Conversely, given a log derivation $(\Di,\Dlog)$, one can reverse this procedure
to get an automorphism.

In the case of a small extension, we notice that $I$ is a $\kk=A'/m_{A'}$-vector
space, and $\O_{\X_{A'}}\otimes_{A'} I=\O_{\X_\kk}\otimes_\kk I$. 
Off of the codimension two 
locus $Z$ where $\X_A^{\dagger}\rightarrow\Spec A^{\dagger}$ fails to be log smooth, 
this module of derivations is $\Theta_{\X_A^{\dagger}/A^{\dagger}}\otimes_{\O_{\X_A}} 
(\O_{\X_\kk}\otimes_\kk I)=\Theta_{\X_\kk^{\dagger}/\kk^{\dagger}}\otimes_\kk I$. By Remark \ref{Thetaextend}, 
this shows the module of derivations on $\X_A^{\dagger}$ 
with values in $\O_{\X_{A'}}\otimes_{A'} I$ is
in fact $\Theta_{\X_\kk^{\dagger}/\kk^{\dagger}}\otimes_\kk I$.
\qed

\medskip

The main theorem for divisorial deformations is then the following.

\begin{theorem}
\label{maindeftheorem}
Let $(B,\P)$ be an integral affine manifold with singularities
and polyhedral decomposition $\P$. Suppose $(B,\P)$ is positive and simple, 
and let
$\X_\kk^{\dagger}$ 
be a toric log Calabi-Yau space with dual intersection
complex $(B,\P)$ and log singular set $\shZ\subseteq\X_\kk$. 
Let $A',A\in Ob(\C_{R})$ with $A'$ a 
small extension of $A$, $A=A'/I$.
Let $\X_A^{\dagger}\rightarrow\Spec A^{\dagger}$ be a divisorial
deformation of $\X_\kk^{\dagger}\rightarrow\Spec
\kk^{\dagger}$. Then
\begin{enumerate}
\item Suppose this deformation lifts to a divisorial deformation $\X_{A'}^{\dagger}
\rightarrow\Spec {A'}^{\dagger}$. Then the set of log automorphisms of
$\X_{A'}^{\dagger}$ over $\Spec {A'}^{\dagger}$ which is the identity on
$\X_A^{\dagger}$ is 
\[
H^0(\X_\kk,\Theta_{\X_\kk^{\dagger}/\kk^{\dagger}}\otimes_\kk I).
\]
\item Under the same assumption as in (1), the set of equivalence classes
of liftings $\X_{A'}^{\dagger}\rightarrow\Spec {A'}^{\dagger}$
of $\X_A^{\dagger}\rightarrow\Spec A^{\dagger}$ is a torsor over
\[
H^1(\X_\kk,\Theta_{\X_\kk^{\dagger}/\kk^{\dagger}}\otimes_\kk I).
\]
Here two liftings are equivalent if there is an isomorphism between them
over $\Spec (A')^{\dagger}$ which is the identity on $\X_A^{\dagger}$.
\item The obstruction to the existence of a lifting as in (1) is in 
\[
H^2(\X_\kk,\Theta_{\X_\kk^{\dagger}/\kk^{\dagger}}\otimes_\kk I).
\]
\end{enumerate}
\end{theorem}

For the proof of Theorem \ref{maindeftheorem}, 
we will need the following lemmas concerning derivations arising in 
Construction \ref{keyP}. The first will be needed to deal with
general points of $\shZ$, and the last two will deal with higher
codimensions.

\begin{lemma}
\label{resolution2}
Let $M'=\ZZ$, $\Delta_0=[0,l]$ for $l>0$ an integer,
and $\Delta_1=[0,1]$ determine data $P$, $\rho\in P$, as in Construction
\ref{keyP}, so that $Y$ is the toric variety defined by the fan in
$M_{\RR}$ given by the cone over the convex hull of $(0,0,1)$, $(1,0,1)$,
$(0,1,0)$ and $(l,1,0)$.
Subdivide this cone as depicted:
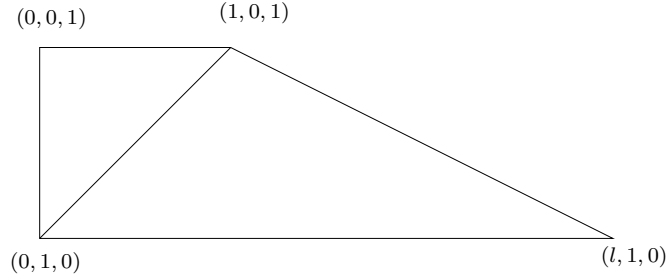
\begin{center}
\begin{picture}(0,0)%
\includegraphics{polygon2.pstex}%
\end{picture}%
\setlength{\unitlength}{1973sp}%
\begingroup\makeatletter\ifx\SetFigFont\undefined%
\gdef\SetFigFont#1#2#3#4#5{%
  \reset@font\fontsize{#1}{#2pt}%
  \fontfamily{#3}\fontseries{#4}\fontshape{#5}%
  \selectfont}%
\fi\endgroup%
\begin{picture}(8271,3397)(3226,-4403)
\put(3301,-1261){\makebox(0,0)[lb]{\smash{{\SetFigFont{8}{9.6}{\rmdefault}{\mddefault}{\updefault}{\color[rgb]{0,0,0}$(0,0,1)$}%
}}}}
\put(5851,-1186){\makebox(0,0)[lb]{\smash{{\SetFigFont{8}{9.6}{\rmdefault}{\mddefault}{\updefault}{\color[rgb]{0,0,0}$(1,0,1)$}%
}}}}
\put(3226,-4336){\makebox(0,0)[lb]{\smash{{\SetFigFont{8}{9.6}{\rmdefault}{\mddefault}{\updefault}{\color[rgb]{0,0,0}$(0,1,0)$}%
}}}}
\put(10651,-4261){\makebox(0,0)[lb]{\smash{{\SetFigFont{8}{9.6}{\rmdefault}{\mddefault}{\updefault}{\color[rgb]{0,0,0}$(l,1,0)$}%
}}}}
\end{picture}%

\end{center}
to define a blow-up $\tilde Y\rightarrow Y$. Let $\pi:\tilde X\rightarrow
X$ be the restriction of this map to the inverse image of $X$, and let
$E\subseteq\tilde X$ be the exceptional locus. Let $\tilde X^{\dagger}$
be the log structure on $\tilde X$ induced by the inclusion $\tilde
X\subseteq\tilde Y$. Then
\begin{enumerate}
\item the composed morphism $\tilde X^{\dagger}\rightarrow X^{\dagger}
\rightarrow \Spec \kk^{\dagger}$ is log-smooth, and $\tilde X^{\dagger}\setminus
E\cong X^{\dagger}\setminus Z$.
\item $\pi_*\Theta_{\tilde X^{\dagger}/\kk^{\dagger}}\cong
\Theta_{X^{\dagger}/\kk^{\dagger}}$.
\item
$R^p\pi_*\Theta_{\tilde X^{\dagger}/\kk^{\dagger}}=0$ for $p>0$.
\end{enumerate}
\end{lemma}

\proof Let $\tilde Y_1$ be the affine open subset of $\tilde Y$ determined 
by the cone generated by \[(0,1,0), (0,0,1), (1,0,1),\] so that 
\[
\tilde  Y_1=\Spec \kk[z^{(1,0,0)}, z^{(0,1,0)}, z^{(-1,0,1)}]=\AA^3_\kk,
\]
and let $\tilde Y_2$ be determined by the cone generated by $(0,1,0), (1,0,1),
(l,1,0)$, so that
\[
\tilde  Y_2=\Spec \kk[z^{(0,0,1)}, z^{(1,0,-1)}, z^{(-1,l,1)}, z^{(0,1,0)}]
=\AA^1_\kk\times\Spec \kk[x,y,t]/(xy-t^l).
\]
Recall $t=z^{e_0^*}=z^{(0,1,0)}$.
If $\tilde X_i=\tilde X\cap \tilde Y_i$, we have
$\tilde X_1=\{z^{(0,1,0)}=0\}\subseteq \AA^3_\kk$
and $\tilde X_2=\AA_\kk^1\times
\{t=0\}\subseteq\tilde Y_2$. If we give $\tilde X_i$ the log structure
induced by the corresponding inclusion, we see $\tilde X_i$ is clearly 
log smooth over $\Spec \kk^{\dagger}$. 
The second statement of (1) is obvious since $\tilde Y\setminus E=
Y\setminus Z$.

For (2) and (3), we use the open cover $\{\tilde X_1,\tilde X_2\}$
to compute \v Cech cohomology of $\Theta_{\tilde X^{\dagger}/\kk^{\dagger}}$.
Note $\Gamma(\tilde X_1, \Theta_{\tilde X^{\dagger}/\kk^{\dagger}})$
is the free $A_1=\kk[z^{(1,0,0)},z^{(-1,0,1)}]$-module generated by
$z^{(1,0,-1)}\partial_{(0,0,1)}$ and $z^{(-1,0,0)}\partial_{(1,0,1)}$
from Proposition \ref{ThetaXl}, while $\Gamma(\tilde X_2,
\Theta_{\tilde X^{\dagger}/\kk^{\dagger}})$ is the free
\[A_2=
\kk[z^{(0,0,1)}, z^{(1,0,-1)}, z^{(-1,l,1)}, z^{(0,1,0)}]/(z^{(0,1,0)})\]
module generated by $z^{(0,0,-1)}\partial_{(1,0,1)}$ and $\partial_{(0,0,1)}$.
In addition $\Gamma(\tilde X_1\cap \tilde X_2,
\Theta_{\tilde X^{\dagger}/\kk^{\dagger}})$ is the free
$A_{12}=\kk[z^{(-1,0,1)},z^{(1,0,-1)},z^{(1,0,0)}]$-module generated by
$\partial_{(0,0,1)}$ and $z^{(-1,0,0)}\partial_{(1,0,1)}$. From this one
calculates easily that the kernel of the \v Cech coboundary map
\[
\Gamma(\tilde X_1,\Theta_{\tilde X^{\dagger}/\kk^{\dagger}})
\oplus
\Gamma(\tilde X_2,\Theta_{\tilde X^{\dagger}/\kk^{\dagger}})
\rightarrow
\Gamma(\tilde X_1\cap\tilde X_2,\Theta_{\tilde X^{\dagger}/\kk^{\dagger}})
\]
is $\Gamma(X,\Theta_{X^{\dagger}/\kk^{\dagger}})$ using the description
of Proposition \ref{ThetaXl} (which tells us $\Gamma(X,
\Theta_{X^{\dagger}/\kk^{\dagger}})$ is generated by $\partial_{(0,0,1)}$, 
$\partial_{(1,0,0)}$, $z^{(1,0,-1)}\partial_{(0,0,1)}$ and $z^{(-1,l,0)}
\partial_{(1,0,1)}$),
and the cokernel of this map is zero,
proving (2) and (3). \qed

\begin{lemma}
\label{lsdef}
In the situation of Construction \ref{keyP},
\begin{enumerate}
\item the (ordinary) sheaf of derivations
$\Theta_{V'(\tau)/\kk}$ is isomorphic to 
the sheaf whose sections over an open subset 
$V\subseteq V'(\tau)$ are
\[
\{(h_p)_{p\in\partial P'}|\hbox{
$h_p\in\Gamma(V,\O_{V(p)})$ and $h_p+h_q=h_{p+q}$ on $V\cap V(p+q)$}\}.
\]
Here $V(p)=\cl(V'(\tau)\setminus\{z^p=0\})$.
\item
Let $\lss_{V'(\tau)}$ denote the cokernel of the inclusion
\[
M'\otimes \O_{V'(\tau)}\rightarrow \Theta_{V'(\tau)/\kk}
\]
given by $M'\ni m\mapsto (\langle m,p\rangle)_{p\in\partial P'}$. 
Then $\Gamma(V,\lss_{V'(\tau)})$ consists of
\[
(f_{\omega})_{\omega}\in \bigoplus_{\dim\omega=1,\omega\subseteq\tau}
\Gamma(V,\O_{V'_{\omega}})
\]
such that for every two-dimensional
face $\tau'$ of $\tau$, 
\[
\sum_{\dim\omega=1} \epsilon_{\tau'}(\omega) f_{\omega}|_{V'_{\tau'}} d_{\omega}
=0\in \Gamma(V,M'\otimes\O_{V'_{\tau'}}).
\]
Here $\epsilon_{\tau'}$ is a choice of sign vector for $\tau'$
(see \cite{PartI}, Definition 3.21), 
$d_{\omega}$ is a choice of primitive generator of the tangent
space of an edge $\omega\subseteq\tau$, and $V'_{\omega}$ is the toric
stratum of $V'(\tau)$ corresponding to $\omega\subseteq\tau$.
We think of $\lss_{V'(\tau)}$ as the Lie algebra of $\shLS_{V'(\tau)}$.
\end{enumerate}
\end{lemma}

\proof
(1) The data $(h_p)_{p\in\partial P'}$ on an open set $V$ defines a derivation
$\Di:\O_V\rightarrow\O_V$ by $\Di(z^p)=h_p z^p$ for $p\in \partial P'$.
This is immediately checked to be a derivation, and it gives an injective
map from the sheaf of such data to $\Theta_{V(\tau)/\kk}$. Conversely, if
$\Di$ is a derivation on an open set $V$, define for each $p\in
\partial P'$ with $V\cap V(p)\not=\emptyset$,
\[
h_p={\Di z^p\over z^p}.
\]
This is clearly a well-defined function on $V'(\tau)\setminus\{z^p=0\}$,
and we need to show it extends to a function on $V(p)$. To show this, we first
note we can restrict the derivation $\Di$ to irreducible components of $V$.
Indeed, assuming $V$ is affine and $W$ is an irreducible component of $V$,
let $f$ be a function on $W$. We can extend it to a function $f'$ on
$V$, hence compute $(\Di f')|_W$. This is independent of the lifting:
if $f', f''$ are two different liftings of $f$, then $\Di(f'-f'')|_W=0$. Now
if we denote by $\partial W$ the toric boundary of $W$ (which is
$\bigcup_{W'} (W'\cap W)$ where $W'$ runs over all irreducible components
of $V$ not equal to $W$), then $\Di$ must preserve the ideal of $\partial W$
in $W$ (extend a function in this ideal by zero off of $W$). Thus
$\Di$ is a log derivation for the pair $(W,\partial W)$ (see Example \ref{preserveideal}).
Thus if $W$ is in particular an irreducible component
of $V(p)\cap V$, it then follows that on $W$, $\Dlog p={\Di z^p\over z^p}$ is
a regular function. Thus $\Di z^p \over z^p$ is a regular function on
$V(p)\cap V$. This gives the data $(h_p)_{p\in\partial P'}$ as desired.

For (2), one repeats the proof of \cite{PartI}, Theorem 3.22, word for word,
but everything is additive rather than multiplicative.
\qed
\medskip

\begin{lemma}
\label{lsdiagram}
In the situation of Construction \ref{keyP}, 
let $\pi_1$, $\pi_2$ be the projections of
$X\cong V'(\tau)\times\AA^q$ to $V'(\tau)$ and $\AA^q$ respectively.
\begin{enumerate}
\item
The locus where $\Theta_{X^{\dagger}/\kk^{\dagger}}$ is not locally
free is contained in $Z=\bigcup_{i=1}^q Z_i$. 
\item
Suppose furthermore the convex hull of 
\[
\bigcup_{i=1}^q (\Delta_i+e_i)
\]
in $M'\oplus\ZZ^q$ is an elementary simplex. If $\shB$ is defined by the exact
sequence
\[
0\rightarrow \Theta_{X^{\dagger}/\kk^{\dagger}}\rightarrow
\Theta_{X/\kk}\rightarrow \shB\rightarrow 0,
\]
then $\Gamma(V,\shB)$ 
for any open subset $V\subseteq X$
consists of those
\[
(f_{\omega})\in \bigoplus_{\dim\omega=1,\omega\subseteq\tau}
\Gamma(V,\shL_{\omega}),
\] 
where
\[
\shL_{\omega}=\begin{cases}
\O_{V_{\omega}}\cdot u_i^{-1}&\hbox{if $\omega\in\Omega_i$,}\\
\O_{V_{\omega}}&\hbox{otherwise,}
\end{cases}
\]
such that $(f_{\omega})$ satisfies the condition that
for every two-dimensional
face $\tau'$ of $\tau$, 
\[
\sum_{\dim\omega=1} \epsilon_{\tau'}(\omega) f_{\omega}|_{V_{\tau'}} d_{\omega}
=0\in M'\otimes K(V_{\tau'}).
\]
Here $K(V_{\tau'})$ denotes the function field of $V_{\tau'}$.
\end{enumerate}
\end{lemma}

\proof
We first show we have a commutative diagram of exact sequences
\[
\xymatrix@C=30pt
{
&0\ar[d]&0\ar[d]&0\ar[d]&\\
0\ar[r]&M'\otimes\O_X\ar[r]\ar[d]&\Theta_{X^{\dagger}/\kk^{\dagger}}
\ar[r]\ar[d]&\bigoplus_{i=1}^q \shI_{Z_i/X} \partial/\partial u_i\ar[r]\ar[d]&
0\\
0\ar[r]&\pi_1^*\Theta_{V'(\tau)/\kk}\ar[d]\ar[r]&\Theta_{X/\kk}\ar[r]\ar[d]
&\pi_2^*\Theta_{\AA^q/\kk}\ar[r]\ar[d]&0\\
0\ar[r]&\pi_1^*\lss_{V'(\tau)}\ar[d]\ar[r]&\shB\ar[r]\ar[d]&
\bigoplus_{i=1}^q\O_{Z_i} \partial/\partial u_i\ar[r]\ar[d]&0\\
&0&0&0&}
\]
As $\Theta_{X/\kk}=\pi_1^*\Theta_{V'(\tau)/\kk}\oplus \pi_2^*
\Theta_{\AA^q/\kk}$, and $\pi_1^*\Theta_{V'(\tau)/\kk}$ consists
of those derivations of $\Theta_{X/\kk}$ which are zero on all functions
pulled back from $\AA^q$, one sees easily from the description of
$\Theta_{X^{\dagger}/\kk^{\dagger}}$ of Proposition \ref{ThetaXl} that
$\Theta_{X^{\dagger}/\kk^{\dagger}}\cap \pi_1^*\Theta_{V'(\tau)/\kk}
=M'\otimes\O_X$. Thus $\coker(M'\otimes\O_X\rightarrow
\Theta_{X^{\dagger}/\kk^{\dagger}})$ injects into $\pi_2^*\Theta_{\AA^q/\kk}$,
and the only thing we need to do is to describe the image of this cokernel.
Now $\pi_2^*\Theta_{\AA^q/\kk}$ is generated by the coordinate vector fields
$\partial/\partial u_i$, $i=1,\ldots,q$. To see an element
of $\shI_{Z_i/X} \partial/\partial u_i$ lifts to a logarithmic derivation,
we consider the various generators of $\shI_{Z_i/X}$. Suppose
$z^p\in\shI_{Z_i/X}$. First consider the case that $p=e_i^*$.
Recalling the notation of a log derivation from \S1,
we have $\partial_{e_i}\in \Theta_{X^{\dagger}/\kk^{\dagger}}$. 
Given a derivation $\Di\in\Theta_{X/\kk}$, we can project it into $\pi_2^*\Theta_{\AA^q/\kk}$
by evaluating it on functions pulled back from $\AA^q$, i.e. 
the image of $\Di$ in $\pi_2^*\Theta_{\AA^q/\kk}$
is
\[
\sum_{i=1}^q (\Di u_i)(\partial/\partial u_i).
\]
Thus in particular the image of $\partial_{e_i}$ is $u_i\partial/\partial u_i$,
so $\partial_{e_i}$ provides a lifting of the latter vector field.
Next consider the case that
$p=n+\sum a_je_j^*\in Q\setminus \{\infty\}$ and there is a unique
vertex $w$ of $\Delta_i$ with $\langle n,w\rangle=-\check\psi_i(n)$. There
are two subcases. If $a_i>\check\psi_i(n)$, then $z^p=z^{p'}u_i$ for some $p'\in Q$,
so we are done as before. If $a_i=\check\psi_i(n)$, then we find
that if $v$ is a vertex of $\Delta_j$ for $j\not=i$, then
\[
\langle p-e_i^*,v+e_j\rangle=\langle p,v+e_j\rangle\ge 0,
\]
and if $v\not=w$ is a vertex of $\Delta_i$, then
\[
\langle p-e_i^*,v+e_i\rangle=\langle n,v\rangle+a_i-1
>-\check\psi_i(n)+\check\psi_i(n)-1=-1,
\]
so $\langle p-e_i^*,v+e_i\rangle\ge 0$. Finally,
\[
\langle p-e_i^*,w+e_i\rangle=-1.
\] 
As $w+e_i$ is a primitive generator of an extremal ray of the cone $K$,
Proposition \ref{ThetaXl} tells us that $z^pu_i^{-1}\partial_{w+e_i}
\in\Theta_{X^{\dagger}/\kk^{\dagger}}$. Now the image of this derivation in
$\pi_2^*\Theta_{\AA^q/\kk}$ is $z^p\partial/\partial u_i$,
and hence provides a lifting of $z^p \partial/\partial u_i$.

Conversely, if $z^p\partial_{m}\in \Gamma(X,\Theta_{X^{\dagger}/\kk^{\dagger}})$
with $m=m'+\sum_{j=1}^q b_je_j$, for $m'\in M'$, we consider its image
$z^p\sum b_j u_j\partial/\partial u_j$ 
in $\pi_2^*\Theta_{\AA^q/\kk}$. If $p\in P$, then clearly the latter is in
$\bigoplus \shI_{Z_i/X}\partial/\partial u_i$. Otherwise there is some
unique $i\not=0$ and unique vertex $w$ of $\Delta_i$ such that
$\langle p,w+e_i\rangle=-1$, $\langle p,v+e_j\rangle\ge 0$ for all
vertices $v\not=w$ of $\Delta_j$. Then $z^pu_i
\in\Gamma(X,\O_X)$, and if $p=n+\sum a_je_j^*$, we have $a_j\ge \check\psi_j(n)$
for all $j\not=i$ and $a_i=\check\psi_i(n)-1$. Furthermore, $w$ is the only
vertex of $\Delta_i$ for which $\langle w,n\rangle=-\check\psi_i(n)$. Thus
$z^pu_i$ is in the ideal of $Z_i$, $m$ is proportional to $w+e_i$,
and $z^p\partial_{w+e_i}$ projects to $z^pu_i\partial/\partial u_i$.
This shows the image of $z^p\partial_m$
is in $\bigoplus \shI_{Z_i/X}\partial/\partial u_i$.
In particular, this shows (1).

To prove (2), first note the assumption that the
$\Delta_i$'s yield an elementary simplex
implies that the sets $\Omega_i$ are disjoint, so the description
of the sheaf $\shL_{\omega}$
in the statement makes sense. Call the sheaf claimed to be
isomorphic to $\shB$ in the statement $\shB'$. To show $\shB'\cong\shB$,
let $j:X\setminus Z\hookrightarrow X$ be the inclusion. Since $Z$
is codimension $2$, $\Gamma(X\setminus Z,\Theta_{X^{\dagger}/\kk^{\dagger}})
=\Gamma(X,\Theta_{X^{\dagger}/\kk^{\dagger}})$, so $H^1_Z(X,\Theta_{X^{\dagger}/
\kk^{\dagger}})=0$, while of course $H^0_Z(X,\Theta_{X/\kk})=0$. Thus $H^0_Z(X,\shB)=0$,
and the adjunction map
\[
\shB\rightarrow j_*j^*\shB=j_*(\pi_1^*\lss_{V'(\tau)})
\]
is injective. 
Now $\shB'$ is clearly a subsheaf of $j_*(\pi_1^*\lss_{V'(\tau)})$,
and we want to show it is the image of $\shB$. We first show the
image of $\shB$ is contained in $\shB'$. Note that $\shB$ is generated by
$\pi_1^*\lss_{V'(\tau)}$ and the images of $\partial/\partial u_1,\ldots,
\partial/\partial u_q$. We can determine the images of $\partial/\partial u_i$
in $j_*(\pi_1^*\lss_{V'(\tau)})$ as follows. After localizing at $u_i$,
consider the log vector field $u_i^{-1}\partial_{e_i}$. This splits as a 
section of
$\Theta_{X/\kk}$ as $\Di_i+\partial/\partial u_i$, where $\Di_i$ is
a section of 
$\pi_1^*\Theta_{V'(\tau)/\kk}$. Thus $\partial/\partial u_i$ is congruent to
$-\Di_i$ modulo $\Theta_{X^{\dagger}/\kk^{\dagger}}$ on $X\setminus\{u_i=0\}$.
To determine $\Di_i$, we evaluate
$u_i^{-1}\partial_{e_i}$ on pull-backs of functions on $V'(\tau)$, i.e.
for $p\in\partial P'$, identified with $N'\cup\{\infty\}$, 
\[
\Di_i(z^p\circ\pi_1)=u_i^{-1}\partial_{e_i}(z^p\circ\pi_1)=u_i^{-1}\partial_{e_i}
(z^{p+\sum_{j=1}^q\check\psi_j(p)e_j^*})=u_i^{-1}\check\psi_i(p)
(z^p\circ \pi_1).
\]
Thus $\Di_i$ corresponds, using the description of $\Theta_{V'(\tau)/\kk}$ of
Lemma \ref{lsdef}, to $(h_p)_{p\in\partial P'}$ with $h_p=u_i^{-1}
\check\psi_i(p)$. The image of $\Di_i$ in $(\pi_1^*\lss_{V'(\tau)})|_{X\setminus
\{u_i=0\}}$ is then given by $(f_{\omega})_{\omega}$ with 
$d_{\omega}\otimes f_{\omega}=u_i^{-1}(\Vert_i(v^-_{\omega})-\Vert_i(v^+_{\omega}))$.
This follows from the additive version of the construction in the proof of 
\cite{PartI}, Theorem 3.22.
By the assumption on the $\Delta_i$'s,
if $\Vert_i(v^-_{\omega})-\Vert_i(v^+_{\omega})\not=0$ then it is equal to $d_{\omega}$.
Thus $f_{\omega}=u_i^{-1}$ if $\omega\in\Omega_i$ and $f_{\omega}=0$ otherwise.
So this section $(f_{\omega})_{\omega}$ of $\pi_1^*(\lss_{V'(\tau)})$ in fact
extends to a section of $\shB'$ on $X$. This section is the image of $-\partial/\partial
u_i$ in $j_*j^*\shB$. Thus the image of $\shB$ in $j_*j^*\shB$
is contained in $\shB'$.

Conversely, let $(f_{\omega})_{\omega}$ be a section of $\shB'$. For a given
$\omega$, if $\omega\not\in\Omega_i$ for any $i$, then $f_{\omega}=:f^0_{\omega}$
is a regular function on $V_{\omega}$. Otherwise, if $\omega\in\Omega_i$, 
then we can write $f_{\omega}=f^0_{\omega}+u_i^{-1}f^1_{\omega}$, where
$f^0_{\omega},f^1_{\omega}$ are regular functions on $V_{\omega}$ with
$f^1_{\omega}$ containing no terms
divisible by $u_i$. Then by testing the additive condition
for each two-face $\tau'$, we find that $(f^0_{\omega})_{\omega}\in
\pi_1^*\lss_{V'(\tau)}$. On the other hand, we claim there
is a function $f_i$ on $X$ such that $f^1_{\omega}=f_i|_{V_{\omega}}$ for
all $\omega\in\Omega_i$. Indeed, let
$\tau'\subseteq\tau$ be a two-dimensional face. Then $\tau'$ determines
a face $\tau'_i$ of $\Delta_i$, the convex hull of vertices of $\Delta_i$
representing $-\check\psi_i$ on the maximal cones of $\check\Sigma_{\tau}$
corresponding to vertices of $\tau'$. Since $\Delta_i$ is an elementary 
simplex and $\dim\tau'_i\le 2$, $\tau'_i$ is a standard simplex. If
$\dim\tau'_i=0$, then $f^1_{\omega}=0$ for all $\omega\subseteq\tau'$;
if $\dim\tau'_i=1$, then $\tau'$ has exactly two edges in $\Omega_i$, say 
$\omega_1$ and $\omega_2$, and then the additive condition tells us
that $f^1_{\omega_1}|_{V_{\tau'}}=f^1_{\omega_2}|_{V_{\tau'}}$. If $\dim
\tau'_i=2$, then $\tau'$ has exactly three edges $\omega_1,\omega_2,\omega_3$
contained in $\Omega_i$,
and the additive condition again says $f^1_{\omega_i}|_{V_{\tau'}}$ is
independent of $i=1,2,3$. Thus the $f^1_{\omega}$'s glue for
$\omega\in\Omega_i$, yielding $f_i$.

Lift $(f^0_{\omega})_{\omega}$
to a derivation $\Di\in \pi_1^*\Theta_{V'(\tau)/\kk}$. Then
$(f_{\omega})_{\omega}$ is the image of 
$\Di-\sum_{i=1}^q f_i\partial/\partial 
u_i$, showing $\shB=\shB'$. \qed

\bigskip

\emph{Proof of Theorem \ref{maindeftheorem}}.  
This result follows from standard methods of deformation theory
(see for example \cite{SGAI}, Expos\'e III)
and Lemma \ref{logauto} once we prove a local uniqueness
statement for divisorial deformations:

\begin{lemma}
\label{claim}
For every closed geometric point $\bar x\in \X_\kk$, there is
a strict \'etale open neighbourhood $\V_\kk^{\dagger}$ of $\bar x$ 
such that for any divisorial deformation $\X_A^{\dagger}$ of $\X_\kk^{\dagger}$ over
$\Spec A^{\dagger}$, the induced deformation $\V_A^{\dagger}$ of
$\V_\kk^{\dagger}$ is
independent of the deformation of $\X_\kk$.
\end{lemma}

\proof If $\bar x\not\in \shZ$, then $\X_\kk\rightarrow
\Spec \kk^{\dagger}$ is log smooth at $\bar x$, and this local uniqueness is
proved for log smooth morphisms in \cite{F.Kato} Proposition 8.3 and
\cite{K.Kato}, (3.14). If $\bar x\in \shZ$, 
let $P$ be the monoid of Theorem \ref{localmodel}, 
yielding $Y^{\dagger}$ and $X^{\dagger}\rightarrow\Spec \kk^{\dagger}$ and
strict \'etale morphisms
$\X_\kk^{\dagger}\mapleft{p}\V_\kk^{\dagger}\mapright{\phi} X^{\dagger}$.
Then we will show that for any divisorial deformation $\X_A^{\dagger}\rightarrow\Spec
A^{\dagger}$ of $\X_\kk^{\dagger}\rightarrow \Spec \kk^{\dagger}$, 
the induced deformation $\V_A^{\dagger}\rightarrow\Spec A^{\dagger}$
is also induced by the deformation
$Y^{\dagger}\times_{\Spec
\kk[t]^{\dagger}}\Spec A^{\dagger}\rightarrow \Spec A^{\dagger}$
of $X^{\dagger}\rightarrow \Spec \kk^{\dagger}$.

Let $\tau,\Delta_1,\ldots,\Delta_q$ be the data provided by Theorem \ref{localmodel}.
We first consider the case that $\dim\tau=1$, so $\bar x\in \X_\kk$ is a double 
point of $\X_\kk$. Then we can assume $\tau=[0,l]$ for some integer
$l>0$, and that either $\Delta_1=\cdots=\Delta_q=\{0\}$, in 
which case $\bar x\not\in \shZ$, a contradiction,
or else $\Delta_1=[0,1]$, $\Delta_2=\cdots=\Delta_q=\{0\}$. This data determines
$Y,X$ as in Construction \ref{keyP}. By replacing $\V_\kk$ by an 
affine \'etale neighbourhood, we can use Condition (2) of
Definition \ref{divisorialdef} to reduce the local uniqueness
statement to the following: given strict \'etale maps $\phi,\phi':
\V_\kk^{\dagger}\rightarrow X^{\dagger}$, let $\V_A^{\dagger}$ and
$(\V_A')^{\dagger}$ be the pull-backs of the thickening $X_A^{\dagger}$ of
$X^{\dagger}$ to $\V_\kk$ via $\phi$ and $\phi'$ respectively. We need to show
there is an isomorphism over $\Spec A^{\dagger}$
\[
\psi:\V_A^{\dagger}\rightarrow (\V_A')^{\dagger}
\]
which is the identity on $\V_\kk^{\dagger}$. Now $Y$ can be written as
$Y'\times \AA_\kk^{q-1}$, $X=X'\times\AA_\kk^{q-1}$, where $Y'$, $X'$ are the $Y$
and $X$ of Lemma \ref{resolution2}. Let $\tilde Y'$, $\tilde X'$ be the 
blow-ups of Lemma \ref{resolution2}, and set $\tilde Y=\tilde Y'\times \AA_\kk^{q-1}$,
$\tilde X=\tilde X'\times\AA^{q-1}_\kk$, and let $\tilde\V_\kk=\V_\kk\times_X \tilde X$.
We note here it is irrelevant whether we take this fibred product
using the map $\phi$ or the map $\phi'$: either way, we are simply taking
the normalization of $\V_\kk$, blowing up one of the two components along
the singular set $p^{-1}(\shZ)$ of the log structure, and then regluing
along the proper transform of the conductor locus.

Let $\tilde \V_A^{\dagger}$, $(\tilde \V_A')^{\dagger}$ be
the pull-backs of the thickenings $\tilde Y^{\dagger}\times_{\Spec\kk[\NN]^{\dagger}}
\Spec A^{\dagger}$ via $\tilde\phi,\tilde\phi':\tilde\V_\kk\rightarrow\tilde X$. 
By Lemma \ref{resolution2}, $\tilde\V_A^{\dagger}$ and $(\tilde\V_A')^{\dagger}$
are log smooth over $\Spec A^{\dagger}$, and these are both log smooth
deformations of the same log smooth scheme $\tilde\V_\kk^{\dagger}$. We
can apply F. Kato's log smooth deformation theory here \cite{F.Kato}, and as
$H^1(\tilde\V_\kk,\Theta_{\tilde\V_\kk^{\dagger}/\kk^{\dagger}})=0$
by Lemma \ref{resolution2}, (3), it follows these deformations are
unique, so there is an isomorphism $\tilde\psi:\tilde\V_A^{\dagger}
\rightarrow(\tilde\V_A')^{\dagger}$ which is the identity on
$\tilde\V_\kk^{\dagger}$. This induces an isomorphism $\psi:\V_A^{\dagger}
\rightarrow (\V_A')^{\dagger}$ which is the identity on $\V_\kk^{\dagger}$:
indeed, we get an isomorphism 
\[
\V_A^{\dagger}\setminus Z
=\tilde\V_A^{\dagger}\setminus E\mapright{\tilde\psi}
(\tilde\V_A')^{\dagger}\setminus E=(\V_A')^{\dagger}\setminus Z
\]
which is the identity on $\V_\kk^{\dagger}\setminus Z$, hence extends to an 
isomorphism $\psi:\V_A^{\dagger}\rightarrow(\V_A')^{\dagger}$ because
these spaces are $S_2$.

All that is needed to complete the proof of Theorem~\ref{maindeftheorem}
is the general case, i.e. $\dim\tau>1$.
In this case, we will in fact show a stronger claim:

\begin{proposition}
\label{extraclaim}
Let $\shZ_3$ be the intersection of $\shZ$ with the union of
toric strata of $\X$ of codimension $\ge 2$, so that $\shZ_3$ is
codimension 3 in $\X$. Let $\bar x\in\shZ_3$, and let $\tilde\shZ_3$
be the pull-back of $\shZ_3$ to $\V_{\kk}$. Let 
$\V_A^{\dagger}\rightarrow\Spec A^{\dagger}$ be a flat deformation
of $\V_{\kk}^{\dagger}\rightarrow\Spec\kk^{\dagger}$ such that
$\V_A^{\dagger}\setminus\tilde\shZ_3\rightarrow\Spec A^{\dagger}$
is a divisorial deformation of $\V_{\kk}^{\dagger}\setminus
\tilde\shZ_3\rightarrow\Spec\kk^{\dagger}$. Then $\V_A^{\dagger}
\rightarrow\Spec A^{\dagger}$ is isomorphic to the deformation
induced by $X_A^{\dagger}\rightarrow \Spec A^{\dagger}$. This shows
in particular that $\V_A^{\dagger}\rightarrow\Spec A^{\dagger}$ is
divisorial, and also shows the local uniqueness in the $\dim\tau>1$
case.
\end{proposition}

\proof
We will proceed inductively. Assume $A'$ is a small extension of $A$
by an ideal $I$, and assume the result is true for 
$\V_A^{\dagger}=\V_{A'}^{\dagger}\times_{\Spec {A'}^{\dagger}}\Spec A^{\dagger}$.
Let $(\V'_{A'})^{\dagger}\rightarrow\Spec A'^{\dagger}$ be the pull-back
of $X_{A'}^{\dagger}\rightarrow\Spec A'^{\dagger}$ via 
$\phi:\V_{\kk}\rightarrow X$. We compare
$(\V'_{A'})^{\dagger}$ and $\V_{A'}^{\dagger}$.
We see $\V_{A'}^{\dagger}\setminus\tilde \shZ_3\rightarrow\Spec A'^{\dagger}$
and $(\V'_{A'})^{\dagger}\setminus\tilde \shZ_3\rightarrow\Spec A'^{\dagger}$
are both liftings of the same deformation
$\V_A^{\dagger}\setminus\tilde\shZ_3\rightarrow\Spec A^{\dagger}$.
As we have already shown the necessary local uniqueness result for points
of $\shZ$ not contained in $\shZ_3$, 
the set of all such deformations is a torsor
over $H^1(\V_\kk\setminus \tilde \shZ_3,\Theta_{\V_\kk^{\dagger}/\kk^{\dagger}}
\otimes I)$, 
so we can write
the difference between these two deformations as an element
$\theta\in H^1(\V_\kk\setminus\tilde \shZ_3,\Theta_{\V_\kk^{\dagger}/\kk^{\dagger}}
\otimes I)
\subseteq
\Ext^1_{\O_{\V_\kk\setminus\tilde \shZ_3}}(\Omega^1_{\V_\kk^{\dagger}/\kk^{\dagger}}|_{\V_\kk
\setminus
\tilde \shZ_3}
\otimes I,\O_{\V_\kk\setminus \tilde \shZ_3})$. 
We wish to show $\theta=0$. We first note that
the natural map $\Omega^1_{\V_\kk/\kk}\rightarrow\Omega^1_{\V_\kk^{\dagger}/\kk^{\dagger}}$
induces a map
\[
\psi:
\Ext^1_{\O_{\V_\kk\setminus \tilde \shZ_3}}(\Omega^1_{\V_\kk^\dagger/\kk^\dagger}|_{\V_\kk
\setminus \tilde \shZ_3}
\otimes I,
\O_{\V_\kk\setminus \tilde \shZ_3})
\rightarrow
\Ext^1_{\O_{\V_\kk\setminus\tilde \shZ_3}}(\Omega^1_{\V_\kk/\kk}|_{\V_\kk\setminus\tilde \shZ_3}
\otimes I,
\O_{\V_\kk\setminus\tilde \shZ_3}).
\]
The set of liftings of $\V_A\setminus\tilde \shZ_3\rightarrow\Spec A$ to flat deformations
of schemes over $\Spec A'$ is a torsor over the latter group, by
traditional deformation theory. In addition, $\psi(\theta)$ is then the
difference between the liftings $\V_{A'}\setminus \tilde \shZ_3$ and 
$\V'_{A'}\setminus\tilde  \shZ_3$
of $\V_A\setminus\tilde \shZ_3$ as schemes. Now since $\V_{A'}$ and $\V'_{A'}$ are both
deformations of $\V_A$ as a scheme, we in fact have 
\[
\psi(\theta)\in\im(\Ext^1_{\O_{\V_\kk}}(\Omega^1_{\V_\kk/\kk}
\otimes I,\O_{\V_\kk})
\rightarrow
\Ext^1_{\O_{\V_\kk\setminus\tilde \shZ_3}}(\Omega^1_{\V_\kk/\kk}|_{\V_\kk\setminus\tilde \shZ_3}
\otimes I,
\O_{\V_\kk\setminus\tilde \shZ_3})).
\]
We also have a commutative diagram with exact rows 
coming from local-global $\Ext$
spectral sequences
{\scriptsize
\[
\xymatrix@C=30pt
{H^1(\V_\kk\setminus\tilde \shZ_3,\Theta_{\V_\kk^{\dagger}/\kk^{\dagger}}
\otimes I)\ar[d]
\ar@{^{(}->}[r]&
\Ext^1_{\O_{\V_\kk\setminus\tilde \shZ_3}}(\Omega^1_{\V_\kk^{\dagger}/\kk^{\dagger}}|_{\V_\kk\setminus\tilde \shZ_3}
\otimes I,
\O_{\V_\kk\setminus\tilde \shZ_3})\ar[d]^{\psi}&\\
H^1(\V_\kk\setminus\tilde \shZ_3,\Theta_{\V_\kk/\kk}
\otimes I) \ar@{^{(}->}[r]
&\Ext^1_{\O_{\V_\kk\setminus\tilde \shZ_3}}(\Omega^1_{\V_\kk/\kk}|_{\V_\kk\setminus\tilde \shZ_3}
\otimes I,
\O_{\V_\kk\setminus\tilde \shZ_3})\ar[r]&H^0(\V_\kk\setminus\tilde \shZ_3,
\shExt^1(\Omega^1_{\V_\kk/\kk}
\otimes I,\O_{\V_\kk}))\\
H^1(\V_\kk,\Theta_{\V_\kk/\kk}\otimes I)\ar[r]
\ar[u]
\ar@{^{(}->}[r]
&\Ext^1_{\O_{\V_\kk}}(\Omega^1_{\V_\kk/\kk}
\otimes I,
\O_{\V_\kk})\ar[r]
\ar[u]
&H^0(\V_\kk,
\shExt^1(\Omega^1_{\V_\kk/\kk}
\otimes I,\O_{\V_\kk}))
\ar@{^{(}->}[u]
}
\]
}
The lower right-hand vertical arrow is an injection as $\V_\kk$ is 
\'etale locally isomorphic to a (reducible) toric variety, and hence the support
of any section of $\shExt^1(\Omega^1_{\V_\kk/\kk}\otimes I,\O_{\V_\kk})$
must be a toric stratum. By a diagram chase, it is then clear that
$\psi(\theta)$ is in fact in the image of
$H^1(\V_\kk,\Theta_{\V_\kk/\kk}
\otimes I)$. Since we can assume 
$\V_\kk$ is affine, the latter group is zero.
Thus $\V_{A'}$ and $\V'_{A'}$ must be isomorphic as schemes over $\Spec A'$.
So we only need to compare the log structures. In particular, this
shows that $\theta\in\ker(H^1(\V_\kk\setminus\tilde \shZ_3,
\Theta_{\V_\kk^{\dagger}/\kk^{\dagger}}
\otimes I)
\rightarrow H^1(\V_\kk\setminus\tilde \shZ_3,\Theta_{\V_\kk/\kk}
\otimes I))$.
We will show this kernel
is zero, thus showing the two liftings are equivalent deformations.
Indeed, this kernel is isomorphic to the cokernel of 
\[
H^0(\V_\kk\setminus\tilde \shZ_3,
\Theta_{\V_\kk/\kk}
\otimes I)\rightarrow H^0(\V_\kk\setminus \tilde \shZ_3,\phi^*\shB
\otimes I),
\]
where
$\shB$ is defined in Lemma \ref{lsdiagram}. By the description of
$\shB$ there, the fact that $\tilde \shZ_3$ is codimension
two in the codimension one strata of $\V_\kk$, and the fact that these strata
are $S_2$, it follows that any element of $H^0(\V_\kk\setminus\tilde \shZ_3,\shB
\otimes I)$
lifts to an element of $H^0(\V_\kk,\shB
\otimes I)$. However, $H^0(\V_\kk,\Theta_{\V_\kk/\kk}
\otimes I)
\rightarrow H^0(\V_\kk,\shB
\otimes I)$ is surjective as $\V_\kk$ is affine, so the cokernel
of the map of sections over $\V_\kk\setminus\tilde \shZ_3$ is zero, as desired.
\qed

\medskip

Having now proved Theorem~\ref{maindeftheorem}, we note

\begin{corollary}
\label{divisorialcodim2}
With the same hypotheses as in Theorem \ref{maindeftheorem},
if $\X_A^{\dagger}\rightarrow \Spec A^{\dagger}$
is a flat deformation of $\X_{\kk}^{\dagger}\rightarrow\Spec\kk^{\dagger}$
which is divisorial off of $\shZ_3$, the intersection of $\shZ$ with
the union of codimension two strata of $\X_{\kk}$, then $\X_A^{\dagger}
\rightarrow\Spec A^{\dagger}$ is divisorial at all points of $\X_A$.
\end{corollary}

\proof This follows from Proposition \ref{extraclaim}.
\qed

\medskip

Recall now that if there exists a strictly convex multi-valued integral piecewise
linear function $\varphi$ on $(B,\P)$, then the main theorem of
\cite{Smoothing} yields a $k$-th order deformation $\X_k^{\dagger}
\rightarrow\Spec\kk[t]/(t^{k+1})^{\dagger}$ of $X_0(B,\P,s)^{\dagger}
\rightarrow\Spec\kk^{\dagger}$ (see \cite{Smoothing}, Theorem~1.29,
Remark~1.28, and Remark~2.40). To know the results of this paper
apply to this deformation, we need to know this deformation is divisorial.

\begin{corollary}
\label{smoothingisdivisorial}
With the same hypotheses as in
Theorem \ref{maindeftheorem}, the $k$-th order
deformation $\X_k^{\dagger}
\rightarrow \Spec \kk[t]/(t^{k+1})^{\dagger}$ of $X_0(B,\P,s)^{\dagger}
\rightarrow\Spec\kk^{\dagger}$ constructed
in \cite{Smoothing} (given a strictly convex multi-valued integral piecewise
linear function $\varphi$ on $(B,\P)$) is divisorial.
\end{corollary}

\proof By Corollary~\ref{divisorialcodim2} it is enough to check divisoriality
off of $\shZ_3$. \cite{Smoothing}, Proposition~2.32 shows log smoothness
away from $\shZ$, while divisoriality at points of $\shZ\setminus\shZ_3$
will follow from a local argument as in the proof of \cite{Smoothing},
Lemma~2.34. Indeed, given $\omega\in\P$ an edge, a thickening of the
(\'etale) open set $V(\omega)$ of $X_0(B,\P,s)$ given by the
construction of \cite{Smoothing} is
\[
V^k(\omega)=\Spec R_{\cup},
\]
where $R_{\cup}=R_-\times_{R_{\cap}} R_+$ and
\begin{eqnarray*}
R_-&=&\kk[\shQ_{\omega}^{\vee}][x,y,t]/(xy-t^l,y^{\beta}t^{\gamma}|\beta l+\gamma
\ge k+1)\\
R_+&=&\kk[\shQ_{\omega}^{\vee}][x,y,t]/(xy-t^l,x^{\alpha}t^{\gamma}|\alpha l+\gamma
\ge k+1)\\
R_{\cap}&=&\big(\kk[\shQ_{\omega}^{\vee}][x,y,t]/(xy-t^l,x^{\alpha}y^{\beta}t^{\gamma}
|\max\{\alpha,\beta\}l+\gamma\ge k+1)\big)_{f_{\omega}}.
\end{eqnarray*}
Here $\shQ_{\omega}=\Lambda_y/\Lambda_{\omega}$ for any $y\in\Int(\omega)
\setminus\Delta$ (it is worth remembering in \cite{Smoothing} we work
with the Legendre dual $B$, hence we use the dimension one $\omega$ instead
of the codimension one $\rho$ of \cite{Smoothing}, Lemma~2.34) and $l$
is the affine length of $\omega$. Furthermore, $f_{\omega}\in 
\kk[\shQ_{\omega}^{\vee}][t]/(t^{k+1})$ is a function with the property that 
$f_{\omega}\mod t$ is the section of $\shLS^+_{\pre,V(\omega)}$ defining
the log structure on $V(\omega)$. (Compare with 
Example~\ref{standarddivisorial}.) The fibred product is defined using
maps $R_+\rightarrow R_{\cap}$, the canonical quotient map
followed by the localization map, while $R_-\rightarrow R_{\cap}$
twists this by 
\[
x\mapsto f_{\omega}x,\quad y\mapsto f_{\omega}^{-1}y.
\]
In particular, as shown in the proof of \cite{Smoothing}, Lemma~2.34,
\[
R_{\cup}=\kk[\shQ_{\omega}^{\vee}][X,Y,t]/(XY-f_{\omega}t^l,t^{k+1}),
\]
with $X=(x,f_{\omega}x)$, $Y=(f_{\omega}y,y)$. 

In addition, the log
structure on $V^k(\omega)$ is defined by gluing together the two standard
log structures on $\Spec R_{\pm}$; these are given by the canonical
charts $P_l\rightarrow R_{\pm}$ where $P_l$ is the monoid generated
by $p_1,p_2$ and $\rho$ with $p_1+p_2=l\rho$ and $p_1\mapsto x$, $p_2\mapsto
y$, $\rho\mapsto t$. Furthermore, these come with canonical log
morphisms to $\Spec \kk[t]/(t^{k+1})^{\dagger}$ (with log structure
given by the chart $\NN\ni n\mapsto t^n$) given by the monoid homomorphism
$\NN\rightarrow P_l$, $1\mapsto \rho$.

The glued log structure on $V^k(\omega)$ can be described using
three charts, on the open sets $U_X,U_Y,U_{f_{\omega}}\subseteq
V^k(\omega)$ defined by localizing at $X,Y$ and $f_{\omega}$ respectively.
Since $U_X\cong \Spec (R_-)_x$ and $U_Y\cong\Spec (R_+)_y$, the charts
for the log structure on these two open sets are just those given
by the above standard charts $P_l\rightarrow R_{\pm}$. 
To see what the chart on $U_{f_{\omega}}$ is, we modify
the chart on $R_-$ to be $p_1\mapsto f_{\omega}^{-1}x$, $p_2\mapsto
f_{\omega}y$, $\rho\mapsto t$. This makes sense after localizing
at $f_{\omega}$, and does not change the log structure or the log
morphism to $\Spec\kk[t]/(t^{k+1})^{\dagger}$. Now the charts $P_l
\rightarrow R_{\pm}$ glue, giving $p_1\mapsto f_{\omega}^{-1}X$,
$p_2\mapsto Y$, $\rho\mapsto t$. This describes the log structure
on $V^k(\omega)\setminus\shZ$, along with a log morphism to
$\Spec\kk[t]/(t^{k+1})^{\dagger}$. The log structure on $V^k(\omega)$
is then obtained by pushing forward the sheaf of monoids on $V^k(\omega)
\setminus\shZ$ to $V^k(\omega)$.

It is then not difficult to check that this morphism $V^k(\omega)^{\dagger}
\rightarrow \Spec\kk[t]/(t^{k+1})^{\dagger}$ coincides with that given by
base-change of
$\Spec \kk[\shQ_{\omega}^{\vee}][X,Y,t]/(XY-f_{\omega}t^l)\rightarrow
\Spec \kk[t]$ as in Example~\ref{standarddivisorial}. As in that Example,
it is easy to check this is a divisorial deformation: the fact that
here $f_{\omega}$ also depends on $t$ does not change the structure
of the deformation \'etale locally.
\qed

\begin{remark} When we began this project in 2001, it was our
original hope that there would be an easy Bogomolov-Tian-Todorov
type argument to show smoothability of log Calabi-Yau spaces. This
approach failed, and instead, if one wants to show unobstructedness
of divisorial deformation theory, one needs to appeal to the explicit
smoothing results of \cite{Smoothing}. We sketch an argument
here, though we do not expect this will be important in the
future for understanding mirror symmetry. 
Let $(B,\P)$ be positive and simple.
We assume that $X_0=X_0(B,\P,s)$
comes along with a polarization, so that we can apply the results of
\cite{Smoothing}. We will also make a number of other assumptions we will
include below.

First, there is a formal versal deformation space for divisorial
deformations of $X_0^\ls$. To make this precise we
have to consider the corresponding functor of log Artin rings as in
\cite{FKato}. 
Let $S_0^\ls:=\Spec\kk^\ls$, (with the obvious $R$-algebra structure
on $\kk$). We then look at the functor from $\C_R$ to sets
associating to $A\in \C_R$
the set of
isomorphism classes of divisorial deformations of $X_0^\ls\to S_0^\ls$
over $S^\ls=\Spec A^{\dagger}$. 
It is shown in \cite{FKato} that if the analogues
($\textrm{H}^\textrm{log}_1$)--($\textrm{H}^\textrm{log}_3$) of
Schlessinger's conditions \cite{Schlessinger} hold for a functor of
log Artin rings then this functor has a hull. In view of
Proposition~\ref{extraclaim} and the arguments for the log smooth deformation
functor in \cite{FKato} this is easy to verify for our functor. It is also
proved in \cite{FKato} that the hull can be chosen to pro-represent the
functor if $H^0(X_0,\Theta_{X_0^\ls/S_0^\ls})=0$, that is, if
$X_0^\ls$ is infinitesimally rigid. Again the proof of this statement
easily translates to our functor of divisorial deformations. This
additional assumption is fulfilled in the most interesting cases of
degenerate Calabi-Yau varieties.
In particular, this holds when the holonomy of $B$ is contained in
$\ZZ^n\rtimes\SL_n(\ZZ)$ and $H^0(B,i_*\Lambda\otimes\kk)=0$, by 
Theorem~\ref{SLnZproperties}.
We restrict the further discussion to the situation where this holds.
We will also need one further restriction to guarantee that all deformations
of $X_0^\ls$ are projective. For this we need $H^2(X_0,\O_{X_0})=0$,
or equivalently, $H^2(B,\kk)=0$.

Thus we now have a formal model $\mathscr{X}^\ls\to \mathscr{M}^\ls$
for a moduli space of log spaces containing $X_0^\ls$. One important feature
of log deformation theory is that by construction,
there is a morphism $\mathscr{M}^\ls\to \Spf\kk\lfor t\rfor^\ls$, and
$H^1(X_0,\Theta_{X_0^\ls/S_0^\ls})$ is actually the \emph{relative}
tangent space at the point corresponding to $X_0^{\dagger}$ 
in $\mathscr{M}^\ls$;
indeed, this relative tangent space is just the set of log morphisms
from $\Spec\kk[\epsilon]^\ls$ to $\mathscr{M}$ over $T^{\ls}$, 
with $\Lambda$-algebra structure
given by $t\mapsto 0$. According to Theorem~\ref{maindeftheorem}, 
this is $H^1(X_0,\Theta_{X_0^\ls/S_0^\ls})$.

The assumption that $H^2(X_0,\O_{X_0})=0$ implies that $\mathscr{X}
\rightarrow\mathscr{M}$ carries a relative polarization, and hence
can be algebraized. If $\mathscr{M}=\Spf R$ for some complete
ring $R$, then setting $\M=\Spec R$, we obtain a morphism of
schemes $\X\rightarrow\M$ whose completion along the closed fibre
is $\mathscr{X}\rightarrow\mathscr{M}$. 
Second, by the main result of \cite{Smoothing} and 
Corollary~\ref{smoothingisdivisorial},
there exists an irreducible component of $\M$ such that
if $\eta\in \M$ is the generic point of this component,
the fiber
$\X_\eta$ is a normal variety over $\kappa(\eta)$ with at
worst codimension four Gorenstein quotient singularities,
by Proposition~\ref{smoothcase}.
For such varieties the deformation
theory is unobstructed by \cite{Ran}, 
thus showing that the fibers of $\cl(\eta)\to
\Spec\kk\lfor t\rfor$ have dimension $h^1(\X_\eta,
\Theta_{\X_\eta/\eta})$. (We note this requires a technical
point of openness of versality for logarithmic deformation spaces,
which we have not checked.)

Third, by Theorem~\ref{basechangetheorem}
and the triviality of $\Omega^n_{\mathscr{X}^\ls/
\mathscr{M}^\ls}$, we know
\[
h^1(X_0,\Theta_{X_0^\ls/S_0^\ls}) =
h^1(X_0,\Omega^{n-1}_{X_0^\ls/S_0^\ls})
=h^1(\X_\eta, \Omega^{n-1}_{\X_\eta/\eta})
=h^1(\X_\eta, \Theta_{\X_\eta/\eta}).
\]

Thus $\M$ contains a subspace of relative dimension
$h^1(X_0,\Theta_{X_0^\ls/S_0^\ls})$ over $\Spec\kk\lfor t\rfor$.
Since this is the relative dimension of the Zariski tangent space of
$\M$ at the closed point we conclude that
$\M\to\Spec \kk\lfor t\rfor$ is smooth of this relative
dimension. In particular, \emph{every} infinitesimal divisorial
deformation can be extended to a divisorial deformation over
$\Spec\kk\lfor t\rfor$. This is the claimed unobstructedness result.
\qed
\end{remark}

\section{Cohomology of log Calabi-Yau spaces}

Our aim in this section is to
calculate the logarithmic de Rham cohomology of
a log Calabi-Yau space, along with its Hodge decomposition.
We will succeed in doing so under certain hypotheses (see Theorem
\ref{bigtheorem}). The condition we require on the log Calabi-Yau
spaces is slightly
stronger than the requirement that the log Calabi-Yau space
be positive and simple. We have partial results when the space
is only positive and simple (Theorem \ref{notquite}).

As usual, we begin with calculations for our local models.

\subsection{Local calculations}
\label{localcalcsection}

In what follows, we want to understand certain operations we will perform
on the sheaves of differentials appearing in Proposition \ref{OmegaXl}.
Initially, we will think of these modules abstractly by considering
modules over monoid rings graded by elements of the monoid, and consider
the operations we will need to make use of in this more general context.

\begin{lemma}
\label{toricmodule}
Let $P$ be a toric monoid, $Q\subseteq P$ a face, $Y=\Spec \kk[P]$, 
$I\subseteq P$ a monoid ideal with radical $P\setminus Q$,
$X=\Spec \kk[P]/I$. Suppose furthermore that $p\in P$,
$q\in Q$, $p+q\in I$ implies $p\in I$.
Consider a $\kk[P]$-module $F$ of the form
\[
F=\bigoplus_{p\in P} z^pF_{\langle p\rangle}
\]
where $\langle p\rangle$ denotes the minimal face of $P$ containing
$p$, $F_{\langle p\rangle}$ a $\kk$-vector space, and $F_{P_1}\subseteq
F_{P_2}$ whenever $P_1\subseteq P_2$. Then
\begin{enumerate}
\item 
\[
F\otimes_{\kk[P]} \kk[P]/I=\bigoplus_{p\in P} z^p\left( {F_{\langle p\rangle}
\over \sum_{p'\in P, q\in I\atop p'+q=p} F_{\langle p'\rangle}}\right).
\]
\item
If $F_Q=F_P$ then 
\[
F_X:=(F\otimes_{\kk[P]} \kk[P]/I)/\Tors=\bigoplus_{p\in P\setminus
I} z^p F_{\langle p\rangle},
\]
where $\Tors$ denotes the submodule of elements of $F\otimes_{\kk[P]}\kk[P]/I$
with support on a proper closed subset of $X$.
\item
Let $J\subseteq P$ be a monoid ideal such that if
$Z\subseteq X$ is the closed
subscheme defined by $(I+J)/I\subseteq\kk[P]/I$, then $Z$ 
is codimension $\ge 2$ in $X$. Let
$\kappa:X\setminus Z\hookrightarrow X$ be the inclusion. If $\shF_X$ is
the sheaf on $X$ associated to $F_X$, then
\[
\Gamma(X,\kappa_*\kappa^*\shF_X)=\bigoplus_{p\in P^{\gp}} z^p
\bigcap_{q\in J\cap Q}
\bigcup_{n\ge 0\atop p+nq\in P\setminus I} F_{\langle p+nq\rangle}
\]
\end{enumerate}
\end{lemma}

\proof
(1) is obvious.
For (2),
we first check that if $p\in I$, then the 
degree $p$ piece of $F\otimes_{\kk[P]} \kk[P]/I$ is annihilated by an element
$z^q$ for some $q\in Q$. Since $\sqrt{I}=P\setminus Q$, $X\not\subseteq
V(z^q)$, and so this will imply the degree $p$ piece is torsion.
Choose an element $q\in Q$
not contained in any proper face of $Q$. Then $F_P=F_Q=F_{\langle q\rangle}
\subseteq F_{\langle p+q\rangle}\subseteq F_P$,
so $F_{\langle p+q\rangle}=F_P$, and
\[
{F_{\langle p+q\rangle}
\over \sum_{p'\in P, q'\in I\atop p'+q'=p+q} F_{\langle p'\rangle}}
=F_P/F_Q=0,
\]
(take $q'=p$, $p'=q$ in the sum). Thus $z^q$ annihilates the degree $p$
piece as desired.

On the other hand, if $p\in P\setminus I$, then the degree 
$p$ part of $F\otimes_{\kk[P]}\kk[P]/I$ is $F_{\langle p\rangle}$
by (1), and by the assumption on $I$, $z^q$ does not annihilate
any element of $F_{\langle p\rangle}$ for any $q\in Q$. This gives (2).

(3) Note that $X\setminus Z$ is covered by open subsets of the
form 
\[
D(z^q)=\{x\in\Spec \kk[P]/I| z^q\not\in x\}
\]
with $q\in J\cap Q$,
as the reduced space $X_{red}$ of $X$ is $\Spec \kk[Q]$ and
$\Spec \kk[Q]/J\cap Q$ has the same support in $X_{red}$ as $Z_{red}$.
Thus we can write
\[
\Gamma(X,\kappa_*\kappa^*\shF_X)=\bigcap_{q\in J\cap Q}\Gamma(D(z^q),\shF_X)
=\bigcap_{q\in J\cap Q} (F_X)_{z^q}.
\]
From (2) and the assumption on $I$, we can write
\[
(F_X)_{z^q}=\bigoplus_{p\in P^{\gp}} \bigcup_{n\ge 0\atop p+nq\in
P\setminus I} z^p F_{\langle p+nq\rangle},
\]
from which follows (3).
\qed

\bigskip

Suppose we are given data $\tau\subseteq M'_{\RR}$, $\Delta_1,\ldots,\Delta_q$
as in Construction \ref{keyP}, yielding a cone $K\subseteq M_{\RR}$,
$P=\dual K\cap N$, $\rho\in P$, $Y=\Spec \kk[P]$, $X=\Spec \kk[P]/(z^{\rho})$,
$\X_k=\Spec \kk[P]/(z^{(k+1)\rho})$.

For every face $\omega$ of $\tau$, we have a stratum $V_{\omega}\subseteq
X$, with
$V_{\omega}=\Spec \kk[P_{\omega}]$ where $P_{\omega}$
is the face of $P$ given by $P\cap (\omega+e_0)^{\perp}$. For every $k\ge 0$,
consider the monoid ideal
\[
I_{\omega}^k=\{p\in P|\hbox{$\langle p,n\rangle >k$ for some
$n\in\omega+e_0$}\}.
\]
This defines a thickening
\[
V^k_{\omega}=\Spec \kk[P]/I_{\omega}^k.
\]
One sees easily that $V_{\omega}^k$ is a closed subscheme of $\X_k$,
and that $I^k_{\omega}$ satisfies the hypotheses of Lemma \ref{toricmodule}
with $Q=P_{\omega}$.
Set $q_{\omega}:V_{\omega}^k\rightarrow\X_k$ the embedding.

Let $Z=\bigcup_i Z_i$ be the subscheme of $X$
defined in Construction \ref{keyP},
with $j:X\setminus Z\hookrightarrow X$
the inclusion. Set $D_{\omega}=\bigcup_{\omega\subsetneq\omega'\subseteq
\tau} V_{\omega'}$. 
This is a subset of the toric boundary of
$V_{\omega}$ consisting of proper intersections of the stratum $V_{\omega}$
with other strata of $X$. Let
\[
\kappa_{\omega}:V_{\omega}^k\setminus (D_{\omega}\cap Z)\rightarrow V^k_{\omega}
\]
be the inclusion. For $\Omega^{r}_k=j_*\Omega^r_{\X_k^{\dagger}/\kk}$ or
$\Omega^{r}_k=j_*\Omega^r_{\X_k^{\dagger}/A_k^{\dagger}}$, let
$\Omega^{r}_{\omega,k}=\kappa_{\omega*}\kappa_{\omega}^*(q_{\omega}^*\Omega^{r}_k
/\Tors)$.

Our goal now is to construct a resolution of the sheaf
$\Omega^r_k$. In the next subsection, we shall use this resolution
for $k=0$ to compute the log Hodge spaces of nice log Calabi-Yau spaces.
The case of $k>0$ will be needed for future period calculations.

\begin{lemma}
\label{localomegatauthickened}
Given $\omega\subseteq\tau$,
let $\omega_i\subseteq\Delta_i$ be the largest face of $\Delta_i$ such
that $\langle n,m\rangle=-\check\psi_i(n)$ for all $n\in\check\omega$,
$m\in \omega_i$. (Here $\check\omega$ is the cone in the normal
fan $\check\Sigma_{\tau}$ of $\tau$ corresponding to $\omega$). Then
\[
\Gamma(V_{\omega}^k,\Omega^{r}_{\omega,k})=\bigoplus_{p\in P_{\omega,k}} z^p
\left({\bigwedge}^r\bigcap_{\{(v,j)|v\in\omega_j,p\in (v+e_j)^{\perp}\}} ((v+e_j)^{\perp}\cap N)
\otimes_{\ZZ} \kk\right)
\]
or
\[
\bigoplus_{p\in P_{\omega,k}} z^p
\left({\bigwedge}^r\bigcap_{\{(v,j)|v\in\omega_j,p\in(v+e_j)^{\perp}\}} (((v+e_j)^{\perp}\cap N)
/\ZZ\rho)\otimes_{\ZZ} \kk\right)
\]
in the $/\kk$ or $/A_k^{\dagger}$ cases respectively, where $v$ runs over
vertices of $\omega_j$ for any $j$ and
\[
P_{\omega,k}:=\left\{p\in P^{\gp}\left|
\begin{array}{l}\hbox{$\langle p,v\rangle\ge 0$
for all $v\in\omega_i+e_i$, $1\le i\le q$}\\
\hbox{$\langle p,v\rangle\le k$ for all $v\in\omega+e_0$}\\
\hbox{$\langle p,v\rangle\ge 0$ for all $v\in \tau+e_0$}
\end{array}
\right.\right\}.
\]
\end{lemma}

\proof We will do the $/\kk$ case; the other case is identical.

Let $Q_1,\ldots,Q_t$ be the maximal proper faces of $P$ containing $\rho$.
Set, for $p\in P$,
\begin{equation}
\label{Omegarp}
\Omega^{r}_p={\bigwedge}^r\bigcap_{\{j| p\in Q_j\}} Q_j^{\gp}\otimes \kk,
\end{equation}
so that $\bigoplus_{p\in P} z^p\Omega^{r}_p$ defines a sheaf
$\Omega^{r}_Y$ on $Y$. Note that $\Omega^{r}_p$ only depends on $\langle p\rangle$,
so set $\Omega^{r}_{\langle p\rangle}:=\Omega^{r}_p$. One checks easily that
Proposition \ref{OmegaXl} implies $\Omega^{r}_Y|_{\X_k}\cong j_*
\Omega^r_{\X_k^{\dagger}/\kk}$. Then
by Lemma \ref{toricmodule}, (2) applied with $I=I_{\omega}^k$ and
$F_{\langle p\rangle}=\Omega^{r}_{\langle p\rangle}$,
\[
\Gamma(V_{\omega}^k,(q_{\omega}^*j_*\Omega^r_{\X_k^{\dagger}/\kk})/\Tors)
=\bigoplus_{p\in P\setminus I_{\omega}^k} z^p\Omega^{r}_p.
\]
Let $\Gamma(V_{\omega}^k,\Omega^{r}_{\omega,k})_p$, 
on the other hand, be the degree $p$ piece of
$\Gamma(V_{\omega}^k,\Omega^{r}_{\omega,k})$. 

First let us determine the structure of
$D_{\omega}\cap Z$. 
The toric strata of $V_{\omega}$ are in one-to-one correspondence
with faces $P'$ of $P_{\omega}$, which in turn are in one-to-one
order reversing correspondence with cones $K'$ with
$K_{\omega}\subseteq K'\subseteq K$, where $K_{\omega}=
C(\omega+e_0)$. Now a stratum corresponding to $K'$
is in $D_{\omega}\cap Z$
if it is contained in $D_{\omega}$ and $Z_i$ for some $i$.
The stratum is contained in $D_{\omega}$ provided $C(\omega'+e_0)
\subseteq K'$ for some $\omega'$ with $\omega\subsetneq\omega'\subseteq
\tau$. On the other hand, it is contained in $Z_i$ if, firstly,
$u_i=0$ on the stratum, i.e. $K'\cap C(\Delta_i+e_i)\not=0$ (otherwise
$K'\subseteq (e_i^*)^{\perp}$); secondly,
the stratum is contained in $V_{\omega''}$ for some $\omega''
\in \Omega_i$, this being equivalent to $\dim\omega_i'>0$. Thus,
let $P_{D_{\omega}\cap Z}$ be the union of faces of $P_{\omega}$
corresponding to cones $K'$ satisfying
\begin{enumerate}
\item
$K'\cap C(\Delta_0+e_0)=C(\omega'+e_0)$ for some $\omega'\supsetneq\omega$;
\item 
$K'\cap C(\Delta_i+e_i)\not=0$ 
and $\dim\omega'_i>0$ for some $1\le i\le q$.
\end{enumerate}
Then 
\[
J:=P\setminus P_{D_{\omega}\cap Z}
\]
is the monoid ideal defining $D_{\omega}\cap Z$.
In particular, by Lemma \ref{toricmodule}, (3),
\begin{equation}
\label{intunioneq}
\Gamma(V_{\omega}^k,\Omega^{r}_{\omega,k})_p
= 
\bigcap_{q\in J\cap P_{\omega}} \bigcup_{n\ge 0\atop
p+nq\in P\setminus I^k_{\omega}} \Omega^{r}_{\langle p+nq\rangle}.
\end{equation}
Let $q\in J\cap P_{\omega}$, and we consider the union in the above expression
for this $q$. Then $Q:=\langle q\rangle\subseteq P_{\omega}$ corresponds
to some $K'$ with $K_{\omega}\subseteq K'\subseteq K$ such that $K'$
fails to satisfy either property (1) or property (2) above. We consider
three cases.

\emph{Case 1}. $K'\cap C(\Delta_0+e_0)=C(\omega+e_0)$. Then 
$K'$ is contained in $K'_{\omega}:=C((\omega+e_0)\cup
\bigcup_{i=1}^q (\omega_i+e_i))$, as this is easily seen to be the largest
face of $K$ satisfying $K'_{\omega}\cap C(\Delta_0+e_0)=C(\omega+e_0)$.
Now as $q\in (K')^{\perp}$ but $\langle q,v\rangle>0$
for all $v\in K\setminus K'$, $p+nq\in P$ for large $n$ if and only if
$\langle p,v\rangle\ge 0$ for all $v\in K'$. This condition becomes
more restrictive the larger $K'$ is, so in order for 
$\Gamma(V_{\omega}^k,\Omega^{r}_{\omega,k})_p$ to be non-zero,
we must certainly
have $\langle p,v\rangle\ge 0$ for all $v\in K_{\omega}'$, i.e.
$\langle p,v\rangle \ge 0$ for $v\in\omega_i+e_i$, for all $i$, which
gives the first condition in the definition of $P_{\omega,k}$, and
also $\langle p,v\rangle \ge 0$ for $v\in\omega+e_0$. Furthermore,
assuming $p+nq\in P$ for large $n$, then
$p+nq\in P\setminus I_{\omega}^k$ for large $n$ if and only if 
$\langle p,v\rangle\le k$ for all $v\in \omega+e_0$. This is the second
condition in the definition of $P_{\omega,k}$.

Now suppose $p+nq\in P\setminus I^k_{\omega}$ for large $n$, and look
at \eqref{intunioneq}. It is clear that for
two choices of $q$, $q'$ with $\langle q\rangle\subseteq\langle q'\rangle$
and both $q$ and $q'$ falling into this first case, then 
$\langle p+nq\rangle\subseteq \langle p+nq'\rangle$ for large $n$, and hence 
$\Omega^{r}_{\langle p+nq\rangle}\subseteq 
\Omega^{r}_{\langle p+nq'\rangle}$ for large $n$. Thus as far
as the intersection is concerned, we can assume that $\langle q\rangle$
is as small as possible, and hence $K'=K'_{\omega}$. 

In this case, $q\in (v+e_i)^{\perp}$ for each vertex $v$ of $\omega_i$,
but $\langle q,v+e_i\rangle>0$ if $v\in\Delta_i\setminus\omega_i$. Thus
for large $n$ and $v$ a vertex of $\omega_i$, 
$p+nq\in (v+e_i)^{\perp}$ if and only if $p\in (v+e_i)^{\perp}$.
On the other hand, if $v\in\Delta_i\setminus\omega_i$,
$p+nq\not\in
(v+e_i)^{\perp}$ for large $n$. Since the $Q_j^{\gp}$'s
are the spaces $(v+e_i)^{\perp}$ for $v$ running over vertices of $\Delta_i$,
for all $i$, we see that for large $n$,
\[
\Omega^{r}_{\langle p+nq\rangle}
={\bigwedge}^r 
\bigcap_{\{(v,j)|v\in\omega_j,p\in (v+e_j)^{\perp}\}}
((v+e_j)^{\perp}\cap N)\otimes_{\ZZ}\kk.
\]

\emph{Case 2.} Here we consider a special case when property (2) doesn't
hold, namely $K'=C(\omega'+e_0)$ for some $\omega'\supseteq\omega$.
In this case, $p+nq\in P$ for large $n$ if and only if 
$\langle p,v\rangle\ge 0$ for all $v\in \omega'+e_0$. In particular,
if we take $K'=C(\tau+e_0)$, we see
the intersection in \eqref{intunioneq} is zero unless
$\langle p,v\rangle\ge 0$ for all $v\in\tau+e_0$, hence the last
condition in the definition of $P_{\omega,k}$. 

\emph{Case 3}. We consider the general case when property (2) doesn't hold,
so $K'\cap C(\Delta_0+e_0)=C(\omega'+e_0)$ for some
$\omega'\supsetneq\omega$,
but $\dim\omega_i'=0$ whenever $K'\cap C(\Delta_i+e_i)\not=0$.
Then $K'$ must be contained in $K'_{\omega'}=
C((\omega'+e_0)+\bigcup_{i=1}^q(\omega'_i+e_i))$, but the extra condition
implies that in fact $K'$ is contained in
$C((\omega'+e_0)+\bigcup_{\{i|\dim\omega'_i=0\}} (\omega'_i+e_i))$.
Since $\omega_i\subseteq\omega'_i$, we have $\omega_i=\omega_i'$
whenever $\dim\omega_i'=0$,
hence $K'\subseteq
C((\omega'+e_0)+\bigcup_{\{i|\dim\omega'_i=0\}}(\omega_i+e_i))$. One 
sees that if $p\in P_{\omega,k}$, then $p+nq\in P\setminus I^k_{\omega}$
for large $n$, and furthermore for large $n$,
$\Omega^{r}_{\langle p+nq\rangle}$
contains $\Omega^{r}_{\langle p+nq'\rangle}$ for $q'$ 
with $\langle q'\rangle$ dual to $K'_{\omega}$.

Putting together these three cases, one obtains the desired result.
\qed

\begin{proposition}
\label{Omegarestricted}
Given faces $\omega\subseteq\omega'\subseteq\tau$, we have
$I_{\omega}^k\subseteq I_{\omega'}^k$, and hence a closed embedding
$V^k_{\omega'}\rightarrow V^k_{\omega}$. Then 
\[
\Gamma(V_{\omega'}^k,\Omega^{r}_{\omega,k}|_{V^k_{\omega'}}/\Tors)
=\bigoplus_{p\in P_{\omega,\omega',k}} z^p
\left({\bigwedge}^r\bigcap_{\{(v,j)|v\in\omega_j,p\in (v+e_j)^{\perp}\}} ((v+e_j)^{\perp}\cap N)
\otimes_{\ZZ} \kk\right)
\]
or
\[
\bigoplus_{p\in P_{\omega,\omega',k}} z^p
\left({\bigwedge}^r\bigcap_{\{(v,j)|v\in\omega_j,p\in(v+e_j)^{\perp}\}} (((v+e_j)^{\perp}\cap N)
/\ZZ\rho)\otimes_{\ZZ} \kk\right)
\]
in the $/\kk$ or $/A_k^{\dagger}$ cases respectively, where $v$
runs over vertices of $\omega_j$ and
\[
P_{\omega,\omega',k}:=\left\{p\in P^{\gp}\left|\begin{array}{l}
\hbox{$\langle p,v\rangle\ge 0$
for all $v\in\omega_i+e_i$, $1\le i\le q$}\\
\hbox{$\langle p,v\rangle\le k$ for all $v\in\omega'+e_0$}\\
\hbox{$\langle p,v\rangle\ge 0$ for all $v\in \tau+e_0$}
\end{array}
\right.\right\}.
\]
Note the only difference between this set and $P_{\omega,k}$ defined in
Lemma \ref{localomegatauthickened} is that in the latter, $\langle p,v
\rangle\le k$ for all $v\in\omega+e_0$ instead of for all $v\in\omega'+e_0$.
\end{proposition}

\proof Let $\tilde P$ be the monoid
\[
\tilde P:=\left\{p\in P^{\gp}\bigg|\begin{array}{l}
\hbox{$\langle p,v\rangle\ge 0$
for all $v\in\omega_i+e_i$, $1\le i\le q$}\\
\hbox{$\langle p,v\rangle\ge 0$ for all $v\in \tau+e_0$}
\end{array}
\right\}
\]
with ideals
\[
\tilde I^k_{\omega}:=\{p\in \tilde P| \hbox{$\langle p,v\rangle>k$
for some $v\in\omega+e_0$}\}
\]
and
\[
\tilde I^k_{\omega'}:=\{p\in \tilde P| \hbox{$\langle p,v\rangle>k$
for some $v\in\omega'+e_0$}\}.
\]
Note $P\subseteq\tilde P$, $I^k_{\omega}=P\cap\tilde I^k_{\omega}$, 
$I^k_{\omega'}=P\cap\tilde I_{\omega'}^k$. If $F$ is the $\kk[\tilde P]$-module
defined by
\[
F=\bigoplus_{p\in\tilde P}z^p
\left({\bigwedge}^r\bigcap_{\{(v,j)|v\in\omega_j,p\in (v+e_j)^{\perp}\}} ((v+e_j)^{\perp}\cap N)
\otimes_{\ZZ} \kk\right)
\]
(or a similar expression in the $/A_k^{\dagger}$ case), then from 
Lemma~\ref{toricmodule}, (2) and Lemma~\ref{localomegatauthickened}, we see that
\[
\Gamma(V^k_{\omega},
\Omega^{r}_{\omega,k})\cong (F\otimes_{\kk[\tilde P]} \kk[\tilde P]/\tilde
I_{\omega}^k)/\Tors.
\]
Now we claim that
\[
\Gamma(V^k_{\omega'},\Omega^{r}_{\omega,k}|_{V^k_{\omega'}}/\Tors)
\cong (F\otimes_{\kk[\tilde P]}\kk[\tilde P]/\tilde I_{\omega'}^k)/\Tors,
\]
from which will follow the result by Lemma \ref{toricmodule}, (2) again.
To prove this claim, first note it is easy to check that $\Tors\subseteq
F\otimes_{\kk[\tilde P]}\kk[\tilde P]/\tilde I^k_{\omega}$ contains no
sections with support containing $V^k_{\omega'}$, from which it follows
that
\[
(\Gamma(V^k_{\omega},\Omega^{r}_{\omega,k})\otimes_{\kk[\tilde P]/\tilde I^k_{\omega}}
\kk[\tilde P]/\tilde I_{\omega'}^k)/\Tors
\cong (F\otimes_{\kk[\tilde P]} \kk[\tilde P]/\tilde I_{\omega'}^k)/\Tors.
\]
So we need to compare
\[
\tilde M=
(\Gamma(V^k_{\omega},\Omega^{r}_{\omega,k})
\otimes_{\kk[\tilde P]/\tilde I^k_{\omega}}
\kk[\tilde P]/\tilde I_{\omega'}^k)/\Tors
\]
with
\[
M=
(\Gamma(V^k_{\omega},\Omega^{r}_{\omega,k})
\otimes_{\kk[P]/I^k_{\omega}}
\kk[P]/I_{\omega'}^k)/\Tors.
\]
Note $\kk[P]/I^k_{\omega'}\subseteq\kk[\tilde P]/\tilde I_{\omega'}^k$,
so $\tilde M$ is also a $\kk[P]/I^k_{\omega'}$-module, and we need
to compare $\tilde M$ and $M$ as $\kk[P]/I_{\omega'}^k$-modules.
We then have a diagram of $\kk[P]/I_{\omega'}^k$-modules with $T$ and
$\tilde T$ the torsion submodules of the modules immediately below them,
(also confusing the sheaf $\Omega^{r}_{\omega,k}$ with its space
of global sections),
\[
\xymatrix@C=30pt
{&&0\ar[d]&0\ar[d]&\\
&&T\ar[r]\ar[d]&\tilde T\ar[d]&\\
0\ar[r]&{\displaystyle\tilde I^k_{\omega'}\Omega^{r}_{\omega,k}
\over \displaystyle I^k_{\omega'}\Omega^{r}_{\omega,k}}\ar[r]&
{\displaystyle\Omega^{r}_{\omega,k}
\over \displaystyle I^k_{\omega'}\Omega^{r}_{\omega,k}}\ar[r]\ar[d]&
{\displaystyle\Omega^{r}_{\omega,k}
\over\displaystyle \tilde I^k_{\omega'}\Omega^{r}_{\omega,k}}\ar[r]\ar[d]&0\\
&&M\ar[d]\ar[r]&\tilde M\ar[d]\ar[r]&0\\
&&0&0&
}
\]
We see that $\tilde I^k_{\omega'}\Omega^{r}_{\omega,k}/
I^k_{\omega'}\Omega^{r}_{\omega,k}\subseteq T$, as for
any $q$ in the interior of $P_{\omega'}$ and $p\in\tilde P$, $p+nq\in P$
for $n\gg 0$, so $z^{nq}$ annihilates
$\tilde I^k_{\omega'}\Omega^{r}_{\omega,k}/
I^k_{\omega'}\Omega^{r}_{\omega,k}$ for $n\gg 0$. Thus by the
snake lemma, $M\cong\tilde M$ provided $T\rightarrow\tilde T$ is
surjective. But if an element of $\tilde T$ is represented by
$\alpha\in \Omega^{r}_{\omega,k}$, we have $z^{nq}\alpha
\in \tilde I^k_{\omega'}\Omega^{r}_{\omega,k}$ for $q$ in the interior 
of $P_{\omega'}$ and $n\gg 0$. But since $z^{nq}\tilde I^k_{\omega'}
\subseteq I^k_{\omega'}$ for $n\gg 0$, we see in fact $\alpha$
represents an element of $T$, proving surjectivity. \qed

\begin{corollary}
\label{inclusionintersectioncor} 
\begin{enumerate}
\item Given faces $\omega_1\subseteq\omega_2\subseteq\omega_3$ of
$\tau$, we have an inclusion 
\[
(\Omega^{r}_{\omega_2,k}|_{V^k_{\omega_3}})/\Tors\subseteq
(\Omega^{r}_{\omega_1,k}|_{V^k_{\omega_3}})/\Tors.
\]
\item
Given $\omega_1\subseteq\omega_2$ faces of $\tau$, 
\[
(\Omega^{r}_{\omega_1,k}|_{V^k_{\omega_2}})/\Tors
=\bigcap_{v\in\omega_1}(\Omega^{r}_{v,k}|_{V^k_{\omega_2}})/\Tors,
\]
where $v$ runs over vertices of $\omega_1$, and the intersection
can be viewed as taking place in 
$j_*(\Omega^{r}_{v,k}|_{V^k_{\omega_2}\setminus Z})$,
which is independent of $v$ since $\Omega^{r}_k$ is locally free
away from $Z$.
\end{enumerate}
\end{corollary}

\proof These statements follow immediately from the explicit formula of
the previous corollary.
\qed

\bigskip

We can now define a resolution of $\Omega^{r}_k$. We define a
barycentric complex (a slight variation of that of \cite{PartI}, Appendix~A)
by
\[
\C^p(\Omega^{r}_k)=\bigoplus_{\omega_0\subsetneq\cdots\subsetneq\omega_p
\subseteq\tau}(\Omega^{r}_{\omega_0,k}|_{V^k_{\omega_p}})/\Tors
\]
and a differential
\[
d_{\bct}:\C^p(\Omega^{r}_k)\rightarrow\C^{p+1}(\Omega^{r}_k)
\]
given by
\[
(d_{\bct}(\alpha))_{\omega_0\subsetneq\cdots\subsetneq\omega_{p+1}}
=\sum_{i=0}^p (-1)^i\alpha_{\omega_0\subsetneq\cdots\hat\omega_i\subsetneq
\cdots\subsetneq\omega_{p+1}}
+(-1)^{p+1}\alpha_{\omega_0\subsetneq\cdots\subsetneq\omega_p}|_{V^k_{\omega_{p+1}}}.
\]
Here we use the inclusions of Corollary \ref{inclusionintersectioncor}, (1)
to identify all these elements with elements of 
$(\Omega^{r}_{\omega_0,k}|_{V^k_{\omega_{p+1}}})/\Tors$.

\begin{theorem}
\label{localresolution}
$\C^{\bullet}(\Omega^{r}_k)$ is a resolution of
$\Omega^{r}_k$.
\end{theorem}

\proof The inclusion $\Omega^{r}_k\rightarrow \C^0(\Omega^{r}_k)$
is defined in the obvious way: an element $\alpha$ of $\Omega^{r}_k$
yields for each face $\omega$ of $\tau$ an element of 
$\Omega^{r}_k|_{V^k_{\omega}}$, and hence an element of
$\Omega^{r}_{\omega,k}$. This gives the map
\[
\epsilon:\Omega^{r}_k\rightarrow\C^0(\Omega^{r}_k)
=\bigoplus_{\omega\subseteq\tau}\Omega^{r}_{\omega,k},
\]
which is clearly injective, as is easily seen by checking on
$\X_k\setminus Z$, where everything is locally free. 
The fact that $\epsilon(\Omega^{r}_k)\subseteq\ker(d_{\bct}:
\C^0(\Omega^{r}_k)\rightarrow\C^1(\Omega^{r}_k))$
is also easily checked off of $Z$, where all these sheaves are
locally free. Conversely, if $\alpha\in\Gamma(U,\C^0(\Omega^{r}_k))$ 
with $d_{\bct}(\alpha)=0$,
then $d_{\bct}(\alpha|_{U\setminus Z})=0$, and again as everything is
locally free on $U\setminus Z$, it is obvious that by gluing we obtain
an element $\beta\in\Gamma(U\setminus Z,\Omega^{r}_k)$ mapping to
$\alpha|_{U\setminus Z}$. But $\Gamma(U\setminus Z,\Omega^{r}_k)
=\Gamma(U,\Omega^{r}_k)$ by definition, so we obtain an
element $\beta\in\Gamma(U,\Omega^{r}_k)$ such that
$\epsilon(\beta)=\alpha$ on $U\setminus Z$. Thus $\epsilon(\beta)=\alpha$
on $U$ as these sheaves have no sections with support on $Z$,
so $\ker(d_{\bct})=\epsilon(\Omega^{r}_k)$.

To check exactness of $\C^{\bullet}$, we can't use \cite{PartI}, \S A.1
directly because our complex takes a slightly different form: 
the module $(\Omega^{r}_{\omega_0,k}|_{V^k_{\omega_p}})/\Tors$
depends on both $\omega_0$ and $\omega_p$, not just $\omega_p$.
However, the result follows from a slight modification of the
argument given in \cite{PartI}, Proposition A.1. The question of exactness then
reduces to a version of an easily checked compatibility condition,
namely if we are given $\omega_0\subsetneq\omega_{p-1}$ and elements 
of $(\Omega_{\omega_0,k}^{r}|_{V^k_{\omega_p}})/\Tors$ 
for all $\omega_{p-1}\rightarrow
\omega_p$ which agree under restriction, then these elements lift to an
element of
$(\Omega^{r}_{\omega_0,k}|_{V^k_{\omega_{p-1}}})/\Tors$. 
However, this is true
again by the explicit description in Proposition~\ref{Omegarestricted}.
\qed

\begin{proposition}
\label{differential}
The differential $d:j_*\Omega^r_{\X_k^{\dagger}/\kk}\rightarrow 
j_*\Omega^{r+1}_{\X_k^{\dagger}/\kk}$ (or 
$d:j_*\Omega^r_{\X_k^{\dagger}/A_k^{\dagger}}\rightarrow 
j_*\Omega^{r+1}_{\X_k^{\dagger}/A_k^{\dagger}}$) is given on the degree 
$p$ piece of $\Gamma(\X_k,j_*\Omega^r_{\X_k^{\dagger}/\kk})$ by $z^p n\mapsto z^p
p\wedge n$. For any pair of faces $\omega_1\subseteq\omega_2\subseteq
\tau$, this induces a map
$d:(\Omega^{r}_{\omega_1,k}|_{V_{\omega_2}^k})/\Tors\rightarrow
(\Omega^{r+1}_{\omega_1,k}|_{V_{\omega_2}^k})/\Tors$.
\end{proposition}

\proof Of course $d(z^p \dlog n)=d(z^p)\wedge\dlog n=
z^p\dlog p\wedge \dlog n$, giving the first formula.
The second statement follows from the fact that the explicit description
in Proposition \ref{Omegarestricted} of 
$(\Omega^{r}_{\omega_1,k}|_{V^k_{\omega_2}})/\Tors$ is clearly closed
under this operation. \qed

\bigskip

\subsection{Global calculations}
\label{globalcalculations}

Let $(B,\P)$ be a positive and simple integral affine manifold with
singularities with toric polyhedral decomposition. Let $s$ be lifted
open gluing data for $(B,\P)$, yielding
$X_0:=X_0(B,\P,s)$. Then $s$ also determines a well-defined 
log structure  on $X_0$
over $\Spec \kk^{\dagger}$ with singular set $Z\subseteq X_0$.
In what follows, we will take
$\Omega^r$ to be the sheaf on $X_0$ 
which is either $j_* \Omega^r_{X_0^{\dagger}/\kk}$ or
$j_*\Omega^r_{X_0^{\dagger}/\kk^{\dagger}}$,
where $j:X_0\setminus Z\rightarrow X_0$ is the inclusion. 
We refer to these as the $/\kk$ and $/\kk^{\dagger}$ cases, respectively.
Our goal is to calculate $H^p(X_0,\Omega^r)$. This section will
be devoted to technical results which essentially lift the local
descriptions of the previous subsection to the global situation.

Our first goal is to obtain 
a nice resolution for $\Omega^r$. We have studied the local form
of this resolution in \S \ref{localcalcsection}. 

Let $q_{\tau}:X_{\tau}\rightarrow X_0$ be the usual inclusion
of strata maps, (\cite{PartI}, Lemma
2.29), $D_{\tau}$ the
toric boundary of $X_{\tau}$, and let
\begin{eqnarray*}
j_{\tau}:X_{\tau}\setminus q_{\tau}^{-1}(Z)&\rightarrow X_{\tau}\\
\kappa_{\tau}:X_{\tau}\setminus (D_{\tau}\cap q_{\tau}^{-1}(Z))
&\rightarrow X_{\tau}
\end{eqnarray*}
be the inclusions. 
We then define 
\[
\Omega^r_{\tau}:=\kappa_{\tau*}\kappa_{\tau}^*(q_{\tau}^*\Omega^r/\Tors).
\]
where $\Tors$ denotes the torsion
subsheaf of $q_{\tau}^*\Omega^r$.

\begin{remark} While we do not need the result here, one can in fact
show that $\Omega^r_{\tau}$ is locally free, by showing the same result
for $\Omega^{r}_{\tau,0}$ in \S \ref{localcalcsection}. This is only
true at $0$th order, not for $k>0$. If one did know this fact, one could
omit all the quotients by $\Tors$ appearing in this section, so the
reader can safely ignore this torsion.
\end{remark}
\bigskip

Recall that $X_0$ can be viewed as the direct limit of a gluing
functor $F_{S,s}$ defined in \cite{PartI}, Definition 2.11. Here we
are taking $S=\Spec\kk$ as the base scheme, while in \S 5 we shall
use a different choice of $S$. Since $S$ and $s$ are given, we shall
write, for $\tau_1\subseteq\tau_2$, 
\[
F_{\tau_1,\tau_2}:X_{\tau_2}\rightarrow X_{\tau_1}
\]
for 
\[
F_{S,s}(\tau_1\rightarrow\tau_2):X_{\tau_2}\rightarrow X_{\tau_1}.
\]
This introduces ambiguity in case there are several morphisms $\tau_1
\rightarrow\tau_2$, but this only happens if the cell $\tau_2$ has
several faces identified inside of $B$. If the reader wants to consider
this case, he can return to the somewhat denser notation of $F_{S,s}(\tau_1
\rightarrow\tau_2)$ instead of $F_{\tau_1,\tau_2}$.

Recall that
\[
q_{\tau_2}=q_{\tau_1}\circ F_{\tau_1,\tau_2}.
\]

Adapting the local results of \S 3.1 to the global situation,

\begin{proposition}
\label{Omegainc}
If $\tau_1\subseteq\tau_2$ with $\tau_1,\tau_2\in\P$, then 
the functorial isomorphism on $X_{\tau_2}\setminus q_{\tau_2}^{-1}(Z)$
\[
\Omega^r_{\tau_2}=q_{\tau_2}^*\Omega^r
\mapright{\cong} F_{\tau_1,\tau_2}^*q_{\tau_1}^*\Omega^r
=F_{\tau_1,\tau_2}^*\Omega^r_{\tau_1}
\]
extends to an inclusion
\[
F_{\tau_1,\tau_2}^*:
\Omega^r_{\tau_2}\rightarrow (F_{\tau_1,\tau_2}^*\Omega^r_{\tau_1})/\Tors.
\]
\end{proposition}

\proof This can be checked in an \'etale neighbourhood of
a point $z\in Z$; by Theorem \ref{localmodel}, this reduces
to the case considered in Corollary \ref{inclusionintersectioncor}, (1).
\qed
\bigskip

We are now able to define our explicit resolution of the sheaf 
$\Omega^r$.

We define a barycentric complex with 
\[
\C^k(\Omega^r)=\bigoplus_{\sigma_0\subsetneq\cdots\subsetneq
\sigma_k} q_{\sigma_k*}((F_{\sigma_0,\sigma_k}^*
\Omega^r_{\sigma_0})/\Tors)
\]
and a differential $d_{\bct}:\C^k(\Omega^r)
\rightarrow\C^{k+1}(\Omega^r)$ given by
\begin{eqnarray*}
(d_{\bct}(\alpha))_{\sigma_0,\cdots,\sigma_{k+1}}
&=&\alpha_{\sigma_1,
\cdots,\sigma_{k+1}}
+\sum_{i=1}^{k} (-1)^i\alpha_{\sigma_0,\cdots,
\check\sigma_i,\cdots,\sigma_{k+1}}
\\
&&+(-1)^{k+1}F_{\sigma_k,\sigma_{k+1}}^*\alpha_{\sigma_0,
\cdots,\sigma_k}.
\end{eqnarray*}
Here $\alpha_{\sigma_1,\cdots,\sigma_{k+1}}
\in (F_{\sigma_1,\sigma_{k+1}}^*\Omega^r_{\sigma_1})/\Tors$
can be viewed, by Proposition \ref{Omegainc}, as an element of $(F_{
\sigma_0,\sigma_{k+1}}^*\Omega^r_{\sigma_0})/\Tors$.

\begin{theorem}
\label{resolution}
$\C^{\bullet}(\Omega^r)$ is a resolution of $\Omega^r$.
\end{theorem}

\proof This follows immediately from the local version, Theorem
\ref{localresolution}. \qed

\begin{corollary}
\label{hyper}
\[
H^p(X_0,\Omega^r)=\HH^p(X_0,\C^{\bullet}(\Omega^r)).
\]
\end{corollary}

Now the differential $d:\Omega^r\rightarrow\Omega^{r+1}$ is defined on
$X_0\setminus Z$, and hence on the pushforward, giving us
a complex $(\Omega^{\bullet},d)$, the \emph{log de Rham complex} of
$X_0$.

By Proposition \ref{differential}, $d$ induces maps, for $e:\tau_1\rightarrow
\tau_2$, 
\[
d:(F_{\tau_1,\tau_2}^*\Omega^r_{\tau_1})/\Tors\rightarrow
(F_{\tau_1,\tau_2}^*\Omega^{r+1}_{\tau_1})/\Tors,
\]
and hence a map of complexes
\[
d:\C^{\bullet}(\Omega^r)\rightarrow\C^{\bullet}(\Omega^{r+1}).
\]
This gives us a double complex $\C^{\bullet}(\Omega^{\bullet})$,
and the obvious

\begin{corollary}
\label{algdeRham}
\[
\HH^r(X_0,\Omega^{\bullet})=\HH^r(X_0,\Tot(\C^{\bullet}
(\Omega^{\bullet}))),
\]
where $\Tot$ denotes the total complex of the double complex.
\end{corollary}

In order to compute these cohomology groups explicitly, we need a 
useful global description for the sheaves $\Omega^r_{\omega}$. We first
describe $\Omega^r_v$ for $v$ a vertex of $\P$.

Pull back the log structure on $X_0$ via $q_v$ to obtain a 
log structure on $X_v\setminus q_v^{-1}(Z)$, with sheaf of monoids
$\M_v$. By \cite{PartI}, Lemma 5.13, we have a split exact sequence
\begin{equation}
\label{Mexact}
0\rightarrow\M^{\gp}_{(X_v,D_v)}\rightarrow \M_v^{\gp}\rightarrow
\ZZ\rho\rightarrow 0,
\end{equation}
where $\M_{(X_v,D_v)}$ is the sheaf of monoids associated to the
divisorial log structure given by $D_v\subseteq X_v$, and $\rho$
as usual is the image of $1\in\NN$ under the map of monoids induced by
the log morphism $X_0^{\dagger}\rightarrow\Spec\kk^{\dagger}$.
Because $q_v^{-1}(Z)\subseteq X_v$ is codimension two,
$j_*\M_v\rightarrow j_*\O_{X_v\setminus q_v^{-1}(Z)}=\O_{X_v}$
determines a log structure on $X_v$, which we write as $X_v^{\dagger}$.
Write $\M_v$ also for $j_*\M_v$. Then the exact sequence (\ref{Mexact})
still holds on $X_v$. From this exact sequence one sees that
$\Omega^1_{X_v^{\dagger}/\kk^{\dagger}}$ coincides with the ordinary sheaf
of log derivations for the pair $(X_v,D_v)$, which is canonically
$\check\Lambda_v\otimes\O_{X_v}$ by \cite{Oda}, Proposition~3.1, while 
$\Omega^1_{X_v^{\dagger}/\kk}$ is $(\check\Lambda_v\oplus\ZZ\rho)
\otimes\O_{X_v}$. We will 
identify this with $\shAff(B,\ZZ)_v\otimes\O_{X_v}$
(see \cite{PartI}, Definition 1.11) in \S \ref{monodromysection}. 

\begin{lemma}
\label{OXD}
Let $v\in\P$ be a vertex. Then $\Omega^r_v$ is naturally isomorphic to
$\Omega^r_{X_v^{\dagger}/\kk}$ or $\Omega^r_{X_v^{\dagger}/\kk^{\dagger}}$
in the $/\kk$ and $/\kk^{\dagger}$ cases respectively.
\end{lemma}

\proof  
We'll do the $/\kk$ case, the $/\kk^{\dagger}$ case being similar.
Functoriality of log differentials
gives a map $q_v^*:q_v^*\Omega^1\rightarrow
\Omega^1_{X_v^{\dagger}/\kk}$ on $X_v\setminus q_{v}^{-1}(Z)$.
This map is injective as it is generically injective and $q_v^*\Omega^1$
is locally free on $X_v\setminus q_v^{-1}(Z)$.
To see it is surjective, we need to recall the technique of
\cite{PartI}, Lemma 5.13. Let $e:v\rightarrow\sigma\in\P_{\max}$.
Without loss of generality, we can view $V(\sigma)=\Spec 
\kk[P_{\sigma}]/(z^{\rho})$ as an open subset of $X_0$. 
(See \cite{PartI}, Definition 2.12 for $P_{\sigma}$.)
Then in the proof
of \cite{PartI}, Lemma 5.13, the log structure on $V(\sigma)\setminus Z$
is given by charts $\varphi_i$
on an open cover $\{U_i\}$ of $V(\sigma)\setminus Z$,
$\varphi_i:P_{\sigma}\rightarrow\O_{U_i}$ a monoid homomorphism. 
Restricting these charts
to $U_i\cap X_v$ gives charts $\varphi_i:P_{\sigma}\rightarrow\O_{U_i
\cap X_v}$, which were shown to be of the form
\[
p\mapsto \begin{cases} 0&p\not\in P_e\\ h_p z^p&p\in P_e\end{cases}
\]
where $P_e$ is the maximal proper face of $P_{\sigma}$ corresponding to
$X_v\cap V(\sigma)$ and
$P_e\ni p\mapsto h_p\in\O_{U_i\cap X_v}^{\times}$ is a monoid
homomorphism. 
Note $P_e^{\gp}\cong\check\Lambda_v$ canonically. This chart lifts
to a monoid homomorphism
$\varphi_i:P_{\sigma}\rightarrow\M_{U_i}$, so for $p\in P_e$,
$\dlog (\varphi_i(p))\in \Gamma(U_i,\Omega^1)$ pulls back via $q_v^*$ to
$\dlog(h_pz^p)={d(h_p)\over h_p} +\dlog (z^p)$ in $\Omega^1_{X_v^{\dagger}/
\kk}$. By extending $h_p$
to $U_i$, we see $dh_p$ is in the image of $q_v^*\Omega^1
\rightarrow\Omega^1_{X_v^{\dagger}/\kk}$, so $\dlog z^p$ is 
also for all $p\in\check\Lambda_v$. 
On the other hand $\dlog\rho$ clearly pulls back to 
$\dlog\rho\in\Omega^1_{X_v^{\dagger}/\kk}$.
Thus $q_v^*$ is surjective
on each $U_i$, hence on $X_{v}\setminus q_{v}^{-1}(Z)$.

Now on 
$X_v\setminus q_v^{-1}(Z)=X_v\setminus (D_v\cap
q_v^{-1}(Z))$, $\Omega^1_v=\kappa_{v*}\kappa_v^*(q_v^*\Omega^1/\Tors)
=\kappa_{v*}\kappa_v^*q_v^*\Omega^1$, 
so we get an isomorphism on $X_v$
\[
\Omega^1_v=\kappa_{v*}\kappa_v^*q_v^*\Omega^1
\rightarrow \kappa_{v*}(\Omega^1_{X_v^{\dagger}/\kk})=\Omega^1_{X_v^{\dagger}/\kk},
\]
the latter equality as $X_v$ is $S_2$ and $q_v^{-1}(Z)$ is codimension
at least two in $X_v$.
Similarly, we obtain $\Omega^r_v\cong\Omega^r_{X_v^{\dagger}/\kk}$.
\qed
\medskip

One way to interpret this Lemma is that one can view
$\Omega^r$ as being obtained by gluing together trivial vector bundles
on the irreducible components of $X_0\setminus Z$. 
Consider the situation
where $\omega\in\P$ is a cell of dimension one, with vertices
$e^{\pm}_{\omega}:v^{\pm}_{\omega}\rightarrow\omega$ arising from
a choice of $d_{\omega}$ a primitive generator of $\Lambda_{\omega}$.
On $X_{\omega}\setminus q_{\omega}^{-1}(Z)$, there are of course
canonical identifications
\[
F_{v^-_{\omega},\omega}^*\Omega_{v^-_{\omega}}^r
= F_{v^-_{\omega},\omega}^*
q_{v^-_{\omega}}^*\Omega^r= F_{v^+_{\omega},\omega}^*
q_{v^+_{\omega}}^*\Omega^r=F_{v^+_{\omega},\omega}^*
\Omega_{v^+_{\omega}}^r.
\]
On the other hand, using the isomorphism of Lemma \ref{OXD} on the left
and right hand sides of the above identifications, we get
on $X_{\omega}\setminus q_{\omega}^{-1}(Z)$ a map
\begin{equation}
\label{Gammao}
\Gamma_{\omega}:F_{v^-_{\omega},\omega}^*
\Omega^r_{X^{\dagger}_{v^-_{\omega}}/\kk}
\mapright{\cong}
F_{v^+_{\omega},\omega}^* \Omega_{X^{\dagger}_{v^+_{\omega}}/\kk}^r
\end{equation}
(or $/\kk^{\dagger}$.) Let's describe $\Gamma_{\omega}$ explicitly.

With $\omega\in\P$ one-dimensional, 
pick $\omega\rightarrow\sigma$ with $\sigma\in\P_{\max}$, so that
we obtain $V(\omega)\subseteq V(\sigma)$.
Recall from \cite{PartI}, Definition 4.21,
that specifying a log smooth structure with singularities
of the desired type on $V(\sigma)$ means giving a section $f_{\sigma}\in
\Gamma(V(\sigma),\shLS^+_{\pre,V(\sigma)})$ with certain
properties. The sheaf
$\shLS^+_{\pre,V(\sigma)}$ restricted to $V(\omega)$ is just the structure
sheaf $\O_{\Sing(V(\omega))}$ of the singular locus of $V(\omega)$, i.e.
the intersection of the two irreducible components $V_{e^{\pm}_{\omega}}$
of $V(\omega)$. So on $V(\omega)$, $f_{\sigma}$ is just a function
on $\Sing(V(\omega))$.

\begin{lemma}
\label{basicgluing}
In the above situation,
identify $\check
\Lambda_{v^+_{\omega}}$ and $\check\Lambda_{v^-_{\omega}}$ via parallel
transport through $\sigma$, and identify these with a lattice $N$. Then
on $\Sing(V(\omega))$, in the $/\kk$ case,
$\Gamma_{\omega}$
is given by, for $n\in\bigwedge^{\bullet}(N\oplus\ZZ\rho)$,
\[
\Gamma_{\omega}(\dlog n)=-\left({df_{\sigma}\over f_{\sigma}}+l_{\omega}\dlog\rho\right)
\wedge \dlog(\iota(d_{\omega})n)+
\dlog n
\]
where $l_{\omega}$ is a positive integer such that there is an integral
affine isomorphism $[0,l_{\omega}]\rightarrow\omega$.
The same formula holds modulo
$\dlog \rho$ in the $/\kk^{\dagger}$ case.
\end{lemma}

\proof
We will follow the notation used in \cite{PartI},
Construction 2.15. We view $\sigma\subseteq M_{\RR}$
as a lattice polytope, with $\omega\subseteq\sigma$ an edge,
and set $P:=\dual{C(\sigma)}\cap (N\oplus\ZZ)$, 
$Q=\dual{C(\omega)}\cap (N\oplus\ZZ)$.
We have the monoid $\partial Q$ we can identify via projection
to $N$ with the set
$N\cup\{\infty\}$, with $V(\omega)=\Spec \kk[\partial Q]$.
If $\check\Sigma$ is the normal fan
to $\sigma$, then $\check\omega^{-1}\check\Sigma$ has two maximal cones,
$\check v^{\pm}_{\omega}$. Addition on $\partial Q=N\cup\{\infty\}$
is given by
\[
p+q=\begin{cases} p+q &\hbox{ if $p,q\in \check v^+_{\omega}$ or
$p,q\in\check v^-_{\omega}$;}\\
\infty&\hbox{otherwise.}\end{cases}
\]

To describe the log structure on $V(\omega)$, we write down a chart
$\varphi:Q\rightarrow \O_U$, for $U\subseteq V(\omega)$ an open
subset. Since any chart must take $q\in Q\setminus\partial Q$ to
$0\in \O_U$, we just give values of the chart on $N\subseteq
\partial Q$:
\[
\partial Q\supseteq N \ni q\mapsto
\begin{cases}
z^q f_{\sigma}^{-\langle d_{\omega},q\rangle}& q\in \check v^+_{\omega};\\
z^q& q\in \check v^-_{\omega}.
\end{cases}
\]
Here, $f_{\sigma}$ is any extension of the section $f_{\sigma}\in
\shLS^+_{\pre,V(\sigma)}|_{V(\omega)}=\O_{\Sing(V(\omega))}$
to a function $\kk[\partial Q]$. It then follows from \cite{PartI}, Theorem 3.22
that the log structure induced by this chart is indeed given by
$f_{\sigma}\in\Gamma(V(\omega),\shLS^+_{\pre,V(\sigma)})$. This chart can also
be considered as giving a map $\varphi:Q\rightarrow\M_U$.

As $\dlog\rho$ is a globally defined section of $\Omega^1$,
$\Gamma_{\omega}(\dlog\rho)=\dlog\rho$. On the other hand,
let $n\in\check\Lambda_{v^-_{\omega}}=N$, and let us calculate
$\Gamma_{\omega}(\dlog n)$. The monoids
$\check v_{\omega}^{\pm}\cap N$ are viewed as submonoids of $N$, and
can be identified with the faces of
$Q$ by the projection $Q\subseteq N\oplus\ZZ\rightarrow N$.
By replacing $n$ by $-n$ if necessary, we can assume $n\in \check v_{\omega}^-
\cap N$, so identifying $n$ with an element of $Q$ means lifting $n$ to
$(v^-_{\omega},1)^{\perp}\subseteq N_{\RR}\oplus\RR$: this lifting is
$n-\langle v^-_{\omega},n\rangle\rho\in Q$, if we view $v^{\pm}_{\omega}$
as the endpoints of $\omega$ in $M$. Under the isomorphism of Lemma
\ref{OXD}, $\dlog n\in \Gamma(V_{e^-_{\omega}},
\Omega^1_{X^{\dagger}_{v^-_{\omega}}/\kk})$ is identified with $\dlog
\varphi(n-\langle v^-_{\omega},n\rangle\rho)\in \Gamma(V_{e^-_{\omega}},
\Omega^1_{v^-_{\omega}})$.
On the other hand, $-n\in \check v^+_{\omega}$, and
\[
n-\langle v^-_{\omega},n\rangle \rho=-(-n-\langle v^+_{\omega},-n\rangle
\rho)+\langle v^+_{\omega}-v^-_{\omega},n\rangle\rho,
\]
so on $V_{e^+_{\omega}}$,
\[\dlog\varphi(n-\langle v^-_{\omega},n\rangle\rho)=
-\dlog\varphi(-n-\langle v^+_{\omega},-n\rangle\rho)
+\langle v^+_{\omega}-v^-_{\omega},n\rangle\dlog\rho
\]
which is identified with
\[
-\dlog(z^{-n}f_{\sigma}^{\langle d_{\omega},n\rangle})
+\langle v^+_{\omega}-v^-_{\omega},n\rangle\dlog\rho
=-\langle d_{\omega},n\rangle {df_{\sigma}\over f_{\sigma}}+\dlog n
+\langle v^+_{\omega}-v^-_{\omega},n\rangle\dlog\rho
\]
in $\Gamma(V_{e^+_{\omega}},\Omega^1_{X^{\dagger}_{e^+_{\omega}}/\kk})$.
Note that $v^+_{\omega}-v^-_{\omega}=-l_{\omega}d_{\omega}$.
Thus for $n=n_1\wedge\cdots\wedge n_r$,
\begin{eqnarray*}
\Gamma_{\omega}(\dlog n)&=&
\bigwedge_{i=1}^r \left(-{df_{\sigma}\over f_{\sigma}}\langle d_{\omega}, n_i\rangle
+\dlog n_i-\langle d_{\omega},n_i\rangle l_{\omega}\dlog\rho\right)\\
&=&-\left({df_{\sigma}\over f_{\sigma}}+l_{\omega}\dlog\rho\right)
\wedge\dlog(\iota(d_{\omega})n)+\dlog n.
\end{eqnarray*}
\qed

Next we describe $\Omega^r_{\tau}$ for $\tau\in\P$ arbitrary.
We shall do this by picking a reference vertex $v\in\P$ with a morphism 
$g:v\rightarrow\tau$. We know by Proposition \ref{Omegainc} that
there is an inclusion of $\Omega^r_{\tau}$ in $F_{v,\tau}^*\Omega^r_v$,
so we only need to describe this subsheaf.

Recall that as we are assuming $(B,\P)$ is simple, for every $\tau\in\P$
with $\dim\tau\not=0,n$ we have as in \cite{PartI}, Definition 1.60,
the following data:
\begin{eqnarray*}
\P_1(\tau)&=&\{\omega\rightarrow\tau\,|\,\dim\omega=1\}\\
\P_{n-1}(\tau)&=&\{\tau\rightarrow\rho\,|\,\dim\rho=n-1\}
\end{eqnarray*}
Simplicity allows us to find disjoint sets
\begin{eqnarray*}
\Omega_1,\ldots,\Omega_q&\subseteq&\P_1(\tau),\\
R_1,\ldots,R_q&\subseteq&\P_{n-1}(\tau),
\end{eqnarray*}
and polytopes
\begin{eqnarray*}
\Delta_1,\ldots,\Delta_q&\subseteq&\Lambda_{\tau,\RR},\\
\check\Delta_1,\ldots,\check\Delta_q&\subseteq&\Lambda_{\tau,\RR}^{\perp}.
\end{eqnarray*}
These have the property that if $\omega\in\Omega_i$, $e:\omega\rightarrow
\tau$,
then the monodromy polytope $\check\Delta_e(\tau)=\check\Delta_i$,
and if $\rho\in R_i$,
$f:\tau\rightarrow\rho$, then $\Delta_f(\tau)=\Delta_i$. (See 
\cite{PartI}, Definition 1.58). Furthermore, the $\Delta_i$'s are
the Newton polytopes of the functions $\check\psi_i$ on $\check\Sigma_{\tau}$,
the normal fan to $\tau$. (See \cite{PartI}, Remark~1.59, where there is  a
typo: $\varphi_{\rho}$ should be $\check\psi_{\rho}$.)
For any $g':v'\rightarrow\tau$, we obtain vertices $\Vert_i(g')$ of $\Delta_i$
as in Construction \ref{keyP}.
The reference vertex $g:v\rightarrow\tau$ then
gives reference vertices $v_i:=\Vert_i(g)\in\Delta_i$. The sets $\Omega_i$
are characterized by $\omega\in\Omega_i$ if and only if $\Vert_i(v^+_{\omega})
\not=\Vert_i(v^-_{\omega})$.

In addition, simplicity includes the condition that
the convex hulls of
\[
\hbox{$\bigcup_{i=1}^q \Delta_i\times\{e_i\}$ and
$\bigcup_{i=1}^q\check\Delta_i\times\{e_i\}$}
\]
in  $\Lambda_{\tau,\RR}\times\RR^q$ and $\Lambda_{\tau,\RR}^{\perp}
\times\RR^q$ respectively are elementary simplices. 
In particular,
$\Delta_1,\ldots,\Delta_q$
and $\check\Delta_1,\ldots,\check\Delta_q$ are themselves elementary
simplices, and their tangent spaces $T_{\Delta_1},\ldots,T_{\Delta_q}$
give a direct sum decomposition of $\sum_{i=1}^q T_{\Delta_i}\subseteq
\Lambda_{\tau,\RR}$ and
$T_{\check\Delta_1},\ldots,T_{\check\Delta_q}$ give a direct sum decomposition
of $\sum_{i=1}^q T_{\check\Delta_i}\subseteq \Lambda_{\tau,\RR}^{\perp}$.

We will often use the obvious

\begin{lemma}
\label{Mink}
If the convex hull of $\bigcup_{i=1}^q\Delta_i\times\{e_i\}$
is an elementary simplex, then
there is a one-to-one correspondence between faces $\sigma$ of
$\Delta_{\tau}:=\Delta_1+\cdots+\Delta_q$ (Minkowski sum) and $q$-tuples
$(\sigma_1,\ldots,\sigma_q)$ with $\sigma_i$ a face of $\Delta_i$,
with $\sigma=\sigma_1+\cdots+\sigma_q$. Furthermore,
\[
\dim\sigma=\sum_{i=1}^q\dim\sigma_i.
\]
\qed
\end{lemma}

According to \cite{PartI}, Corollary 5.8, 
$q_{\tau}^{-1}(Z)=Z_1^{\tau}\cup\cdots\cup Z_q^{\tau}\cup
Z'$ where $Z'\subseteq D_{\tau}$ is of codimension at least two in 
$X_{\tau}$ and $Z_i^{\tau}$ is a hypersurface in $X_{\tau}$,
with Newton polytope $\check\Delta_i$. Furthermore, from
the proof of \cite{PartI}, Corollary 5.8, $Z_i^{\tau}=
F_{\omega,\tau}^{-1}
(Z_{\omega})$, for any $\omega\in\Omega_i$, 
where $Z_{\omega}$ is the irreducible component of $Z$ contained in
the codimension one stratum $X_{\omega}$ of $X$.

For an index set $I\subseteq \{1,\ldots,q\}$, set $Z^{\tau}_I:=
\bigcap_{i\in I} Z^{\tau}_i$.
Pull back the log structure $X_v^{\dagger}$ on $X_v$ to $X_{\tau}$
via $F_{v,\tau}$, and
then restrict further to $Z_I^{\tau}$, for any $I$. We write these 
structures as $X_{\tau}^{\dagger}$ and $(Z_I^{\tau})^{\dagger}$, but keep
in mind these are not intrinsic and depend on the choice of vertex
$g:v\rightarrow\tau$. Note that these are all defined over
$\Spec \kk^{\dagger}$, by composing the inclusions into $X_v$ 
with $X_v^{\dagger}\mapright{q_v}X_0^{\dagger}\rightarrow
\Spec \kk^{\dagger}$.

Viewing $Z_I^{\tau}\subseteq X_v$ via the inclusion $F_{v,\tau}:X_{\tau}
\rightarrow X_v$, we have 

\begin{lemma}
\label{OmegaZ}
\begin{enumerate}
\item
There are exact sequences
\[
0\rightarrow\bigoplus_{i\in I}\O_{Z_I^{\tau}}(-Z_i^{\tau})
\rightarrow \Omega^1_{X_v^{\dagger}/\kk}|_{Z^{\tau}_I}\rightarrow
\Omega^1_{(Z_I^{\tau})^{\dagger}/\kk}\rightarrow 0
\]
and
\[
0\rightarrow\bigoplus_{i\in I}\O_{Z^{\tau}_I}(-Z^{\tau}_i)\rightarrow 
\Omega^1_{X_v^{\dagger}/\kk^{\dagger}}|_{Z^{\tau}_I}\rightarrow
\Omega^1_{(Z^{\tau}_I)^{\dagger}/\kk^{\dagger}}\rightarrow 0.
\]
Here $\O_{Z_I^{\tau}}(-Z_i^{\tau})$ denotes the restriction of 
the line bundle $\O_{X_{\tau}}(-Z^{\tau}_i)$ to $Z^{\tau}_I$.
In addition, $\Omega^1_{(Z_I^{\tau})^{\dagger}/\kk}$ and 
$\Omega^1_{(Z_I^{\tau})^{\dagger}/\kk^{\dagger}}$ are locally free
$\O_{Z_I^{\tau}}$-modules.
\item
If $Y\subseteq X_{\tau}$ is a toric stratum, then
$\Omega^r_{(Z_I^{\tau})^{\dagger}/\kk}|_Y=\Omega^r_{(Z_I^{\tau}\cap Y)^{\dagger}/\kk}$
and 
\[
\Tor^{\O_{X_{\tau}}}_j(\Omega^r_{(Z_I^{\tau})^{\dagger}/\kk},\O_Y)=0
\]
for $j>0$.
Here the log structure on $Z_I^{\tau}\cap Y$ is the pull-back of the one on
$Z_I^{\tau}$. The same holds for the $/\kk^{\dagger}$ case.
\end{enumerate}
\end{lemma}

\proof Again we do the $/\kk$ case. There is always a functorial map
$\Omega^1_{X_v^{\dagger}/\kk}|_{Z^{\tau}_I}\rightarrow \Omega^1_{(Z^{\tau}_I)^{\dagger}/\kk}$ as
$(Z^{\tau}_I)^{\dagger}\rightarrow X_v^{\dagger}$ is a morphism of log schemes. We
also have a map $\shI_{Z_I^{\tau}/X_v}\rightarrow\Omega^1_{X_v^{\dagger}/\kk}|_{Z^{\tau}_I}$
by $f\mapsto df$. Then we first check exactness of
\[
\shI_{Z^{\tau}_I/X_v}\rightarrow \Omega^1_{X_v^{\dagger}/\kk}|_{Z^{\tau}_I}\rightarrow
\Omega^1_{(Z^{\tau}_I)^{\dagger}/\kk}\rightarrow 0
\]
on affine pieces. Choosing $h:\tau\rightarrow\sigma\in\P_{\max}$,
we obtain an affine open set $V_{h\circ g}=\Spec \kk[P_{h\circ g}]$ 
of $X_v$, 
where $P_{h\circ g}$ is the maximal proper face of $P_{\sigma}$ corresponding
to $h\circ g:v\rightarrow\sigma$. 
The log structure
$X_v^{\dagger}$ is given by the chart $P_{\sigma}\rightarrow \kk[P_{h\circ g}]$
that maps $p\mapsto z^p$ if $p\in P_{h\circ g}$ and $p\mapsto 0$ otherwise.
Then
\[
\Omega^1_{X_v^{\dagger}/\kk}=(\Omega^1_{X_v/\kk}\oplus (\O_{X_v}\otimes P_{\sigma}^{\gp}))/\shK_{X_v}
\]
where $\shK_{X_v}$ is the submodule generated by $(d(z^p),-z^p\cdot p)$ for 
$p\in P_{\sigma}$,
and 
\[
\Omega^1_{(Z^{\tau}_I)^{\dagger}/\kk}=(\Omega^1_{Z^{\tau}_I/\kk}\oplus (\O_{Z^{\tau}_I}
\otimes P_{\sigma}^{\gp}))/\shK_{Z^{\tau}_I}
\]
where $\shK_{Z^{\tau}_I}$ is defined analogously. Then we have a diagram
\[
\xymatrix@C=30pt
{&&\shI_{Z^{\tau}_I/X_v}\ar[d]^{(d,0)}&&\\
&\shK_{X_v}|_{Z^{\tau}_I}\ar[r]\ar[d]&\Omega^1_{X_v/\kk}|_{Z^{\tau}_I}\oplus (\O_{Z^{\tau}_I}\otimes P_{\sigma}^{\gp})\ar[d]\ar[r]&
\Omega^1_{X_v^{\dagger}/\kk}|_{Z^{\tau}_I}\ar[r]\ar[d]&0\\
0\ar[r]&\shK_{Z^{\tau}_I}\ar[r]&\Omega^1_{Z^{\tau}_I/\kk}\oplus (\O_{Z^{\tau}_I}\otimes P_{\sigma}^{\gp})\ar[r]\ar[d]&
\Omega^1_{(Z^{\tau}_I)^{\dagger}/\kk}\ar[r]\ar[d]&0\\
&&0&0&}
\]
where the map $\shK_{X_v}|_{Z^{\tau}_I}\rightarrow\shK_{Z^{\tau}_I}$ is clearly surjective,
and the middle column is exact from the standard result for K\"ahler
differentials. So we see
\[
\shI_{Z^{\tau}_I/X_v}\mapright{d}\Omega^1_{X_v^{\dagger}/\kk}|_{Z^{\tau}_I}\mapright{}
\Omega^1_{(Z^{\tau}_I)^{\dagger}/\kk}\mapright{} 0
\]
is exact by the snake lemma. 

Now view $Z^{\tau}_I\subseteq X_{\tau}\subseteq X_v$.
Looking at the composition $\shI_{X_{\tau}/X_v}\hookrightarrow
\shI_{Z^{\tau}_I/X_v}\mapright{d} \Omega^1_{X_v^{\dagger}/\kk}|_{Z^{\tau}_I}$,
we see that as the former ideal is generated by monomials $z^p$ for
various $p$, and $d(z^p)=z^p\dlog(p)$, it follows that $\shI_{X_{\tau}/X_v}$
maps to zero in $\Omega^1_{X_v^{\dagger}/\kk}|_{Z^{\tau}_I}$, and thus 
the map 
\[
\shI_{Z^{\tau}_I/X_v}\mapright{d} \Omega^1_{X_v^{\dagger}/\kk}|_{Z^{\tau}_I}
\]
factors through $\shI_{Z^{\tau}_I/X_v}\rightarrow \shI_{Z^{\tau}_I/X_{\tau}}$, so
we have an exact sequence
\[
\shI_{Z^{\tau}_I/X_{\tau}}/\shI_{Z^{\tau}_I/X_{\tau}}^2
\rightarrow
\Omega^1_{X_v^{\dagger}/\kk}|_{Z^{\tau}_I}\rightarrow \Omega^1_{(Z^{\tau}_I)^{\dagger}/\kk}
\rightarrow 0.
\]
The conormal bundle is of course $\bigoplus_{i\in I}
\O_{Z^{\tau}_I}(-Z^{\tau}_i)$, and
we just need to check injectivity. This can be done locally,
where $X_v=\Spec \kk[P]$, $X_{\tau}=\Spec \kk[Q]$ for $Q$ a face of $P$,
and $Z^{\tau}_I$ is defined by equations $\{f_i=0|i\in I\}$ in $X_{\tau}$,
with $f_i$ having, up to translation, Newton polytope $\check\Delta_i$.
We can write
$f_i=1+\sum_{j=1}^{l_i} a_{ij} z^{q_{ij}}$ where the set $\{q_{ij}\}$
is linearly independent in $Q^{\gp}\otimes_{\ZZ}\QQ$ by simplicity. 
Now $df_i=\sum_{j=1}^{l_i} a_{ij} z^{q_{ij}}\dlog q_{ij}$.
We can't have all $z^{q_{ij}}$'s vanishing at a point of $Z^{\tau}_i$,
so $df_i$ is non-vanishing, and
by the linear independence of the $q_{ij}$'s, 
$\{df_i|i\in I\}$ are an everywhere
linearly independent set of sections
of $\Omega^1_{X_v^{\dagger}/\kk}|_{Z^{\tau}_I}$.
Thus the map $\shI_{Z^{\tau}_I/X_{\tau}}/\shI_{Z^{\tau}_I/X_{\tau}}^2\mapright{d}
\Omega^1_{X_v^{\dagger}/\kk}|_{Z^{\tau}_I}$ is injective, and the cokernel
is locally free. (Recall that $\Omega^1_{X_v^{\dagger}/\kk}$ is locally
free.) This shows (1).

For (2), just restrict the sequences of (1) of locally free sheaves to
$Z_I^{\tau}\cap Y$. 
This gives the sequence of (1) 
for $\Omega^1_{(Z_I^{\tau}\cap Y)^{\dagger}/\kk}$, proving
$\Omega^1_{(Z^{\tau}_I\cap Y)^{\dagger}/\kk}=\Omega^1_{(Z_I^{\tau})^{\dagger}/
\kk}|_Y$
(or $/\kk^{\dagger}$).
To show the vanishing of the $\Tor$'s, by the local freeness statement of
(1), it is enough to show that $\Tor_j^{\O_{X_{\tau}}}(\O_{Z_I^{\tau}},
\O_Y)=0$
for $j>0$. To do this, note using simplicity that $Z_I^{\tau}$ is a complete
intersection, so we have the Koszul complex
\begin{equation}
\label{koszul}
0\rightarrow {\bigwedge}^{\# I}\bigoplus_{i\in I}\O_{X_{\tau}}(-Z_i^{\tau})
\rightarrow\cdots\rightarrow {\bigwedge}^1 \bigoplus_{i\in I}\O_{X_{\tau}}
(-Z_i^{\tau})\rightarrow \O_{X_{\tau}}\rightarrow\O_{Z_I^{\tau}}
\rightarrow 0
\end{equation}
giving a resolution for $\O_{Z_I^{\tau}}$. Tensoring with $\O_Y$ gives the
Koszul resolution for $Y\cap Z_I^{\tau}$, 
which is still a complete intersection in $Y$. This gives the result.
\qed

\begin{remark}
\label{poleremark}
We shall use on a number of occasions variations on
 the following simple observation: given a form $\alpha\in \Omega^r_{X_v^{\dagger}/
\kk}$ or $\Omega^r_{X_v^{\dagger}/\kk^{\dagger}}$, if $f_j=0$ locally defines
$Z_j$ in an affine subset of $X_v$, then ${df_j\over f_j}\wedge\alpha$
has no pole along $Z_j$ if and only if $\alpha|_{Z_j}=0$. This follows
from linear algebra: if $V$ is a vector space, $\beta\in V^*$, $\alpha\in \bigwedge^r
V^*$, then $\beta\wedge\alpha=0$ if and only if $\alpha|_{\ker\beta}=0$.
\end{remark}

\begin{proposition}
\label{delta0kernel}
Given $v\rightarrow\tau_1\rightarrow\tau_2$, the image of the inclusion
$(F_{\tau_1,\tau_2}^*\Omega^r_{\tau_1})/\Tors$ in $F_{v,\tau_2}^*
\Omega^r_v$ is
\[
\ker\bigg(
F_{v,\tau_2}^*\Omega^r_v\mapright{\delta_0}
\bigoplus_{i=1,\ldots,q\atop
w_i\not=v_i} \Omega^{r-1}_{(Z_i^{\tau_2})^{\dagger}/\kk}\bigg)
\]
or
\[
\ker\bigg(
F_{v,\tau_2}^*\Omega^r_v\mapright{\delta_0}
\bigoplus_{i=1,\ldots,q\atop
w_i\not=v_i} \Omega^{r-1}_{(Z_i^{\tau_2})^{\dagger}/\kk^{\dagger}}\bigg)
\]
in the $/\kk$ and $/\kk^{\dagger}$ cases respectively,
where: 
\begin{enumerate}
\item 
The direct sum is over all $i$ and all vertices $w_i$
of $\Delta_i$, $w_i\not=v_i$, and
$\Delta_1,\ldots,\Delta_q$ are part of the simplicity data for
$\tau_1$.
\item $Z_i^{\tau_2}=F_{\tau_1,\tau_2}^{-1}(Z_i^{\tau_1})$ where $Z_1^{\tau_1},
\ldots,Z_q^{\tau_1}$ are as usual the codimension one irreducible components
of $q_{\tau_1}^{-1}(Z)$ with Newton polytopes $\check\Delta_1,\ldots,\check
\Delta_q$.
\item
For $\alpha\in F_{v,\tau_2}^*\Omega^r_v$, the component of
$\delta_0(\alpha)$ in the direct 
summand $\Omega^{r-1}_{(Z_i^{\tau_2})^{\dagger}/\kk}$
or $\Omega^{r-1}_{(Z_i^{\tau_2})^{\dagger}/\kk^{\dagger}}$
corresponding to some $w_i$ is given by
$\iota(\partial_{w_i-v_i})\alpha|_{(Z^{\tau_2}_i)^{\dagger}}$. 
\end{enumerate}
\end{proposition}

\proof 
Again we'll do the $/\kk$ case. The result will follow by using
the characterization of Corollary \ref{inclusionintersectioncor}, (2)
of $(F_{\tau_1,\tau_2}^*\Omega^r_{\tau_1})/\Tors$ as
\[
\bigcap_{g':v'\rightarrow\tau_1} F_{v',\tau_2}^*\Omega^r_{v'},
\]
using Lemma \ref{basicgluing} for the explicit identification
of this intersection with a subsheaf of $F_{v,\tau_2}^*\Omega^r_v$.
Let $\alpha$ be a section of $F_{v,\tau_2}^*\Omega^r_v$. Then
for any $j$ and vertex $w_j\not=v_j$ of $\Delta_j$, we can
find a sequence of edges $h_i:\omega_i\rightarrow\tau_1$, $i=1,\ldots,m$
of $\tau_1$, with $d_{\omega_i}$ chosen appropriately, so that
\begin{itemize}
\item $v^-_{\omega_1}=v$;
\item $v^+_{\omega_i}=v^-_{\omega_{i+1}}$ for $i<m$;
\item $\Vert_l(v^+_{\omega_i})=v_l$ for $i<m$, for all $l$;
\item
$\Vert_l(v^+_{\omega_m})=\begin{cases} v_l&l\not=j,\\ w_j&l=j.\end{cases}$
\end{itemize}
Choose a maximal cell $\sigma$ containing $\tau_2$ for reference,
and let $f_1,\ldots,f_q$ be the equations defining $Z_1^{\tau_2},
\ldots,Z_q^{\tau_2}$ in the affine chart $V(\sigma)\cap X_{\tau_2}$.
Using Lemma \ref{basicgluing}, apply $\Gamma_{\omega_1},\ldots,
\Gamma_{\omega_m}$
successively to $\alpha$. Then in fact, for $1\le i\le m-1$,
\[
\Gamma_{\omega_i}(\dlog n)=-l_{\omega_i}\dlog\rho\wedge\dlog(\iota(d_{\omega_i})n)
+\dlog n.
\]
(Here and throughout this proof we write $d_{\omega_i}$ rather than
$\partial_{d_{\omega_i}}$ etc.\ for typographical convenience.)
Indeed, $\Vert_l(v_{\omega_i}^+)=\Vert_l(v^-_{\omega_i})$ for each
$1\le i\le m-1$, so $\omega_1,\ldots,\omega_{m-1}$ are not in any
$\Omega_i$. Thus $F_{\omega_i,\tau_2}^{-1}(Z_{\omega_i})=\emptyset$
(see discussion after Lemma \ref{Mink} or \cite{PartI}, Cor. 5.8), so the
function $f_{\sigma}$ appearing in Lemma \ref{basicgluing} is constant and hence
does not play a role. On the other hand, $\omega_m\in\Omega_j$, so
\[
\Gamma_{\omega_m}(\dlog n)=-\left({df_j\over f_j}+l_{\omega_m}\dlog\rho\right)\wedge
\dlog(\iota(d_{\omega_m})n)+\dlog n.
\]
Putting this together, we see that 
\[
\Gamma_{\omega_{m-1}}\circ\cdots\circ\Gamma_{\omega_1}(\alpha)
=\dlog\rho\wedge\iota(v_{\omega_{m-1}}^+-v)\alpha+\alpha
\]
and
\[
\Gamma_{\omega_{m}}\circ\cdots\circ\Gamma_{\omega_1}(\alpha)
=
-{df_j\over f_j}\wedge\iota(d_{\omega_m})[\dlog\rho\wedge
\iota(v^+_{\omega_{m-1}}-v)\alpha+\alpha]
+\dlog\rho\wedge\iota(v_{\omega_m}^+-v)\alpha+\alpha.
\]
Now ${df_j\over f_j}\wedge\iota(d_{\omega_m})\alpha$ has no pole along
$f_j=0$ if and only if $(\iota(d_{\omega_m})\alpha)|_{Z_j^{\tau_2}}=0$
in $\Omega^{r-1}_{(Z_j^{\tau_2})^{\dagger}/\kk}$, by Remark~\ref{poleremark}. 
If this latter
condition holds, then of course $\iota(d_{\omega_m}\wedge
(v^+_{\omega_{m-1}}-v))\alpha|_{Z_j^{\tau_2}}=0$ also, so $\Gamma_{\omega_m}
\circ
\cdots\circ\Gamma_{\omega_1}(\alpha)$ has no pole if and only if 
$\iota(v_j-w_j)\alpha|_{Z_j^{\tau_2}}=0$. Hence we see that
\[
\bigcap_{v'\rightarrow\tau_1} F_{v',\tau_2}^*\Omega^r_{v'}
\subseteq\ker\delta_0.
\]

Conversely, if $\alpha\in\ker\delta_0$, let $v'$ be any vertex of $\tau_1$.
Then we can find a sequence of edges $\omega_i\rightarrow\tau_1$,
$i=1,\ldots,m$ of $\tau_1$, with $d_{\omega_i}$ chosen appropriately, so that
\begin{itemize}
\item $v^-_{\omega_1}=v$;
\item $v^+_{\omega_i}=v^-_{\omega_{i+1}}$ for $i<m$;
\item $v^+_{\omega_m}=v'$;
\item For each $1\le l\le q$, there is at most one $i$ such that
$\Vert_l(v^-_{\omega_i})\not= \Vert_l(v^+_{\omega_i})$, and for this $i$,
$\Vert_l(v^-_{\omega_i})=v_l$, $\Vert_l(v^+_{\omega_i})=v_l'=\Vert_l(v')$.
\end{itemize}
Then again using Lemma \ref{basicgluing} repeatedly along each $\omega_i$
to identify $\alpha$ with a rational section of $F_{v',\tau_2}^*
\Omega^r_{v'}$, one sees that $\alpha$ is identified with $\alpha+\cdots$,
where the terms in $\cdots$ which may not be regular sections of
$F_{v',\tau_2}^*\Omega^r_{v'}$ are of the form
\[
\bigg(\bigwedge_{i\in I}{df_i\over f_i}\bigg)\wedge\iota\bigg(
\bigwedge_{i\in I}(v_i-v_i')\bigg)\alpha
\]
or
\[
\bigg(\bigwedge_{i\in I}{df_i\over f_i}\bigg)\wedge\dlog\rho\wedge
\iota\bigg(w\wedge
\bigwedge_{i\in I}(v_i-v_i')\bigg)\alpha
\]
for $I\subseteq \{1,\ldots,q\}$ some index set and some $w\in\Lambda_{\tau_1}$.
Here we use the fact that $v_i-v_i'\in\Lambda_{\tau_1}$, while the monomials
appearing in the $f_i$'s are in $\Lambda_{\tau_2}^{\perp}$, so
$\iota(v_i-v_i')df_k=0$. However, for $i\in I$, this has no pole along
$Z_i^{\tau_2}$ again by Remark~\ref{poleremark}
because by assumption $\iota(v_i-v_i')\alpha|_{Z_i^{\tau_2}}$
is zero. 
Thus we see that $\alpha$ is identified with a regular
section of $F_{v',\tau_2}^*\Omega^r_{v'}$, and hence $\alpha$
is in $\bigcap_{g':v'\rightarrow\tau_1}
F_{v',\tau_2}^*\Omega^r_{v'}$. \qed

\bigskip

For $e:\tau_1\rightarrow\tau_2$,
we will now calculate the cohomology of $(F_{\tau_1,\tau_2}^*\Omega^r_{\tau_1})/
\Tors$ by
building a convenient resolution of this sheaf. The first two terms of
this resolution are given by Proposition~\ref{delta0kernel}; we need to
extend this two-term complex.

For $V\subseteq\Lambda_{\tau,\RR}$
a subspace, we have a subsheaf
$\Omega^r_v|_{X_{\tau}}\cap V^{\perp}$ of
$\Omega^r_v|_{X_{\tau}}$
given by forms $\alpha$
with $\iota(\partial_m)\alpha=0$ for all $m\in V$. 
We define $\Omega^r_{(Z_I^{\tau})^{\dagger}/\kk}\cap V^{\perp}$
or 
$\Omega^r_{(Z_I^{\tau})^{\dagger}/\kk^{\dagger}}\cap V^{\perp}$
to be the image of $\Omega^r_v|_{X_{\tau}}\cap V^{\perp}$ in 
$\Omega^r_{(Z_I^{\tau})^{\dagger}/\kk}$ (or $/\kk^{\dagger}$).
For $m\in\Lambda_{\tau}$, note that 
\[
\iota(\partial_m)\left(\im\left(\bigoplus_{i\in I}\O_{Z_I^{\tau}}
(-Z_i^{\tau})\mapright{d}
\Omega^1_v|_{Z_I^{\tau}}\right)\right)=0,
\]
as all monomials occuring in the equations for the $Z_i^{\tau}$'s are in 
$\Lambda_{\tau}^{\perp}$. 
We thus in particular have from Lemma \ref{OmegaZ} an exact sequence
\[
0\rightarrow \bigoplus_{i\in I}\O_{Z_I^{\tau}}(-Z_i^{\tau})\rightarrow
\Omega^1_v|_{Z_I^{\tau}}\cap V^{\perp}\rightarrow
\Omega^1_{(Z_I^{\tau})^{\dagger}/\kk}\cap V^{\perp}\rightarrow 0
\]
and a similar exact sequence for the $/\kk^{\dagger}$ case.

Given $g:v\rightarrow\tau_1$ as usual, we 
can now define a complex $\shF_v^{r,\bullet}$ by
\[
\shF^{r,p}_v=\bigoplus_{\sigma\subseteq\Delta_{\tau_1}:\bar v\in\sigma
\atop \dim\sigma=p}
\Omega^{r-p}_{(Z^{\tau_1}_{I(\sigma)})^{\dagger}/\kk}\cap T_{\sigma}^{\perp},
\]
in the $/\kk$ case and
\[
\shF^{r,p}_v=\bigoplus_{\sigma\subseteq\Delta_{\tau_1}:\bar v\in\sigma
\atop \dim\sigma=p}
\Omega^{r-p}_{(Z^{\tau_1}_{I(\sigma)})^{\dagger}/\kk^{\dagger}}
\cap T_{\sigma}^{\perp},
\]
in the $/\kk^{\dagger}$ case,
where the sum is over all $p$-dimensional $\sigma=\sigma_1+\cdots+
\sigma_q$ where $\sigma_i$ is a face of $\Delta_i$ containing $v_i$,
and
\begin{eqnarray*}
&&\hbox{$\Delta_{\tau_1}=\Delta_1+\cdots+\Delta_q$;}\\
&&\hbox{$\bar v=v_1+\ldots+v_q$;}\\
&&\hbox{$T_{\sigma}$ is the tangent space to $\sigma$ in $\Lambda_{\tau,\RR}$;}\\
&&\hbox{$I(\sigma)=\{i|\sigma_i\not=\{v_i\}\}$.}
\end{eqnarray*}
We use the convention that if $I(\sigma)=\emptyset$ then 
$\Omega^r_{(Z_{I(\sigma)}^{\tau})^{\dagger}/\kk}$ (or $/\kk^{\dagger}$)
is $\Omega^r_v|_{X_{\tau}}$.

We define differentials
$\delta_p:\shF^{r,p}_v\rightarrow\shF^{r,p+1}_v$ by
\[
(\delta_p\alpha)_{\sigma'}=\sum_{\sigma\subseteq\sigma':\bar v\in\sigma\atop
\dim\sigma=p} \iota(\partial_{w_j-v_j})
\alpha_{\sigma}|_{Z_{I(\sigma')}^{\tau_1}}.
\]
Here $\sigma'$ is a face of $\Delta_{\tau_1}$ of dimension $p+1$, and
we sum over all faces $\sigma$ of $\sigma'$ of dimension $p$ containing $v$. For
each such $\sigma'$, 
by Lemma \ref{Mink} there is a unique $j$ such that $\sigma'_j\not=\sigma_j$,
and $w_j$ is the unique vertex of $\sigma_j'$ not contained in $\sigma_j$.
By Proposition \ref{delta0kernel}, 
\[
\Omega_{\tau_1}^r=
\ker(\delta_0:\shF^{r,0}_v\rightarrow \shF^{r,1}_v).
\]

\begin{lemma}
\label{keylemma}
For any $\tau_1\subseteq\tau_2$, 
\[
F_{\tau_1,\tau_2}^* \shF_v^{r,\bullet}
\]
is a resolution of $(F_{\tau_1,\tau_2}^*\Omega^r_{\tau_1})/\Tors$.
\end{lemma}

\proof
We show this for $\tau_1=\tau_2=\tau$,
and then the complex remains a resolution under pull-back by 
Lemma \ref{OmegaZ}, (2). We will do this for the $/\kk$ case, the
$/\kk^{\dagger}$ case being essentially identical.
We will proceed by induction. Consider faces $v\in\omega\subseteq\omega'
\subseteq\Delta_{\tau}$, and consider the complex
$\shF^{\bullet}_{\omega,\omega'}$ defined by 
\[
\shF^{p}_{\omega,\omega'}=\bigoplus_{\omega\subseteq\sigma\subseteq\omega'
\atop \dim\sigma=p} \Omega^{r-p}_{(Z^{\tau}_{I(\sigma)})^{\dagger}/\kk}\cap 
T_{\sigma}^{\perp},
\]
with differential $\delta_p$ defined as before.
Note that if $\omega=v$, $\omega'=\Delta_{\tau}$, this is $\shF^{r,p}_v$.
We will show 
$H^i(\shF^{\bullet}_{\omega,\omega'})=0$ for $i>\dim\omega$ inductively
on $\dim\omega'-\dim\omega$. (Here $H^i$ denotes cohomology of the complex.)
If $\omega=\omega'$, the statement
is trivial. If $\omega\not=\omega'$, let
$w_j$ be a vertex of $\omega_j'$ not in $\omega_j$, and write
$\omega\cup\{w_j\}:=\omega_1+\cdots+\Conv\{\omega_j,w_j\}+\cdots+\omega_p$,
and $\omega'\setminus\{w_j\}:=\omega_1'+\cdots+\omega_j'\setminus\{w_j\}+
\cdots+\omega_p'$, where $\omega_j'\setminus\{w_j\}$ is the convex hull
of those vertices of $\omega_j'$ not equal to $w_j$. 
Then we have an exact sequence of complexes
\[
0\rightarrow \shF^{\bullet}_{\omega\cup\{w_j\},\omega'}
\rightarrow\shF^{\bullet}_{\omega,\omega'}\rightarrow 
\shF^{\bullet}_{\omega,\omega'\setminus\{w_j\}}\rightarrow
0,
\]
and hence an exact sequence
\[
0\mapright{} H^{\dim\omega}(\shF^{\bullet}_{\omega,\omega'})\mapright{}
H^{\dim\omega}(\shF^{\bullet}_{\omega,\omega'\setminus\{w_j\}})
\mapright{\tilde \delta}H^{\dim\omega+1}(\shF^{\bullet}_{\omega\cup\{w_j\},
\omega'})
\mapright{} H^{\dim\omega+1}(\shF^{\bullet}_{\omega,\omega'})\mapright{}
0
\]
and $H^j(\shF^{\bullet}_{\omega,\omega'})=0$ for $j>\dim\omega+1$,
by the induction hypothesis. To complete the result, we need to show the
map $\tilde\delta$ is surjective.

There are two cases. First, suppose $\omega_j\not=\{v_j\}$. Then 
$j\in I(\omega)=I(\omega\cup\{w_j\})$, and we get a diagram
with exact rows
{\scriptsize
\[
\xymatrix@C=30pt
{0\ar[r]&H^{\dim\omega}(\shF^{\bullet}_{\omega,\omega'\setminus\{w_j\}})
\ar[r]\ar[d]^{\tilde\delta}&
\Omega^{r-\dim\omega}_{(Z^{\tau}_{I(\omega)})^{\dagger}/\kk}
\cap T_{\omega}^{\perp}\ar[r]^<<<<<<{\delta_{\dim\omega}}\ar[d]&
\displaystyle{\bigoplus_{\omega\subseteq\sigma\subseteq\omega'\setminus\{w_j\}
\atop\dim\sigma=\dim\omega+1}} 
\Omega^{r-\dim\omega-1}_{(Z^{\tau}_{I(\sigma)})^{\dagger}/\kk}\cap T_{\sigma}^{\perp}
\ar[d]\\
0\ar[r]&H^{\dim\omega+1}(\shF^{\bullet}_{\omega\cup\{w_j\},\omega'})
\ar[r]&\Omega^{r-\dim\omega-1}_{(Z^{\tau}_{I(\omega)})^{\dagger}/\kk}
\cap T_{\omega\cup\{w_j\}}^{\perp}
\ar[r]^<<<<{\delta_{\dim\omega+1}}&
\displaystyle{\bigoplus_{\omega\subseteq\sigma\subseteq\omega'\setminus\{w_j\}
\atop\dim\sigma=\dim\omega+1}} 
\Omega^{r-\dim\omega-2}_{(Z^{\tau}_{I(\sigma)})^{\dagger}/\kk}\cap T_{\sigma\cup
\{w_j\}}^{\perp}
}
\]
}
The vertical maps are given by contraction by $\partial_{w_j-v_j}$. 
By simplicity, we can find
an $n\in \check\Lambda_{v,\RR}$ such that 
$\langle n,w_j-v_j\rangle =1$ and $\langle n,w_i'-v_i\rangle=0$ for all
$w_i'$ a vertex of $\Delta_i$ not equal to $w_j$.
Let $\alpha\in H^{\dim\omega+1}(\shF^{\bullet}_{\omega
\cup\{w_j\},\omega'})$. We claim that $\dlog n\wedge
\alpha\in H^{\dim\omega}(\shF^{\bullet}_{\omega,\omega'\setminus\{w_j\}})$
and $\tilde\delta(\dlog n\wedge\alpha)=\alpha$. First,
for any $\omega\subseteq\sigma\subseteq
\omega'\setminus\{w_j\}$ with $\dim\sigma=\dim\omega+1$, $\sigma_i=\omega_i$
for $i\not=l$ and
$\sigma_l=\omega_l\cup\{w_l\}$ for some unique $l$ and vertex
$w_l$ of $\omega'_l\setminus\{w_j\}$. Computing the component of
$\delta_{\dim\omega}(\dlog n\wedge\alpha)$ labelled by $\sigma$ gives 
\[
\iota(\partial_{w_l-v_l})(\dlog n\wedge\alpha)
=-\dlog n\wedge (\iota(\partial_{w_l-v_l})\alpha),
\]
which restricts to zero on $(Z^{\tau}_{I(\sigma)})^{\dagger}$ because
$\iota(\partial_{w_l-v_l})\alpha$ does. Thus 
\[
\dlog n\wedge\alpha
\in
H^{\dim\omega}(\shF^{\bullet}_{\omega,\omega'\setminus\{v_j\}}).
\]
On the other hand 
$\tilde\delta(\dlog n\wedge\alpha)=\alpha$,
as $\partial_{w_j-v_j}\in T_{\omega\cup\{w_j\}}$, so $\iota(
\partial_{w_j-v_j})\alpha=0$.
This shows surjectivity of
$\tilde\delta$ in this case.

In the second case, suppose $\omega_j=\{v_j\}$. Then $j\not\in I(\omega)$,
and we get a diagram
{\scriptsize
\[
\xymatrix@C=30pt
{0\ar[r]&H^{\dim\omega}(\shF^{\bullet}_{\omega,\omega'\setminus\{w_j\}})
\ar[r]\ar[d]^{\tilde\delta}&
\Omega^{r-\dim\omega}_{(Z^{\tau}_{I(\omega)})^{\dagger}/\kk}
\cap T_{\omega}^{\perp}\ar[r]\ar[d]&
\displaystyle{\bigoplus_{\omega\subseteq\sigma\subseteq\omega'\setminus\{w_j\}
\atop\dim\sigma=\dim\omega+1}} 
\Omega^{r-\dim\omega-1}_{(Z^{\tau}_{I(\sigma)})^{\dagger}/\kk}
\cap T_{\sigma}^{\perp}
\ar[d]\\
0\ar[r]&H^{\dim\omega+1}(\shF^{\bullet}_{\omega\cup\{w_j\},\omega'})
\ar[r]&\Omega^{r-\dim\omega-1}_{(Z^{\tau}_{I(\omega)}\cap Z^{\tau}_j)^{\dagger}
/\kk} \cap T_{\omega \cup\{w_j\}}^{\perp}
\ar[r]&
{\bigoplus_{\omega\subseteq\sigma\subseteq\omega'\setminus\{w_j\}
\atop\dim\sigma=\dim\omega+1}}
\Omega^{r-\dim\omega-2}_{(Z^{\tau}_{I(\sigma)}\cap Z^{\tau}_j)^{\dagger}/\kk}
\cap T_{\sigma\cup \{w_j\}}^{\perp}
}
\]
}
Note the differences between this and the previous diagram: the sheaves
of log differentials in the lower line are now over subschemes of those
on the upper line.
We need a more explicit description of the kernels. 
In particular, we will describe a set of generators of
$H^{\dim\omega+1}(\shF^{\bullet}_{\omega\cup\{w_j\},\omega'})$
as a sheaf of modules and give liftings of each of these generators.
To do so,
choose dual bases for
$\Lambda_{v,\RR}\oplus(\RR\rho)^*$ and $\check\Lambda_{v,\RR}\oplus\RR\rho$
as follows. Start with $e_0$ a generator of $(\RR\rho)^*$, throw in
a basis for $\sum_{i=1}^q T_{\Delta_i}
\subseteq \Lambda_{\tau,\RR}\subseteq\Lambda_{v,\RR}$ 
with elements $\partial_{w-v_i}$, where $w$ is a vertex of $\Delta_i$,
$w\not=v_i$, and $i$ runs from $1$ to $q$.
Extend to a basis for $\Lambda_{v,\RR}\oplus(\RR\rho)^*$. 
Call this basis $e_0,\ldots,
e_n$, with dual basis $e_0^*,\ldots,e_n^*$. Consider for $i=1,\ldots,q$
the following subsets of $\{0,\ldots,n\}$:
\begin{eqnarray*}
I_i&=&\{l|\hbox{$e_l+v_i$ is a vertex of $\omega_i'$, but not $\omega_i$,
and $e_l+v_i\not=w_j$}\}\\
K_i&=&\{l|\hbox{$e_l+v_i$ is a vertex of $\omega_i$, or $e_l+v_i=w_j$}\}.
\end{eqnarray*}
Now if $\alpha\in H^{\dim\omega+1}(\shF^{
\bullet}_{\omega\cup\{w_j\},\omega'})$, let's examine the property
$\delta_{\dim\omega+1}(\alpha)=0$ component by component. We have
a component for each $l\in I_i$, $i=1,\ldots,q$. Then
either $i\in I(\omega\cup\{w_j\})$ and
$\iota(e_l)\alpha=0$ because $Z^{\tau}_{I(\omega\cup\{e_l+v_i\})}\cap Z^{\tau}_j
=Z^{\tau}_{I(\omega)}\cap Z^{\tau}_j$, 
or else $i\not\in I(\omega\cup\{w_j\})$ and 
$\iota(e_l)\alpha|_{Z_{I(\omega)}^{\tau}\cap Z^{\tau}_i\cap Z^{\tau}_j}=0$. 
From this observation
we will describe a set of generators for the sheaf 
$H^{\dim\omega+1}(\shF^{\bullet}_{\omega\cup\{w_j\},\omega'})$
as an $\O_{Z^{\tau}_{I(\omega)}\cap Z^{\tau}_j}$-module. For an index set 
$J=\{j_1,\ldots,j_r\}\subseteq\{0,\ldots,n\}$, 
write $e^*_J=e^*_{j_1}\wedge\cdots\wedge e^*_{j_r}$.
Let $I=\bigcup_{i\not\in I(\omega\cup\{w_j\})} I_i$, $K=\bigcup_{i=1}^q K_i$.
Then a basis for $T_{\omega\cup\{w_j\}}^{\perp}$ in $\bigwedge^{\bullet}
(\check\Lambda_{v,\RR}\oplus\RR\rho)$ is 
\[
\{e^*_L| \hbox{$L$ an index set with $L\cap K=\emptyset$}\}.
\]
Working on affine charts, let $f_1,\ldots,f_q$ be the functions on
an affine chart of 
$X_{\tau}$ defining $Z^{\tau}_1,\ldots,Z^{\tau}_q$. 
Then we will show locally a generating set for
$H^{\dim\omega+1}(\shF^{\bullet}_{\omega\cup\{w_j\},\omega'})$
consists of forms of the shape
\[
\bigg(\prod_{i\in M} f_i\bigg)
\dlog(e^*_L)\wedge \dlog(e^*_N)\wedge \bigwedge_{i\in
O\setminus M} df_i,
\]
where: 
\begin{eqnarray*}
&&\hbox{$L\subseteq I$;}\\
&&\hbox{$N\subseteq\{0,\ldots,n\}\setminus (K\cup\bigcup_{i=1}^q I_i);$}\\
&&\hbox{$O=\{i|L\cap I_i\not=\emptyset\}$};\\
&&\hbox{$M\subseteq O$.}
\end{eqnarray*}
Indeed, first note such an element is in $\ker(\delta_{\dim\omega+1})$:
given $l\in I_i$ for some $i$, if $i\in I(\omega\cup\{w_j\})$
then $e_l^*$ does not occur in $e_L^*$ or $e_N^*$ or any $df_j$,
so the contraction with $e_l$ is zero. If $i\not\in I(\omega\cup
\{w_j\})$, then either $l\not\in L$ and the contraction is zero, or else
$l\in L$ and the contraction restricts to zero on $Z^{\tau}_{I(\omega)}
\cap Z_i^{\tau}\cap Z_j^{\tau}$, either because $i\in M$
or $i\in O\setminus M$.

Conversely, we need to show this set generates $\ker(\delta_{\dim\omega+1})$.
If $\alpha\in\ker(\delta_{\dim\omega+1})$, then any term involving
$\dlog e_{L'}^*$ appearing in $\alpha$ must satisfy $L'\cap K=\emptyset$.
Furthermore, as the contraction of $\alpha$ with $e_l$ for $l\in I_i$,
$i\in I(\omega\cup\{w_j\})$ is zero, such an $l$ cannot appear in $L'$.
Thus we can decompose $L'$ as $L\cup N$ as claimed. Furthermore,
the condition that $\iota(e_l)\alpha|_{Z^{\tau}_{I(\omega)}\cap Z_i^{\tau}
\cap Z_j^{\tau}}=0$ for $l\in I_i$, $i\not\in I(\omega\cup\{w_j\})$ implies
that $\alpha$ must be a linear combination of the forms of the type given, with 
coefficients in $\O_{Z^{\tau}_{I(\omega)}\cap Z_j^{\tau}}$.

We now have to lift these generators.
So assume furthermore that $e_1=w_j-v_j$. Then 
the above generator
can be lifted to $\Omega_{(Z^{\tau}_{I(\omega)})^{\dagger}/\kk}^{r-\dim\omega}
\cap T_{\omega}^{\perp}$ as $\big(\prod_{i\in M}f_i\big)\dlog(e_1^*)\wedge\dlog(e^*_L)
\wedge\dlog(e^*_N)\wedge\bigwedge_{i\in O\setminus M}df_i$, which as
before is easily checked to be in $\ker(\delta_{\dim\omega})$. 
This shows surjectivity. \qed

\subsection{The Hodge decomposition}
\label{Hodgedecompsection}

We continue with the notation of the previous section, and will use
the technical results of that section to finish the calculation
of the Hodge decomposition.
Having now constructed the resolution of the sheaves
$(F_{\tau_1,\tau_2}^*\Omega^r_{\tau_1})/\Tors$, we wish to use these resolutions
to compute the cohomology of these sheaves. As a first step, we calculate
the cohomology of the individual sheaves appearing in the resolution.
Under important additional hypotheses, these sheaves have no cohomology
in degree $\ge 1$, and their global sections are easily expressed in
terms of data on $B$.

\begin{lemma}
\label{vanishing}
Suppose that for the cell $\tau\in\P$, 
$\Conv(\bigcup_{i=1}^q\check\Delta_i\times\{e_i\})$
is a standard simplex (as opposed to just an elementary simplex). Then
\begin{enumerate}
\item For $\sigma\subseteq\Delta_{\tau}$ a face, 
\[
\Gamma(X_{\tau},\Omega^r_{(Z_I^{\tau})^{\dagger}/\kk^{\dagger}}
\cap T_{\sigma}^{\perp})
= {{\bigwedge}^r T_{\sigma}^{\perp}\over
\Top(I)_r}\otimes \kk,
\]
for $T_{\sigma}^{\perp}\subseteq \check\Lambda_{v,\RR}$,
$\Top(I)_r$ the degree $r$ part of the ideal in the exterior
algebra of $T_{\sigma}^{\perp}$ generated by 
\[
\bigcup_{i\in I} {\bigwedge}^{\top} T_{\check\Delta_i}.
\]
\item
$H^j(X_{\tau}, \Omega^{r}_{(Z_I^{\tau})^{\dagger}/\kk^{\dagger}}
\cap T_{\sigma}^{\perp})=0$
for $j>0$.
\end{enumerate}
\end{lemma}

\proof
Let $W$ be a complementary subspace to $\sum_{i\in I} T_{\check\Delta_i}
\subseteq T_{\sigma}^{\perp}$. Then we can split $\Omega^1_v|_{X_{\tau}}
\cap T_{\sigma}^\perp$ as $(\O_{X_{\tau}}\otimes W)\oplus
\bigoplus_{i\in I}(\O_{X_{\tau}}\otimes T_{\check\Delta_i})$, and in addition
$d(\O(-Z^{\tau}_i))\subseteq \O_{X_{\tau}}\otimes T_{\check\Delta_i}$, 
as the polynomial
defining $Z^{\tau}_i$ only involves monomials in $\check\Delta_i$. 
Let $d_i=\dim\check\Delta_i$.
Then we obtain a splitting of the exact sequence of Lemma \ref{OmegaZ}
\[
0\rightarrow \bigoplus_{i\in I}\O_{Z^{\tau}_I}(-Z^{\tau}_i)
\rightarrow
\Omega^1_v|_{Z^{\tau}_I}\cap T_{\sigma}^{\perp}\rightarrow
\Omega^1_{(Z^{\tau}_I)^{\dagger}/\kk^{\dagger}}\cap T_{\sigma}^{\perp}
\rightarrow 0
\]
into exact sequences, for $i\in I$,
\begin{equation}
\label{splitexact}
0\rightarrow \O_{Z^{\tau}_I}(-Z^{\tau}_i)
\rightarrow \O_{Z^{\tau}_I}\otimes T_{\check\Delta_i}
\rightarrow\Omega^1_i\rightarrow 0,
\end{equation}
where this sequence defines locally free rank $d_i-1$ sheaves $\Omega^1_i$.
In addition, we have one remaining direct summand
of the original exact sequence,
\[
0\rightarrow 0\rightarrow \O_{Z_I^{\tau}}\otimes W\rightarrow
\O_{Z_I^{\tau}}\otimes W\rightarrow 0.
\]
If we show that  $H^j(Z_I^{\tau},\bigotimes_{i\in I}\Omega_i^{r_i})=0$ for $j>0$
and any collection of $r_i$'s with $0\le r_i\le d_i-1$, then that
will show (2) of the Lemma. We will in fact show
that 
\[
\hbox{$
H^j\left(Z^{\tau}_I,\left(\bigotimes_{i\in I}\Omega_i^{r_i}\right)
\left(-\sum_{i\in I} a_iZ^{\tau}_i\right)\right)=0$ for 
 $0\le a_i\le d_i-1-r_i$,} 
\]
by induction on $\sum_{i\in I} r_i$, starting with $\sum_{i\in I} r_i=0$.

{\it Base case:} $\sum_{i\in I} r_i=0$. 
We need to show $H^j\left(Z^{\tau}_I,\O_{Z^{\tau}_I}\left(-\sum_{i\in I}
a_iZ^{\tau}_i\right)\right)=0$ for $0\le a_i\le d_i-1$. For this we use toric methods.
Let $\psi_i$ be the piecewise linear function on the fan $\Sigma_{\tau}$
defining the line bundle $\O_{X_{\tau}}(Z^{\tau}_i)$. These are convex
functions and $\check\Delta_i$ is the Newton polytope of $\psi_i$ (\cite{PartI}, 
Remark~1.59). Then $-\sum_{i\in I}a_i\psi_i$ corresponds to
$\O_{X_{\tau}}(-\sum_{i\in I} a_iZ^{\tau}_i)$. 
With $\shQ_{\tau,\RR}=(\Lambda^{\perp}_{\tau})^*$ as in \cite{PartI},
Definition 1.33,
by \cite{Oda}, \S 2.2,
\[
H^j(\O_{X_{\tau}}(-\sum_{i\in I} a_iZ^{\tau}_i))=\bigoplus_{n\in\Lambda_{\tau}^{\perp}}
H^j_{Z(n)}(\shQ_{\tau,\RR},\kk),
\]
where 
\[
Z(n):=\{m\in\shQ_{\tau,\RR}|\langle m,n\rangle\ge
\sum_{i\in I} a_i\psi_i(m)\}.
\]
(\cite{Oda}, p.\ 74. Note that our sign convention for the piecewise linear
function corresponding to a divisor is the opposite of that used
in \cite{Oda}.) Also $H^j_{Z(n)}
(\shQ_{\tau,\RR},\kk)\cong H^{j-1}(\shQ_{\tau,\RR}\setminus
Z(n),\kk)$ for $j\ge 2$, and we have an exact sequence
\[
0\rightarrow H^0_{Z(n)}(\shQ_{\tau,\RR},\kk)\rightarrow \kk\rightarrow
H^0(\shQ_{\tau,\RR}\setminus Z(n),\kk)\rightarrow H^1_{Z(n)}(\shQ_{\tau,\RR},\kk)
\rightarrow 0.
\]
(\cite{Oda}, p.\ 74.) Because $\sum_{i\in I} a_i\psi_i$ is convex for 
$a_i\ge 0$, it is easy to see $\shQ_{\tau,\RR}\setminus Z(n)$ is
contractible unless $Z(n)$ is a vector subspace of $\shQ_{\tau,\RR}$,
and this happens if and only if $\langle m,n\rangle\le\sum_{i\in I}
a_i\psi_i(m)$ for all $m\in\shQ_{\tau,\RR}$ with equality if and only if
$m$ is in the subspace $Z(n)$. 
This in turn is equivalent to $n$ being in the relative
interior of $-\sum_{i\in I} a_i\check\Delta_i$, as 
\[
\check\Delta_i=\{n\in\Lambda_{\tau}^{\perp}|\langle m,n\rangle\ge -\psi_i(m)
\quad \forall m\in\shQ_{\tau,\RR}\}.
\]
Now each $\check\Delta_i$ is assumed to be a standard simplex, and hence
$\Int(a_i\check\Delta_i)
\cap \Lambda_{\tau}^{\perp}$ contains no integral points
for $a_i\le d_i$. By simplicity,
$\sum_{i\in I} T_{\check\Delta_i}$ is an interior direct sum, and by the
additional assumption that $\Conv(\bigcup_{i\in I} \check\Delta_i\times
\{e_i\})$ is a standard simplex, this direct sum decomposition also works
over $\ZZ$. Thus if $n\in \Int(-\sum_{i\in I} a_i\check\Delta_i)
\cap\Lambda_{\tau}^{\perp}$, we can write $n=\sum_{i\in I} n_i$
with $n_i\in \Int(-a_i\check\Delta_i)\cap \Lambda_{\tau}^{\perp}$. Thus there
is no such point $n$ if $0\le a_i\le d_i$. In particular, $H^j(\O_{X_{\tau}}
(-\sum_{i\in I} a_iZ^{\tau}_i))=0$ for $j>0$, $0\le a_i\le d_i$.

Now use the Koszul resolution of $\O_{Z^{\tau}_I}$. Tensoring (\ref{koszul})
with $\O_{X_{\tau}}(-\sum_{i\in I}a_i Z^{\tau}_i)$ 
for $0\le a_i\le d_i-1$ we obtain
the base case using the above vanishing.

{\it Induction case.} Suppose $0\le a_i\le d_i-1-r_i$ for all $i\in I$. Pick
$i_1\in I$ with $r_{i_1}>0$. The exact sequence (\ref{splitexact})
yields an exact sequence (see \cite{Hartshorne}, Exercise~II~5.16)
\begin{equation}
\label{wedgesplit}
0\rightarrow \Omega_{i_1}^{r_{i_1}-1}(-Z^{\tau}_i)
\rightarrow \O_{Z^{\tau}_I}\otimes {\bigwedge}^{r_{i_1}} T_{\check\Delta_{i_1}}
\rightarrow \Omega^{r_{i_1}}_{i_1}\rightarrow 0,
\end{equation}
which we tensor with $\left(\bigotimes_{i\in I\setminus\{i_1\}} \Omega_i^{r_i}
\right) \left(-\sum_{i\in I} a_iZ^{\tau}_i\right)$. Then
\[
H^j\bigg(Z^{\tau}_I,\bigg(\bigotimes_{i\in I\setminus\{i_1\}} \Omega_i^{r_i}
\bigg) \big(-\sum_{i\in I} a_iZ^{\tau}_i\big)\bigg)
\]
 and
\[
H^{j+1}\bigg(Z^{\tau}_I,\bigg(\bigotimes_{i\in I\setminus\{i_1\}} \Omega_i^{r_i}
\bigg)\otimes\Omega_{i_1}^{r_{i_1}-1} \big(-Z^{\tau}_{i_1}-\sum_{i\in I} 
a_iZ^{\tau}_i\big)\bigg)
\]
both vanish for $j>0$ by the induction hypothesis, hence 
\[
H^j\left(Z^{\tau}_I,\left(\bigotimes_{i\in I} \Omega_i^{r_i}
\right) \left(-\sum_{i\in I} a_iZ^{\tau}_i\right)\right)
\]
vanishes, completing
the induction step.

This proves (2).

To prove (1), we use 
\[
\Omega^r_{(Z^{\tau}_I)^{\dagger}/\kk^{\dagger}}\cap T_{\sigma}^{\perp}
=\bigoplus_{t,r_i\atop
t+\sum_{i\in I} r_i=r}
\bigg({\bigwedge}^t W\otimes \bigotimes_{i\in I}\Omega^{r_i}_i\bigg).
\]
Since $\rank \Omega_i=d_i-1$, the direct summand is zero unless $r_i
\le d_i-1$ for each $i$. Thus, if we show
\[
\Gamma(Z^{\tau}_I,{\bigwedge}^t W\otimes\bigotimes_{i\in I}\Omega_i^{r_i})
={\bigwedge}^t W\otimes \bigotimes_{i\in I}{\bigwedge}^{r_i} T_{\check\Delta_i}
\]
whenever $r_i\le d_i-1$ for all $i\in I$, the result will follow.
We will show inductively that
\[
H^0(Z^{\tau}_I,(\bigotimes_{i\in I}\Omega_i^{r_i})(-\sum_{i\in I} a_iZ^{\tau}_i))=0
\]
for $0\le a_i\le d_i-1-r_i$ if at least one $a_i>0$, and 
\[
H^0(Z^{\tau}_I,\bigotimes_{i\in I}\Omega_i^{r_i})=\bigotimes_{i\in I}{\bigwedge}^{r_i} 
T_{\check\Delta_i}.
\] 
Again, induction is on $\sum_{i\in I} r_i$. The base case follows immediately
from 
\[
H^0(X_{\tau},\O_{X_{\tau}}(-\sum a_iZ^{\tau}_i))=0
\]
for $a_i\ge 0$ and at
least one $a_i>0$, and the Koszul resolution (\ref{koszul}). For the induction 
step, we tensor (\ref{wedgesplit}) with $(\bigotimes_{i\in I\setminus\{i_1\}}
\Omega_i^{r_i})(-\sum_{i\in I} a_iZ^{\tau}_i)$ as before; if at least one $a_i>0$,
then the 
vanishing of $H^0$ follows from the induction hypothesis, while if $a_i=0$
for all $i$, we get
\begin{eqnarray*}
0&\rightarrow& H^0(Z^{\tau}_I,\Omega_{i_1}^{r_{i_1}-1}(-Z^{\tau}_i)\otimes \bigotimes_{i\in I
\setminus\{i_1\}}\Omega_i^{r_i})\rightarrow H^0(Z^{\tau}_I,{\bigwedge}^{r_{i_1}}
T_{\check\Delta_{i_1}}\otimes \bigotimes_{i\in I\setminus\{i_1\}} \Omega_i^{r_i})\\
&\rightarrow& H^0(Z^{\tau}_I,\bigotimes_{i\in I}\Omega^{r_i}_i)\rightarrow
H^1(Z^{\tau}_I,\Omega_{i_1}^{r_{i_1}-1}(-Z^{\tau}_i)\otimes \bigotimes_{i\in I
\setminus\{i_1\}}\Omega_i^{r_i})\rightarrow 0.
\end{eqnarray*}
From the induction hypothesis and part (2) already proved, 
we get $H^0(Z^{\tau}_I,\bigotimes_{i\in I}
\Omega_i^{r_i})\cong \bigotimes_{i\in I}\bigwedge^{r_i} T_{\check\Delta_i}\otimes \kk$,
as desired.
\qed

\begin{lemma}
\label{almostthere}
With the same hypotheses as in Lemma \ref{vanishing}, in the
$/\kk^{\dagger}$ case, we have for any morphism
$e:\tau_1\rightarrow\tau_2$, $W_e\subseteq B$ the open subset
defined in \cite{PartI}, Lemma 2.9,
\[
\Gamma(W_e,i_*{\bigwedge}^{r}\check\Lambda\otimes \kk)
\cong H^0(X_{\tau_2},(F_{\tau_1,\tau_2}^*\Omega^{r}_{\tau_1})/\Tors)
\]
and
\[
H^j(X_{\tau_2},(F_{\tau_1,\tau_2})^*\Omega^{r}_{\tau_1})/\Tors)=0
\]
for $j>0$.
Without the additional hypotheses on $\bigcup_{i=1}^q \check\Delta_i\times
\{e_i\}$, there is still an injective map 
\[
\Gamma(W_e,i_*{\bigwedge}^{r}\check\Lambda\otimes \kk)
\rightarrow H^0(X_{\tau_2},(F_{\tau_1,\tau_2}^*\Omega^{r}_{\tau_1})/\Tors).
\]
\end{lemma}

\proof
Pick a vertex $g:v\rightarrow\tau_1$.
Then
\begin{eqnarray*}
H^j(X_{\tau_2},(F_{\tau_1,\tau_2}^*
\Omega^r_{\tau_1})/\Tors)&\cong&\HH^j(X_{\tau_2},
F_{\tau_1,\tau_2}^*\shF_v^{r,\bullet})\\
&=&H^j(\Gamma(X_{\tau_2},F_{\tau_1,\tau_2}^*\shF_v^{r,\bullet}))
\end{eqnarray*}
by Lemma \ref{keylemma} and Lemma \ref{vanishing}, (2). In addition,
by Lemma \ref{vanishing}, (1), the complex $\Gamma(X_{\tau_2},
F_{\tau_1,\tau_2}^*\shF^{r,\bullet}_v)$ coincides with the complex of
$\kk$-vector spaces
$F^{\bullet}$, where $F^{\bullet}$ is
defined by, if $\Omega_i,R_i,\Delta_i,\check\Delta_i$ is the simplicity 
data for $\tau_1$,
\[
F^{s}=\bigoplus_{\sigma\subseteq\Delta_{\tau_1}:\bar v\in\sigma
\atop \dim\sigma=s} \left({\bigwedge}^{r-s}T_{\sigma}^{\perp}\right)/
\Top(e,I(\sigma))_{r-s},
\]
where $\sigma=\sigma_1+\cdots+\sigma_q$, $I(\sigma)=\{i|\sigma_i
\not=\{v_i\}\}$ as before, and $\Top(e,I(\sigma))_{r-s}$ is the 
degree $r-s$ part of the ideal of the exterior algebra of
$T_{\sigma}^{\perp}$ generated by
\[
\bigcup_{i\in I(\sigma)}{\bigwedge}^{\top} (T_{\check \Delta_i}\cap 
\Lambda^{\perp}_{\tau_2}).
\]
The term $\Lambda_{\tau_2}^{\perp}$
appears because $\check\Delta_i$ is the Newton polytope for $Z_i^{\tau_1}
\subseteq X_{\tau_1}$, while in Lemma~\ref{vanishing},(1), one uses the tangent
space to the Newton polytope for the corresponding divisor in $X_{\tau_2}$,
i.e. $F_{\tau_1,\tau_2}^{-1}(Z_i^{\tau_1})$. This Newton polytope is a face
of $\check\Delta_i$, and its tangent space is precisely 
$T_{\check\Delta_i}\cap\Lambda^{\perp}_{\tau_2}$.
We interpret the ideal to be the unit ideal if $T_{\check\Delta_i}\cap
\Lambda^{\perp}_{\tau_2}=\{0\}$ for some $i\in I(\sigma)$, i.e. if
$F_{\tau_1,\tau_2}^{-1}(Z^{\tau_1}_i)=\emptyset$.

Furthermore, the differential $\delta_s:F^{s}\rightarrow
F^{s+1}$ is defined by 
\[
(\delta_s\alpha)_{\sigma'}=\sum_{\sigma\subseteq\sigma':\bar v\in\sigma\atop
\dim\sigma=s}\iota(\partial_{w_j-v_j})\alpha_{\sigma}
\]
as before. We can then show $H^j(F^{\bullet})=0$ for $j>0$
by repeating the argument of Lemma~\ref{keylemma}, defining analogous
complexes $F^{\bullet}_{\omega,\omega'}$ and proceeding by induction.
The key stage looks at the diagram
{\scriptsize
\[
\xymatrix@C=30pt
{0\ar[r]&H^{\dim\omega}(F^{\bullet}_{\omega,\omega'\setminus\{w_j\}})
\ar[r]\ar[d]^{\tilde\delta}&{{\bigwedge}^{r-\dim\omega}T_{\omega}^{\perp}\over
\Top(e,I(\omega))_{r-\dim\omega}}\ar[r]\ar[d]&
\displaystyle{\bigoplus_{\omega\subseteq\sigma\subseteq\omega'\setminus\{w_j\}
\atop\dim\sigma=\dim\omega+1}}
{{\bigwedge}^{r-\dim\omega-1}T_{\sigma}^{\perp}\over 
\Top(e,I(\sigma))_{r-\dim\omega-1}}
\ar[d]\\
0\ar[r]&H^{\dim\omega+1}(F^{\bullet}_{\omega\cup\{w_j\},\omega'})
\ar[r]&{{\bigwedge}^{r-\dim\omega-1}T_{\omega
\cup\{w_j\}}^{\perp}\over \Top(e,I(\omega\cup\{w_j\}))_{r-\dim\omega-1}}
\ar[r]&
\displaystyle{\bigoplus_{\omega\subseteq\sigma\subseteq\omega'\setminus\{w_j\}
\atop\dim\sigma=\dim\omega+1}}
{\bigwedge^{r-\dim\omega-2}T_{\sigma\cup
\{w_j\}}^{\perp}\over \Top(e,I(\sigma\cup\{w_j\}))_{r-\dim\omega-2}}
}
\]
}
and we need to show surjectivity of $\tilde \delta$. We explicitly write
down a basis for 
\[
H^{\dim\omega+1}(F^{\bullet}_{\omega\cup\{w_j\},
\omega'}).
\]
We can choose a basis $e_1,\ldots,e_n$ of
$\Lambda_{v,\RR}$, $e_1^*,\ldots,e_n^*$ the dual basis for 
$\check\Lambda_{v,\RR}$, with the properties listed in the proof of
Lemma \ref{keylemma}, but in addition make a choice so that 
there are index sets $J_i$, for $i=1,\ldots,q$ with $\{e_k^*|k\in J_i\}$
a basis for $T_{\check
\Delta_i}\cap \Lambda_{\tau_2}^{\perp}$. Then a basis for
$\bigwedge^{r-\dim\omega-1}
T_{\omega\cup\{w_j\}}^{\perp}/\Top(e,I(\omega\cup
\{w_j\}))_{r-\dim\omega-1}$ is
\[
\{e_L^*|\hbox{$L$ satisfies $L\cap K=\emptyset$ and $J_i\not\subseteq
L$ for any $i\in I(\omega\cup\{w_j\})$}\}.
\]
A basis for $H^{\dim\omega+1}(F^{\bullet}_{\omega\cup\{w_j\},
\omega'})$ then is 
\[
\{e_L^*|\hbox{$L\cap K=\emptyset$, $J_i\not\subseteq L$ for any $i
\in I(\omega\cup\{w_j\})$, and for each $i$, $L\cap I_i\not=\emptyset$ implies
$J_i\subseteq L$}\}.
\]
One then checks as in the proof of Lemma \ref{keylemma}
that each basis vector lifts.

Finally, we calculate $H^0(F^{\bullet})$, and compare this with
$\Gamma(W_e,i_*\bigwedge^{r}\check\Lambda\otimes_{\ZZ} \kk)$. We identify
this with monodromy invariant elements of $i_*\bigwedge^{r}
\check\Lambda_v\otimes_{\ZZ} \kk$ for loops based at $v$ whose interior is in
$W_e$. The monodromy action is then generated by transformations of the
form $T_{feg,feg'}^{\rho}:\Lambda_v\rightarrow\Lambda_v$,
(see \cite{PartI}, \S 1.5)
where we have $f:\tau_2\rightarrow\rho$ with $\rho$ codimension one,
and $g':v'\rightarrow\tau_1$ a vertex. Then as in \cite{PartI}, \S 1.5,
\[
T^{\rho}_{feg,feg'}(m)=m+\langle \check d_{\rho},m\rangle
m^{\rho}_{feg,feg'},
\]
and hence the action on $\bigwedge^{r}\check\Lambda_v\otimes_{\ZZ} \kk$
is
\[
T^{\rho}_{feg,feg'}(n)=n+\check d_{\rho}\wedge\iota(m^{\rho}_{feg,feg'})n.
\]
Thus $n\in\bigwedge^{r}\check\Lambda_v\otimes_{\ZZ} \kk$ is invariant under
all such monodromy operations if and only if $\check d_{\rho}
\wedge\iota(m^{\rho}_{feg,feg'})n=0$ for all choices of $f$ and $g'$. 
Note that as $f\circ e$ runs through elements of $R_i$ which factor
through $e$, $\check d_{\rho}$ runs through a generating set for
$T_{\check\Delta_i}\cap \Lambda_{\tau_2}^{\perp}$, and for
any given $f$ with $f\circ e\in R_i$, as $g'$ varies over
all vertices of $\tau_1$, $m^{\rho}_{feg,feg'}$ runs over
$\{v_i'-v_i|\hbox{$v_i':=\Vert_i(g')$ a vertex of $\Delta_i$}\}$.
From this description, it is then clear that $\check d_{\rho}\wedge
\iota(m^{\rho}_{feg,feg'})n=0$ for all $f,g'$ if and only if
$n\in H^0(F^{\bullet})$. 

Note that in any event, without the extra hypotheses, if $n\in
\Gamma(W_e,i_*\bigwedge^r\check\Lambda\otimes_{\ZZ} \kk)$ is viewed as an element
$n\in\bigwedge^r\check\Lambda_v\otimes_{\ZZ} \kk$, hence defining $\dlog n$
in $F_{v,\tau_2}^*\Omega^r_v$, it is easy to check from the
above formulae that monodromy invariance of $n$ implies $\dlog n
\in \ker\delta_0$, i.e., $\dlog n\in \Gamma(X_{\tau_2},(F_{\tau_1,\tau_2}^*
\Omega^r_{\tau_1})/\Tors)$, giving the map in general, which
is clearly injective, as $\dlog n=0$ if and only if $\dlog n=0$ in
$F_{v,\tau_2}^*\Omega^r_v$.
\qed

\medskip

We can now prove the main theorem of this section: the identification
of the logarithmic Dolbeault groups with the expected cohomology
groups on $B$.

\begin{theorem}
\label{bigtheorem}
Let $B$ be an integral affine manifold with singularities,
with polyhedral decomposition $\P$, and suppose $(B,\P)$
is positive and simple. Assume furthermore that for all $\tau\in\P$,
$\Conv(\bigcup_{i=1}^q\check\Delta_i\times\{e_i\})$ is a standard simplex.
Let $s$ be lifted gluing data, with $X_0=X_0(B,\P,s)$. Then
there is a canonical isomorphism
\[
H^p(X_0,j_*\Omega^r_{X_0^{\dagger}/\kk^{\dagger}})
\cong H^p(B,i_*{\bigwedge}^r\check\Lambda\otimes \kk).
\]
\end{theorem}

\proof
By Corollary \ref{hyper} and Lemma \ref{almostthere}, 
\[
H^p(X_0,j_*\Omega^r_{X_0^{\dagger}/\kk^{\dagger}})
=H^p(\Gamma(X_0,\C^{\bullet}(\Omega^r))),
\]
and 
\[
\Gamma(X_0,\C^p(\Omega^r))=\bigoplus_{\sigma_0\subsetneq\cdots
\subsetneq\sigma_p}
\Gamma(W_{\sigma_0\rightarrow\sigma_p}, i_*{\bigwedge}^r
\check\Lambda\otimes \kk).
\]
However, $\Gamma(W_{\sigma_0\rightarrow\sigma_p},i_*\bigwedge^r
\check\Lambda\otimes \kk)=\Gamma(W_{\sigma_0\rightarrow\cdots\rightarrow\sigma_p},
i_*\bigwedge^r \check\Lambda\otimes \kk)$ 
(where $W_{\sigma_0\rightarrow\cdots\rightarrow\sigma_p}$
is the connected component of $W_{\sigma_0}\cap\dots\cap W_{\sigma_p}$
indexed by $\sigma_0\rightarrow\cdots\rightarrow\sigma_p$; 
if $\P$ has no self-intersecting cells, then $W_{\sigma_0}\cap
\cdots\cap W_{\sigma_p}$ only has one connected component anyway)
because the relevant
monodromy operators, as considered in the proof of Lemma \ref{almostthere},
only depend on $\sigma_0\rightarrow\sigma_p$. Under this identification,
the differential then agrees with the \v Cech differential for
$i_*\bigwedge^r\check\Lambda\otimes \kk$ with respect to the standard
open covering $\{W_{\sigma}\}$. This proves the theorem. 
\qed

\medskip

We can still prove more limited results without the additional
hypotheses.

\begin{theorem}
\label{notquite}
Let $B$ be an integral affine manifold with singularities
with polyhedral decomposition $\P$ and $\dim B=n$, with $(B,\P)$ positive and
simple.
Let $s$ be lifted gluing data, with $X_0=X_0(B,\P,s)$. Then
there is a canonical isomorphism
\[
\hbox{
$H^p(X_0,j_*\Omega^r_{X_0^{\dagger}/\kk^{\dagger}})\cong
H^p(B,i_*\bigwedge^r\check\Lambda\otimes \kk)$
for $0\le p \le n$, $r=0,1,n-1$ or $n$.}
\]
In addition there is always a canonical homomorphism
\[
H^p(B, i_*{\bigwedge}^r\check\Lambda\otimes \kk)
\rightarrow 
H^p(X_0,\Omega^r_{X_0^{\dagger}/\kk^{\dagger}}).
\]
\end{theorem}

\proof 
The case when $r=0$ is already \cite{PartI}, Proposition 2.37.

We need to show the proof of Theorem \ref{bigtheorem} goes
through for forms of degree one, $n-1$, and $n$. For $\Omega^1$, we just
prove Lemma \ref{almostthere} in degree one. 
The resolution $\shF_v^{1,\bullet}$ gives a resolution
\[
0\rightarrow \Omega^1_{\tau}\rightarrow F_{v,\tau}^*\Omega^1_v
\rightarrow \bigoplus_{\bar v\in\sigma\subseteq\Delta_{\tau}\atop
\dim\sigma=1} \O_{Z_{I(\sigma)}}\rightarrow 0.
\]
Note $\# I(\sigma)=1$ if $\dim\sigma=1$. Also $H^p(X_{\tau},\O_{Z_i})=0$
for $p>0$, as follows from $H^p(X_{\tau},\O_{X_{\tau}}(-Z_i))=0$ for all
$p$, which comes from the same toric argument given in the base case
of the proof of Lemma \ref{vanishing}, (2), as $-\check\Delta_i$
contains no interior points. From this it immediately follows that
Lemma \ref{almostthere} holds in general in degree one. Thus
Theorem \ref{bigtheorem} holds for $r=1$, as desired.

Now consider the resolution for $\Omega_{\tau}^{n-1}$. The $p$-th term
of this resolution is a direct sum of terms of the form 
$\Omega^{n-1-p}_{Z^{\dagger}_{I(\sigma)}}\cap T_{\sigma}^{\perp}$, for
$v\in \sigma\subseteq\Delta_{\tau}$, with $\dim\sigma=p$. Now 
$\Omega^1_{Z^{\dagger}_{I(\sigma)}}\cap T_{\sigma}^{\perp}$ is a locally
free sheaf of rank $(n-p)-\# I(\sigma)$, so 
$\Omega^{n-1-p}_{Z^{\dagger}_{I(\sigma)}}\cap T_{\sigma}^{\perp}
=\bigwedge^{n-1-p}(\Omega^1_{Z_{I(\sigma)}^{\dagger}}\cap T_{\sigma}^{\perp})
=0$ unless $\#I(\sigma)=1$. In this case, if $I(\sigma)=\{i\}$, split
$T_{\sigma}^{\perp}=T_{\check\Delta_i}\oplus W$. Then as in the proof of Lemma
\ref{vanishing}, with $d_i=\dim T_{\check\Delta_i}$, $\Omega^1_{Z_i^{\dagger}}\cap
T_{\sigma}^{\perp}\cong \Omega^1_i\oplus(\O_{X_{\tau}}\otimes W)$
and $\Omega^{n-1-p}_{Z_i^{\dagger}}\cong \Omega^{d_i-1}_i\otimes
\bigwedge^{\top}W$. On the other hand, from the definition of $\Omega_i$,
(\ref{splitexact}), we see
$\Omega_i^{d_i-1}\cong\O_{Z_i}(Z_i)$. As $H^j(X_{\tau},\O_{X_\tau})=
H^j(X_{\tau},\O_{X_{\tau}}(Z_i))=0$ for $j>0$, $H^j(X_{\tau},
\Omega^{n-1-p}_{Z_i^{\dagger}}\cap T_{\sigma}^{\perp})=0$ for $j>0$.
In addition, $H^0(X_{\tau},\O_{Z_i}(Z_i))=\coker(H^0(X_{\tau},\O_{X_{\tau}})
\rightarrow H^0(X_{\tau},\O_{X_{\tau}}(Z_i))$, and $\dim H^0(X_{\tau},
\O_{X_{\tau}}(Z_i))=d_i+1$: there is a basis of sections corresponding
to the vertices of $\check\Delta_i$. Thus $\dim H^0(X_{\tau},\O_{Z_i}(Z_i))=d_i$.
On the other hand, the map $\bigwedge^{d_i-1}T_{\check\Delta_i}\otimes \kk
\rightarrow H^0(X_{\tau},\Omega^{d_i-1}_i)$ is injective. To see this, 
take a basis $e_1^*,\ldots,e_{d_i}^*$ of $T_{\check\Delta_i}$, chosen so that
the vertices of $\check\Delta_i$ are $v_i,v_i+e_1^*,\ldots,v_i+e_{d_i}^*$.
Then a basis for $\bigwedge^{d_i-1}T_{\check\Delta_i}$ is $\{
\hat e_j^*:=(-1)^{j-1}e_1^*\wedge\cdots\wedge\widehat{e_j^*}\wedge\cdots
\wedge e_{d_i}^*\}$. If $f_i=0$ locally defines $Z_i$, 
$f_i=1+\sum_{j=1}^{d_i} a_jz^{e_j^*}$, then
for $n\in \bigwedge^{d_i-1}T_{\check\Delta_i}$,
$\dlog n$ maps to zero in $H^0(X_{\tau},\Omega^{d_i-1}_i)$ if and only if
$df_i\wedge \dlog n=0$ along $Z_i$. But writing $n=\sum b_j\hat e_j^*$,
$df_i\wedge \dlog n=\sum_{j=1}^{d_i} a_jb_j z^{e_j^*} \dlog (e_1^*\wedge
\cdots\wedge e_{d_i}^*)$, which is never zero everywhere on $Z_i$ unless
$n=0$ since all the $a_j$'s are non-zero. 
Thus by comparing dimensions, we see $H^0(X_{\tau},\Omega_i^{d_i-1})
\cong\bigwedge^{d_i-1} T_{\check\Delta_i}\otimes \kk$, and thus
$H^0(X_{\tau},\Omega^{n-1-p}_{Z_i^{\dagger}}\cap T_{\sigma}^{\perp})
={\bigwedge^{n-1-p}T_{\sigma}^{\perp}\over \Top(\{i\})}\otimes \kk$. So
Lemma \ref{vanishing} holds for the appropriate degrees and 
Lemma \ref{almostthere} follows for degree $n-1$ forms. This proves
the result for $r=n-1$.

The case $q=n$ is easier: in this case $\Omega^n_{\tau}\cong \O_{\tau}$
for each $\tau$, and the proof goes through immediately.

In any event, by Lemma \ref{almostthere}, we always have maps, for
$e:\tau_1\rightarrow\tau_2$,
\[
\Gamma(W_e,i_*\bigwedge^r\check\Lambda\otimes \kk)
\rightarrow H^0(X_{\tau_2},(F_{\tau_1,\tau_2}^*\Omega^r_{\tau_1})/\Tors)
\]
giving a morphism of complexes from the \v Cech complex of 
$i_*\bigwedge^r\check\Lambda\otimes \kk$ to the complex
$\Gamma(X_0,\C^{\bullet}(\Omega^r))$, which in turn maps to
the complex $\RR\Gamma(X_0,\C^{\bullet}(\Omega^r))$.
After taking the cohomology of this composition, we get a morphism
$H^p(B,i_*\bigwedge^r\check\Lambda\otimes_{\ZZ} \kk)\rightarrow
H^p(X_0,j_*\Omega^r_{X_0^{\dagger}/\kk^{\dagger}})$,
as desired.
\qed

\medskip

Relating this to the cohomology of the tangent bundle, we get

\begin{theorem}
\label{SLnZproperties}
Let $(B,\P)$ be positive and simple,
and suppose the holonomy of $B$ is contained in
$\ZZ^n\rtimes \SL_n(\ZZ)$ (rather than $\Aff(\ZZ^n)=\ZZ^n
\rtimes \GL_n(\ZZ)$). Then with $X_0=X_0(B,\P,s)$,
\begin{enumerate}
\item $j_*\Omega^n_{X_0^{\dagger}/\kk^{\dagger}}
\cong\O_X$, so that $j_*\Omega_{X_0^{\dagger}/\kk^{\dagger}}^{n-r}\cong
j_*\bigwedge^r\Theta_{X_0^{\dagger}/\kk^{\dagger}}$.
\item $i_*\bigwedge^n\check\Lambda\cong\ZZ$ so that $i_*\bigwedge^{n-r}\check
\Lambda\cong i_*\bigwedge^r\Lambda$.
\item There is a canonical isomorphism
$H^p(X,\Theta_{X_0^{\dagger}/\kk^{\dagger}})\cong 
H^p(B,i_*\Lambda\otimes \kk).$
\end{enumerate}
\end{theorem}

\proof The statement about $i_*\bigwedge^n\check\Lambda$ is obvious,
since the monodromy of $\Lambda$ is the linear part of the
holonomy. Thus by Theorem \ref{notquite}, $H^0(X_0,
\Omega^n_{X_0^{\dagger}/\kk^{\dagger}})=
\Gamma(B,i_*\bigwedge^n\check\Lambda\otimes \kk)=\kk$, so
$\Omega^n_{X_0^{\dagger}/\kk^{\dagger}}$ has a section which is nowhere zero, as $\Omega^n_v=\O_{X_v}$
for each $v$. Thus on $X_0\setminus Z$, 
$j_*\Omega^n_{X_0^{\dagger}/\kk^{\dagger}}$ is trivial. 
The remaining statements are then obvious from previous results.
\qed
\bigskip

We also obtain the proverbial interchange of Hodge numbers of mirror 
symmetry:

\begin{corollary}
\label{Hodgeexchange}
Let $(B,\P)$ satisfy the hypotheses of Theorem \ref{bigtheorem}, 
and suppose the holonomy of $B$ is contained in 
$\ZZ^n\rtimes \SL_n(\ZZ)$. Suppose furthermore $\varphi$ is a multi-valued
strictly convex integral piecewise linear function on
$B$ (see \cite{PartI}, Definition~1.47). Let $(\check B,\check \P,
\check\varphi)$ be the discrete Legendre transform of $(B,\P,\varphi)$
(see \cite{PartI}, \S1.4). Suppose $(\check B,\check\P)$ also satisfies
the hypotheses of Theorem \ref{bigtheorem}. Let $s$ be lifted open gluing
data for $(B,\P)$ and $\check s$ be lifted open gluing data for
$(\check B,\check\P)$. Then there is a canonical isomorphism
\[
H^p(X_0(B,\P,s),j_*\Omega^q_{X_0(B,\P,s)^{\dagger}/
\kk^{\dagger}})
\cong
H^p(X_0(\check B,\check \P,\check s),j_*
\Omega^{n-q}_{X_0(\check B,\check \P,\check s)^{\dagger}/
\kk^{\dagger}})
\]
\end{corollary}

\proof Using successively Theorem~\ref{bigtheorem}, 
\cite{PartI}, Proposition~1.50,(1), and 
Theorem~\ref{SLnZproperties},(2),
\begin{eqnarray*}
H^p(X_0(B,\P,s),j_*\Omega^q_{X_0(B,\P,s)^{\dagger}/\kk^{\dagger}})
&\cong&H^p(B,i_*{\bigwedge}^q\check\Lambda^B\otimes\kk)\\
&\cong&H^p(\check B,i_*{\bigwedge}^q\Lambda^{\check B}\otimes\kk)\\
&\cong&H^p(\check B,i_*{\bigwedge}^{n-q}\check\Lambda^{\check B}\otimes\kk)\\
&\cong&
H^p(X_0(\check B,\check \P,\check s),j_*\Omega^{n-q}_{X_0(\check B,\check 
\P,\check s)^{\dagger}/ \kk^{\dagger}}).
\end{eqnarray*}
Here the superscript $B$ or $\check B$ on $\Lambda$ indicates which affine
structure we are using to define $\Lambda$.
\qed

\begin{remark}
\label{stringycomment}
The discrete Legendre transform interchanges the role of the $\Delta_i$'s
and $\check\Delta_i$'s in simplicity data, as follows easily from the
definition in \cite{PartI}, \S 1.5.
Thus by Proposition~\ref{smoothcase} and Corollary~\ref{smoothingisdivisorial},
the deformation of $X_0(\check B,\check\P,\check s)$ constructed
in \cite{Smoothing} over $\Spec\kk\lfor t\rfor$ has non-singular generic
fibre if $(B,\P)$ satisfies the hypotheses of Theorem~\ref{bigtheorem}.
Thus if both $(B,\P)$ and $(\check B,\check\P)$ satisfies
the hypotheses of Theorem~\ref{bigtheorem},  one obtains a
non-singular mirror pair.
If $(B,\P)$ is only simple, then this generic fibre may have orbifold
singularities. As is well-known \cite{Bat}, in dimensions $\ge 4$, 
orbifold singularities are often necessary, in which case stringy
Hodge numbers are needed \cite{BatStringy}. It is 
not yet clear what the relationship
between stringy Hodge numbers and the affine geometry of $B$ is, but
see \cite{Ruddat}.
\end{remark}

\medskip

We now obtain the Hodge decomposition:

\begin{theorem}
\label{hodgedecomp}
With the hypotheses of Theorem \ref{bigtheorem}, there is a canonical
isomorphism
\[
\HH^r(X_0,j_*\Omega^{\bullet}_{X_0^{\dagger}/\kk^{\dagger}})
\cong\bigoplus_{p+q=r} H^p(X_0,j_*\Omega^q_{X_0^{\dagger}/
\kk^{\dagger}}).
\]
\end{theorem}

\proof 
By Corollary \ref{algdeRham} and Lemma \ref{almostthere},
\[
\HH^r(X_0,\Omega^{\bullet})
=H^r(\Gamma(X_0,\Tot(\C^{\bullet}(\Omega^{\bullet})))).
\]
But as
$\Gamma(X_0,(F_{\tau,\sigma}^*\Omega^{\bullet}_{\tau})/\Tors)$
consists entirely of differentials of
the form $\dlog n$, $d$ is in fact zero in 
$\Gamma(X_0,\C^{\bullet}(\Omega^{\bullet}))$, and thus the global
sections of the total complex split as a direct sum
$\bigoplus_q \Gamma(X_0,\C^{\bullet}(\Omega^q)[-q])$, hence the result.
\qed

\section{Basechange and the cohomology of smoothings}

\begin{theorem}\label{hypercohomology base change}
Let $A$ be a local Artinian $\kk[t]$-algebra with residue class field $\kk$
and $\Spec A^\ls$ the scheme $\Spec A$ with log structure induced by
$\NN\to A$, $1\mapsto t$. Assume that
\[
\pi:\X^\ls=(\X,\M_\X)\lra \Spec A^\ls
\]
is a divisorial deformation of a positive and simple toric log Calabi-Yau space
$X_0^{\dagger}\rightarrow
\Spec\kk^{\dagger}$. As before denote by $\shZ\subseteq\X$ 
the singular set of the
log structure of relative codimension two, $j:\X\setminus \shZ\to \X$ the
inclusion of the complement and write $\Omega^\bullet_\X:=
j_*\Omega_{\X^\ls/A^\ls}^\bullet$. Then
$\mathbb{H}^p(\X,\Omega_\X^\bullet)$ is a free $A$-module, and it
commutes with base change.
\end{theorem}

\proof
We follow \cite{steenbrink}, \cite{KN}. By the
cohomology and base change theorem it suffices to prove the
surjectivity of the restriction map
\[
\mathbb{H}^p(\X,\Omega_\X^\bullet)\lra
\mathbb{H}^p(X_0,\Omega_{X_0}^\bullet).
\]
Here $\Omega_{X_0}^\bullet= j_*\Omega^{\bullet}_{X_0^\ls/\kk^\ls}$. As in
\cite{KN}, p.404 it suffices to prove this for
$A=\kk[t]/(t^{k+1})$ with the obvious $\kk[t]$-algebra structure. For
structural clarity we nevertheless keep the notation $A$ for the
base ring.

Consider the complex of $\O_\X$-modules
\[
\mathcal{L}^\bullet= j_*\Omega^\bullet_{\X^\ls/\kk}[u]
=\bigoplus_{s=0}^\infty j_*\Omega^\bullet_{\X^\ls/\kk}\cdot u^s
\]
with differential
\begin{eqnarray*}
\di\Big(\sum_{s=0}^N \alpha_s u^s \Big)&=&
\sum_{s=0}^N \di\alpha_s\cdot u^s + s \dlog \rho\wedge\alpha_s\cdot
u^{s-1}\\
&=& \di\alpha_N\cdot  u^N
+\sum_{s=0}^{N-1} (\di\alpha_s + (s+1) \dlog \rho\wedge\alpha_{s+1})
\cdot u^s,
\end{eqnarray*}
where
$\rho\in\Gamma(\M_\X)$ is the pull-back of the section of $\M_A$
induced by $t$. 
Note that these are differentials relative $\Spec \kk$ rather than
relative $\Spec A^\ls$, so $\dlog\rho\neq0$ unlike in $\Omega_\X$.  
In this complex the dummy variable $u$
formally behaves like $\log t$, and the use of considering this
complex is to trade powers of $\dlog \rho$ with powers of $u$.

Now projection $\sum\alpha_s u^s\mapsto \alpha_0$ defines a map
\[
\shL^\bullet \lra \Omega^\bullet_\X.
\]
To finish the proof it suffices to show that the composition
\[
\varphi^\bullet: \shL^\bullet\lra \Omega_\X^\bullet
\lra \Omega_{X_0}^\bullet
\]
is a quasi-isomorphism, that is, induces isomorphisms of cohomology
sheaves $H^p(\shL^\bullet)\to H^p(\Omega_{X_0}^\bullet)$. In fact,
then the induced composed map of hypercohomology
\[
\mathbb{H}^p(\shL^\bullet)\lra \mathbb{H}^p(\Omega_\X^\bullet)
\lra \mathbb{H}^p(\Omega_{X_0}^\bullet)
\]
is an isomorphism, and hence the second map is surjective as needed.

By this argument and since $\X^\ls\to \Spec A^\ls$ is a divisorial
deformation of ${X_0}^\ls\to \Spec \kk^\ls$, for the rest of the proof we
consider the following local situation. There is a toric variety
$Y=\Spec \kk[P]$ containing $X_0$ as a toric Cartier divisor $V(z^\rho)$
such that the deformation $\X^\ls\to \Spec A^\ls$ is given by
\[
\pi:\Spec \kk[P]/(z^{(k+1)\cdot \rho})\lra \Spec \kk[t]/(t^{k+1}),\quad
\pi^*(t)=z^\rho.
\]
Since $\varphi^r:\shL^r \to \Omega_{X_0}^r$ is surjective for any $r$ we
obtain a short exact sequence
\[
0\lra\shK^\bullet\lra \shL^\bullet\stackrel{\varphi^\bullet}{\lra}
\Omega_{X_0}^\bullet\lra 0
\]
of complexes by defining $\shK^\bullet=\ker \varphi^\bullet$. Now
$\varphi^\bullet$ is a quasi-isomorphism if and only if $\shK^\bullet$
is acyclic, and this is what we are going to show.

For an explicit description of $\shK^r$ let $\sum_{s=0}^N\alpha_s
u^s\in \shL^r$, that is, $\alpha_s\in j_*\Omega_{\X^{\dagger}/\kk}^r$ for all $s$. Then
$\sum_{s=0}^N \alpha_s u^s\in\ker \varphi^r$ iff $\alpha_0|_{X_0}=0$. On
the other hand, the closedness equation $\di (\sum\alpha_s u^s)=0$ is
equivalent to the system of equations
\begin{eqnarray}\label{closedness eqn}
\begin{array}{rcl}
\di\alpha_N&=&0\\
\di\alpha_s+(s+1)\dlog \rho\wedge \alpha_{s+1} &=&0,
\quad s<N.\end{array}
\end{eqnarray}
It is easy to solve these equations after decomposing the coefficients
$\alpha_s$ according to weights, that is, according to the
$P$-grading. First, Proposition~\ref{OmegaXl} gives a decomposition of
$\Gamma(\X,j_*\Omega_{\X^{\dagger}/\kk}^r)$ into homogeneous pieces as follows:
\[
\Gamma(\X,j_*\Omega_{\X^{\dagger}/\kk}^r)=
\bigoplus_{p\in P\setminus((k+1)\rho+P)}
z^p\cdot{\bigwedge}^r\Big(\bigcap_{\{j\,|\, p\in Q_j\}} Q_j^\gp \Big)
\otimes_\ZZ \kk.
\]
The $Q_j$ in this formula are the submonoids of $P$ defined before
Proposition~\ref{Thetatoric}. Thus we obtain $P$-gradings on
$\shL^\bullet$ by imposing the $P$-grading on each direct summand
$j_*\Omega_{\X^{\dagger}/\kk}^r\cdot u^s\subset \shL^\bullet$, and on each
$\Omega_{X_0}^\bullet$ by plugging $k=0$ into the formula above and dividing
by $\ZZ\rho$. Second,
the differentials on $\shL^\bullet$ and on $\Omega_{X_0}^\bullet$ commute
with the respective $P$-gradings, and so does $\varphi^\bullet$. In
fact, this follows from
\begin{eqnarray}\label{P-graded differential}
\di\big(z^p \dlog\omega\cdot u^s\big)&=&
z^p \big(\dlog(p\wedge \omega)\cdot
u^s+s\dlog(\rho\wedge\omega)\cdot u^{s-1}\big).
\end{eqnarray}
Third, all sheaves involved are pull-backs under the morphism
$\mathrm{Et}:\X^\mathrm{et}\to \X^\mathrm{Zar}$ relating the
Zariski site on $\X$ to the \'etale site.

We may thus assume that the $\alpha_s$ in (\ref{closedness eqn}) are
of the form
\[
\alpha_s=z^p \dlog\omega_s,\quad s=0,\ldots,N,
\]
with $\omega_s\in\bigwedge^r V_p\otimes_\ZZ \kk$,
$V_p:=\bigcap_{\{j\,|\, p\in Q_j\}} Q_j^\gp$. Taking into account
(\ref{P-graded differential}) the closedness
condition~(\ref{closedness eqn}) now reads
\begin{eqnarray}\label{P-homogeneous closedness eqn}
\begin{array}{rcl}
p\wedge\omega_N&=&0\\
p\wedge\omega_s+(s+1)\rho\wedge\omega_{s+1}&=&0,
\quad s=0,\ldots,N-1.
\end{array}
\end{eqnarray}
Let us first assume $p\neq 0$. We claim that $\sum\alpha_s u^s=
d\big(\sum z^p\dlog \tau_s\cdot u^s \big)$ for some
$\sum z^p\dlog\tau_s\cdot u^s\in\shL^{r-1}$ if and only if
$d(\sum\alpha_su^s)=0$. To show this, given that
$\sum\alpha_s u^s$ is closed, we need to find
$\tau_0,\ldots,\tau_N$ such that
\begin{eqnarray}\label{P-homogeneous solution}
\begin{array}{rcl}
p\wedge \tau_N&=&\omega_N\\
p\wedge\tau_s+(s+1)\rho\wedge\tau_{s+1}&=&\omega_s,
\quad s=0,\ldots,N-1.
\end{array}
\end{eqnarray}
Recall for a vector space $V$ and $\omega\in\bigwedge^r V$,
$p\in V\setminus\{0\}$, a necessary and sufficient condition for
solvability of the equation $p\wedge\tau=\omega$ is $p\wedge \omega=0$
(integrability condition). By the first line of (\ref{P-homogeneous
closedness eqn}) we can therefore find $\tau_N\in\bigwedge^{r-1}V_p$
with $p\wedge\tau_N=\omega_N$. Assuming inductively that
$\tau_N,\ldots,\tau_{s+1}$ fulfilling (\ref{P-homogeneous solution})
have already been found the integrability condition for $\tau_s$ reads
\[
p\wedge\big( \omega_s-(s+1)\rho\wedge\tau_{s+1}\big)=0.
\]
This follows from (\ref{P-homogeneous closedness eqn}) for $s$ and
from (\ref{P-homogeneous solution}) for $s+1$:
\[
p\wedge\omega_s=-(s+1)\rho\wedge \omega_{s+1}
= -(s+1)\rho\wedge p\wedge\tau_{s+1} =
(s+1)p\wedge\rho\wedge\tau_{s+1}.
\]
Thus there exists a $\tau_s$ satisfying the second line of 
\eqref{P-homogeneous solution}. 
Note that $\tau_s\in\bigwedge^{r-1} V_p$ and hence $\sum
z^p\dlog\tau_s\cdot u^s\in\shL^{r-1}$, proving the claim.

Now suppose furthermore that
\[
\sum z^p\dlog\omega_s\cdot u^s\in \ker\varphi^r.
\]
This is the case if and only if the image of $z^p\dlog\omega_0$
in $\Omega^r_{X_0}$ is zero, which by Proposition~\ref{OmegaXl} and 
Corollary~\ref{OmegaXldagger} holds if and only if $p\in \rho+P$ or
$\rho\wedge\omega_0=0$. If $p\in\rho+P$, then also
$\sum z^p\dlog\tau_s\cdot u^s\in\ker\varphi^{r-1}$. On the other hand, if
$\rho\wedge\omega_0=0$, then by the second line of \eqref{P-homogeneous solution},
\[
\rho\wedge p\wedge\tau_0=\rho\wedge\omega_0=0,
\]
so $\tau_0=\rho\wedge\tau_0'+p\wedge\tau_0''$. We can replace
$\tau_0$ by $\rho\wedge\tau_0'$ without affecting \eqref{P-homogeneous solution},
and hence $\sum z^p\dlog\tau_s\cdot u^s\in\ker\varphi^{r-1}$ also.
This shows acyclicity of the $p$-graded component of $\shK^\bullet$
for $p\neq 0$.

For $p=0$ Equation~(\ref{P-homogeneous closedness eqn}) simply says
$\rho\wedge \omega_s=0$ for $s=1,\ldots, N$. Now $V_0=\bigcap_j
Q_j\otimes_\ZZ \kk$, but as $\rho\in Q_j$ for all $j$, for $s>0$ there
exists $\tau_{s+1}\in \bigwedge^{r-1}V_0$ 
with $\omega_s=(s+1)\rho\wedge\tau_{s+1}$.
Finally, $\varphi^r\big(\sum z^0\dlog\omega_s\cdot u^s \big)=0$
implies $\rho\wedge\omega_0=0$. Hence 
there exists a $\tau_1$ such that $\rho\wedge\tau_1=\omega_0$, and
$\di\big(\sum_{s=1}^N \dlog\tau_s\cdot u^s
\big)=\sum\dlog\omega_s\cdot u^s$.
\qed

\begin{theorem}
\label{basechangetheorem}
With the same hypotheses as Theorem \ref{hypercohomology base change},
with $X_0=X_0(B,\P,s)$, suppose $(B,\P)$ satisfies the hypotheses of
Theorem \ref{bigtheorem}. Then $H^p(\X,j_*\Omega^q_{\X^{\dagger}/A^{\dagger}})$
is a locally free $A$-module, and it commutes with base change.
\end{theorem}

\proof This follows from Theorems \ref{hodgedecomp} 
and \ref{hypercohomology base change}
in a standard way, see \cite{Deligne}, \S 5. \qed

\begin{remark}
Combining this base-change result with Corollary~\ref{Hodgeexchange},
we obtain the interchange of Hodge numbers for non-singular Calabi-Yau varieties
obtained as smoothings of mirror pairs of log Calabi-Yau spaces, of course
subject to the hypotheses of Corollary~\ref{Hodgeexchange}: see the discussion
of Remark~\ref{stringycomment}. For example, as it was shown that
the Batyrev construction of \cite{Bat} is a special case of our mirror
symmetry construction using the discrete Legendre transform on affine
manifolds, the exchange of Hodge numbers in the Batyrev construction
follows from the results in this paper. Of course, our proof is much
more involved, but it covers a much more general
construction and will ultimately give deeper insight into the role that periods
play in calculating Gromov-Witten invariants of the mirror.
\end{remark}

\section{Monodromy and the logarithmic Gauss-Manin connection}

\subsection{Monodromy}
\label{monodromysection}

We begin with a few observations about $B$, if $B$ is an integral
affine manifold with singularities and toric polyhedral decomposition
$\P$. Recall on $B_0$ the extension class of the exact sequence
\[
0\rightarrow\ZZ\rightarrow \shAff(B_0,\ZZ)\rightarrow\check\Lambda
\rightarrow 0
\] 
is the radiance obstruction (\cite{PartI}, Proposition 1.12) 
$c_{B_0}\in H^1(B_0,\Lambda)
=\Ext^1(\check\Lambda,\ZZ)$.
Here $\shAff(B_0,\ZZ)$ is the sheaf of affine linear functions on $B_0$
with integral slope and integral constant term, see \cite{PartI},
Definitions 1.13 and 1.39. As a consequence, the extension class
of the natural exact sequence
\[
0\rightarrow \ZZ\otimes{\bigwedge}^{r-1}\check\Lambda
\rightarrow{\bigwedge}^r\shAff(B_0,\ZZ)\rightarrow{\bigwedge}^r\check
\Lambda\rightarrow 0
\]
in $\Ext^1(\bigwedge^r\check\Lambda,\ZZ\otimes\bigwedge^{r-1}\check\Lambda)
=H^1(B_0,\shHom(\bigwedge^r\check\Lambda,\ZZ\otimes\bigwedge^{r-1}\check\Lambda))$
is by direct check
the image of $c_{B_0}$ under the natural map $\Lambda\rightarrow
\shHom(\bigwedge^r\check\Lambda,\ZZ\otimes\bigwedge^{r-1}\check\Lambda)$
given by 
\[
m\mapsto \bigg(\bigwedge^r\check\Lambda\ni n\mapsto 1\otimes\iota(m)n\bigg).
\]
Pushing forward by $i$, we obtain the exact sequence
\begin{equation}
\label{wedgeaff}
0\rightarrow \ZZ\otimes i_*{\bigwedge}^{r-1}\check\Lambda
\rightarrow i_*{\bigwedge}^r\shAff(B_0,\ZZ)\rightarrow i_*{\bigwedge}^r\check
\Lambda\rightarrow 0.
\end{equation}
We have surjectivity on the right because by local triviality 
of the radiance obstruction on $B$ (\cite{PartI}, Prop. 1.29), one can
find an open cover $\{U_j\}$ of $B$ such that $\shAff(B_0,\ZZ)
\cong \ZZ\oplus\Lambda$ on $U_j\cap B_0$. This also implies the extension
class of (\ref{wedgeaff}) lies in $H^1(B,\shHom(i_*\bigwedge^r\check
\Lambda,\ZZ\otimes i_*\bigwedge^{r-1}\check\Lambda))\subseteq 
\Ext^1(i_*\bigwedge^r\check\Lambda,
\ZZ\otimes i_*\bigwedge^{r-1}\check\Lambda)$. Since $c_{B_0}\in H^1(B_0,\Lambda)$ is in 
$H^1(B,i_*\Lambda)\subseteq H^1(B_0,\Lambda)$ (by local triviality of $c_{B_0}$)
we see the image of $c_{B_0}$ in \[
H^1(B,\shHom(i_*{\bigwedge}^r\check
\Lambda,\ZZ\otimes i_*{\bigwedge}^{r-1}\check\Lambda))
\]
induced by the natural map
$i_*\Lambda\rightarrow \shHom(i_*\bigwedge^r\check
\Lambda,\ZZ\otimes i_*\bigwedge^{r-1}\check\Lambda)$ as above is the extension class of 
(\ref{wedgeaff}). The connecting homomorphisms for (\ref{wedgeaff})
\[
H^p(B,i_*{\bigwedge}^r\check\Lambda)\rightarrow H^{p+1}(B,
\ZZ\otimes i_*{\bigwedge}^{r-1}\check\Lambda)
\]
is then given by $c_{B_0}\cup$ using the natural cup product
\[
H^1(B,i_*\Lambda)\otimes
H^p(B,i_*{\bigwedge}^r\check\Lambda)
\rightarrow H^{p+1}(B,\ZZ\otimes i_*{\bigwedge}^{r-1}\check\Lambda).
\]

\begin{theorem}
Let $(B,\P)$ be a positive and 
simple integral affine manifold with singularities
and polyhedral decomposition $\P$. Let $s$ be lifted gluing data.
Let $X_0^{\dagger}:=X_0(B,\P,s)^{\dagger}$.
Then there is an exact sequence of complexes
\[
0\rightarrow 
f^*\Omega^1_{\kk^{\dagger}/\kk}\otimes
j_*\Omega_{X_0^{\dagger}/\kk^{\dagger}}^{\bullet}[-1]
\rightarrow j_*\Omega_{X_0^{\dagger}/\kk}^{\bullet}\rightarrow
j_*\Omega_{X_0^{\dagger}/\kk^{\dagger}}^{\bullet}
\rightarrow 0.
\]
If $(B,\P)$ satisfies the hypotheses of Theorem \ref{bigtheorem}, then
there are canonical isomorphisms and equalities
\begin{enumerate}
\item
\[
H^p(X_0,j_*\Omega^r_{X_0^{\dagger}/\kk})\cong 
H^p(B,i_*{\bigwedge}^r\shAff(B_0,\ZZ)\otimes \kk).
\]
\item 
\[
\HH^r(X_0,j_*\Omega^{\bullet}_{X_0^{\dagger}/\kk})=
\bigoplus_{p+q=r} H^p(X_0,
j_*\Omega_{X_0^{\dagger}/\kk}^q).
\]
\item 
\[
H^p(X_0,
f^*\Omega^1_{\kk^{\dagger}/\kk}\otimes 
\Omega_{X_0^{\dagger}/\kk^{\dagger}}^{r-1})
\cong 
\Omega^1_{\kk^{\dagger}/\kk}\otimes
H^p(X_0,\Omega^{r-1}_{X_0^{\dagger}/\kk^{\dagger}}).
\]
\item 
Let $\Res_0\nabla$ denote the coboundary map
\begin{eqnarray*}
\HH^r(X_0,
j_*\Omega_{X_0^{\dagger}/\kk^{\dagger}}^{\bullet})&\rightarrow&
\Omega^1_{\kk^{\dagger}/\kk}\otimes
\HH^{r+1}(X_0,j_*\Omega_{X_0^{\dagger}/\kk^{\dagger}}^{\bullet}
[-1])
\\
&=&
\Omega^1_{\kk^{\dagger}/\kk}\otimes\HH^{r}(X_0,
j_*\Omega_{X_0^{\dagger}/\kk^{\dagger}}^{\bullet})
\end{eqnarray*}
followed by evaluation on $\partial_t$, where $\partial_t$ is
the log derivation on $\Spec \kk^{\dagger}$ with $\partial_t\log:\M_{\kk}
=\kk^{\times}\times\NN\rightarrow \kk$ given by $(\alpha,n)\mapsto n$. 
Then under the identifications
of Theorem~\ref{bigtheorem} and Theorem \ref{hodgedecomp},
\[
\Res_0\nabla(\alpha)=c_B \cup \alpha.
\]
\end{enumerate}
\end{theorem}

\begin{remark}
Before proving this theorem, let us explain its importance. Suppose that we 
have a flat deformation $\X\rightarrow D$ over a disk $D$ of $\X_0
=X_0(B,\P,s)$ such that the log structure on $\X_0$ induced by the
divisorial log structure given by $\X_0\subseteq\X$ coincides
with $X_0(B,\P,s)^{\dagger}$. Then there is a monodromy operator
acting on the de Rham cohomology of a fibre $\X_t$, $t\not=0$, and
this monodromy operator is determined by the behaviour of the
Gauss-Manin connection of this family at $0$. In particular, the residue
of the Gauss-Manin connection at $0$, $\Res_0(\nabla)$, coincides with
the map defined above (see \cite{steenbrink}, (2.19)--(2.21) for the
normal crossings case) and the
monodromy operator satisfies
\[
T=e^{-2\pi i \Res_0(\nabla)}.
\]
This description of the monodromy operator should be compared to that of
\cite{SlagI}, Theorem~4.1, which calculated the effect of monodromy assuming
the topological monodromy is described by translation by a section of
an SYZ fibration. In fact, in the current context, there is a canonical
section to consider: one expects the topological SYZ fibration on $\X_t$ to
be a compactification of $X(B_0):=\T_{B_0}/\Lambda$. The torus
fibration $f_0:X(B_0)\rightarrow B_0$ has a canonical section
which can be described locally as the graph of the developing map
$\tilde B_0\rightarrow\RR^n$ of the affine structure (see \cite{PartI},
\S 1.1). As in \cite{SlagI}, this section defines a class in $H^1(B_0,
R^{n-1}f_{0*}\ZZ)\cong H^1(B_0,\Lambda)$. One can show this class 
coincides with $c_{B_0}$. Hence, up to insignificant differences of
sign and factors of $2\pi i$, the description of monodromy given above
and in \cite{SlagI}, Theorem 4.1, agree.
\end{remark}

\proof 
Off of $Z$, we always have an exact sequence
\[
0\rightarrow f^*\Omega^1_{\kk^{\dagger}/\kk}\rightarrow
\Omega^1_{X_0^{\dagger}/\kk}
\rightarrow\Omega^1_{X_0^{\dagger}/\kk^{\dagger}}\rightarrow
0;
\]
hence, as $\Omega^1_{\kk^{\dagger}/\kk}$ is rank 1,
\[
0\rightarrow 
f^*\Omega^1_{\kk^{\dagger}/\kk}\otimes 
j_*\Omega^{r-1}_{X_0^{\dagger}/\kk^{\dagger}}
\rightarrow
j_*\Omega^r_{X_0^{\dagger}/\kk}\rightarrow
j_*\Omega^r_{X_0^{\dagger}/\kk^{\dagger}}
\]
is exact.
Surjectivity on the right can be checked \'etale locally, which
can be done immediately from Proposition \ref{OmegaXl} and Corollary
\ref{OmegaXldagger}. This gives the exact sequence of the theorem.
This sequence induces, for $\tau\in\P$, homomorphisms
$\Omega^1_{\kk^{\dagger}/\kk}\otimes
\Omega^{r-1}_{\tau/\kk^{\dagger}}
\rightarrow \Omega^r_{\tau/\kk}$
and $\Omega^r_{\tau/\kk}\rightarrow\Omega^r_{\tau/\kk^{\dagger}}$,
where we now distinguish between the two different types of $\Omega^r_{\tau}$'s
by adding the $/\kk$ or $/\kk^{\dagger}$. One sees immediately
from Lemma \ref{localomegatauthickened} that these give an exact sequence
\begin{equation}
\label{sestau}
0\rightarrow
\Omega^1_{\kk^{\dagger}/\kk}\otimes
\Omega^{r-1}_{\tau/\kk^{\dagger}}
\rightarrow\Omega^r_{\tau/\kk}\rightarrow\Omega^r_{\tau/\kk^{\dagger}}
\rightarrow 0
\end{equation}
for each $\tau$, and hence an exact sequence of complexes
\[
0\rightarrow\Omega^1_{\kk^{\dagger}/\kk}\otimes 
\C^{\bullet}(j_*\Omega^{r-1}_{X_0^{\dagger}/\kk^{\dagger}})
\rightarrow
\C^{\bullet}(j_*\Omega^r_{X_0^{\dagger}/\kk})\rightarrow
\C^{\bullet}(j_*\Omega^r_{X_0^{\dagger}/\kk^{\dagger}})
\rightarrow 0.
\]
We wish to understand the long exact cohomology sequence of
this sequence. Consider $g:v\rightarrow\tau$.
We can write $\Omega^r_{v/\kk}=\Omega^r_{v/\kk^{\dagger}}\oplus
\dlog\rho\wedge\Omega^{r-1}_{v/\kk^{\dagger}}$. Under this splitting,
it follows from Proposition \ref{delta0kernel} that
$\Omega^r_{\tau/\kk}$ as a subsheaf of $F_{v,\tau}^*\Omega^r_{v/\kk}$
splits as $\Omega^r_{\tau/\kk^{\dagger}}\oplus
\dlog\rho\wedge\Omega^{r-1}_{\tau/\kk^{\dagger}}$. 
This splitting is not canonical, but depends on the
choice of $g$.

Now assume $(B,\P)$ satisfies the hypotheses of Theorem \ref{bigtheorem}.
To show (1), we need, for $e:\tau_1\rightarrow\tau_2$, a canonical
isomorphism
\[
H^0(X_{\tau_2},(F_{\tau_1,\tau_2}^*\Omega_{\tau_1/\kk}^r)/\Tors)
\cong\Gamma(W_e,i_*{\bigwedge}^r\shAff(B_0,\ZZ)\otimes \kk).
\]
Now the sequence (\ref{wedgeaff}) splits on $W_e$ by local vanishing
of the radiance obstruction (\cite{PartI}, Proposition~1.29), and given any
splitting $\varphi:i_*\check\Lambda\rightarrow\shAff(B_0,\ZZ)$ on
$W_e$, a choice of vertex $v\rightarrow\tau_1$ gives a well-defined
splitting $\varphi_v$ of $\Gamma(W_e,\shAff(B_0,\ZZ))$ by
$\varphi_v(n)=\varphi(n)-\varphi(n)(v)$ for $n\in \Gamma(W_e,
i_*\check\Lambda)$. Thus we lift sections of $i_*\check\Lambda$
to affine linear
functions vanishing on $v$; the resulting splitting is independent
of the original choice $\varphi$. Thus we get an
isomorphism, depending only on $v$,
\[
\Gamma(W_e,i_*{\bigwedge}^r\shAff(B_0,\ZZ))
\cong\Gamma(W_e,i_*{\bigwedge}^r\check\Lambda)
\oplus\big(\ZZ\otimes\Gamma(W_e,i_*{\bigwedge}^{r-1}\check\Lambda)\big).
\]
Now let $h:\omega\rightarrow\tau_1$ be an edge. Then for $n\in
\Gamma(W_e,i_*\check\Lambda)$,
\begin{eqnarray*}
\varphi_{v^+_{\omega}}(n)&=&\varphi_{v^-_{\omega}}(n)
-\varphi_{v^-_{\omega}}(n)(v^+_{\omega})\\
&=&\varphi_{v^-_{\omega}}(n)-\langle n,v^+_{\omega}-v^-_{\omega}\rangle\\
&=&\varphi_{v^-_{\omega}}(n)+\langle n,l_{\omega}d_{\omega}\rangle,
\end{eqnarray*}
where the notation is as in Lemma \ref{basicgluing}.
Thus we have two isomorphisms
\[
\Gamma(W_e,i_*{\bigwedge}^{r}\check\Lambda)
\oplus(\ZZ\otimes\Gamma(W_e,i_*{\bigwedge}^{r-1}\check\Lambda))
\mapright{\cong}
\Gamma(W_e,i_*{\bigwedge}^{r}\shAff(B_0,\ZZ))
\]
given by the splittings $\varphi_{v^-_{\omega}}$ and
$\varphi_{v^+_{\omega}}$ respectively. The composition of the inverse of
the second isomorphism with the first isomorphism is
\[
\Gamma_h':(\alpha_1,\alpha_2)\mapsto (\alpha_1,-l_{\omega}\iota(d_{\omega})
\alpha_1+\alpha_2).
\]
We now have by Lemma \ref{almostthere} and \eqref{sestau} a diagram
\[
\xymatrix@C=30pt
{
\Gamma(X_{\tau_2},(F_{\tau_1,\tau_2}^*\Omega^r_{\tau_1/\kk}))/\Tors
\ar[d]^{\Gamma_h}\ar[r]^-{\cong}&
\Gamma(W_e,i_*{\bigwedge}^{r}\check\Lambda\otimes k)
\oplus \dlog\rho\wedge
\Gamma(W_e,i_*{\bigwedge}^{r-1}\check\Lambda\otimes k)
\ar[d]_{\Gamma_h'}\\
\Gamma(X_{\tau_2},(F_{\tau_1,\tau_2}^*\Omega^r_{\tau_1/\kk}))/\Tors
\ar[r]_-{\cong}&
\Gamma(W_e,i_*{\bigwedge}^{r}\check\Lambda\otimes k)
\oplus \dlog\rho\wedge
\Gamma(W_e,i_*{\bigwedge}^{r-1}\check\Lambda\otimes k)\\
}
\]
where we think of the two occurrences of 
$\Gamma(X_{\tau_2},(F_{\tau_1,\tau_2}^*\Omega^r_{\tau_1/\kk})/\Tors)$ 
as being contained in
\[
\hbox{$\Gamma(X_{\tau_2},F_{v^-_{\omega},\tau_2}^*\Omega^r_{v^-_{\omega}/\kk})$
and
$\Gamma(X_{\tau_2},F_{v^+_{\omega},\tau_2}^*\Omega^r_{v^+_{\omega}/\kk})$}
\]
respectively. This diagram is then commutative by Lemma 
\ref{basicgluing}.
Thus we can canonically identify
$\Gamma(X_{\tau_2},(F_{\tau_1,\tau_2}^*\Omega_{\tau_1/\kk}^r)/\Tors)$
with $\Gamma(W_e,i_*\bigwedge^r\shAff(B_0,\ZZ)\otimes \kk)$.
Item (1) then follows because $\Gamma(X_0,\C^{\bullet}(j_*
\Omega_{X_0^{\dagger}/\kk}^r))$
now coincides with the \v Cech complex for $i_*\bigwedge^r
\shAff(B_0,\ZZ)\otimes \kk$ with respect to the standard open cover $\{W_{\tau}\}$.
Item (2) follows as in Theorem \ref{hodgedecomp}, item (3) is obvious, and
item (4) follows from the discussion of $c_B$ and (1)--(3). \qed

\subsection{The connection on moduli}

Let $(B,\P)$ be positive and simple. Recall from \cite{PartI}, Theorem 5.4,
that the set of positive log Calabi-Yau spaces with dual intersection
complex $(B,\P)$, modulo isomorphisms preserving $B$, is canonically
$H^1(B,i_*\Lambda\otimes\kk^{\times})$. By the universal coefficient theorem,
there is an exact sequence
\[
0\rightarrow H^1(B,i_*\Lambda)\otimes \kk^{\times}\rightarrow 
H^1(B,i_*\Lambda\otimes\kk^{\times})
\rightarrow H^2(B,i_*\Lambda)_{tors}\rightarrow 0,
\]
so $H^1(B,i_*\Lambda\otimes\kk^{\times})$ can be viewed as a disjoint union 
of torsors over $H^1(B,i_*\Lambda)\otimes\kk^{\times}$. This latter
algebraic
torus can be viewed as $S:=\Spec \kk[H^1(B,i_*\Lambda)^*_f]$,
where $G_f:=G/G_{\tors}$ denotes the torsion free part of the finitely
generated abelian group $G$. Observe
$\Gm(S)=\kk^{\times}\times H^1(B,i_*\Lambda)^*_f$, and
\[
H^1(B,i_*\Lambda\otimes\Gm(S))=H^1(B,i_*\Lambda\otimes\kk^{\times})
\times (H^1(B,i_*\Lambda)\otimes H^1(B,i_*\Lambda)^*_f).
\]
For the remainder of this section fix an element
\[
{\bf s}\in H^1(B,i_*\Lambda\otimes\Gm(S))
\]
of the form
\[
{\bf s}=(s_0,s_{\id}),
\]
where $s_0\in H^1(B,i_*\Lambda\otimes\kk^{\times})$ is arbitrary
and $s_{\id}\in H^1(B,i_*\Lambda)\otimes H^1(B,i_*\Lambda)^*_f$
corresponds to a choice of splitting $H^1(B,i_*\Lambda)_f
\rightarrow H^1(B,i_*\Lambda)$. One can view $s_0$ as selecting
the component of $H^1(B,i_*\Lambda\otimes\kk^{\times})$
containing it, and $(s_0,s_{\id})$ can then be viewed as a ``universal 
element'' on this component,
in the following sense. We can view ${\bf s}$ as an
$H^1(B,i_*\Lambda\otimes\kk^{\times})$-valued 
function on $S$. At a closed point 
$s$ of $S$ corresponding to a group
homomorphism $\chi_s:H^1(B,i_*\Lambda)_f^*\rightarrow\kk^{\times}$,
the value of ${\bf s}$ is $s_0\cdot (1\otimes\chi_s)(s_{\id})
=s_0\cdot s$.

Choose a \v Cech representative $({\bf s}_e)\in\bigoplus_e
\Gamma(W_e,i_*\Lambda\otimes
\Gm(S))$ for ${\bf s}$ 
by choosing \v Cech representatives $(s_{0,e})$ and $(s_{\id,e})$
for $s_0$ and $s_{id}$
respectively. As in \cite{PartI}, 
Definition 5.1, $({\bf s}_e)$ can be viewed as 
open gluing data for $(B,\P)$, and hence we obtain a family of algebraic
spaces $\pi:X_0(B,\P,{\bf s})\rightarrow S$ by \cite{PartI}, Definition 2.28.

\begin{remark} The choice of splitting $s_{\id}$ does not affect the isomorphism
class of each fibre of $\pi:X_0(B,\P,{\bf s})^{\dagger}\rightarrow S^{\dagger}$.
However, it does affect the total space $X_0(B,\P,{\bf s})^{\dagger}$.
The results we give here apply to any choice.
\end{remark}

We need to put a log structure on this space. Unfortunately, in \cite{PartI},
we did not deal with log structures for families; however, since the spaces
are reduced, the methods of \cite{PartI}, \S\S 3.2, 3.3 still apply. 
Let $X_0=X_0(B,\P,{\bf s})$. Then
$X_0$ is covered by open sets of the form $V(\sigma)\times
S$ for $\sigma\in\P_{\max}$, and the sheaf $\shLS_{V(\sigma)\times S}$
on $V(\sigma)\times S$ is determined by \cite{PartI}, Theorem 3.24; in 
particular, $\shLS^+_{\pre,V(\sigma)\times S}$ is just the
pull-back of $\shLS^+_{\pre,V(\sigma)}$ under the projection to 
$V(\sigma)$. Sections of $\shLS^+_{\pre,V(\sigma)\times S}$ transform
under the gluing map exactly as in \cite{PartI}, Theorem 3.27,
keeping in mind quantities involving ${\bf s}$ are now functions on $S$.
This defines the sheaf $\shLS^+_{\pre,X_0}$, and as in
\cite{PartI}, Theorem 5.2, (2), ${\bf s}$ determines a unique
normalized section $f\in\Gamma(X_0,\shLS^+_{\pre,X_0})$ 
which in fact is a section of $\shLS_{X_0}$
off of some singular set of relative codimension $\ge 2$ over $S$
not containing any toric stratum of $X_0$. This defines
a log structure $X_0^{\dagger}$ on $X_0$.
Furthermore, if $S^{\dagger}$ denotes the log structure on $S$ given
by the chart 
\[
\NN\rightarrow\O_S, \quad n\mapsto \begin{cases} 1&n=0\\
0&n>0,\end{cases}
\]
then $\pi$ lifts to a log morphism
$\pi:X_0^{\dagger}\rightarrow S^{\dagger}$, smooth away from
$Z$. Note in addition there is an obvious canonical map $S^{\dagger}\rightarrow
\Spec\kk^{\dagger}$.

\medskip

We have the following generalisation of Theorem \ref{localmodel}, whose
proof is identical.

\begin{theorem}
\label{localmodel2}
In the above situation,
let $\bar x\rightarrow Z\subseteq X_0$ 
be a closed geometric point. 
Then there exists data $\tau$, $\check\psi_1,\ldots,
\check\psi_q$ as in Construction \ref{keyP}
 defining a monoid $P$, hence log spaces $Y^{\dagger}$, $X^{\dagger}\rightarrow
\Spec \kk^{\dagger}$ as in \S 1, 
such that there is a diagram over $S^{\dagger}$ 
\begin{equation*}
\xymatrix{
&V^{\dagger}\ar[ld]\ar[rd]^{\phi}&\\
X_0^{\dagger}&&X^{\dagger}\times_{\Spec\kk^{\dagger}} S^{\dagger}
}
\end{equation*}
with both maps strict \'etale and $V^{\dagger}$
an \'etale neighbourhood of $\bar x$.
\end{theorem}

Our goal now is to compute the Gauss-Manin connection for the family $\pi$, 
\[
\nabla_{GM}:\RR^r\pi_*(j_*\Omega^{\bullet}_{X_0^{\dagger}/
S^{\dagger}})\rightarrow\Omega^1_{S^{\dagger}/\kk^{\dagger}}
\otimes\RR^r\pi_*(j_*\Omega^{\bullet}_{X_0^{\dagger}/
S^{\dagger}}),
\]
where $j:X_0\setminus Z
\rightarrow X_0$ is the inclusion
as usual.
Just as in \cite{KatzOda}, this connection can be defined as follows.
We have on $X_0\setminus Z$ an exact sequence of locally
free sheaves
\begin{equation}
\label{GMexactseq}
0\rightarrow \pi^*\Omega^1_{S^{\dagger}/\kk^{\dagger}}
\rightarrow \Omega^1_{X_0^{\dagger}/\kk^{\dagger}}
\rightarrow \Omega^1_{X_0^{\dagger}/S^{\dagger}}
\rightarrow 0,
\end{equation}
hence a filtration of complexes
\[
\Omega^{\bullet}_{X_0^{\dagger}/\kk^{\dagger}}
=F^0(\Omega^{\bullet}_{X_0^{\dagger}/\kk^{\dagger}})
\supseteq
F^1(\Omega^{\bullet}_{X_0^{\dagger}/\kk^{\dagger}})
\supseteq\cdots
\]
where
\begin{eqnarray*}
F^i&:=&F^i(
\Omega^{\bullet}_{X_0^{\dagger}/\kk^{\dagger}})\\
&=&\im(
\pi^*(\Omega^i_{S^{\dagger}/\kk^{\dagger}})
\otimes
\Omega^{\bullet}_{X_0^{\dagger}/\kk^{\dagger}}[-i]
\rightarrow\Omega^{\bullet}_{X_0^{\dagger}/\kk^{\dagger}})
\end{eqnarray*}
with
\[
F^i/F^{i+1}=\pi^*(\Omega^i_{S^{\dagger}/\kk^{\dagger}})
\otimes\Omega^{\bullet}_{X_0^{\dagger}/S^{\dagger}}[-i].
\]
Noting that by Theorem \ref{localmodel2}, \eqref{GMexactseq} is locally split
after applying $j_*$,
we get an exact sequence on $X_0$
\[
0\rightarrow \pi^*\Omega^1_{S^{\dagger}/\kk^{\dagger}}
\rightarrow j_*\Omega^1_{X_0^{\dagger}/\kk^{\dagger}}
\rightarrow j_*\Omega^1_{X_0^{\dagger}/S^{\dagger}}
\rightarrow 0,
\]
as well as
\[
j_*F^i/j_*F^{i+1}=\pi^*(\Omega^i_{S^{\dagger}/\kk^{\dagger}})
\otimes j_*\Omega^{\bullet}_{X_0^{\dagger}/S^{\dagger}}[-i].
\]
Then the Gauss-Manin connection is the boundary map coming from the
exact sequence
\begin{equation}
\label{GMexact}
0\rightarrow j_*F^1/j_*F^2\rightarrow j_*F^0/j_*F^2\rightarrow
j_*F^0/j_*F^1\rightarrow 0,
\end{equation}
i.e.
\begin{eqnarray*}
\nabla_{GM}:\RR^r\pi_*(j_*\Omega^{\bullet}_{X_0^{\dagger}
/S^{\dagger}})&\rightarrow&
\RR^{r+1}\pi_*(\pi^*(\Omega^1_{S^{\dagger}/\kk^{\dagger}})
\otimes j_*\Omega^{\bullet}_{X_0^{\dagger}/S^{\dagger}}[-1])\\
&=&\Omega^1_{S^{\dagger}/\kk^{\dagger}}
\otimes\RR^r\pi_*(j_*\Omega^{\bullet}_{X_0^{\dagger}/
S^{\dagger}}).
\end{eqnarray*}

\begin{theorem}
Let $(B,\P)$ be positive and 
simple, satisfying the hypotheses of Theorem \ref{bigtheorem}.
Let ${\bf s}$ be as above. Let $X_0=X_0(B,\P,{\bf s})$.
Then
\begin{enumerate}
\item 
there exists a canonical isomorphism
\[
R^q\pi_*(j_*\Omega^p_{X_0^{\dagger}/S^{\dagger}})
\cong H^q(B,i_*{\bigwedge}^p\check\Lambda)\otimes\O_S.
\]
\item 
\[
\RR^r\pi_*(j_*\Omega^{\bullet}_{X_0^{\dagger}/S^{\dagger}})
\cong\bigoplus_{p+q=r} R^q\pi_*(j_*\Omega^p_{X_0^{\dagger}/
S^{\dagger}}).
\]
\item (Griffiths Transversality)
\[
\nabla_{GM}(R^q\pi_*(j_*\Omega^p_{X_0^{\dagger}/S^{\dagger}}))
\subseteq \Omega^1_{S^{\dagger}/\kk^{\dagger}}
\otimes R^{q+1}\pi_*(j_*\Omega^{p-1}_{X_0^{\dagger}/S^{\dagger}}).
\]
\item Identifying $\Omega^1_{S^{\dagger}/\kk^{\dagger}}=\Omega^1_{S/\kk}$ with
$H^1(B,i_*\Lambda)^*_f\otimes\O_S$, via $\alpha\in H^1(B,i_*\Lambda)^*_f$
yielding the differential $\dlog\alpha$,
then $-\nabla_{GM}$ restricted to the constant sections (under the
isomorphism of (1)) of
$R^q\pi_*(j_*\Omega^p_{X_0^{\dagger}/S^{\dagger}})$ is the map
\[
-\nabla_{GM}:H^q(B,i_*{\bigwedge}^p\check\Lambda\otimes\kk)
\rightarrow H^1(B,i_*\Lambda\otimes\kk)^*\otimes H^{q+1}(B,
i_*{\bigwedge}^{p-1}\check\Lambda\otimes\kk)
\]
induced by the cup product composed with contraction:
\[
H^1(B,i_*\Lambda\otimes\kk)\otimes H^q(B,i_*{\bigwedge}^p\check\Lambda
\otimes\kk)\rightarrow H^{q+1}(B,i_*{\bigwedge}^{p-1}\check\Lambda\otimes
\kk).
\]
We denote this latter map by $\eta^*\otimes n\mapsto \iota(\eta^*)n$,
for $\eta^*\in H^1(B,i_*\Lambda\otimes\kk)$ and 
$n\in H^q(B,i_*{\bigwedge}^p\check\Lambda
\otimes\kk)$.
\end{enumerate}
\end{theorem}

\proof \emph{Step 1.} \emph{Review of how the lifted
gluing data $({\bf s}_e)$ determines the section $f\in
\Gamma(X_0,\shLS^+_{\pre,X_0})$.}

For $\sigma\in\P_{\max}$, denote by $f_{\sigma}$ the pull-back of $f$
to $V(\sigma)\times S$. We follow the proof of \cite{PartI}, Prop. 4.25.
We can write
\[
\shLS^+_{\pre,V(\sigma)\times S}=\bigoplus_{e\in\coprod\Hom(\omega,\sigma)
\atop\dim\omega=1} \O_{V_e\times S},
\]
which allows us to decompose $f_{\sigma}$ into components $f_{\sigma,e}
\in \Gamma(\O_{V_e\times S})$. 
We can then write, for $e:\omega\rightarrow\sigma$,
\[
f_{\sigma,e}=\sum_{p\in\check\Delta(\omega)\cap\check\Lambda_y}
f_{\sigma,e,p} z^p,
\]
where $y\in\omega$ is a point near $v^+_{\omega}$ and $\check\Delta(\omega)$
is the convex hull in $\Lambda^{\perp}_{\omega,\RR}\subseteq\check
\Lambda_{\RR,y}$ of
\[
\{n_{\omega}^{\sigma\sigma'}|e':\omega\rightarrow\sigma'\in\P_{\max}\}.
\]
(See \cite{PartI}, \S 1.5).
By simplicity this is an elementary simplex. Given $e':\omega\rightarrow
\sigma'$, we have a diagram
\[
\xymatrix@C=30pt
{&&\sigma\\
\omega\ar@/^/[rru]^{e}\ar[r]^h\ar@/_/[drr]_{e'}&\tau
\ar[ru]_{g}\ar[rd]^{g'}&\\
&&\sigma'
}
\]
with $(g,g')$ maximal (\cite{PartI}, Definition 2.21). Then by \cite{PartI}, (4.4),
\[
f_{\sigma,e,n_{\omega}^{\sigma\sigma'}}={D({\bf s}_g,\omega,\sigma)\over
D({\bf s}_{g'},\omega,\sigma')}{\bf s}_g(n_{\omega}^{\sigma\sigma'}).
\]
However, since ${\bf s}_g$, ${\bf s}_{g'}$ are multiplicative maps
on $\check\Lambda_y$, rather
than piecewise multiplicative
(see \cite{PartI}, Definition 5.1), $D({\bf s}_g,\omega,\sigma)
=D({\bf s}_{g'},\omega,\sigma')=1$. On the other hand,
by the \v Cech cocycle condition, ${\bf s}_h\cdot {\bf s}_g=
{\bf s}_e$. Since ${\bf s}_h\in\Gamma(W_h,i_*\Lambda\otimes\Gm(S))$,
${\bf s}_h$ is invariant under the monodromy operator 
$T_{\omega}^{\sigma\sigma'}$
(see \cite{PartI}, \S 1.5, where we write $T_{\omega}^{ee'}$), 
which takes the form (expressed multiplicatively
on $\Lambda_y\otimes\Gm(S)$)
\[
T_{\omega}^{\sigma\sigma'}({\bf s}_h)={\bf s_h}\cdot (d_{\omega}\otimes {\bf s}_h
(n_{\omega}^{\sigma\sigma'})).
\]
Thus 
\begin{equation}
\label{shmonpoly}
{\bf s}_h(n_{\omega}^{\sigma\sigma'})=1,
\end{equation}
and we see
\[
f_{\sigma,e,n_{\omega}^{\sigma\sigma'}}={\bf s}_e(n_{\omega}^{\sigma\sigma'}).
\]
Since as $e'$ varies, $n_{\omega}^{\sigma\sigma'}$ runs over all integral points
of $\check\Delta(\omega)$, we in fact see
\begin{equation}
\label{fsigmae}
f_{\sigma,e}=\sum_{p\in\check\Delta(\omega)\cap\check\Lambda_y}
{\bf s}_e(p)z^p.
\end{equation}

\smallskip

\emph{Step 2}. \emph{The resolutions of $j_*\Omega^r_{X_0^{\dagger}/
\kk^{\dagger}}$ and $j_*\Omega^r_{X_0^{\dagger}/S^{\dagger}}$
and the proof of (1) and (2).}

We have maps $q_{\tau}:X_{\tau}\times S\rightarrow X_0$ from
the construction of $X_0$ in \cite{PartI}, \S 2.2, as usual,
and using this, we can define the sheaves $\Omega^r_{\tau/\kk^{\dagger}}$
and $\Omega^r_{\tau/S^{\dagger}}$ on $X_{\tau}$ using
$\Omega^r_{X_0^{\dagger}/\kk^{\dagger}}$ and
$\Omega^r_{X_0^{\dagger}/S^{\dagger}}$ just as in
\S \ref{globalcalculations}. As usual these yield resolutions 
$\C^{\bullet}(j_*\Omega^r_{X_0^{\dagger}/\kk^{\dagger}})$
and
$\C^{\bullet}(j_*\Omega^r_{X_0^{\dagger}/S^{\dagger}})$
of $j_*\Omega^r_{X_0^{\dagger}/\kk^{\dagger}}$ and
$j_*\Omega^r_{X_0^{\dagger}/S^{\dagger}}$ respectively.
This follows from Theorem \ref{localmodel2} and the local form of these
results (Theorem \ref{localresolution}).

Now if $s\in S$, then $\Omega^r_{X_0(B,\P,{\bf s})^{\dagger}
/S^{\dagger}}|_{\pi^{-1}(s)}
=\Omega^r_{X_0(B,\P,s_0\cdot s)^{\dagger}/\kk^{\dagger}}$ on
$X_0(B,\P,s_0\cdot s)\setminus Z$. One checks locally that the 
same continues to hold after pushing forward by $j$, i.e.
\[
(j_*\Omega^r_{X_0(B,\P,{\bf s})^{\dagger}/S^{\dagger}})|_{\pi^{-1}(s)}
=j_*\Omega^r_{X_0(B,\P,s_0\cdot s)^{\dagger}/\kk^{\dagger}}.
\]
Thus $R^q\pi_*(j_*\Omega^p_{X_0^{\dagger}/S^{\dagger}})$
is locally free. Given a \v Cech cocycle
\[
(W_{\sigma_0\rightarrow\cdots\rightarrow\sigma_q},n_{\sigma_0\rightarrow\cdots
\rightarrow\sigma_q})
\]
representing an element
$n\in H^q(B,i_*\check\Lambda\otimes\kk)$, $n_{\sigma_0\rightarrow\cdots
\rightarrow\sigma_q}$ determines a section of
\[
\Gamma(X_{\sigma_q}\times S,F_{\sigma_0,\sigma_q}^*
\Omega^p_{\sigma_0/S^{\dagger}}),
\]
namely $\dlog n_{\sigma_0\rightarrow
\cdots\rightarrow\sigma_q}$, exactly as in Lemma 
\ref{almostthere}, and hence $n$ determines an element of 
\[
H^q(\Gamma(X_0,\C^{\bullet}(j_*
\Omega^p_{X_0^{\dagger}/S^{\dagger}})))
\]
which 
restricts to the class $n$ on each fibre under the isomorphism of Theorem
\ref{bigtheorem}. This shows item (1), while item (2) follows exactly
as in the proof of Theorem \ref{hodgedecomp}.

\smallskip

\emph{Step 3}. \emph{ Describing
$\Omega^r_{\tau/\kk^{\dagger}}$.}

The methods of \S \ref{globalcalculations}
apply to calculate $\Omega^r_{\tau/\kk^{\dagger}}$,
as nowhere did we assume properness while calculating these sheaves.
However, there is an important subtlety.
Consider first an irreducible component $X_v\times S$ of $X_0$.
Lemma \ref{OXD} still applies, so $\Omega^r_{v/\kk^{\dagger}}$ is a trivial
vector bundle, but we need to be careful about how we identify this bundle.
The reason is that $X_v\times S$ is obtained by gluing together affine schemes
of the form $V_e\times S$, with $e:v\rightarrow\sigma\in\P_{\max}$. However,
they are not glued in the canonical way, but rather with a twist induced
by ${\bf s}$, and this needs to be taken into account. As a result,
this bundle is abstractly 
\[
\Omega^r_{X_v^{\dagger}\times S^{\dagger}/\kk^{\dagger}}
={\bigwedge}^r(\check\Lambda_v\oplus
H^1(B,i_*\Lambda)^*_f)\otimes\O_{X_v\times S};
\]
however there is only a canonical identification after choosing a maximal
cell $\sigma\in\P_{\max}$ containing $v$. We shall now
see how this representation
depends on $e:v\rightarrow\sigma
\in\P_{\max}$. We have $V_e\times S$ is an open affine subset of $X_v\times S$,
and given a diagram
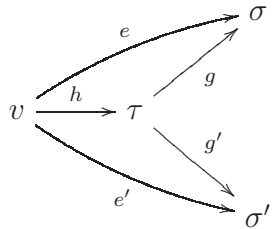
\[
\xymatrix@C=30pt
{&&\sigma\\
v\ar@/^/[rru]^{e}\ar[r]^h\ar@/_/[drr]_{e'}&\tau
\ar[ru]_{g}\ar[rd]^{g'}&\\
&&\sigma'
}
\]
with $(g,g')$ a maximal pair, $V_e\times S$ and $V_{e'}\times S$ 
are glued along the open subset $V_h\times S$ using
$\Phi_{gg'}({\bf s})$. This is given on the level of rings as follows,
as explained in \cite{PartI}, \S 2.2.
The vertex $v$ determines (different) cones $\check v$ in the normal
fans $\check\Sigma_{\sigma}$ and $\check\Sigma_{\sigma'}$
of $\sigma$ and $\sigma'$ in $\Lambda^*_{\sigma,\RR}$ and $\Lambda^*_{\sigma',
\RR}$ respectively. Furthermore, $V_e=\Spec\kk[\check v\cap\Lambda_{\sigma}^*]$
and $V_{e'}=\Spec \kk[\check v\cap\Lambda^*_{\sigma'}]$. Then $V_h
\subset V_e, V_{e'}$ is described via the localizations
$\Spec \kk[(\check v+\RR\check\tau)\cap\Lambda_{\sigma}^*]$
and $\Spec \kk[(\check v+\RR\check\tau)\cap\Lambda^*_{\sigma'}]$
respectively. Identifying $\Lambda_{\sigma}$ and $\Lambda_{\sigma'}$
via parallel transport through $v$ yields an isomorphism
\[
\Phi_{gg'}:\Spec \kk[(\check v+\RR\check\tau)\cap\Lambda^*_{\sigma}]
\rightarrow \Spec\kk[(\check v+\RR\check\tau)\cap\Lambda^*_{\sigma'}],
\]
hence defining an identification of $V_h\subseteq V_{e'}$ with
$V_h\subseteq V_e$, which we also denote by $\Phi_{gg'}$.
Then to construct $X_0$, we glue $V_h\times S
\subset V_{e'}\times S$ and $V_h\times S\subset V_e\times S$
via 
\begin{eqnarray*}
\Phi_{gg'}({\bf s})&:=&{\bf s}^{-1}_g\circ (\Phi_{gg'}\times\id)\circ
{\bf s}_{g'}\\
&=&{\bf s}_e^{-1}\circ (\Phi_{gg'}\times\id)\circ {\bf s}_{e'}.
\end{eqnarray*}
Here ${\bf s}_g$ acts on the coordinate ring of $V_h\times S\subseteq
V_e\times S$ via $z^n\mapsto {\bf s}_g(n)z^n$, where as usual we
view ${\bf s}_g$ as a multiplicative function on $\Lambda_{\sigma}^*$ with
values in $\Gm(S)$.

Now if $n\in\Lambda_{\sigma}^*$ represents a logarithmic 1-form $\dlog n$
on $V_e$, let us consider the pulled-back 1-form $\Phi_{gg'}({\bf s})^*
(\dlog n)$. Clearly $(\Phi_{gg'}\times \id)^*(\dlog n)$ just has the effect
of parallel transporting $n$ from $\Lambda_{\sigma}^*$ to 
$\Lambda_{\sigma'}^*$ through $v$, 
so identify both those spaces with $\check\Lambda_v$. We note
\begin{equation}
\label{transform0}
{\bf s}_e^*(\dlog n)=\dlog {\bf s}_e(n)z^n
=\dlog {\bf s}_e(n)+\dlog n=\dlog s_{id,e}(n)+\dlog n,
\end{equation}
as ${\bf s}_e(n)=s_{0,e}(n)s_{\id,e}(n)$ and 
$\dlog s_{0,e}(n)=0$. Thus if $n\in \bigwedge^r\check\Lambda_v$,
\begin{equation}
\label{transform2}
{\bf s}_e^*(\dlog n)=
\dlog \iota(s_{id,e})n+\dlog n\mod {}{\bigwedge}^2
H^1(B,i_*\Lambda)_f^*,
\end{equation}
where $\iota(s_{id,e})n\in H^1(B,i_*\Lambda)_f^*\otimes
\bigwedge^{r-1}\check\Lambda_{\sigma}$.
In any event, we can write, for $n\in\bigwedge^r\check\Lambda_{\sigma}$,
\begin{equation}
\label{transform1}
\Phi_{gg'}({\bf s})^*(\dlog n)={\bf s}_{e'}^*({\bf s}_e^{-1})^*(\dlog n).
\end{equation}

The upshot of this is that whenever we specify a logarithmic form
on $X_v\times S$ defined over $\kk^{\dagger}$ as $\dlog n$
for some $n\in\bigwedge^r(H^1(B,i_*\Lambda)^*_f\oplus\check\Lambda_v)$,
we need to specify $\sigma\in\P_{\max}$ containing $v$ to indicate in
which affine piece of $X_v\times S$ we are representing this form.
Furthermore, if we wish to represent an element of $\Omega^r_{\tau/
\kk^{\dagger}}$ we should pick $v\mapright{e}\tau\mapright{g}\sigma
\in\P_{\max}$, describe $\Omega^r_{\tau/\kk^{\dagger}}$ as a subsheaf
of $F_{v,\tau}^*\Omega^r_{v/\kk^{\dagger}}$, but describe the
latter sheaf in the affine open subset $V_g\times S$ of $X_{\tau}\times S$.

\smallskip

\emph{Step 4.} \emph{The Gauss-Manin connection
in terms of our resolutions.}

We note that $\Omega^r_{\tau/\kk^{\dagger}}$ has a filtration
\[
F^i(\Omega^r_{\tau/\kk^{\dagger}})=\im 
(\pi^*\Omega^i_{S^{\dagger}/\kk^{\dagger}}\otimes
\Omega^{r-i}_{\tau/\kk^{\dagger}}\rightarrow
\Omega^r_{\tau/\kk^{\dagger}})
\]
and similarly for $e:\tau_1\rightarrow\tau_2$,
\[
F^i((F_{\tau_1,\tau_2}^*\Omega^r_{\tau_1/\kk^{\dagger}})/\Tors)=\im 
(\pi^*\Omega^i_{S^{\dagger}/\kk^{\dagger}}\otimes
(F_{\tau_1,\tau_2}^*\Omega^{r-i}_{\tau_1/\kk^{\dagger}})/\Tors\rightarrow
(F_{\tau_1,\tau_2}^*\Omega^r_{\tau_1/\kk^{\dagger}})/\Tors)
\]
yields a filtration of $(F_{\tau_1,\tau_2}^*
\Omega^r_{\tau_1/\kk^{\dagger}})/\Tors$.
We have
\[
F^i((F_{\tau_1,\tau_2}^*\Omega^r_{\tau_1/\kk^{\dagger}})/\Tors)/
F^{i+1}((F_{\tau_1,\tau_2}^*\Omega^r_{\tau_1/\kk^{\dagger}})/\Tors)
\cong \pi^*\Omega^i_{S^{\dagger}/\kk^{\dagger}}
\otimes (F_{\tau_1,\tau_2}^*\Omega^{r-i}_{\tau_1/S^{\dagger}})/\Tors.
\]
This follows from Theorem \ref{localmodel2}, which tells us \'etale locally
$\Omega^1_{\tau/\kk^{\dagger}}$ splits as $\Omega^1_{S^{\dagger}/\kk^{\dagger}}
\oplus\Omega^1_{\tau/S^{\dagger}}$.
Similarly, this gives a resolution of the exact sequence \eqref{GMexact}:
\[
\xymatrix@C=30pt
{&0\ar[d]&0\ar[d]\\
0\ar[r]&j_*F^1/j_*F^2\ar[r]\ar[d]&\pi^*\Omega^1_{S^{\dagger}/\kk^{\dagger}}
\otimes \C^{\bullet}(j_*\Omega^{\bullet}_{X_0^{\dagger}
/S^{\dagger}}[-1])\ar[d]
\\
0\ar[r]&j_*F^0/j_*F^2\ar[r]\ar[d]&
{\displaystyle\C^{\bullet}
(j_*\Omega^{\bullet}_{X_0^{\dagger}/\kk^{\dagger}})
\over
\displaystyle F^2(\C^{\bullet}
(j_*\Omega^{\bullet}_{X_0^{\dagger}/\kk^{\dagger}}))}\ar[d]\\
0\ar[r]&j_*F^0/j^*F^1\ar[r]\ar[d]&
\C^{\bullet}(j_*\Omega^{\bullet}_{X_0^{\dagger}
/S^{\dagger}})\ar[d]\\
&0&0
}
\]
where $F^2(\C^{\bullet}
(j_*\Omega^{\bullet}_{X_0^{\dagger}/\kk^{\dagger}}))$
is obtained by replacing each $\Omega^r_{\tau/\kk^{\dagger}}$ with
$F^2(\Omega^r_{\tau/\kk^{\dagger}})$.

\emph{Step 5.} \emph{Calculating $\nabla_{GM}$ and the proof of (3) and
(4).}

Consider $n\in H^q(B,i_*\bigwedge^p\check\Lambda\otimes\kk)$,
represented as before by a \v Cech cocycle
\[
(W_{\sigma_0\rightarrow\cdots\rightarrow\sigma_q},n_{\sigma_0\rightarrow
\cdots\rightarrow\sigma_q}).
\]
Then as before, $\dlog n_{\sigma_0\rightarrow\cdots\rightarrow\sigma_q}
\in \Gamma(X_{\sigma_q}\times S,(F_{\sigma_0,\sigma_q}^*
\Omega^p_{\sigma_0/S^{\dagger}})/\Tors)$, and these elements represent a section
of $\RR^{p+q}\pi_*(j_*\Omega^{\bullet}_{X_0/S^{\dagger}})$
over $S$ precisely because they yield a $(p+q)$-cocycle in
$\Gamma(X_0,\Tot(\C^{\bullet}(j_*
\Omega^{\bullet}_{X_0^{\dagger}/S^{\dagger}})))$, keeping
in mind 
\[
d(\dlog n_{\sigma_0\rightarrow\cdots\rightarrow\sigma_q})=0.
\]
To compute the value of the boundary homomorphism on this class, we first
lift this to a $(p+q)$-cochain in 
\[
\Gamma\bigg(X_0, 
\Tot\bigg(
{\C^{\bullet}
(j_*\Omega^{\bullet}_{X_0^{\dagger}/\kk^{\dagger}})
\over
F^2(\C^{\bullet}
(j_*\Omega^{\bullet}_{X_0^{\dagger}/\kk^{\dagger}}))}
\bigg)\bigg),
\]
and then apply the total differential to this cochain to get a $(p+q+1)$-cocycle
for the complex
$\pi^*\Omega^1_{S^{\dagger}/\kk^{\dagger}}
\otimes \Tot(\C^{\bullet}(j_*\Omega^{\bullet}_{X_0^{\dagger}/S^{\dagger}}[-1]))$.

We now carry out this procedure. We need to lift $\dlog n_{\sigma_0\rightarrow
\cdots\rightarrow\sigma_q}$. For convenience, we drop the subscript
for the moment. Pick $e:v\rightarrow\sigma_0$, $g:\sigma_q\rightarrow
\sigma\in\P_{\max}$, let $g_0$ be the composition $\sigma_0\mapright{}
\sigma_q\mapright{g}\sigma$, and set
\[
\tilde n={\bf s}_{g_0}^*(\dlog n)
\]
viewed as an element of $F_{v,\sigma_q}^*
\Omega^p_{v/\kk^{\dagger}}$. 
From \eqref{transform2}, $\tilde n$ is a lifting of $\dlog n$
in 
\[
F_{v,\sigma_q}^*
\Omega^p_{v/\kk^{\dagger}}.
\]
We claim that in fact $\tilde n$
is a section of $(F_{\sigma_0,\sigma_q}^*
\Omega^p_{\sigma_0/\kk^{\dagger}})/\Tors$.
We use Proposition \ref{delta0kernel} to do this. 
Let $\Delta_i,\check\Delta_i$ be the simplicity data for $\sigma_0$,
yielding $Z_i^{\sigma_q}\subseteq X_{\sigma_q}\times S$.
Let $f_i$ be the equation
determining $Z_i^{\sigma_q}$ in the affine subset $V_g\times S$
of $X_{\sigma_q}\times S$
determined by $g$. Each $f_i$ is determined as follows: there
is some $\omega_i\rightarrow\sigma_0$, $\dim\omega_i=1$, $e_i:\omega_i
\rightarrow\sigma$ such that $f_{\sigma,e_i}|_{V_g\times S}=f_i$, and
by \eqref{fsigmae}, 
\begin{eqnarray*}
f_i&=&\sum_{p\in\check\Delta(\omega_i)\cap\Lambda_{\sigma_q}^{\perp}}
{\bf s}_{e_i}(p)z^p\\
&=&\sum_{p\in\check\Delta_i\cap\Lambda_{\sigma_q}^{\perp}}
{\bf s}_{e_i}(p)z^p.
\end{eqnarray*}
But $e_i$ factors as $\omega_i\mapright{h}\sigma_0\mapright{g_0}\sigma$,
and ${\bf s}_h(p)=1$ for $p\in\check\Delta_i\cap\Lambda_{\sigma_q}^{\perp}$
by \eqref{shmonpoly},
so we can write
\[
f_i=\sum_{p\in\check\Delta_i\cap\Lambda^{\perp}_{\sigma_q}}
{\bf s}_{g_0}(p)z^p={\bf s}_{g_0}^*(f^0_i),
\]
where
\[
f_i^0=\sum_{p\in\check\Delta_i\cap\Lambda_{\sigma_q}^{\perp}}z^p.
\]
Since $\dlog n\in \Gamma(X_{\sigma_q}\times S, F_{\sigma_0,
\sigma_q}^*\Omega^p_{\sigma_0/S^{\dagger}}/\Tors)$,
it follows from Proposition~\ref{delta0kernel}, Remark~\ref{poleremark}
and Lemma~\ref{OmegaZ}
that for each $v_i'\not=v_i$ a vertex of $\Delta_i$, 
${df_i\over f_i}\wedge\dlog
\iota(v_i'-v_i)n$ has no pole along the locus $f_i=0$,
as a form defined over $S^{\dagger}$. However, over $S^{\dagger}$, this
is equivalent to the statement that ${df_i^0\over f_i^0}\wedge
\dlog\iota(v_i'-v_i)n$ has no pole along the locus $f_i^0=0$ (i.e. the
coefficients don't affect whether or not this holds). As forms 
over $\kk^{\dagger}$, there is some invertible function $g$ pulled back
from $S$ such that
\[
{\bf s}_{g_0}^*\left({df_i^0\over f_i^0}\wedge\dlog(\iota(v_i'-v_i)n)\right)
=g{df_i\over f_i}\wedge\dlog\iota(v_i'-v_i){\bf s}_{g_0}^*(\dlog n).
\]
Thus this latter form has no poles along $Z_i^{\sigma_q}$. So
\[
\dlog\tilde n\in\Gamma(X_{\sigma_q}\times S,
(F_{\sigma_0,\sigma_q}^*\Omega^p_{\sigma_q/\kk^{\dagger}})
/\Tors),
\]
by Proposition \ref{delta0kernel}.
Note in fact this lifting is independent of the choice of $\sigma$:
given a diagram
\[
\xymatrix@C=30pt
{&&\sigma\\
\sigma_0\ar@/^/[rru]^{g_0}\ar[r]\ar@/_/[drr]_{g_0'}&\sigma_q
\ar[ru]_{g}\ar[rd]^{g'}&\\
&&\sigma'
}
\]
we note that we get liftings ${\bf s}_{g_0}^*(\dlog n)$, ${\bf s}_{g'_0}^*(\dlog
n)$ on
$V_g\times S$ and $V_{g'}\times S$ respectively. These are in fact
identified under the transformation \eqref{transform1}. 

Thus we have a well-defined lifting of our $(p+q)$-cocycle, to get
the $(p+q)$-cochain $(\tilde n_{\sigma_0\rightarrow\cdots\rightarrow\sigma_q})$.
We now apply the total differential to this. Since the liftings $\tilde n$
are themselves $d$-closed, this just means applying $d_{\bct}$. 
In particular, (3) is now clear. For (4) to see
what the contribution to $\sigma_0\rightarrow\cdots\rightarrow\sigma_{q+1}$
is,
choose $g:\sigma_{q+1}\rightarrow\sigma$, and let $g_i:\sigma_i\mapright{}
\sigma_{q+1}\mapright{g}\sigma$ be the composition. Then using the
representation determined by $\sigma$, the contribution is
\[
{\bf s}_{g_1}^*(\dlog n_{\sigma_1\rightarrow\cdots\rightarrow\sigma_{q+1}})
+
\sum_{i=1}^{q+1}(-1)^i{\bf s}_{g_0}^*(\dlog n_{\sigma_0\rightarrow\cdots
\rightarrow\hat\sigma_i\rightarrow\cdots\rightarrow\sigma_{q+1}}).
\]
We view this modulo $\Omega^2_{S^{\dagger}/\kk^{\dagger}}$, noting
from \eqref{transform2} that 
\[
{\bf s}_g^*(\dlog n)=\dlog \iota(s_{\id,g})n+\dlog n\mod \pi^*
\Omega^2_{S^{\dagger}/\kk^{\dagger}},
\]
so we obtain
\begin{align*}
-\dlog\iota(s_{\id,\sigma_0\rightarrow\sigma_1})
n_{\sigma_1\rightarrow\cdots\rightarrow\sigma_{q+1}}&
+\sum_{i=0}^{q+1} (-1)^i\dlog n_{\sigma_0\rightarrow\cdots\rightarrow
\hat\sigma_i\rightarrow\cdots\rightarrow\sigma_{q+1}}\\
&+\sum_{i=0}^{q+1}(-1)^i\dlog\iota(s_{\id,g_0})n_{\sigma_0\rightarrow
\cdots\rightarrow\hat\sigma_i\rightarrow\cdots\rightarrow\sigma_{q+1}}.
\end{align*}
Now as $(\sigma_0\rightarrow\cdots\rightarrow\sigma_q,
\dlog n_{\sigma_0\rightarrow\cdots\rightarrow\sigma_q})$ is already
$d_{\bct}$-closed over $S^{\dagger}$, the last two terms in
fact vanish, and we are left with the $p+q+1$-cocycle
\[
(\sigma_0\rightarrow\cdots\rightarrow\sigma_{q+1},
-\dlog\iota(s_{\id,\sigma_0\rightarrow\sigma_1})n_{\sigma_1\rightarrow\cdots
\rightarrow\sigma_{q+1}}).
\]
By \cite{Bredon}, III, 4.15, this represents in \v Cech cohomology the negative
of
the image of $s_{\id}\otimes n$ under the composition
{\scriptsize
\begin{eqnarray*}
H^1(B,i_*\Lambda\otimes H^1(B,i_*\Lambda)_f^*)\otimes
H^q(B,i_*{\bigwedge}^p\check\Lambda\otimes\kk)
&\rightarrow& H^{q+1}(B,(i_*\Lambda\otimes H^1(B,i_*\Lambda)_f^*)
\otimes i_*{\bigwedge}^p\check\Lambda\otimes\kk)\\
&\rightarrow& 
H^{q+1}(B,H^1(B,i_*\Lambda)^*_f\otimes i_*{\bigwedge}^{p-1}
\check\Lambda\otimes\kk)\\
&=&H^1(B,i_*\Lambda\otimes\kk)^*\otimes H^{q+1}(B,i_*{\bigwedge}^{p-1} 
\check\Lambda\otimes\kk),
\end{eqnarray*}
}\noindent
where the first map is cup product and the second is contraction.
However, this map is the same as the map
\[
H^1(B,i_*\Lambda\otimes\kk)\otimes 
H^q(B,i_*{\bigwedge}^p\check\Lambda\otimes\kk)
\rightarrow H^{q+1}(B,i_*{\bigwedge}^{p-1}\otimes\kk),
\]
tensored with $H^1(B,i_*\Lambda\otimes\kk)^*$. 
Now in general, if $U,V$ and $W$ are vector spaces, the canonical
identification $\Hom(U\otimes V,W)=\Hom(V,U^*\otimes W)$
can be realised by, given a map $\varphi:U\otimes V\rightarrow
W$, tensoring with $U^*$ to get $\id_{U^*}\otimes\varphi:
U^*\otimes U\otimes V\rightarrow U^*\otimes W$, which induces
a map $V\rightarrow U^*\otimes W$ via $v\mapsto (\id_{U^*}\otimes
\varphi)(\id_U\otimes v)$, thinking of $\id_U\in U^*\otimes U$.
Since $s_{\id}$
is just the identity in $H^1(B,i_*\Lambda\otimes\kk)\otimes
H^1(B,i_*\Lambda\otimes\kk)^*$, the result stated in (4) follows.
\qed

\begin{corollary}
Choose a basis $\eta_1,\ldots,\eta_h$ of $H^1(B,i_*\Lambda)_f^*$,
and let $\eta_1^*,\ldots,\eta_h^*$ be the dual basis. Let $q_1,
\ldots,q_h$ be the corresponding monomials in $\kk[H^1(B,i_*\Lambda)_f^*]$.
If 
\[
n\in H^q(B,i_*{\bigwedge}^p\check\Lambda\otimes\kk),
\]
then
\begin{eqnarray*}
n+\sum_{i_1=1}^h (\iota(\eta_{i_1}^*)n)\log q_i&+&
\sum_{i_1,i_2=1}^h (\iota(\eta_{i_1}^*)\iota(\eta_{i_2}^*)n)
{\log q_i\log q_j\over 2}\\
\cdots&+&
\sum_{i_1,\ldots,i_p=1}^h (\iota(\eta^*_{i_1})\ldots\iota(\eta^*_{i_p})n)
{\log q_{i_1}\cdots\log q_{i_p}\over p!}
\end{eqnarray*}
is flat with respect to $\nabla_{GM}$. Here $\log q_i$ can be viewed
formally via the property $d(\log q_i)=\dlog q_i$.
\end{corollary}

\proof Simply apply $\nabla_{GM}$ to this section.
\qed

\begin{remark} We have now in fact described a variation of mixed
Hodge structures over $S$: indeed, 
$\RR^r\pi_*(j_*\Omega^{\bullet}_{X_0^{\dagger}/
S^{\dagger}})$ is a vector
bundle on $S$ with $\nabla_{GM}$ a flat connection defining a local
system of flat sections of this vector bundle. We have the
Hodge filtration
\[
\shF^i=\bigoplus_{p\ge i} R^{r-p}\pi_*(j_*\Omega^p_{X_0^{\dagger}
/S^{\dagger}})
\]
and weight filtration
\[
W_{2i}=W_{2i+1}
=\bigoplus_{p\le i} R^{r-p}\pi_*(j_*\Omega^p_{X_0^{\dagger}
/S^{\dagger}}).
\]
What is missing from this description is the \emph{integral} structure
on the local system. This will be addressed in \cite{TopologyPaper}.
\end{remark}


\end{document}